\documentclass[11pt,a4paper]{article}
\RequirePackage[latin1]{inputenc}
\RequirePackage{amsmath,amssymb,graphics}
\RequirePackage[xdvi]{epsfig} \RequirePackage{eepic}
\RequirePackage{epic}
\newcommand{\epspic}[4]{\begin{figure}[!ht]
\begin{center}
\epsfig{file=#1,height=#2,angle=#3}
\end{center}
\caption{\emph{#4}}
\label{eps:#1}
\end{figure}}

\usepackage{amsfonts,amssymb,epsfig}
\usepackage[latin1]{inputenc}
\usepackage{amsmath,amsthm,amssymb}
\usepackage{graphicx,graphics,epsfig}

\def\proof{{\it Proof.}\ }

\def\eq#1{(\ref{#1})}

\def\neweq#1{\begin{equation}\label{#1}}
\def\endeq{\end{equation}}
\def\weak{\rightharpoonup}
\def\ep{\varepsilon}
\def\la{\lambda}

\def\RR{{\mathbb R} }

\def\di{\displaystyle}
\def\ri{\rightarrow}

\def\un{\underline}
\def\ov{\overline}
\def\la{\lambda}
\newtheorem{thm}{Theorem}[section]
 \newtheorem{cor}[thm]{Corollary}
 \newtheorem{lem}[thm]{Lemma}
 \newtheorem{prop}[thm]{Proposition}
\newtheorem{defin}[thm]{Definition}
\newtheorem{rem}{Remark}
\title{ \huge Singular phenomena in nonlinear elliptic problems\\
\Large From blow-up boundary solutions to equations with singular nonlinearities}
\author{
Vicen\c tiu D. R\u ADULESCU\thanks{The author is partially
supported by Grant CEEX 05-D11-36 {\it Analysis and
Control of Differential Systems}.}\\
\small Department of Mathematics, University of Craiova, 200585
Craiova, Romania\\ \small {\tt http://inf.ucv.ro/$\sim$radulescu}
\qquad E-mail: {\tt vicentiu.radulescu@math.cnrs.fr}}

\date{}

\begin{document}
\baselineskip16pt \maketitle
\renewcommand{\theequation}{\arabic{section}.\arabic{equation}}

{\small {\bf Abstract.} In this survey we report on some recent results
related to various singular phenomena arising in the study of some classes
of nonlinear elliptic equations. We establish qualitative results on the
existence, nonexistence or the uniqueness
of solutions and we focus on the following types of problems:
(i) blow-up boundary solutions of logistic equations; (ii) Lane-Emden-Fowler equations
with singular nonlinearities and subquadratic convection term.
We study the combined effects of various terms involved in these problems:
sublinear or superlinear nonlinearities, singular nonlinear terms, convection nonlinearities,
as well as sign-changing potentials.
We also take into account bifurcation nonlinear problems
 and we establish the precise rate decay of the solution
 in some concrete situations.  Our approach
combines standard techniques
based on the maximum principle with
non-standard arguments, such as the Karamata regular variation theory.


 }

\section{Motivation and Previous Results}
Let $\Omega$ be a bounded domain with smooth boundary in $\RR^N$,
$N\geq 2$. We are concerned in this paper with the following types
of stationary singular problems:

{\bf I. The logistic equation}
\begin{equation}\label{problem1} \left\{\begin{tabular}{ll}
$\Delta u=\Phi (x,u,\nabla u)$ \quad & $\mbox{\rm in}\ \Omega,\ $\\
$u>0$ \quad & $\mbox{\rm in}\ \Omega,$\\
$u=+\infty$ \quad & $\mbox{\rm on}\ \partial\Omega$.
\end{tabular} \right.\end{equation}

{\bf II. The Lane-Emden-Fowler equation}
\begin{equation}\label{problem2} \left\{\begin{tabular}{ll}
$-\Delta u=\Psi (x,u,\nabla u)$ \quad & $\mbox{\rm in}\ \Omega,\ $\\
$u>0$ \quad & $\mbox{\rm in}\ \Omega,$\\
$u=0$ \quad & $\mbox{\rm on}\ \partial\Omega,$
\end{tabular} \right.\end{equation}
where $\Phi$ is a smooth nonlinear function, while $\Psi$ has one or more singularities.
The solutions of \eq{problem1} are called {\it large}
(or {\it blow-up}) {\it solutions}.

In this work we focus on Problems \eq{problem1} and \eq{problem2}
and we establish several recent contributions in the study of
these equations. In order to illustrate the link between these
problems, consider the most natural case where $\Phi (u,\nabla
u)=u^p$, where $p>1$. Then the function $v=u^{-1}$ satisfies
\eq{problem2} for $\Psi (u,\nabla v)=v^{2-p}-  2v^{-1}\,|\nabla
v|^2$.

The study of large solutions has been initiated in 1916 by
Bieberbach \cite{bi} for the particular case $\Phi (x,u,\nabla
u)=\exp(u)$ and $N=2$. He showed that there exists a unique
solution of \eq{problem1} such that $u(x)-\log (d(x)^{-2})$ is
bounded as $x\to\partial\Omega$, where
$d(x):=\mbox{dist}\,(x,\partial\Omega)$. Problems of this type
arise in Riemannian geometry: if a Riemannian metric of the form
$|ds|^2=\exp(2u(x))|dx|^2$ has constant Gaussian curvature $-c^2$
then $\Delta u=c^2 \exp(2u)$. Motivated by a problem in
mathematical physics, Rademacher \cite{rad} continued the study of
Bieberbach on smooth bounded domains in $ {\mathbb R}^3 $. Lazer
and McKenna \cite{lm} extended the results of Bieberbach and
Rademacher for bounded domains in ${\mathbb R}^N$ satisfying a
uniform external sphere condition and for nonlinearities $\Phi
(x,u,\nabla u)=b(x)\exp(u)$, where $b$ is continuous and strictly
positive on $\overline{\Omega}$. Let $\Phi (x,u,\nabla u)=f(u)$
where $f\in C^1[0,\infty)$, $f'(s)\geq 0$ for $ s\geq 0$, $f(0)=0$
and $ f(s)>0$ for $ s>0$. In this case, Keller \cite{ke} and
Osserman \cite{os} proved that large solutions of \eq{problem1}
exist if and only if
$$
\int_1^\infty  \frac{dt}{\sqrt{F(t)}} <\infty,\quad \mbox{where}\
F(t)=\int_0^t f(s)\,ds.
$$
In a celebrated paper, Loewner and Nirenberg \cite{ln}
linked the uniqueness of the blow-up solution to
the growth rate at the boundary.
Motivated by certain geometric problems, they
established the uniqueness for the case $f(u)=u^{(N+2)/(N-2)}$, $N>2$.
Bandle and Marcus \cite{bm} give results on asymptotic behaviour and uniqueness
of the large solution for more general nonlinearities including
$f(u)=u^p$ for any $p>1$. We refer to Bandle \cite{bandle},
Bandle and M. Ess\`{e}n \cite{be}, Bandle and Marcus \cite{bm1}, Du and Huang \cite{dh},
Garc\'{i}a-Meli\'{a}n, Letelier-Albornoz, and Sabina de Lis \cite{gls},
Lazer and McKenna \cite{lm2}, Le Gall \cite{gall}, Marcus and V\'eron \cite{mv1,mv2},
Ratto, Rigoli and V\'eron \cite{rrv}
and the references therein for several results on
large solutions extended to $N$-dimensional
domains and for other classes of nonlinearities.

Singular problems like \eq{problem2} have been intensively studied in the last decades.
Stationary problems involving singular nonlinearities,
 as well as the associated evolution equations,
describe naturally several physical phenomena. At our best
knowledge, the first study in this direction is due to Fulks and
Maybee \cite{ful}, who proved existence and uniqueness results by
using a fixed point argument; moreover, they showed that solutions
of the associated parabolic problem tend to the unique solution of
the corresponding elliptic equation. A different approach (see
Coclite and Palimieri \cite{cp}, Crandall, Rabinowitz, and Tartar
\cite{crt}, Stuart \cite{stu}) consists in approximating the
singular equation with a regular problem, where the standard
techniques (e.g., monotonicity methods) can be applied and then
passing to the limit to obtain the solution of the original
equation. Nonlinear singular boundary value problems arise in the
context of chemical heterogeneous catalysts and chemical catalyst
kinetics, in the theory of heat conduction in electrically
conducting materials, singular minimal surfaces, as well as in the
study of non-Newtonian fluids, boundary layer phenomena for
viscous fluids (we refer for more details to Caffarelli, Hardt,
and L. Simon \cite{caf}, Callegari and Nachman \cite{cn1,cn},
D\'iaz \cite{diaz}, D\'iaz,  Morel, and Oswald \cite{dmo} and the
more recent papers by Haitao \cite{hai}, Hern\'andez, Mancebo, and
 Vega \cite{her,her1}, Meadows \cite{mea}, Shi and Yao \cite{shi,sy2}).
We also point out that, due to the meaning of the unknowns
(concentrations, populations, etc.), only the positive solutions
are  relevant in most cases. For instance, problems of this type
characterize some reaction-diffusion processes where $u\geq 0$ is
viewed as the density of a reactant and the region where $u=0$ is
called the {\it dead core}, where no reaction takes place (see
Aris \cite{a} for the study of a single, irreversible steady-state
reaction). Nonlinear singular elliptic equations are also
encountered in glacial advance, in transport of coal slurries down
conveyor belts and in several other geophysical and industrial
contents (see Callegari and Nachman \cite{cn}  for the case of the
incompressible flow of a uniform stream past a semi-infinite flat
plate at zero incidence).

In \cite{crt}, Crandall, Rabinowitz and Tartar established that the boundary value problem
$$
 \left\{\begin{tabular}{ll}
$-\Delta u-u^{-\alpha}= -u$ \quad & ${\rm in}\
\Omega,$\\
$u>0$ \quad & ${\rm in}\ \Omega,$\\
$u=0$ \quad & ${\rm on}\ \partial\Omega$
\end{tabular} \right.
$$
has a solution, for any $\alpha>0$. The importance of the linear and nonlinear
terms is crucial for the existence of solutions. For instance, Coclite and Palmieri
studied in \cite{cp}  the
problem
\begin{equation} \label{unuunu}
 \left\{\begin{tabular}{ll}
$-\Delta u-u^{-\alpha}=\lambda u^p$ \quad & ${\rm in}\
\Omega,$\\
$u>0$ \quad & ${\rm in}\ \Omega,$\\
$u=0$ \quad & ${\rm on}\ \partial\Omega,$\\
\end{tabular} \right.
\end{equation}
where $\lambda\geq 0$ and $\alpha,p\in(0,1).$ In \cite{cp} it is
proved that problem \eqref{unuunu} has at least one solution for
all $\lambda\geq 0$ and $0<p<1$. Moreover, if $p\geq 1,$ then
there exists $\lambda^*$ such that problem \eqref{unuunu} has a
solution for $\lambda\in[0,\lambda^*)$ and no solution for
$\lambda>\lambda^*$. In \cite{cp} it is also proved a related
non--existence result. More exactly, the  problem
$$
 \left\{\begin{tabular}{ll}
$-\Delta u+u^{-\alpha}= u$ \quad & ${\rm in}\
\Omega,$\\
$u>0$ \quad & ${\rm in}\ \Omega,$\\
$u=0$ \quad & ${\rm on}\ \partial\Omega$\\
\end{tabular} \right.
$$
has no solution, provided that $0<\alpha <1$ and $\lambda_1\geq 1$
(that is, if $\Omega$ is ``small''), where $\lambda_1$ denotes the
first eigenvalue of $(-\Delta)$ in $H^1_0(\Omega)$.

Problems related to multiplicity and uniqueness become difficult
even in simple cases. Shi and Yao studied in \cite{shi} the existence of
radial symmetric solutions of the problem
$$
 \left\{\begin{tabular}{ll}
$\displaystyle \Delta u+\lambda(u^p-u^{-\alpha})=0$ \quad & ${\rm
in}\ B_1,$\\
$u>0$ \quad & ${\rm in}\ B_1,$\\
$u=0$ \quad & ${\rm on}\ \partial B_1,$\\
\end{tabular} \right.
$$
where $\alpha>0$, $0<p<1,$ $\lambda>0$, and $B_1$ is the unit ball
in $\RR^N$. Using a bifurcation theorem of Crandall and
Rabinowitz, it has been shown in \cite{shi}  that there exists
$\lambda_1>\lambda_0>0$ such that the above problem has no
solutions for $\lambda<\lambda_0,$ exactly one solution for
$\lambda=\lambda_0$ or $\lambda>\lambda_1,$  and two solutions for
$\lambda_0 < \lambda\leq\lambda_1$.

The author's interest for the study of singular problems is
motivated by several stimulating discussions with Professor Haim
Brezis in Spring 2001. I would like to use this opportunity to
thank once again Professor Brezis for his constant scientific
support during the years.

This work is organized as follows. Sections 2--5 are mainly
devoted to the study of blow-up boundary solutions of logistic
type equations with absorption.  In the second part of this work
(Sections 6--8), in connection with the previous results, we are
concerned with the study of the Dirichlet boundary value problem
for the singular Lane-Emden-Fowler equation. Our framework
includes the presence of a convection term.

\section{Large solutions of elliptic equations with absorption and subquadratic convection term}
Consider the problem
\neweq{unu}
 \left\{\begin{tabular}{ll}
$\Delta u+q(x)|\nabla u|^a=p(x)f(u)$ \quad & ${\rm in}\ \Omega\,,$\\
$u\geq 0,\ u\not\equiv 0$ \quad & ${\rm in}\ \Omega\,,$\\
\end{tabular} \right.
\endeq
where
$\Omega\subset\RR^N$ ($N\geq 3$) is a smooth domain (bounded or
possibly unbounded) with compact (possibly empty) boundary.
We assume that $a\leq 2$ is a positive real
number, $p,q$ are non-negative function such that $p\not\equiv 0$,
$p,q\in C^{0,\alpha}(\overline\Omega )$ if $\Omega$ is bounded, and
$p,q\in C^{0,\alpha}_{\rm loc}(\Omega )$, otherwise. Throughout this section we
assume that the nonlinearity $f$ fulfills the following conditions

\smallskip
\noindent $(f1)\qquad f\in C^{1}[0,\infty),\ f'\geq 0,
\ f(0)=0$ and $ f>0$ on $(0,\infty )$.

\smallskip
\noindent $(f2)\qquad \di\int\limits_{1}^{\infty}
[F(t)]^{-1/2}\,dt<\infty\,,
 \quad \mbox{where}\quad
F(t)=\di\int\limits_{0}^{t}f(s)\,ds. $

\smallskip
\noindent$(f3)\qquad \di\frac{F(t)}{f^{2/a}(t)} \quad\to 0$\quad as
 \quad $t\to \infty$.

\smallskip
Cf. V\'eron \cite{ver}, $f$ is called an
absorption term.
 The above conditions hold provided that $f(t)=t^k$,
$k>1$
and $0<a<\frac{2r}{r+1}(<2)$, or $f(t)=e^t-1$, or $f(t)=e^t-t$ and $a<2$.
We observe that by $(f1)$ and $(f3)$ it follows that
$f/F^{a/2}\geq \beta>0$ for  $t$ large
enough, that is, $(F^{1-a/2})'\geq \beta>0$ for $t$ large enough
which yields $0<a\leq 2$.
We also deduce that conditions $(f2)$ and $(f3)$ imply
$\di\int\limits_{1}^{\infty}f^{-1/a}(t)dt<\infty\,.$

We are mainly interested in finding properties of {\it large (explosive) solutions}
of \eq{unu}, that is solutions $u$
satisfying $u(x)\to
\infty$ as $ {\rm dist}\, (x,\partial\Omega)\to 0 $ (if
$\Omega\not\equiv \RR^N$), or
$u(x)\to
\infty$ as $|x|\to\infty$ (if $\Omega = \RR^N$). In the latter case the
solution is called an {\it entire large (explosive)} solution.

Problems of this type appear in stochastic control theory and have
been first study by Lasry and Lions \cite{ll}. The corresponding
parabolic equation was considered in Quittner \cite{qu} and in
Galaktionov and V\'azquez \cite{gal}. In terms of the dynamic
programming approach, an explosive solution of \eq{unu}
corresponds to a value function (or Bellman function) associated
to an infinite exit cost (see Lasry and Lions \cite{ll}).

Bandle and Giarrusso \cite{bangia} studied the
existence of a large solution of
problem \eq{unu} in the case $p\equiv 1$, $q\equiv 1$ and  $\Omega$
bounded. Lair and Wood \cite{lw} studied the sublinear case corresponding to $p\equiv 1$,
while C\^\i rstea and R\u adulescu \cite{crna} proved the existence of  large
solutions to \eq{unu}
in the case $q\equiv 0$.

As observed by Bandle and Giarrusso \cite{bangia}, the simplest
case is $a=2$, which can be reduced to a problem without gradient
term. Indeed, if
 $u$ is a solution of \eq{unu} for $q\equiv 1$, then
the function $v=e^u$ (Gelfand transformation) satisfies
$$
 \left\{\begin{tabular}{ll}
$\Delta v=p(x)vf(\ln v)$ \quad & ${\rm in}\ \Omega\,,$\\
$v(x)\ri +\infty$ \quad & ${\rm if}\  {\rm dist}\, (x,\partial\Omega)\ri 0.$\\
\end{tabular} \right.
$$
We shall therefore mainly consider the case where $0<a<2$.

The main results in this Section are due to Ghergu, Niculescu, and
R\u adulescu \cite{gnr}. These results generalize those obtained
by C\^{\i}rstea and R\u adulescu \cite{crna} in the case of the
presence of a convection (gradient) term.

Our first result concerns the existence of a large solution to problem
\eq{unu} when $\Omega$ is bounded.

\begin{thm}\label{t1} Suppose that $\Omega$ is bounded and assume that $p$ satisfies

\noindent $(p1)$\
for every $x_0\in \Omega$ with $p(x_0)=0$, there
exists a domain $\Omega_0\ni x_0$ such that
$\overline{\Omega_0}\subset
\Omega$ and $p>0$ on $\partial\Omega_0$.

 Then problem
\eq{unu} has a positive large solution.\end{thm}

A crucial role in the proof of the above result is played by the
following  auxiliary result (see Ghergu, Niculescu, and R\u
adulescu \cite{gnr}).

\begin{lem} \label{l1} Let  $\Omega$ be a bounded domain. Assume that
$p,q\in C^{0,\alpha}(\overline\Omega )$ are non-negative functions,
$0<a<2$
is a real number,
$f$ satisfies $(f1)$
and $g:\partial\Omega\ri (0,\infty )$
is continuous. Then the boundary value problem
\neweq{trei}
 \left\{\begin{tabular}{ll}
$ \di \Delta u+q(x)|\nabla u|^a=p(x)f(u) $, & \quad $\mbox{in}\ \Omega
$\\
$ \di u=g $, & \quad $ \mbox{on}\ \partial\Omega $\\
$\di u\geq  0,\ u\not\equiv 0$, & \quad $\mbox{in}\ \Omega $\\
\end{tabular}\right.
\endeq
has a classical solution. If $p$ is positive, then the solution is
unique.
\end{lem}

{\it Sketch of the proof of Theorem \ref{t1}.} By Lemma \ref{l1}, the boundary value problem
 $$\left\{\begin{tabular}{ll}
$ \di \Delta v_{n}+q(x)|\nabla v_n|^a=\left(p(x)+\frac{1}{n}
\right)f(v_{n}) $, & \quad $ \mbox{in}\ \Omega $\\
$ \di v_{n}=n $, & \quad $ \mbox{on}\ \partial\Omega $\\
$v_n\geq 0,\ v_n\not\equiv 0$, & \quad $ \mbox{in}\ \Omega $\\
\end{tabular}\right.$$
has a unique positive solution, for any $n\geq 1$.
Next, by the maximum principle,
the sequence $(v_n)$ is non-decreasing and is bounded from below in $\Omega$
by a positive function.

To conclude the proof, it is sufficient to show that

$(a)$  for all $x_0 \in \Omega$ there exists an
open set ${\cal O}\subset\subset\Omega $ which
contains $x_0$
 and $M_0=M_0(x_0)>0$ such
that $v_n\leq M_0 $
in ${\cal O}$ for all $n\geq 1$

$(b)$  $\lim_{x\to \partial\Omega}v(x)=\infty$,
where $v(x)=\lim_{n\to \infty}v_n(x).$

We observe that the statement $(a)$ shows that the sequence
$(v_n)$
is uniformly bounded on every compact subset of $\Omega$.
Standard elliptic regularity arguments (see Gilbarg and Trudinger \cite{gt})
show that $v$ is a solution
of problem \eq{unu}. Then, by
(b), it follows that $v$ is a large solution of problem \eq{unu}.

To prove $(a)$ we distinguish two cases :

\noindent
{\sc Case $p(x_0)>0$.} By the continuity of $p$,
there exists a ball
$B=B(x_0,r)\subset\subset\Omega$ such that
$$m_0:=\min\,\{p(x);\ x\in\overline{B}\}>0.$$
Let $w$ be a positive solution of the problem
 $$\left\{\begin{tabular}{lll}
$ \di \Delta w+q(x)|\nabla w|^a=m_0 f(w) $,  & \quad $\mbox{in}\ B$\\
$ \di w(x)\to \infty$,  &\quad  $ \mbox{as}\
x\to \partial B$.\\
\end{tabular}\right. $$
The existence of $w$ follows by considering the problem
$$\left\{\begin{tabular}{lll}
$ \di \Delta w_n+q(x)|\nabla w_n|^a=m_0 f(w_n) $,  & \quad $\mbox{in}\ B$\\
$ \di w_n=n$,  &\quad  $ \mbox{on}\ \partial B$.\\
\end{tabular}\right. $$
The maximum principle implies $w_n\leq w_{n+1}\leq\theta$, where
$$\left\{\begin{tabular}{lll}
$ \di \Delta \theta +\| q\|_{L^\infty}|\nabla \theta|^a=m_0 f(\theta ) $,  & \quad $\mbox{in}\ B$\\
$ \di \theta (x)\to \infty$,  &\quad  $ \mbox{as}\
x\to \partial B$.
\end{tabular}\right. $$

Standard arguments show that $v_n\leq w$ in $B$.
Furthermore,
$w$ is bounded in $\overline{B(x_0,r/2)}$. Setting
$M_0=\sup\limits_{{\cal O}}\ w$, where ${\cal O}=B(x_0,r/2)$,
we obtain (a).

\smallskip
\noindent {\sc Case $p(x_0)=0$.} Our hypothesis $(p1)$ and the
boundedness
of $\Omega$ imply
the existence of a domain ${\cal O}\subset\subset\Omega$ which contains
$x_0$ such that $p>0$ on $\partial {\cal O}$. The above
case shows that for any
$x\in \partial{\cal O}$ there exist a ball $B(x,r_x)$ strictly
contained in $\Omega$ and a constant $M_x>0$ such that $v_n\leq M_x$ on
$B(x,r_x/2)$, for any $n\geq
1$.
Since
$\partial{\cal O}$ is compact, it follows that it
may be covered by a finite number of such balls, say
$B(x_i,r_{x_i}/2)$, $i=1,\cdots,k_0 $.
Setting $M_0=\max\,\{M_{x_1},\cdots,
M_{x_{k_0}}\}$ we have $v_n\leq M_0$ on $\partial{\cal O}$, for any
$n\geq 1$.
Applying  the maximum principle
 we obtain $v_n\leq M_0$ in ${\cal O}$ and $(a)$ follows.

\smallskip
Let $z$ be the unique function satisfying
$-\Delta z=p(x)$ in $\Omega$ and
$z=0$,  on $\partial\Omega.$
Moreover, by the maximum principle, we have $z>0$ in $\Omega$.
We first observe that for proving $(b)$ it is sufficient to show that
\neweq{sapte} \int\limits_{v(x)}^{\infty}
\frac{dt}{f(t)}\leq z(x) \quad \mbox{for any}\ x\in\Omega. \endeq
By \cite [Lemma 1]{crna}, the left hand-side of
\eq{sapte} is well defined
in $\Omega$. We choose $R>0$ so that $\overline\Omega
\subset B(0,R)$ and fix $\ep>0$.
Since $v_n=n$ on $\partial\Omega$, let
$n_1=n_1(\ep)$ be such that
\neweq{opt}
n_1\, >\frac{1}{\ep(N-3)(1+R^2)^{-1/2}+3\ep(1+R^2)^{-5/2}}\,,
\endeq
and
\neweq{noua}
  \di \int\limits_{v_n(x)}^{\infty}\frac{dt}{f(t)}
  \leq z(x)+\ep(1+|x|^2)^{-1/2}
\qquad \di \forall \;x\in\partial\Omega\,, \forall\; n\geq n_1\,.\\
 \endeq

In order to prove \eq{sapte}, it is enough to show that
\neweq{zece} \int\limits_{v_n(x)}^{\infty}\frac{dt}{f(t)}\leq
z(x)+\ep(1+|x|^2)^{-1/2}\qquad \forall \;x\in\Omega\,,\ \forall n\geq
n_1.\endeq
Indeed, taking $n\to \infty$ in \eq{zece} we deduce \eq{sapte}, since
$\ep >0$ is arbitrarily chosen.
Assume now, by contradiction, that \eq{zece} fails. Then
$$\max\limits_{x\in\overline
\Omega}\left\{\int\limits_{v_n(x)}^{\infty}
\frac{dt}{f(t)}-z(x)-\ep(1+|x|^2)^{-1/2}\right\}>0.$$
Using \eq{noua} we see that the point where the maximum is achieved
must
lie in $\Omega$.
A straightforward computation shows that at this point,
say $x_0$, we have
$$0\geq \Delta\left(\int\limits_{v_n(x)}^{\infty}
\frac{dt}{f(t)}-z(x)-\ep(1+|x|^2)^{-1/2}\right)_{|x=x_0} >0.$$
This contradiction shows that inequality \eq{noua} holds and the proof
of Theorem~1 is complete.
\qed

Similar arguments based on the maximum principle and the approximation of large balls $B(0,n)$
imply the following existence result.

\begin{thm}\label{theott}
Assume that $\Omega =\RR^N$ and that problem \eq{unu} has at least one
solution.
Suppose that $p$ satisfies the condition

\noindent $(p1)'\quad$ There exists a sequence of smooth bounded
domains $(\Omega_n)_{n\geq 1}$ such that $\overline{\Omega_n}
\subset \Omega_{n+1}$, $\RR^N=\cup_{n=1}^{\infty}
\Omega_n$,
and $(p1)$ holds in  $\Omega_n$, for any $n\geq
1$.

Then there exists a classical solution $U$ of \eq{unu} which
is a maximal
solution if $p$ is positive.

Assume that $p$ verifies the additional condition

\noindent $(p2)\qquad \di\int\limits_{0}^{\infty} r
\Phi(r)\,dr<\infty\,,
\quad \mbox{where}\ \Phi(r)=\max\,\{p(x):\ |x|=r\}$.

Then $U$ is an entire large solution of \eq{unu}. \end{thm}

We now consider the case in which $\Omega\not= \RR^{N}$ and $\Omega$ is
unbounded. We say that
a large solution $u$ of \eq{unu} is {\it regular}
if $u$ tends to zero at infinity. In \cite[Theorem 3.1]{marcus7} Marcus
proved for this case (and if $q=0$) the existence of regular
large solutions to problem \eq{unu}
by assuming that
there exist $\gamma>1$ and $\beta>0$
such that
$$ \liminf_{t\to 0}f(t)t^{-\gamma}>0 \qquad\mbox{and}\qquad
\liminf_{|x|\to \infty}p(x)|x|^{\beta}>0.$$
The large solution constructed in Marcus \cite{marcus7} is the {\it smallest}
large solution of problem \eq{unu}. In the next result we show that
problem \eq{unu} admits a {\it maximal} classical solution $U$ and that $U$
blows-up at infinity if $\Omega=
\RR^N\setminus \overline{B(0,R)}$.

\begin{thm}\label{theor} Suppose that $\Omega\not=\RR^N$ is
unbounded and that problem \eq{unu} has at least a solution. Assume
that
$p$ satisfies condition $(p1)'$ in $\Omega$.
Then there exists a classical
solution $U$ of problem \eq{unu} which is maximal solution if $p$ is
positive.

If $\Omega=
\RR^N\setminus \overline{B(0,R)}$ and $p$ satisfies the additional
condition $(p2)$,
with $\Phi(r)=0$ for $r\in [0,R]$, then the solution $U$ of \eq{unu}
is a large solution that blows-up at infinity.
\end{thm}

We refer to Ghergu, Niculescu and R\u adulescu \cite{gnr} for complete proofs
of Theorems \ref{theott} and \ref{theor}.

A useful observation is given in the following

\begin{rem}\label{remarka} Assume that $p\in C(\RR^{N})$ is a non-negative
and
non-trivial function which satisfies $(p2)$. Let $f$ be a
function satisfying assumption $(f1)$. Then condition
\neweq{doi} \int\limits_{1}^{\infty}\frac{dt}{f(t)}<\infty\endeq
is necessary for the existence of
entire large solutions to \eq{unu}.\end{rem}

Indeed,
let $u$ be an entire large solution of problem \eq{unu}.
Define
$$ \bar u(r)=\frac{1}{\omega_Nr^{N-1}}\int\limits_{|x|=r}
\left(\int\limits_{a_0}^{u(x)}\frac{dt}{f(t)}\right)\,dS=\frac{1}{\omega_N}
\int\limits_{|\xi|=1}\left(\int\limits_{a_0}^{u(r\xi)}\frac{dt}{f(t)}\right)\,dS,$$
where
$\omega_N$
denotes the surface area of the unit sphere in $\RR^{N}$ and $a_0$ is
chosen
such that $a_0\in (0,u_0)$, where $u_0=\inf_{\RR^N}u>0$.
By the divergence theorem we have
$$ \bar
u'(r)=\frac{1}{\omega_Nr^{N-1}}\int\limits_{B(0,r)}\Delta\left(
\int\limits_{a_0}^{u(x)}\frac{dt}{f(t)}\right)\, dx.
$$
Since $u$ is a positive classical solution it follows that
$$|\bar u'(r)|\leq Cr\rightarrow 0\qquad\mbox{as}\ r\rightarrow 0\,.$$
On the other hand
$$
\di\omega_N\left(R^{N-1}\bar u'(R)-r^{N-1}\bar u'(r)\right)=
 \int\limits_{r}^{R}\left(\;\int\limits_{|x|=z}\Delta\left(
\int\limits_{a_0}^{u(x)}\frac{dt}{f(t)}\right)\,dS\right)\,dz.
$$
Dividing by $R-r$ and taking $R\ri r$ we
find
$$ \begin{array}{lll}
 \omega_N(r^{N-1}\bar u'(r))'&\di=\int\limits_{|x|=r}\Delta\left(
\int\limits_{a_0}^{u(x)}\frac{dt}{f(t)}\right)\,dS=\int\limits_{|x|=r}
{\rm div}\,\left(\frac{1}{f(u(x))}\nabla u(x)\right)\,dS\\
& \di
=\int\limits_{|x|=r}\left[\left(\frac{1}{f}\right)'(u(x))\cdot
|\nabla
u(x)|^{2}+\frac{1}{f(u(x))}\Delta u(x)\right]\,dS\\
&\di\leq\int\limits_{|x|=r}
\frac{p(x)f(u(x))}{f(u(x))}\,dS\leq
\omega_Nr^{N-1}\Phi(r).\end{array}$$
The above inequality yields by integration
\neweq{treispe}
\bar u(r)\leq \bar u(0)+\int\limits_{0}^{r}\sigma^{1-N}\left(\int
\limits_{0}^{\sigma}\tau^{N-1}\Phi(\tau)\,d\tau\right)\,d\sigma \qquad
\forall
r\geq 0. \endeq
On the other hand, according to $(p2)$, for all $r>0$ we have
$$\begin{array}{lll}
\di\int\limits_{0}^{r}\sigma^{1-N}\left(\,\int\limits_{0}^{\sigma}\tau^{N-1}
\Phi(\tau)\,d\tau\right)\,d\sigma
&\di
=\frac{1}{2-N}r^{2-N}\int\limits_{0}^{r}\tau^{N-1}\Phi(\tau)
\,d\tau-\frac{1}{2-N}\int\limits_{0}^{r}\sigma\Phi(\sigma)\,d\sigma\\
&\di \leq \frac{1}{N-2}
\int\limits_{0}^{\infty}r\Phi(r)\,dr<\infty. \end{array} $$
So, by \eq{treispe},
$ \bar u(r)\leq \bar u(0)+K$, for all $r\geq 0.$
The last inequality implies that $\bar u$ is
bounded and assuming that \eq{doi} is not fulfilled it follows that $u$
cannot be a large solution. \qed

We point out that the hypothesis $(p2)$ on $p$ is essential in the statement of
Remark \ref{remarka}.
Indeed, let us consider $f(t)=t$, $p\equiv 1$, $\alpha\in(0,1)$,
$q(x)=2^{\alpha-2}\cdot|x|^\alpha$, $a=2-\alpha\in(1,2)$.
Then the corresponding problem has the entire large solution $u(x)=|x|^2+2N$, but
\eq{doi} is not fulfilled.

\section{Singular solutions with lack of the Keller-Osserman condition}
We have already seen that if $f$ is smooth and increasing on
$[0,\infty)$ such that $f(0)=0$ and $ f>0$ in $(0,\infty)$, then
the problem
$$ \left\{\begin{tabular}{ll}
$\Delta u=f (u)$ \quad & $\mbox{\rm in}\ \Omega,\ $\\
$u>0$ \quad & $\mbox{\rm in}\ \Omega,$\\
$u=+\infty$ \quad & $\mbox{\rm on}\ \partial\Omega$
\end{tabular} \right.$$
has a solution if and only if the Keller-Osserman condition
$\int_1^\infty  \left[F(t)\right]^{-1/2}dt <\infty$ is fulfilled, where
$F(t)=\int_0^t f(s)\,ds$. In particular, this implies that $f$ must have a superlinear growth.
In this section we are concerned with the problem
 \begin{equation}\label{unucpaa}
                         \left\{\begin{tabular}{ll}
                         $\Delta u+|\nabla u|=p(x)f(u)$ \qquad & $\mbox{ in }\
                         \Omega,$\\
                         $u\geq 0$ \qquad & $\mbox{ in }\ \Omega,$\\
                         \end{tabular} \right.
                         \end{equation}
                         where
                         $\Omega\subset\RR^N$ ($N\geq 3$) is either a smooth
                         bounded domain or the whole space. Our main assumptions on
$f$ is that it has a {\it sublinear} growth, so we cannot expect that Problem \eq{unucpaa}
admits a blow-up boundary solution. Our main purpose in this section is to establish
a necessary and sufficient
condition on the variable potential $p(x)$ for the existence of an entire large solution.

Throughout this section we assume that $p$ is a
                         non-negative function such that
                         $p\in C^{0,\alpha}(\overline\Omega )\,(0<\alpha<1)$ if
                         $\Omega$ is bounded, and
                         $p\in C^{0,\alpha}_{\rm loc}(\RR^N)$, otherwise. The
                         non-decreasing
                         non-linearity $f\in C^{0,\alpha}_{\rm
                         loc}[0,\infty)$ fulfills
$f(0)=0$ and $f>0$ on $(0,\infty )$.
             We also assume that $f$ is sublinear at infinity, in the sense that
$\Lambda := \sup _{s\geq 1}\frac{f(s)}{s}<
                         \infty.$

The main results in this section have been established by Ghergu
and R\u adulescu \cite{grcpaa}.

       If $\Omega$ is bounded we prove the following non-existence result.

 \begin{thm}\label{th1cpaa}  Suppose that $\Omega\subset\RR^N$ is a smooth
bounded domain. Then problem~(\ref{unucpaa}) has no positive large
solution in   $\Omega$.
\end{thm}

\proof Suppose by contradiction that problem~(\ref{unucpaa}) has a positive
                         large solution $u$ and define
                         $v(x)=\ln(1+u(x)),\;x\in\Omega.$ It follows that $v$
                         is positive
                         and $v(x)\rightarrow \infty$ as $\,{\rm dist}\,
                         (x,\partial\Omega)\to 0$.
                         We have
                         $$\displaystyle \Delta v=\frac{1}{1+u}\Delta
                         u-\frac{1}{(1+u)^2}|\nabla
                         u|^2\qquad\mbox{in}\;\Omega$$
                         and so
             $$\displaystyle\Delta v\leq
                         p(x)\frac{f(u)}{1+u}\leq\|p\,\|_{\infty}\frac{f(u)}{1+u}\leq
                         A\qquad
                         \mbox{ in }\Omega,$$
                         for some  constant $A>0$.
                         Therefore
                         $$\displaystyle \Delta(v(x)-A|x|^2)<0,\qquad\mbox{for
                         all}\;\:x\in\Omega.$$
                         Let $w(x)=v(x)-A|x|^2,\;x\in\Omega$. Then $\displaystyle\Delta
                         w<0$ in $\Omega.$ Moreover,
                         since $\Omega$ is bounded, it follows that
                         $w(x)\rightarrow \infty$ as ${\rm dist}(x,\partial\Omega)\rightarrow
                         0.$

                         Let $M>0$ be arbitrary. We claim that $w\geq M$ in
                         $\Omega$.
                         For all $\delta>0,$ we set
                         $$\displaystyle \Omega_{\delta}=\{x\in\Omega\,;\;{\rm
                         dist}(x,\partial\Omega)>\delta\}.$$
                         Since $w(x)\rightarrow \infty$ as ${\rm
                         dist}(x,\partial\Omega)\rightarrow  0,$ we
                         can choose $\delta>0$ such that
                         \begin{equation}\label{awm}\displaystyle w(x)\geq M \qquad\mbox{
                         for all }\;
                         x\in\Omega\setminus\Omega_{\delta}.
                         \end{equation}
                         On the other hand,
                         $$\begin{tabular}{ll}
                         $\displaystyle -\Delta (w(x)-M)>0$ \qquad & ${\rm in}\
                         \Omega_{\delta},$\\
                         $\displaystyle\qquad\qquad w(x)-M\geq 0$ \qquad & ${\rm on}\
                         \partial\Omega_{\delta}.$\\
                         \end{tabular}$$
                         By the maximum principle we get
                         $w(x)-M\geq 0$ in $\Omega_{\delta}$. So, by~(\ref{awm}),
                          $w\geq M$ in $\Omega.$
                         Since $M>0$ is arbitrary, it follows that
                         $w\geq n$ in $\Omega$,
                         for all $n\geq 1$. Obviously, this is a
                         contradiction and the proof is now complete.
                         \qed

                         Next, we consider the problem~(\ref{unucpaa})
                         when $\Omega=\RR^{N}$.
                         For all $r\geq 0$ we set
                         $$\phi(r)=\max \limits _{|x|=r}p(x),\qquad\psi(r)=\min
                         \limits
                         _{|x|=r}p(x),\qquad
                         \mbox{and}\qquad h(r)=\phi(r)-\psi(r).$$
                         We suppose that
                         \begin{equation}
                         \int\limits_{0}^{\infty}rh(r)\Psi(r)dr<\infty,\label{doicpaa}
                         \end{equation}
                         where
                         $$\Psi(r)=\exp\left(\Lambda_N\int\limits_{0}^{r}s\psi(s)ds\right),
                         \qquad  \Lambda_N=\frac{\Lambda}{N-2}.\\$$
                          Obviously, if $p$ is radial then $h\equiv 0$ and
                         ~(\ref{doicpaa}) occurs.
                         Assumption~(\ref{doicpaa}) shows that the variable potential
                         $p(x)$ has a slow variation.
                         An example of non-radial potential for which~(\ref{doicpaa})
                         holds is
                         $\displaystyle p(x)=\frac{1+|x_1|^2}{(1+|x_1|^2)(1+|x|^2)+1}.$
                         In this case
                         $\displaystyle\phi(r)=\frac{r^2+1}{(r^2+1)^2+1}$ and
                         $\displaystyle\psi(r)=\frac{1}{r^2+2}.$
                         If $\Lambda_N=1,$ by direct computation we get
                         $rh(r)\Psi(r)=O\left(r^{-2}\right)$ as $r\rightarrow \infty$ and
                         so~(\ref{doicpaa}) holds.

    \begin{thm}\label{th2cpaa} Assume $\Omega =\RR^N$ and
                         $p$ satisfies~(\ref{doicpaa}).
                         Then problem~(\ref{unucpaa}) has a positive entire large solution if
                         and only if
                         \begin{equation}
                         \int\limits_{1}^{\infty}e^{-t}t^{1-N}
                         \int\limits_{0}^{t}e^ss^{N-1}\psi(s)dsdt=\infty.\label{apatru}
                         \end{equation}
 \end{thm}

\proof Several times in the proof of Theorem~\ref{th2cpaa} we shall apply
                         the  following elementary inequality:
                         \begin{equation}
                         \int\limits_{0}^{r}e^{-t}t^{1-N}
                         \int\limits_{0}^{t}e^ss^{N-1}g(s)dsdt\leq
                         \frac{1}{N-2}
                         \int\limits_{0}^{r}tg(t)dt,
                         \qquad\forall\;r>0,\label{doispe}
                         \end{equation}
              for any continuous function
                          $g:[0,\infty)\rightarrow  [0,\infty)$.
                       The proof follows easily by integration by parts.

                         \noindent{\sc Necessary condition}.
                         Suppose that ~(\ref{doicpaa}) fails and the equation ~(\ref{unucpaa})
                         has a positive entire
                         large solution $u$. We claim that
                         \begin{equation}
                         \displaystyle \int\limits_{1}^{\infty}e^{-t}t^{1-N}
                         \int\limits_{0}^{t}e^ss^{N-1}\phi(s)dsdt<\infty.\label{phi}
                         \end{equation}

                         \noindent We first recall that $\phi=h+\psi.$ Thus
             $$\begin{tabular}{ll}
                         $\displaystyle \int\limits_{1}^{\infty}e^{-t}t^{1-N}
                         \int\limits_{0}^{t}e^ss^{N-1}\phi(s)dsdt$
                         &$\displaystyle
                         =\int\limits_{1}^{\infty}e^{-t}t^{1-N}
                         \int\limits_{0}^{t}e^ss^{N-1}\psi(s)dsdt$\\
                         &$\displaystyle \;\;\;+\int\limits_{1}^{\infty}e^{-t}t^{1-N}
                         \int\limits_{0}^{t}e^ss^{N-1}h(s)dsdt.$\\
                         \end{tabular}$$
                         By virtue of ~(\ref{doispe}) we find
                         $$\begin{tabular}{ll}
                         $\displaystyle \int\limits_{1}^{\infty}e^{-t}t^{1-N}
                         \int\limits_{0}^{t}e^ss^{N-1}\phi(s)dsdt$&$\displaystyle 
\leq\int\limits_{1}^{\infty}e^{-t}t^{1-N}           \int\limits_{0}^{t}e^ss^{N-1}\psi(s)dsdt+
\frac{1}{N-2}\int\limits_{0}^{\infty}th(t)dt$\\
                         &$\displaystyle \leq\int\limits_{1}^{\infty}e^{-t}t^{1-N}
\int\limits_{0}^{t}e^ss^{N-1}\psi(s)dsdt+
\frac{1}{N-2}\int\limits_{0}^{\infty}th(t)\Psi(t)dt.$
                         \end{tabular}$$
                         Since $\displaystyle \int\limits_{1}^{\infty}e^{-t}t^{1-N}
                         \int\limits_{0}^{t}e^ss^{N-1}\psi(s)dsdt<\infty,$ by
                         ~(\ref{doicpaa})
                         we deduce that ~(\ref{phi}) follows.

                         Now, let $\bar u$ be the spherical average of $u,$
                         i.e.,
                         $$\bar
                         u(r)=\frac{1}{\omega_Nr^{N-1}}\int\limits_{|x|=r}
                         u(x)d\sigma_x,\qquad r\geq 0,$$
                         where $\omega_N$ is the surface area of the unit
                         sphere in $\RR ^N$. Since $u$
                         is a positive entire large solution of ~(\ref{unu}) it
                         follows that $\bar u$
                         is positive and $\bar u(r)\rightarrow \infty$ as $r\rightarrow \infty.$
                         With the change of variable $\,x\rightarrow  ry,$ we have
                         $$\displaystyle \bar u(r)=\frac{1}{\omega_N}\int\limits_{|y|=1}
                         u(ry)\,d\sigma_y,\qquad r\geq 0$$
                         and
                         \begin{equation}\label{acinci}
                         \displaystyle \bar u'(r)=\frac{1}{\omega_N}\int\limits_{|y|=1}
                         \nabla u(ry)\cdot y\,d\sigma_y,\qquad r\geq 0.
                         \end{equation}
                         Hence
                         $$\displaystyle  \bar u'(r)=
                         \frac{1}{\omega_N}\int\limits_{|y|=1}\frac{\partial
                         u}{\partial r}
                         (ry)\,d\sigma_y=
                         \frac{1}{\omega_Nr^{N-1}}\int\limits_{|x|=r}\frac{\partial
                         u}{\partial r}
                         (x)\,d\sigma_x,$$
                         that is
                         \begin{equation}\label{asase}
                         \bar
                         u'(r)=\frac{1}{\omega_Nr^{N-1}}\int\limits_{B(0,R)}
                         \Delta u(x)\,dx,\qquad \mbox{for all}\;\;r\geq 0.
                         \end{equation}

                         Due to the gradient term $|\nabla u|$ in ~(\ref{unu}), we
                         cannot infer that
                         $\Delta u\geq 0$ in $\RR ^N$ and so we cannot expect that
                         $\bar u'\geq 0$
                         in $[0,\infty)$. We define the auxiliary function
                         \begin{equation}\label{aopt}
                         \displaystyle  U(r)=\max_{0\leq t\leq r}\bar u(t),\qquad r\geq 0.
                         \end{equation}
                         Then $U$ is positive and
                         non-decreasing.
                         Moreover, $U\geq\bar u$ and $U(r)\rightarrow \infty$ as
                         $r\rightarrow \infty$.

                         The assumptions $(f1)$ and $(f2)$ yield
                         $\displaystyle  f(t)\leq \Lambda(1+t),$ for all $t\geq 0.$ So, by
                          ~(\ref{acinci})
                         and ~(\ref{asase}),
                         $$\begin{tabular}{ll}
                         $\displaystyle  \bar u''+\frac{N-1}{r}\,\bar u'+\bar
                         u'\!\!\!$&$\displaystyle \leq\;\frac{1}{\omega_Nr^{N-1}}
                         \int\limits_{|x|=r}\left[\Delta u(x)+|\nabla
                         u|(x)\right]d\sigma_x=\displaystyle  \frac{1}{\omega_Nr^{N-1}}
                         \int\limits_{|x|=r}p(r)f(u(x))d\sigma_x$\\
                         &$\displaystyle \leq \;\Lambda \phi(r)\frac{1}{\omega_Nr^{N-1}}
                         \int\limits_{|x|=r}\left(1+u(x)\right)d\sigma_x
                         =\;\Lambda \phi(r)\left(1+\bar u(r)\right) \leq\;\Lambda \phi(r)\left(1+U(r)\right),$\\
                         \end{tabular}$$
                         for all $\,r\geq 0$. It follows that
                         $$\displaystyle  \left(r^{N-1}e^r\bar u'\right)'\leq\;\Lambda
                         e^rr^{N-1}\phi(r)\left(1+U(r)\right),
                         \qquad\mbox{for all}\;\,r\geq 0.$$
                         So, for all $\,r\geq r_0>0\,$,
                         $$\displaystyle \bar u(r)\leq\bar u(r_0)+\Lambda
                         \int_{r_0}^re^{-t}t^{1-N}\int_0^te^ss^{N-1}\phi(s)(1+U(s))dsdt.$$
                         The monotonicity of $U$ implies
                         \begin{equation}\label{azece}
                         \displaystyle \bar u(r)\leq \bar
                         u(r_0)+\Lambda(1+U(r))\int_{r_0}^r
                         e^{-t}t^{1-N}\int_0^te^ss^{N-1}\phi(s)dsdt,
                         \end{equation}
                         for all $r\geq r_0\geq 0.$
                          By ~(\ref{phi}) we can choose $r_0\geq 1$ such that
                         \begin{equation}\label{aunspe}
                         \displaystyle \int_{r_0}^{\infty}e^{-t}t^{1-N}\int_0^te^ss^{N-1}\phi(s)dsdt
                         <\frac{1}{2\Lambda}.
                         \end{equation}
                         Thus ~(\ref{azece}) and ~(\ref{aunspe}) yield
                         \begin{equation}\label{ff}
                         \displaystyle  \bar u(r)\leq \bar u(r_0)+\frac{1}{2}(1+U(r)),
                         \qquad\mbox{for all}\;\;r\geq r_0.
                         \end{equation}
                         By the definition of $U$ and
                         $\displaystyle \lim_{r\rightarrow \infty}\bar u(r)=\infty,$
                         we find $r_1\geq r_0$ such that
                         \begin{equation}\label{adoispe}
                         \displaystyle  U(r)=\max_{r_0\leq t\leq r}\bar
                         u(r),\qquad\mbox{for all}\;\;r\geq r_1.
                         \end{equation}
                         Considering now ~(\ref{ff}) and ~(\ref{adoispe}) we obtain
                         $$\displaystyle  U(r)\leq \bar
                         u(r_0)+\frac{1}{2}(1+U(r)),\qquad\mbox{for all}\;\;r\geq
                         r_1. $$
                         Hence
                         $$\displaystyle  U(r)\leq 2\bar u(r_0)+1,\qquad\mbox{for
                         all}\;\;r\geq r_1.$$
                         This means that $U$ is bounded, so $u$ is also
                         bounded, a contradiction. It follows that ~(\ref{unu}) has
                         no
                         positive entire large solutions.

             \medskip
             \noindent{\sc Sufficient condition.}
                         We need the following auxiliary comparison result.

  \begin{lem}\label{l1cpaa}       Assume that ~(\ref{doicpaa}) and ~(\ref{apatru})
                         hold.
                         Then the equations
                         \begin{equation}
                         \Delta v+|\nabla v|=\phi(|x|)f(v) \qquad
                         \Delta w+|\nabla w|=\psi(|x|)f(w)\\\label{cinci}
                         \end{equation}
                         have positive entire large solution such that
                         \begin{equation}
                         v\leq w  \qquad in\;\; \RR ^N.\\ \label{sase}
                         \end{equation}
                         \end{lem}

                        \proof Radial solutions of ~(\ref{cinci}) satisfy
                         $$v''+\frac{N-1}{r}v'+|v'|=\phi(r)f(v)$$ and
                         $$w''+\frac{N-1}{r}w'+|w'|=\psi(r)f(w).$$
                         Assuming that $v'$ and $w'$ are non-negative, we
                         deduce
                         $$\displaystyle
                         \left(e^rr^{N-1}v'\right)'=e^rr^{N-1}\phi(r)f(v)$$
                         and
                         $$\displaystyle \left(e^rr^{N-1}w'\right)'=e^rr^{N-1}\psi(r)f(w).$$
                         Thus any positive solutions $v$ and $w$ of the
                         integral equations
                         \begin{equation}
                         v(r)=1+\int\limits_{0}^{r}e^{-t}t^{1-N}
                         \int\limits_{0}^{t}e^ss^{N-1}\phi(s)f(v(s))dsdt,\qquad
                         r\geq 0,\\ \label{saptecpaa}
                         \end{equation}
                         \begin{equation}
                         w(r)=b+\int\limits_{0}^{r}e^{-t}t^{1-N}
                         \int\limits_{0}^{t}e^ss^{N-1}\psi(s)f(w(s))dsdt,\qquad
                         r\geq 0,\\ \label{optcpaa}
                         \end{equation}
                         provide a solution of ~(\ref{cinci}), for any $b>0$.
                         Since $w\geq b$, it follows that $f(w)\geq f(b)>0$
                         which yields
                         $$w(r)\geq b+f(b)\int\limits_{0}^{r}e^{-t}t^{1-N}
                         \int\limits_{0}^{t}e^ss^{N-1}\psi(s)dsdt,\qquad r\geq
                         0.\\$$
                         By ~(\ref{apatru}), the right hand side of this
                         inequality goes to $+\infty$ as
                         $r\rightarrow \infty$. Thus $w(r)\rightarrow \infty$ as
                         $r\rightarrow \infty.$
                         With a similar argument we find $v(r)\rightarrow \infty$ as
                         $r\rightarrow \infty.$

                           Let $b>1$ be fixed.
                         We first show that ~(\ref{optcpaa}) has a positive solution.
                         Similarly,
                         ~(\ref{saptecpaa}) has a positive solution.
                         \medskip

                         Let $\{w_k\}$ be the sequence defined by $w_1=b$ and
                         \begin{equation}
                         w_{k+1}(r)=b+\int\limits_{0}^{r}e^{-t}t^{1-N}
                         \int\limits_{0}^{t}e^ss^{N-1}\psi(s)f(w_k(s))dsdt,\qquad
                         k\geq 1.\\ \label{nouacpaa}
                         \end{equation}

                         We remark that $\{w_k\}$ is a non-decreasing
                         sequence.
                         To get the convergence of $\{w_k\}$ we will show
                         that $\{w_k\}$ is bounded from above on bounded
                         subsets.
                         To this aim, we fix $R>0$ and we prove that
              \begin{equation}
                         w_k(r)\leq be^{Mr},\qquad\mbox{for any } 0\leq r\leq R,\;\mbox{ and for all
                         } k\geq 1,\label{zececpaa}
                         \end{equation}
                         where $\displaystyle  M\equiv\Lambda_N
                         \max_{t\in[0,R]}\,t\psi(t).$

                         We achieve ~(\ref{zececpaa}) by induction. We first notice
                         that ~(\ref{zececpaa})
                         is true for $k=1$. Furthermore, the assumption $(f2)$
                         and the fact that $w_k\geq 1$
                         lead us to $f(w_k)\leq\Lambda w_k$, for all $k\geq
                         1$.
                         So, by ~(\ref{nouacpaa}),
                         $$w_{k+1}(r)\leq
                         b+\Lambda\int\limits_{0}^{r}e^{-t}t^{1-N}
                         \int\limits_{0}^{t}e^ss^{N-1}\psi(s)w_k(s)dsdt, \qquad
                         r\geq 0.$$\\
                         Using now ~(\ref{doispe}) (for
                         $g(t)=\psi(t)w_k(t)$) we deduce
                         $$w_{k+1}(r)\leq b+\Lambda
                         _N\int\limits_{0}^{r}t\psi(t)w_k(t)dt,
                         \qquad\forall\;r\in[0,R].$$
             The induction hypothesis yields
                         $$w_{k+1}(r)\leq
                         b+bM\int\limits_{0}^{r}e^{Mt}dt=be^{Mr},\qquad\forall\;r\in[0,R].$$
             Hence, by induction, the sequence
                         $\{w_k\}$ is
                         bounded in $[0,R]$, for any $R>0$.
                         It follows that $\displaystyle  w(r)=\lim_{k\rightarrow \infty}w_k(r)$ is
                         a positive solution of
                          ~(\ref{optcpaa}). In a similar way we conclude that
                         ~(\ref{saptecpaa}) has a
                         positive solution on $[0,\infty)$.

                         The next step is to show that the constant $b$ may be
                         chosen sufficiently large
                         so that ~(\ref{sase}) holds. More exactly, if
                         \begin{equation}
                         b>1+K\Lambda_N\int\limits_{0}^{\infty}sh(s)\Psi(s)ds,\label{treispecpaa}
                         \end{equation}
                         where
                         $K=\exp\left(\Lambda_N\int\limits_{0}^{\infty}th(t)dt\right),$
                         then ~(\ref{sase}) occurs.

                         We first prove that the solution $v$ of ~(\ref{saptecpaa})
                         satisfies
                         \begin{equation}
                         v(r)\leq K\Psi(r),\qquad \forall\;r\geq 0.\label{paispe}
                         \end{equation}
                         Since $v\geq 1$, from $(f2)$ we have $f(v)\leq \Lambda
                         v$.
                         We use this fact in ~(\ref{saptecpaa}) and then we apply the
                         estimate
                         ~(\ref{doispe}) for $g=\phi.$ It follows that
                         \begin{equation}
                         v(r)\leq
                         1+\Lambda_N\int\limits_{0}^{r}s\phi(s)v(s)ds,\qquad\forall\;r\geq
                         0.\label{cincispe}
                         \end{equation}
                         By Gronwall's inequality we obtain
                         $$v(r)\leq
                         \exp\left(\Lambda_N\int\limits_{0}^{r}s\phi(s)ds\right),
                         \qquad\forall\;r\geq 0,$$
                         and, by ~(\ref{cincispe}),
                         $$v(r)\leq 1+\Lambda_N\int\limits_{0}^{r}s\phi(s)\exp
                         \left(\Lambda_N\int\limits_{0}^{s}t\phi(t)dt\right)ds,
                         \qquad\forall\;r\geq 0.$$
                         Hence
                         $$v(r)\leq 1+\int\limits_{0}^{r}\left(\exp
                         \left(\Lambda_N\int\limits_{0}^{s}t\phi(t)dt\right)\right)'ds,
                         \qquad\forall\;r\geq 0,$$
                         that is
                         \begin{equation}
                         \displaystyle
                         v(r)\leq\exp\left(\Lambda_N\int\limits_{0}^{r}t\phi(t)dt\right),
                         \qquad\forall\;r\geq 0.\label{58}
                         \end{equation}
                         Inserting $\phi=h+\psi$ in ~(\ref{58}) we have
                         $$\displaystyle  v(r)\leq
                         e^{\Lambda_N\int\limits_{0}^{r}th(t)dt}\Psi(r)\leq
                         K\Psi(r),\qquad\forall\;r\geq 0,$$
                         so ~(\ref{paispe}) follows.

                         Since $b>1$ it follows that $v(0)<w(0).$ Then there
                         exists $R>0$
                         such that $v(r)<w(r),$ for any $0\leq r\leq R$. Set
                         $$R_{\infty}=\sup\{\ R>0 \,|\,v(r)<w(r),\;\,
                         \forall\,r\in[\,0,R] \,\}.$$
                         In order to conclude our proof, it remains to
                         show that $R_{\infty}=\infty$.
                         Suppose the contrary. Since
                         $v\leq w$ on $[\,0,R_{\infty}]$ and $\phi=h+\psi,$
                         from ~(\ref{saptecpaa}) we deduce
                         $$  v(R_{\infty})=\displaystyle
                         1+\int\limits_{0}^{\,R_{\infty}}e^{-t}t^{1-N}
                         \int\limits_{0}^{t}e^ss^{N-1}h(s)f(v(s))dsdt
+\int\limits_{0}^{R_{\infty}}e^{-t}t^{1-N}
                         \int\limits_{0}^{t}e^ss^{N-1}\psi(s)f(v(s))dsdt.$$
                         So, by ~(\ref{doispe}),
                         $$\displaystyle  v(R_{\infty})
                         \leq
                         1+\frac{1}{N-2}\int\limits_{0}^{R_{\infty}}th(t)f(v(t))dt
                         +\int\limits_{0}^{R_{\infty}}e^{-t}t^{1-N}
                         \int\limits_{0}^{t}e^ss^{N-1}\psi(s)f(w(s))dsdt.$$
                         Taking into account that $v\geq 1$ and  the
                         assumption $(f2),$
                         it follows that
                         $$ v(R_{\infty})
                         \leq
                         1+K\Lambda_N\int\limits_{0}^{R_{\infty}}th(t)\Psi(t)dt
                         +\int\limits_{0}^{R_{\infty}}e^{-t}t^{1-N}
                         \int\limits_{0}^{t}e^ss^{N-1}\psi(s)f(w(s))dsdt.$$
                         Now, using ~(\ref{treispecpaa}) we obtain
                         $$\displaystyle
                         v(R_{\infty})<b+\int\limits_{0}^{R_{\infty}}e^{-t}t^{1-N}
                         \int\limits_{0}^{t}e^ss^{N-1}\psi(s)f(w(s))dsdt=w(R_{\infty}).$$
                         Hence $v(R_{\infty})<w(R_{\infty}).$ Therefore,
                         there exists $R>R_{\infty}$ such that $v<w$ on
                         $[\,0,R]$, which
                         contradicts the maximality of $R_{\infty}$.
                         This contradiction shows that inequality ~(\ref{sase})
                         holds and the
                         proof of Lemma~\ref{l1} is now complete.
                         \qed

\smallskip
                         {\it Proof of Theorem~\ref{th2cpaa} completed.}
                         Suppose that ~(\ref{apatru}) holds. For all $k\geq 1$ we
                         consider the problem
                         \begin{equation}
                         \left\{\begin{tabular}{ll}
                         $\displaystyle \Delta u_k+|\nabla u_k|=p(x)f(u_k)$&$\displaystyle  \mbox{ in
                         }\,\, B(0,k),$\\
                         $\displaystyle  u_k(x)=w(k)$&$\displaystyle \mbox{ on }\,\,\partial
                         B(0,k).$\\ \label{anouaspe}
                         \end{tabular}\right.
                         \end{equation}
                         Then $v$ and $w$ defined by
                         (\ref{saptecpaa})
                         and (\ref{optcpaa}) are positive sub and super-solutions of
                        (\ref{anouaspe}).
                         So this problem has at least a positive solution
                         $u_k$ and
                         $$\displaystyle  v(|x|)\leq u_k(x)\leq w(|x|) \qquad\mbox{ in
                         }\;B(0,k), \mbox{ for all }\;k\geq 1.$$
                         By Theorem~14.3 in Gilbarg and Trudinger \cite{gt}, the
                         sequence
                         $\{\nabla u_k\}$ is bounded on every compact set in
                         $\RR ^N$.
                         Hence the sequence $\{u_k\}$ is bounded and
                         equicontinuous on compact subsets of
                         $\RR ^N.$
                         So, by the Arzela-Ascoli Theorem, the sequence
                         $\{u_k\}$
                         has a uniform convergent subsequence, $\{u_k^1\}$ on
                         the ball
                         $B(0,1).$
                         Let $u^1=\lim_{k\rightarrow\infty}u_k^1$. Then
                         $\{f(u_k^1)\}$ converges
                         uniformly to $f(u^1)$ on  $B(0,1)$ and, by
                         (\ref{anouaspe}), the sequence
                          $\{\Delta u_k^1+|\nabla u_k^1|\}$ converges uniformly
                         to
                         $pf(u^1).$
                         Since the sum of the Laplace and Gradient operators is a closed
                         operator, we deduce
                         that $u^1$ satisfies (\ref{unu}) on $B(0,1).$

                         Now, the sequence $\{u_k^1\}$ is bounded and
                         equicontinuous on the ball $B(0,2)$,  so it  has a
                         convergent subsequence $\{u_k^2\}.$
                         Let $u^2=\lim \limits _{k\rightarrow\infty}u_k^2\,$ on
                         $\,B(0,2)\,$
                         and $u^2$ satisfies (\ref{unu}) on $B(0,2).$ Proceeding
                         in the same way, we
                         construct a sequence $\{u^n\}$ so that $u^n$ satisfies
                         (\ref{unu}) on
                         $B(0,n)$
                         and $u^{n+1}=u^n$ on $B(0,n)$ for all $n$. Moreover,
                         the
                         sequence $\{u^n\}$
                         converges  in $L^\infty_{\rm loc}(\RR ^N)$ to the
                         function $u$
                         defined by
                          $$u(x)=u^m(x), \qquad\mbox{for }\,\,x\in B(0,m).$$
                          Since $v\leq u^n\leq w$ on $B(0,n)$ it follows that
                         $v\leq u\leq
                         w$ on $\RR ^N,$ and $\,u\,$ satisfies (\ref{unu}). From
                         $v\leq u$ we
                         deduce that $u$ is a positive entire large solution of
                         (\ref{unu}).
                         This completes the proof.
                         \qed

\section{Blow-up boundary solutions of the logistic equation}
Consider the semilinear elliptic equation
\neweq{sep} \Delta u+au=b(x)f(u)\qquad\mbox{in}\ \Omega , \endeq
where $\Omega$ is a smooth bounded domain in  $\RR^N$, $N\geq 3$.
Let $a$ be a real parameter and $b\in C^{0,\mu}(\overline\Omega)$,
$0<\mu<1$, such that $b\geq 0$ and $b\not\equiv 0$ in $\Omega$.
Set
$$ \Omega_0={\rm int}\,\{x\in \Omega:\ b(x)=0\} $$
and suppose, throughout, that
$\ov{\Omega}_0\subset \Omega$ and $b>0$ on
$\Omega\setminus\ov{\Omega}_0$.
Assume that $f\in C^1[0,\infty)$ satisfies

\medskip
\noindent $(A_1)\quad  f\geq 0$
and $f(u)/u$ is increasing on $(0,\infty)$.

\medskip
Following Alama and Tarantello \cite{al_t1}, define
by $H_{\infty}$
the Dirichlet Laplacian on  $\Omega_0$
as the unique self-adjoint operator associated to the quadratic
form $\psi(u)=\int_{\Omega}|\nabla u|^2\,dx$
with form domain
$$ H^1_D(\Omega_0)=\{u\in H_0^1(\Omega):\ u(x)=0\quad
\mbox{for a.e.}\ x\in \Omega\setminus\Omega_0\}. $$
If $\partial\Omega_0$ satisfies the exterior
cone condition then, according to \cite{al_t1}, $H^1_D(\Omega_0)$
coincides with $H^{1}_{0}(\Omega_0)$ and $H_{\infty}$
is the classical Laplace operator with Dirichlet
condition on $\partial\Omega_0$.

Let $\lambda_{\infty,1}$ be the first Dirichlet eigenvalue of
$H_{\infty}$
in $\Omega_0$.
We understand $\lambda_{\infty,1}=\infty$
if $\Omega_0=\emptyset$.

Set
$ \mu_{0}:=\lim_{u\searrow 0}\frac{f(u)}{u}$,
$\mu_{\infty}:=\lim_{u\to \infty}\frac{f(u)}{u}$, and
denote by $\lambda_1(\mu_0)$ (resp., $\lambda_1(\mu_{\infty})$)
the first eigenvalue of the operator $H_{\mu_0}=
-\Delta+\mu_0 b$ (resp., $H_{\mu_{\infty}}=-\Delta+\mu_{\infty}b$)
in $H_{0}^{1}(\Omega)$. Recall that
$\lambda_1(+\infty)=\lambda_{\infty,1}$.

Alama and Tarantello \cite{al_t1} proved that problem
\eq{sep} subject to the Dirichlet boundary condition
\neweq{bc}
u=0 \quad \mbox{on}\ \partial\Omega
\endeq
has a positive solution $u_a$ if and only if
$a\in (\lambda_{1}(\mu_0),\lambda_{1}(\mu_\infty))$. Moreover,
$u_a$ is the unique positive solution for \eq{sep}+\eq{bc}
(see \cite[Theorem A (bis)]{al_t1}). We shall refer to the combination
of
\eq{sep}+\eq{bc} as problem $(E_a)$.

Our first aim in this section is to give a corresponding necessary
and sufficient condition, but for the existence of {\it large} (or
{\it explosive}) solutions of \eq{sep}. An elementary argument
based on the maximum principle shows that if such a solution
exists, then it is {\it positive} even if $f$ satisfies a weaker
condition than $(A_1)$, namely
\medskip

\noindent $(A_1)' \quad  f(0)=0$,\ $f'\geq 0$ and
$f>0$ on $(0,\infty)$.
\medskip

We recall that Keller \cite{ke} and Osserman \cite{os} supplied a necessary and
sufficient condition on $f$ for the existence of large solutions
to (1) when $a\equiv 0$, $b\equiv 1$ and $f$ is assumed to fulfill
$(A_1)'$.
More precisely, $f$ must satisfy the Keller-Osserman condition
(see \cite{ke,os}),

\smallskip

\noindent $(A_2)\quad \di \int_{1}^{\infty}\frac{dt}{\sqrt{F(t)}}<
\infty\,, \ \ \mbox{where}\
F(t)=\di\int_{0}^{t} f(s)\,ds.$
\smallskip

Typical examples of non-linearities satisfying
$(A_{1})$ and $(A_{2})$ are:
$$ {\rm (i)}\ f(u)=e^u-1;\quad {\rm (ii)}\ f(u)=u^p,\ p>1;\quad
{\rm (iii)}\ f(u)=u[\ln\,(u+1)]^p,\ p>2.$$

Our first result
gives  the maximal interval for the
parameter $a$ that ensures the existence of large solutions to
problem \eq{sep}. More precisely, we prove

\begin{thm}\label{teo0} Assume that $f$ satisfies conditions $(A_1)$
and $(A_2)$. Then problem \eq{sep} has a
large solution if and only if $a\in (-\infty,\lambda_{\infty,1})$.
\end{thm}

We point out that our framework in the above result includes the
case when $b$ vanishes at some points on $\partial\Omega$, or even
if $b\equiv 0$ on $\partial\Omega$. This later case includes the ``competition" $0\cdot\infty$
on $\partial\Omega$.
We also point out that, under our hypotheses, $\mu_{\infty}:=\lim_{u\to \infty}f(u)/u=
\lim_{u\to \infty}f'(u)=\infty$. Indeed, by l'Hospital's rule, $\lim_{u\to
\infty}F(u)/u^2=\mu_{\infty}/2$. But,
by $(A_2)$, we deduce that $\mu_{\infty}=\infty$. Then, by $(A_1)$ we
find that
$f'(u)\geq f(u)/u$ for any $u>0$, which shows that $\lim_{u\to
\infty}f'(u)=\infty$.

Before giving the proof of Theorem \ref{teo0} we claim that
assuming $(A_1)$, then problem \eq{sep} can have large
solutions
only if $f$ satisfies the Keller-Osserman condition $(A_2)$.
Indeed, suppose that problem \eq{sep} has a large solution
$u_\infty$. Set $\tilde f(u)=|a|u+\|b\|_{\infty}f(u)$ for $u\geq 0$.
Notice that
$\tilde f\in C^1[0,\infty)$ satisfies $(A_1)'$. For any $n\geq 1$,
consider the problem
$$ \left\{\begin{array}{lll}
& \di \Delta u=\tilde f(u) \quad & \mbox{in}\ \Omega\,, \\
& \di u=n \quad & \mbox{on}\ \partial\Omega\,,\\
& \di u\geq 0 \quad & \mbox{in}\ \Omega\,.
\end{array} \right. $$
A standard argument based on the maximum principle shows that this problem has a unique
solution, say $u_n$,
which, moreover, is positive in $\overline\Omega$. Applying again the maximum principle we deduce
that
$ 0<u_n\leq u_{n+1}\leq u_{\infty}$, in $\Omega$,
for all $n\geq 1.$
Thus, for every $x\in \Omega$, we can define $\bar u(x)=\lim_{n\to
\infty}u_n(x)$.
Moreover, since $(u_n)$ is uniformly bounded on every compact subset of
$\Omega$,
standard elliptic regularity arguments show that $\bar u$ is a positive
large solution
of the problem $\Delta u=\tilde f(u)$. It follows that $\tilde f$
satisfies
the Keller-Osserman condition $(A_2)$. Then, by $(A_1)$,
$\mu_{\infty}:=\lim_{u\to \infty}f(u)/u>0$ which yields
$\lim_{u\to \infty}\tilde
f(u)/f(u)=|a|/\mu_{\infty}+\|b\|_{\infty}<\infty$. Consequently,
our claim follows.

\medskip
{\it Proof of Theorem \ref{teo0}}.
{\sc A. Necessary condition.} Let $u_\infty$ be a large
solution
of problem \eq{sep}. Then, by the maximum principle,
$u_\infty$ is positive. Suppose $\lambda_{\infty,1}$ is finite. Arguing
by
contradiction, let us assume $a\geq \lambda_{\infty,1}$.
Set $\lambda\in (\lambda_{1}(\mu_0),\lambda_{\infty,1})$
and denote by $u_\lambda$ the unique positive solution
of problem $(E_a)$ with $a=\lambda$. We have
$$ \left\{\begin{array}{lll}
& \di \Delta (Mu_{\infty})+
\lambda_{\infty,1}(M u_{\infty})\leq b(x)f(Mu_{\infty}) \quad &
\mbox{in}\ \Omega\,, \\
& \di Mu_{\infty}=\infty \quad & \mbox{on}\ \partial\Omega\,,\\
& \di Mu_{\infty}\geq u_{\lambda} \quad & \mbox{in}\ \Omega\,,
\end{array} \right. $$
where $M:=\max\left\{\max_{\ov{\Omega}}\,u_\lambda/
\min_{\Omega}\,u_{\infty};\,1\right\}$.
By the sub-super solution method
we conclude that problem $(E_a)$ with $a=\lambda_{\infty,1}$
has at least a positive solution (between $u_{\lambda}$ and
$Mu_{\infty}$). But this is a contradiction. So, necessarily, $a\in
(-\infty,\lambda_{\infty,1})$.
\medskip

\noindent B. {\sc Sufficient condition.}
This will be proved with the aid of several results.

\begin{lem}\label{lb7}
Let $\omega$ be a smooth
bounded domain in $\RR^N$. Assume $p$, $q$, $r$ are
$C^{0,\mu}$-functions on $\ov{\omega}$ such that $r\geq 0$ and
$p>0$ in $\ov{\omega}$. Then for any non-negative function
$0\not\equiv \Phi\in C^{0,\mu}(\partial\omega)$ the boundary value
problem
\neweq{cr3} \left\{ \begin{array}{lll}
& \di \Delta u+q(x)u=p(x)f(u)-r(x) \quad & \mbox{in}\ \omega, \\
& \di u>0 \quad & \mbox{in}\ \omega,\\
& \di u=\Phi \quad & \mbox{on}\ \partial\omega,
\end{array} \right. \endeq
has a unique solution.\end{lem}

We refer to  C\^{\i}rstea and R\u adulescu \cite[Lemma 3.1]{crccm}
for the proof of the above result.

Under the assumptions of Lemma~\ref{lb7} we obtain the following result
which generalizes \cite[Lemma~1.3]{mv1}.

\begin{cor}\label{corM} There exists a positive
large solution of the problem
\neweq{cr2} \Delta u+q(x) u=p(x)f(u)-r(x) \qquad \mbox{in}\
\omega.\endeq \end{cor}
\proof
Set $\Phi=n$ and let $u_n$ be the unique solution
of \eq{cr3}. By the maximum principle, $u_n\leq u_{n+1}\leq
\ov{u}$ in $\omega$, where $\ov{u}$ denotes a large solution
of
$$\Delta u+\|q\|_{\infty}u=p_0 f(u)-\bar r \quad  \mbox{in}\
\omega.$$
Thus $\lim_{n\to \infty}u_n(x)=u_\infty(x)$
exists and is a positive large solution of \eq{cr2}.
Furthermore, every positive large solution of \eq{cr2} dominates
$u_\infty$, i.e., the solution $u_\infty$ is the {\it minimal large
solution}.
This follows from the definition of $u_\infty$ and the maximum principle.
\qed

\begin{lem}\label{lb8} If  $0\not\equiv \Phi\in C^{0,\mu}(\partial
\Omega)$
is a non-negative function and $b>0$ on $\partial\Omega$,
then the boundary value problem
\neweq{cr8} \left\{\begin{array}{lll}
& \di \Delta u+au=b(x)f(u) \quad & \mbox{in}\ \Omega, \\
& \di u>0 \quad & \mbox{in}\ \Omega,\\
& \di u=\Phi \quad & \mbox{on}\ \partial\Omega,
\end{array} \right. \endeq
has a solution if and only if $a\in (-\infty,\lambda_{\infty,1})$.
Moreover, in this case, the solution is unique. \end{lem}
\proof The first part follows exactly in the same way as the proof
of Theorem~\ref{teo0} (necessary condition).

For the sufficient condition, fix $a<\lambda_{\infty,1}$
and let $\lambda_{\infty,1}>\lambda_*>\max\,\{a,\lambda_1(\mu_0)\}$.
Let $u_*$ be the unique positive solution of $(E_a)$
with $a=\lambda_*$.

Let $\Omega\,_i$ ($i=1,2$) be subdomains of $\Omega$
such that $\Omega_0 \subset\subset \Omega_1
\subset\subset \Omega_2 \subset\subset \Omega$ and
$\Omega\setminus\ov{\Omega}_1$ is smooth.

\noindent We define $u_+\in C^2(\Omega)$ as a positive function in
$\Omega$
such that $u_+\equiv u_\infty $ on $\Omega\setminus\Omega_2$ and
$u_+\equiv u_*$ on $\Omega_1$.
Here $u_{\infty}$ denotes a positive large
solution of \eq{cr2} for $p(x)=b(x)$, $r(x)=0$, $q(x)=a$
and $\omega=\Omega\setminus \ov{\Omega}_1$. So, since
$b_0:=\inf_{\Omega_2\setminus\Omega_1}b$ is positive, it
is easy to check that if $C>0$ is large enough then
$\ov{v}_{\Phi}=Cu_+$ satisfies
$$ \left\{\begin{array}{lll}
& \di \Delta \ov{v}_\Phi+ a\ov{v}_\Phi\leq b(x)f(\ov{v}_\Phi) &\quad
\mbox{in}\ \Omega\,,\\
& \di \ov{v}_\Phi=\infty & \quad \mbox{on}\ \partial\Omega\,.\\
& \di \ov{v}_\Phi\geq \max_{\partial\Omega}\Phi & \quad \mbox{in}\
\Omega\,.
\end{array} \right. $$
Let $\un{v}_{\Phi}$ be the unique
classical solution of the problem
$$ \left\{\begin{array}{lll}
& \di \Delta \un{v}_{\Phi}=|a| \un{v}_{\Phi}+\|b\|_{\infty}
f(\un{v}_{\Phi}) \quad & \mbox{in}\ \Omega, \\
& \di \un{v}_{\Phi}>0 \quad & \mbox{in}\ \Omega,\\
& \di \un{v}_{\Phi}=\Phi \quad & \mbox{on}\ \partial\Omega\,.
\end{array} \right. $$
It is clear that $\un{v}_{\Phi}$ is a positive sub-solution of
\eq{cr8} and $\un{v}_{\Phi}\leq \max_{\partial\Omega}\Phi\leq
\ov{v}_\Phi$ in $\Omega$. Therefore, by the sub-super solution
method, problem \eq{cr8} has at least a solution $v_{\Phi}$
between $\un{v}_{\Phi}$ and $\ov{v}_{\Phi}$. Next, the uniqueness
of solution to \eq{cr8} can be obtained by using essentially the
same technique as in \cite[Theorem 1]{bros} or \cite[Appendix
II]{bk}. \qed

\medskip
{\it Proof of Theorem \ref{teo0} completed}. Fix $a\in (-\infty
,\lambda_{\infty ,1})$. Two cases may occur:
\medskip

{\sc Case 1:}\quad $b>0$ on $\partial\Omega$. Denote by $v_n$ the
unique
solution of \eq{cr8} with $\Phi\equiv n$.
For $\Phi\equiv 1$,
set $v:=\un{v}_{\Phi}$ and $V:=\ov{v}_{\Phi}$,
where $\un{v}_{\Phi}$ and $\ov{v}_{\Phi}$
are defined in the proof of Lemma \ref{lb8}.
The sub and super-solutions method combined with
the uniqueness of solution
of \eq{cr8} shows that $v\leq v_n\leq v_{n+1}\leq V$ in $\Omega$.
Hence $v_{\infty}(x):=
\lim_{n\to \infty}v_n(x)$ exists and is a positive large solution
of \eq{sep}.
\medskip

{\sc Case 2:}\quad $b\geq 0$
on $\partial\Omega$.
Let $z_n$ ($n\geq 1$) be the unique solution of \eq{cr3} for
$p\equiv b+1/n$, $r\equiv 0$, $q\equiv a$, $\Phi\equiv n$
and $\omega=\Omega$. By the maximum principle, $(z_n)$ is non-decreasing. Moreover,
$(z_n)$ is uniformly bounded on every compact subdomain of $\Omega$.
Indeed, if $K\subset \Omega$ is an arbitrary compact set, then
$d:={\rm dist}\, (K,\partial\Omega)>0$. Choose $\delta \in (0,d)$ small
enough so that $\ov{\Omega}_0\subset C_\delta$, where
$C_\delta=\{x\in \Omega:\ {\rm dist}\, (x,\partial\Omega)>\delta\}$.
Since $b>0$ on $\partial C_\delta$, Case~1 allows us to define
$z_+$ as a positive large solution of \eq{sep} for $\Omega=C_\delta$.
Using
A standard argument based on the maximum principle implies that
 $z_n\leq z_+$ in $C_\delta$, for all $n\geq 1$. So, $(z_n)$
is uniformly bounded on $K$. By the monotonicity of $(z_n)$, we
conclude that
$z_n\to \underline{z}$ in $L_{\rm loc}^{\infty}(\Omega)$.
Finally, standard elliptic regularity arguments lead to $z_n\to
\underline{z}$ in $C^{2,\mu}(\Omega)$.
This completes the proof of Theorem \ref{teo0}. \qed

\medskip
Denote by ${\cal D}$ and ${\cal R}$ the boundary operators
$$ {\cal D}u:=u\qquad \mbox{and}\qquad {\cal R}u:=\partial_{\nu}u+
\beta(x)u, $$
where $\nu$ is the unit outward normal to $\partial \Omega$,
and $\beta\in C^{1,\mu}(\partial\Omega)$ is non-negative.
Hence, ${\cal D}$ is the {\it Dirichlet} boundary operator
and ${\cal R}$ is either the {\it Neumann} boundary operator,
if $\beta\equiv 0$, or the {\it Robin} boundary operator, if
$\beta\not\equiv 0$. Throughout this work, ${\cal B}$
can define any of these boundary operators.

Note that the Robin condition ${\cal R}=0$ relies
essentially to heat flow problems in a body with constant temperature
in the surrounding medium. More generally, if $\alpha$ and $\beta$ are
smooth
functions on $\partial\Omega$ such that $\alpha ,\beta\geq 0$, $\alpha
+\beta>0$,
then the boundary condition $Bu=\alpha \partial_{\nu}u+\beta u=0$
represents the exchange
of heat at the surface of the reactant by Newtonian cooling. Moreover,
the boundary condition $Bu=0$ is called isothermal (Dirichlet)
condition if
$\alpha\equiv 0$, and it becomes an adiabatic (Neumann) condition if
$\beta\equiv 0$.
An intuitive meaning of the condition $\alpha +\beta >0$ on
$\partial\Omega$ is that,
for the diffusion process described by problem \eq{sep}, either the
reflection
phenomenon or the absorption phenomenon may occur at each point of the
boundary.

We are now concerned with the following
boundary blow-up problem
\neweq{b1}
\left\{ \begin{array}{lll}
& \di \Delta u+au=b(x)f(u) \quad &\mbox{in}\ \Omega\setminus
\overline{\Omega}_0\,,\\
& \di {\cal B}u=0 \quad &\mbox{on}\ \partial \Omega\,,\\
& \di u=\infty \quad &\mbox{on}\ \partial\Omega_0\,,\\
\end{array} \right.\endeq
where $b>0$ on $\partial\Omega$, while
$\ov{\Omega}_0$ is non-empty, connected and with smooth boundary. Here,
$u=\infty$ on $\partial\Omega_0$ means that
$u(x)\to \infty$ as $x\in \Omega\setminus \ov{\Omega}_0 $ and
$d(x):={\rm dist}\,(x,\Omega_0)\to 0$.

The question of existence and uniqueness of positive solutions for
problem \eq{b1} in the case of pure superlinear power in
the non-linearity is treated by Du-Huang \cite{dh}.
Our next results extend their previous paper to the case
of much more general non-linearities of Keller-Osserman type.

In the following,
by $(\tilde A_1)$ we mean that
$(A_1)$ is fulfilled and
there exists $\lim_{u\to \infty}\left(F/f\right)'(u):=\gamma$.
Then, $\gamma\geq 0$.
\medskip

We prove

\begin{thm}\label{teo1ccm} Let $(\tilde A_1)$ and $(A_2)$ hold.
Then, for any $a\in \RR$, problem \eq{b1} has a
minimal (resp., maximal) positive solution $\un{U}_a$
(resp., $\ov{U}_a$).\end{thm}

\proof
In proving Theorem~\ref{teo1ccm} we rely on an appropriate
comparison principle which allows us
to prove that $(u_n)_{n\geq 1}$
is non-decreasing, where $u_n$ is the unique
positive solution of problem \eq{b11}
with $\Phi\equiv n$. The minimal positive solution of
\eq{b1} will be obtained as the limit of the sequence
$(u_n)_{n\geq 1}$. Note that, since
$b=0$ on $\partial\Omega_0$, the main difficulty is related
to the construction of an upper bound of this sequence
which must fit to our general framework. Next, we deduce the maximal positive
solution of \eq{b1} as the limit of the non-increasing sequence
$(v_m)_{m\geq m_1}$ provided $m_1$ is large so that
$\Omega_{m_1}\subset\subset \Omega$. We denoted by
$v_m$ the minimal positive solution
of \eq{b1} with $\Omega_0$ replaced by
\neweq{om} \Omega_{m}:=\{x\in \Omega:\ d(x)<1/m\},
\qquad
m\geq m_1.\endeq

We start with the following
auxiliary result (see C\^{\i}rstea and R\u adulescu \cite{crccm}).

\begin{lem}\label{bl3} Assume $b>0$ on $\partial\Omega$.
If $(A_1)$ and $(A_2)$ hold, then
for any positive function
$\Phi\in C^{2,\mu}(\partial\Omega_0)$ and $a\in\RR$
the problem
\neweq{b11} \left\{\begin{array}{lll}
&\di \Delta u+au=b(x)f(u) \quad & \mbox{in}\ \Omega\setminus
\overline{\Omega}_0\,,\\
&\di {\cal B}u=0 \quad & \mbox{on}\ \partial\Omega\,,\\
&\di u=\Phi \quad & \mbox{on}\ \partial\Omega_0\,,
\end{array}\right.\endeq has a unique positive solution. \end{lem}

We now come back
to the proof of Theorem~\ref{teo1ccm}, that will be divided
into two steps:

{\it Step 1. Existence of the minimal positive solution for problem
\eq{b1}}.

For any $n\geq 1$, let $u_n$ be the unique
positive solution of problem \eq{b11}
with $\Phi\equiv n$. By
the maximum principle, $u_n(x)$ increases with $n$
for all $x\in\ov{\Omega}\setminus\ov\Omega_0$. Moreover, we prove

\begin{lem}\label{mb3} The sequence $(u_n(x))_n$
is bounded from above by some function
$V(x)$ which is uniformly bounded on all compact
subsets of $\ov{\Omega}\setminus \ov{\Omega}_0$. \end{lem}

 \proof
Let $b^*$ be a $C^2$-function on $\ov\Omega\setminus\Omega_0$ such that
$$ 0<b^*(x)\leq b(x)\quad \forall x\in
\ov{\Omega}\setminus\ov{\Omega}_0.$$
For $x$ bounded away from $\partial\Omega_0$ is not a problem to
find such a function $b^*$. For $x$ satisfying
$0<d(x)<\delta$ with $\delta>0$ small such that
$x\to d(x)$ is a $C^2$-function, we can take
$$
b^*(x)=\int_{0}^{d(x)}\int_{0}^{t}[\min_{d(z)\geq
s}b(z)]\,ds\,dt.$$

Let $g\in {\cal G}$ be a function such that $(A_g)$ holds.
Since $b^*(x)\to 0$ as $d(x)\searrow 0$,
we deduce, by $(A_1)$,
the existence of some $\delta>0$
such that for all $x\in \Omega$ with $0<d(x)<\delta$
and $\xi>1$
$$ \frac{b^*(x)f(g(b^*(x))\xi)}{g''(b^*(x))\,\xi}>
\sup_{\ov{\Omega}\setminus \Omega_0}|\nabla b^*|^2
+\frac{g'(b^*(x))}{g''(b^*(x))}\,\inf_{\ov{\Omega}\setminus \Omega_0}
(\Delta b^*)+a\,\frac{g(b^*(x))}{g''(b^*(x))}.$$
Here, $\delta>0$ is taken sufficiently small
so that $g'(b^*(x))<0$
and $g''(b^*(x))>0$ for all $x$ with
$0<d(x)<\delta$.

For $n_0\geq 1$ fixed, define $V^*$ as follows

${\rm (i)}\quad \ \,V^*(x)=u_{n_0}(x)+1\quad \mbox{for}\ x\in
\ov{\Omega}\
\mbox{and near}\ \partial\Omega\,;$

${\rm (ii)}\quad \, V^*(x)=g(b^*(x))\quad \mbox{for}\ x\
\mbox{satisfying}\
0<d(x)<\delta\,;$

${\rm (iii)}\quad V^*\in C^2(\ov{\Omega}\setminus\ov{\Omega}_0)
\ \mbox{is positive on}\
\ov{\Omega}\setminus\ov{\Omega}_0\,.$

We show that for $\xi>1$ large enough the upper bound
of the sequence $(u_n(x))_{n}$ can be taken as $V(x)=\xi V^*(x)$.
Since $$ {\cal B}V(x)=\xi\, {\cal B} V^*(x)
\geq \xi \min\,\{1,\beta(x)\}\geq 0,
\quad \forall x\in \partial\Omega
\quad\mbox{and}\quad \lim_{d(x)\searrow 0}[u_n(x)-V(x)]=-\infty<0,
$$
to conclude that $u_n(x)\leq V(x)$ for all $x\in \ov{\Omega}\setminus
\ov{\Omega}_0$ it is sufficient to show that
\neweq{b12} -\Delta V(x)\geq aV(x)-b(x)f(V(x)),\qquad \forall x\in
\Omega\setminus\ov
{\Omega}_{0}.\endeq
For $x\in \Omega$ satisfying $0<d(x)<\delta$ and $\xi>1$ we have
$$ \begin{array}{lll}
& \di -\Delta V(x)-aV(x)+b(x)f(V(x))=-\xi\Delta g(b^*(x))-a\,\xi
g(b^*(x))+b(x)f(g(b^*(x))\xi)\\
&\di \quad \quad \geq\xi  g''(b^*(x))
\left(-\frac{g'(b^*(x))}{g''(b^*(x))}\,\Delta b^*(x)-
|\nabla b^*(x)|^2-a\,\frac{g(b^*(x))}{g''(b^*(x))}+
b^*(x)\,\frac{f(g(b^*(x))\xi)}{g''(b^*(x))\,\xi}\right)>0.
\end{array} $$
For $x\in \Omega$ satisfying $d(x)\geq \delta$,
$$ -\Delta V(x)-aV(x)+b(x)f(V(x))=
\xi\left(-\Delta V^*(x)-aV^*(x)+b(x)\,\frac{f(\xi
V^*(x))}{\xi}\right)\geq 0$$
for $\xi$ sufficiently large.  It follows that \eq{b12} is
fulfilled provided $\xi$ is large enough. This finishes the proof
of the lemma.\qed
\smallskip

 By Lemma~\ref{mb3},
$\un{U}_a(x)\equiv \lim_{n\to \infty}u_n(x)$ exists,
for any $x\in \ov{\Omega}\setminus \ov{\Omega}_0$.
Moreover,
$\un{U}_a$ is a positive solution of \eq{b1}.
Using the maximum principle once more, we find that
any positive solution $u$ of \eq{b1} satisfies
$u\geq u_n$ on $\ov{\Omega}\setminus \ov{\Omega}_0$, for all $n\geq 1$.
Hence $\un{U}_a$ is the minimal positive
solution of \eq{b1}.
 \medskip

{\it Proof of Theorem~\ref{teo1ccm} completed.}

{\it Step 2. Existence of the maximal positive solution for problem
\eq{b1}.}

\begin{lem}\label{*} If $\Omega_0$ is replaced by $\Omega_m$ defined
in
\eq{om}, then problem \eq{b1} has a minimal positive solution provided
that
$(A_1)$ and $(A_2)$ are fulfilled. \end{lem}

\proof
The argument used here (more easier, since $b>0$
on $\ov{\Omega}\setminus \Omega_m$) is similar to that in Step~1.
The only difference which appears in the proof (except the replacement
of
$\Omega_0$ by $\Omega_m$)
is related to the construction of $V^*(x)$ for $x$ near
$\partial\Omega_m$.
Here, we use
our Theorem~\ref{teo0} which says that,
for any $a\in\RR$, there exists a positive large solution
$u_{a,\infty}$
of problem \eq{sep} in the domain $\Omega\setminus\ov{\Omega}_m$.
We define $V^*(x)=u_{a,\infty}(x)$ for $x\in
\Omega\setminus\ov{\Omega}_m$
and near $\partial\Omega_m$. For $\xi>1$ and $x\in \Omega\setminus
\ov{\Omega}_m$ near $\partial\Omega_m$ we have
$$ \begin{array}{lll}
& \di -\Delta V(x)-aV(x)+b(x)f(V(x))=
-\xi \Delta V^*(x)-a\xi V^*(x)+b(x)f(\xi V^*(x))\\
& \quad \quad \di= b(x)[f(\xi V^*(x))-\xi f(V^*(x)]\geq 0.
\end{array} $$
This completes the proof. \qed
\medskip

Let $v_m$ be the minimal positive solution for the problem considered
in
the statement of Lemma~\ref{*}.
By the maximum principle, $v_m\geq v_{m+1}\geq u$
on $\ov{\Omega}\setminus \ov{\Omega}_m$, where
$u$ is any positive solution of \eq{b1}. Hence $\ov{U}_a(x):=
\lim_{m\to \infty}v_m(x)\geq u(x)$.
A regularity and compactness argument shows that $\ov{U}_a$
is a positive solution of \eq{b1}. Consequently,
$\ov{U}_a$ is the maximal positive solution. This concludes the proof
of Theorem~\ref{teo1ccm}.\qed

\medskip
The next question is whether one can conclude
the uniqueness of positive solutions
of problem \eq{b1}. We recall first what
is already known in this direction.
When $f(u)=u^p$, $p>1$, Du-Huang \cite{dh}
proved the uniqueness of solution to problem \eq{b1} and established
its
behaviour near $\partial\Omega_0$, under the assumption
\neweq{nh} \lim_{d(x)\searrow 0}\frac{b(x)}{[d(x)]^\tau}=c\quad
\mbox{for some positive constants}\ \tau,c>0.\endeq

We shall give a general uniqueness result provided that
$b$ and $f$ satisfy the following assumptions:
\smallskip

\noindent $(B_1)\quad \di\lim_{d(x)\searrow 0}
\frac{b(x)}{k(d(x))}=c \quad $ for some constant $c>0$,
where $0<k\in C^1(0,\delta_0)$ is increasing and satisfies
\medskip

\noindent $(B_2)\quad
\di K(t)=\frac{\int_0^{t} \sqrt{k(s)}\,ds}{\sqrt{k(t)}}
\in C^1[0,\delta_0)$, for some $\delta_0>0$.
\medskip

Assume there exist $\zeta>0$ and $t_0\geq 1$ such that
\medskip

\noindent $(A_3)\quad f(\xi t)\leq
\xi^{1+\zeta}f(t),\quad \forall \xi\in (0,1),\quad \forall t\geq
t_0/\xi$

\noindent $(A_4)\quad \mbox{the mapping}\
(0,1]\ni \xi\longmapsto A(\xi)=\lim_{u\to \infty}
\di\frac{f(\xi u)}{\xi f(u)}$ is a continuous positive function.
\medskip

Our uniqueness result is

\begin{thm}\label{teo2ccm} Assume the conditions $(\tilde A_1)$
with $\gamma\not =0$, $(A_3)$, $(A_4)$, $(B_1)$ and $(B_2)$ hold.
Then, for any $a\in \RR$, problem \eq{b1} has a unique positive
solution
$U_a$. Moreover,
$$ \lim_{d(x)\searrow 0}\frac{U_a(x)}{h(d(x))}=\xi_0,$$
where $h$ is defined by
\neweq{mb} \int_{h(t)}^{\infty}\frac{ds}{\sqrt{2F(s)}}=
\int_{0}^{t}\sqrt{k(s)}\,ds,\quad \forall t\in (0,\delta_0) \endeq
and $\xi_0$ is the unique positive solution of $A(\xi)=\di
\frac{K'(0)(1-2\gamma)+2\gamma}{c} $.\end{thm}

\begin{rem}\label{cy}
(a)\quad $(A_1)+(A_3) \Rightarrow (A_2)$.
{\rm Indeed,
$\lim_{u\to \infty}\frac{f(u)}{u^{1+\zeta}}>0$
since
$\frac{f(t)}{t^{1+\zeta}}$ is non-decreasing for
$t\geq t_0$. }

\noindent (b)\quad $K'(0)(1-2\gamma)+2\gamma\in (0,1]$
when $(\tilde A_1)$ with $\gamma\not =0$, $(A_2)$, $(B_1)$ and $(B_2)$
hold.

\noindent (c)\quad
The function $(0,\infty)\ni \xi
\longmapsto A(\xi)\in (0,\infty)$ is bijective when
$(A_3)$ and $(A_4)$ hold (see Lemma~\ref{vi}).
\end{rem}

Among the non-linearities $f$ that satisfy the assumptions
of Theorem~\ref{teo2ccm} we note:
(i) $f(u)=u^p,\ p>1$;
(ii) $f(u)=u^p \ln(u+1)$, $p>1$;
(iii) $f(u)=u^p \arctan u$, $ p>1$.

\medskip
{\it Proof of Theorem \ref{teo2ccm}}.
By $(A_4)$ we deduce that the mapping $(0,\infty)\ni \xi\longmapsto
A(\xi)=\di\lim_{u\to \infty}\frac{f(\xi u)}{\xi f(u)}$
is a continuous positive function, since
$A(1/\xi)=1/A(\xi)$ for any $\xi\in (0,1)$.
Moreover, we claim

\begin{lem}\label{vi} The function $A:(0,\infty)\to (0,\infty)$
is bijective, provided that
$(A_3)$ and $(A_4)$ are fulfilled.
\end{lem}

 \proof By the continuity of $A$, we see that the surjectivity
of $A$ follows if we prove that $\lim_{\xi\searrow 0}A(\xi)=0$. To this
aim,
let $\xi\in (0,1)$ be fixed. Using $(A_3)$ we find
$$ \frac{f(\xi u)}{\xi f(u)}\leq \xi^\zeta, \quad \forall u\geq
\frac{t_0}{\xi} $$
which yields $A(\xi)\leq \xi^\zeta$. Since $\xi\in (0,1)$
is arbitrary, it follows that $\lim_{\xi\searrow 0}A(\xi)=0$.

We now prove that the function $\xi\longmapsto A(\xi)$ is increasing
on $(0,\infty)$ which concludes our lemma. Let $0<\xi_1<\xi_2<\infty$
be chosen arbitrarily. Using assumption $(A_3)$ once more, we obtain
$$ f(\xi_1 u)=f\left(\frac{\xi_1}{\xi_2}\xi_2 u\right)\leq
\left(\frac{\xi_1}{\xi_2}\right)^{1+\zeta}f(\xi_2 u),\quad \forall
u\geq t_0\, \frac{\xi_2}{\xi_1}. $$
It follows that
$$ \frac{f(\xi_1 u)}{\xi_1 f(u)}\leq
\left(\frac{\xi_1}{\xi_2}\right)^{\zeta}
\frac{f(\xi_2 u)}{\xi_2 f(u)},
\quad \forall u\geq t_0\,\frac{\xi_2}{\xi_1}.$$
Passing to the limit as $u\to \infty$ we find
$$ A(\xi_1)\leq
\left(\frac{\xi_1}{\xi_2}\right)^{\zeta}A(\xi_2)<A(\xi_2),$$
which finishes the proof. \qed
\medskip

{\it Proof of Theorem~\ref{teo2ccm} completed.}
 Set
$ \Pi(\xi)=\di
\lim_{d(x)\searrow 0}b(x)\,\frac{f(h(d(x))\xi)}{h''(d(x))\,\xi}$,
for any $\xi> 0$. Using $(B_1)$  we find
$$ \begin{array}{lll}
\Pi(\xi)& =\di\lim_{d(x)\searrow 0}\frac{b(x)}{k(d(x))}\,
\frac{k(d(x))f(h(d(x)))}{h''(d(x))}\,\frac{f( h(d(x))\xi )}{\xi
f(h(d(x)))}
=c\,\lim_{t\searrow 0}\frac{k(t)f(h(t))}{h''(t)}\lim_{u\to \infty}
\frac{f(\xi u)}{\xi f(u)}\\
& \di =\frac{c}{K'(0)(1-2\gamma)+2\gamma}\,A(\xi).
\end{array}$$
This and Lemma~\ref{vi} imply that the function $\Pi:(0,\infty)\to
(0,\infty)$ is bijective.
Let $\xi_0$ be the unique positive solution of $\Pi(\xi)=1$, that
is $A(\xi_0)=\di\frac{K'(0)(1-2\gamma)+2\gamma}{c}$.

For $\ep\in (0,1/4)$ arbitrary,
we denote $\xi_1=\Pi^{-1}(1-4\ep)$, respectively
$\xi_2=\Pi^{-1}(1+4\ep)$.

We choose $\delta>0$
small enough such that

${\rm (i)}\ \ \,\di
{\rm dist}\, (x,\partial\Omega_0)\ \mbox{is a}\ C^2\ \mbox{function on
the set}\
\{x\in \Omega:\ {\rm dist}\,(x,\partial\Omega_0)\leq 2\delta\};$

${\rm (ii)}\ \di\left| \frac{h'(s)}{h''(s)}\,\Delta d(x)+
a\,\frac{h(s)}{h''(s)}\right|<\ep\ \mbox{and}\ h''(s)>0
\quad \mbox{for all}\ s\in (0,2\delta)
\ \mbox{and}\ x\ \mbox{satisfying}\ 0<d(x)<2\delta;$

${\rm (iii)}\
\di(\Pi(\xi_2)-\ep)\,\frac{h''(d(x))\,\xi_2}{f(h(d(x))\xi_2)}\leq
b(x)\leq
(\Pi(\xi_1)+\ep)\,\frac{h''(d(x))\,\xi_1}{f(h(d(x))\xi_1)},\quad
\mbox{for every}\ x\ \mbox{with}\ 0<d(x)<2\delta.$

${\rm (iv)}\
\di b(y)<(1+\ep)b(x), \quad \mbox{for every}\ x,y\ \mbox{with}\
0<d(y)<d(x)<2\delta$.
\medskip

Let $\sigma\in (0,\delta)$ be arbitrary. We define
$\un{v}_{\sigma}(x)=h(d(x)+\sigma)\xi_1$,
for any $x$ with $d(x)+\sigma<2\delta $,
respectively $\ov{v}_{\sigma}(x)=h(d(x)-\sigma)\xi_2$
for any $x$ with $\sigma<d(x)<2\delta$.

Using (ii), (iv) and the first inequality in (iii),
when $\sigma<d(x)<2\delta$, we obtain
(since $|\nabla d(x)|\equiv 1$)
$$\begin{array}{lll}
&\di -\Delta
\ov{v}_{\sigma}(x)-a\ov{v}_{\sigma}(x)+b(x)f(\ov{v}_{\sigma}(x))\\
& \di \quad\quad =\xi_2\left(-h'(d(x)-\sigma)\,\Delta
d(x)-h''(d(x)-\sigma)-
a\,h(d(x)-\sigma)+\frac{b(x)f(h(d(x)-\sigma)\xi_2)}{\xi_2}\right)\\
&\di \quad\quad =
\xi_2\,h''(d(x)-\sigma)\left(-\frac{h'(d(x)-\sigma)}{h''(d(x)-\sigma)}
\,\Delta d(x)-a\,\frac{h(d(x)-\sigma)}{h''(d(x)-\sigma)}-1+
\frac{b(x)f(h(d(x)-\sigma)\xi_2)}{h''(d(x)-\sigma)\,\xi_2}\right)\\
&\di \quad \quad \geq
\xi_2\,h''(d(x)-\sigma)\left(-\frac{h'(d(x)-\sigma)}{h''(d(x)-\sigma)}
\,\Delta d(x)-a\,\frac{h(d(x)-\sigma)}{h''(d(x)-\sigma)}-1+
\frac{\Pi(\xi_2)-\ep}{1+\ep}\right)\geq 0 \end{array} $$
for all $x$ satisfying $\sigma<d(x)<2\delta$.

Similarly, using (ii), (iv) and the second inequality in (iii), when
$d(x)+\sigma<2\delta$ we find
$$\begin{array}{lll}
&\di -\Delta
\un{v}_{\sigma}(x)-a\un{v}_{\sigma}(x)+b(x)f(\un{v}_{\sigma}(x))\\
&\di \quad\quad =
\xi_1\,h''(d(x)+\sigma)\left(-\frac{h'(d(x)+\sigma)}{h''(d(x)+\sigma)}
\,\Delta d(x)-a\,\frac{h(d(x)+\sigma)}{h''(d(x)+\sigma)}-1+
\frac{b(x)f(h(d(x)+\sigma)\xi_1)}{h''(d(x)+\sigma)\,\xi_1}\right)\\
&\di \quad \quad \leq
\xi_1\,h''(d(x)+\sigma)\left(-\frac{h'(d(x)+\sigma)}{h''(d(x)+\sigma)}
\,\Delta d(x)-a\,\frac{h(d(x)+\sigma)}{h''(d(x)+\sigma)}-1+
(1+\ep)(\Pi(\xi_1)+\ep)\right)\leq 0, \end{array} $$
for all $x$ satisfying $d(x)+\sigma<2\delta$.

Define $ \Omega_{\delta}\equiv \{x\in\Omega:\ d(x)<\delta\}$.
Let $\omega\subset\subset \Omega_0$ be such that the first Dirichlet
eigenvalue of $(-\Delta)$ in the smooth domain
$\Omega_0\setminus\ov{\omega}$ is strictly greater than $a$. Denote by
$w$ a positive large solution to the following problem
$$ -\Delta w=aw-p(x)f(w)\qquad {\rm in}\ \Omega_{\delta}\,,$$
where $p\in C^{0,\mu}(\ov{\Omega}_{\delta})$
satisfies $0<p(x)\leq b(x)$ for $x\in \ov{\Omega}_{\delta}\setminus
\ov{\Omega}_0$, $p(x)=0$ on $\ov{\Omega}_0\setminus \omega$
and $p(x)>0$ for $x\in\omega$.
The existence of $w$ is guaranteed by our Theorem~\ref{teo0}.

Suppose that $u$ is an arbitrary solution of \eq{b1} and let
$v:=u+w$. Then $v$ satisfies
$$ -\Delta v\geq av-b(x)f(v)\quad \mbox{in}\ \Omega_{\delta}\setminus
\overline{\Omega}_{0}.$$
Since
$$ v|_{\partial \Omega_0}=\infty>\un{v}_{\sigma}|_{\partial
\Omega_0}\quad
\mbox{and}\quad
v|_{\partial \Omega_\delta}=\infty>\un{v}_{\sigma}|_{\partial
\Omega_\delta},$$
 we find
\neweq{b7} u+w\geq \un{v}_{\sigma}\ \ \mbox{on}\ \Omega_\delta\setminus
\overline{\Omega}_0.
\endeq Similarly
\neweq{b8} \ov{v}_{\sigma}+w\geq u\ \ \mbox{on}\ \Omega_\delta\setminus
\overline{\Omega}_\sigma.\endeq
Letting $\sigma\to 0$ in \eq{b7} and \eq{b8}, we deduce
$$ h(d(x))\,\xi_2+2w\geq u+w\geq h(d(x))\,\xi_1,\quad \forall x\in
\Omega_{\delta}\setminus \overline{\Omega}_0.$$
Since $w$  is uniformly bounded on $\partial \Omega_0$, it follows that
\neweq{b9} \xi_1\leq \liminf_{d(x)\searrow 0}\,\frac{u(x)}{h(d(x))}\leq
\limsup_{d(x)\searrow 0}\,\frac{u(x)}{h(d(x))}\leq \xi_2.\endeq
Letting $\ep\to 0$ in \eq{b9} and looking at the definition
of $\xi_1$ respectively $\xi_2$ we find
\neweq{b10} \lim_{d(x)\searrow 0}\frac{u(x)}{h(d(x))}=\xi_0.\endeq
This behavior of the solution will be speculated
in order to prove that problem
\eq{b1} has a unique solution. Indeed, let $u_1$, $u_2$ be two positive
solutions of \eq{b1}. For any $\ep>0$, denote $\tilde u_i=
(1+\ep)\,u_i$, $i=1,2$. By virtue of \eq{b10} we get
$$ \lim_{d(x)\searrow 0}\frac{u_1(x)-\tilde u_2(x)}{h(d(x))}=
\lim_{d(x)\searrow 0}\frac{u_2(x)-\tilde u_1(x)}{h(d(x))}=-\ep\,\xi_0<0
$$
which implies
$$ \lim_{d(x)\searrow 0}[u_1(x)-\tilde u_2(x)]=\lim_{d(x)\searrow
0}[u_2(x)-
\tilde u_1(x)]=-\infty. $$
On the other hand, since $\frac{f(u)}{u}$ is increasing
for $u>0$, we obtain
$$ \begin{array}{lll}
&\di -\Delta \tilde u_i=-(1+\ep)\,\Delta
u_i=(1+\ep)\left(a\,u_i-b(x)f(u_i)\right)\geq
a\,\tilde u_i-b(x)f(\tilde u_i)\quad \mbox{in}\ \Omega\setminus
\ov{\Omega}_0,\\
&\di {\cal B}\tilde u_i={\cal B}u_i=0\quad \mbox{on}\ \partial\Omega.
\end{array}$$
So, by the maximum principle,
$$ u_1(x)\leq \tilde u_2(x),\quad u_2(x)\leq \tilde u_1(x),\quad
\forall x\in \Omega\setminus \ov{\Omega}_0.$$
Letting $\ep\to 0$, we obtain $u_1\equiv u_2$. The proof of
Theorem~\ref{teo2ccm}
is complete. \qed

The above results have been established by C\^{\i}rstea and R\u
adulescu \cite{crccm,crhouston}.

\subsection{Uniqueness and asymptotic behaviour of the large solution. A Karamata
regular variation theory approach}

The major purpose in this section is to advance innovative methods
to study the uniqueness and asymptotic behavior of large solutions
of \eq{sep}.  This approach is due to C\^{\i}rstea and R\u
adulescu \cite{crasun,crasas,crcras04,crasan,crtams} and it relies
essentially on the {\it regular variation theory} introduces by
Karamata (see Bingham, Goldie, and
 Teugels \cite{bgt}, Karamata \cite{karamata}), not only in the statement but in
the proof as well. This enables us to obtain significant
information about the qualitative behavior of the large solution
to \eq{sep} in a general framework that removes previous
restrictions in the literature.

\begin{defin}\label{def1cras}  A positive
measurable function $R$ defined on $[D,\infty)$, for
some $D>0$, is called regularly varying (at infinity)
with index  $q\in \RR$ (written $R\in RV_q$) if
for all $\xi>0$
$$ \lim_{u\to \infty}R(\xi u)/R(u)=\xi^q. $$
When the index of regular variation $q$ is zero, we say that the
function is slowly varying. \end{defin}

We remark that
any function $R\in RV_q$ can be written
in terms of a slowly varying function. Indeed, set
$R(u)=u^qL(u)$. From the above definition we easily deduce that $L$ varies slowly.

The canonical $q$-varying function is $u^q$. The functions $\ln (1+u)$,
$\ln \ln(e+u)$, ${\rm exp}\,\{(\ln u)^\alpha\}$, $\alpha\in (0,1)$
vary slowly, as well as any measurable function
on $[D,\infty)$ with positive limit at infinity.

In what follows $L$ denotes an arbitrary slowly
varying function
and $D>0$ a positive number. For details
on the below properties, we refer to Seneta \cite{seneta}.

\begin{prop}\label{prop1cras} (i) For any $m>0$, $u^m L(u)\to
\infty$,\quad $u^{-m} L(u)\to 0$ as $u\to \infty$.

(ii) Any positive $C^1$-function on $[D, \infty)$
satisfying $ u L_1'(u)/L_1(u)\to 0$ as $u\to \infty $
is slowly varying. Moreover, if the above limit is
$q\in \RR$, then $L_1\in RV_q$.

(iii) Assume $R:[D,\infty)\to (0,\infty)$
is measurable and Lebesgue integrable on each finite
subinterval of
$[D,\infty)$. Then $R$ varies regularly iff
there exists $j\in \RR$ such that
\neweq{lt} \lim_{u\to
\infty}\frac{u^{j+1}R(u)}{\int_D^u x^j R(x)\,dx}
\endeq
exists and is a positive number, say $a_j+1$.
In this case, $R\in RV_q$ with $q=a_j-j$.

(iv) (Karamata Theorem, 1933).
If $R\in RV_q$ is Lebesgue integrable
on each finite subinterval of $[D,\infty)$, then the
limit defined by
\eq{lt} is $q+j+1$, for every $j>-q-1$.
\end{prop}

\begin{lem}\label{lema1cras} Assume $(A_1)$ holds. Then we have
the equivalence
$$ a)\ f'\in RV_\rho \Longleftrightarrow \
b)\ \lim_{u\to \infty}uf'(u)/f(u):=\vartheta<\infty
\Longleftrightarrow \ c)\
\lim_{u\to \infty}\left(F/f\right)'(u):=\gamma>0. $$
\end{lem}

\begin{rem}\label{rem2cras} Let $a)$ of Lemma \ref{lema1cras} be
fulfilled. Then the following assertions hold

 (i) $\rho$ is non-negative;

 (ii) $ \gamma=1/(\rho+2)=1/(\vartheta+1)$;

 $(iii)$ If $\rho\not=0$, then $(A_2)$ holds
({\rm use
$\lim_{u\to \infty}f(u)/u^p=\infty$, $\forall p\in
(1,1+\rho)$)}.
The converse implication is not
necessarily true (take $f(u)=u \ln^4 (u+1)$).
However, there are cases when
$\rho=0$ and $(A_2)$ fails so that
\eq{sep} has {\bf no} large solutions. This is
illustrated by
$f(u)=u$ or $f(u)=u \ln (u+1)$. \end{rem}

Inspired by the definition of $\gamma$, we denote by
${\cal K}$ the set of all positive, increasing
$C^1$-functions $k$
defined on $(0,\nu)$, for some $\nu>0$, which satisfy
$\lim_{t\to 0^+}\di\left(\frac{\int_0^t
k(s)\,ds}{k(t)}\right)^{(i)}:=\ell_i,\
i=\overline{0,1}$.

It is easy to see that
$\ell_0=0$ and $\ell_1\in [0,1]$, for every $k\in {\cal K}$.
Our next result gives examples of functions $k\in
{\cal K}$ with
$\lim_{t\to 0^+}k(t)=0$, for {\it every} $\ell_1\in
[0,1]$.

\begin{lem}\label{lema2cras}
Let $S\in C^1[D,\infty)$ be such that
$S'\in RV_q$ with $q>-1$. Hence the following hold:

$a)$\quad If $\di k(t)={\rm exp}\,\{-S(1/t)\}\quad
\forall t\leq 1/D$,
then
$k\in {\cal K}$ with $\ell_1=0$.

 $b)$\quad If $\di k(t)=
1/S(1/t) \quad \forall t\leq 1/D$, then $k\in {\cal
K}$ with
$\ell_1=1/(q+2)\in (0,1)$.

 $c)$\quad If $\di k(t)=
1/\ln S(1/t)\quad \forall t\leq 1/D$, then $k\in {\cal
K}$ with
$\ell_1=1$.
\end{lem}

\begin{rem}\label{rem3cras} If $S\in C^1[D,\infty)$, then
$S'\in RV_q$ with $q>-1$ iff for some $m>0$, $C>0$
and $B>D$ we have
$ S(u)= Cu^m {\rm exp}\left\{ \int_B^u\frac{y(t)}{t}\,dt\right\}$,
$\forall u\geq B $, where $y\in C[B,\infty)$ satisfies $\lim_{u\to \infty}
y(u)=0$. In this case, $S'\in RV_q$ with $q=m-1$. (This is a consequence
of Property \ref{prop1cras} (iii) and (iv).\end{rem}

Our main result is

\begin{thm}\label{1karamata}
Let $(A_1)$ hold and $f'\in RV_\rho$ with $\rho>0$.
Assume $b\equiv 0$ on $\partial\Omega$ satisfies
\medskip

\noindent $(B)\quad b(x)=c\,k^2(d(x))+o(k^2(d(x)))$ as
$d(x)\to 0$,
for some constant $c>0$ and $k\in {\cal K}$\,.
\medskip

Then, for any $a\in (-\infty,\lambda_{\infty,1})$,
Eq. \eq{sep} admits a unique large
solution $u_a$. Moreover,
\neweq{po}
\lim_{d(x)\to 0}\frac{u_a(x)}{h(d(x))}=\xi_0,\endeq
where $\xi_0=\di\left(\frac{2+\ell_1
\rho}{c(2+\rho)}\right)^{1/\rho}$
and
$h$ is defined by
\neweq{mbcras}
\int_{h(t)}^{\infty}\frac{ds}{\sqrt{2F(s)}}=
\int_{0}^{t}k(s)\,ds,\quad \forall t\in (0,\nu).
\endeq
\end{thm}

By Remark \ref{rem3cras}, the assumption $f'\in RV_\rho$ with $\rho>0$ holds
if and only if
there exist $p>1$ and $B>0$ such that
$ f(u)=C u^p {\rm exp}\left\{\int_B^u\frac{y(t)}{t}\,dt\right\}$, for all $u\geq B $
($y$ as before and $p=\rho+1$). If $B$ is large enough ($y>-\rho$ on $[B,\infty)$), then
$f(u)/u$ is increasing on $[B,\infty)$.
Thus, to get the whole range of functions $f$ for which our Theorem \ref{1karamata}
applies we have only to ``paste" a suitable smooth function on $[0,B]$
in accordance with $(A_1)$. A simple way to do this is to define
$ f(u)=u^p{\rm exp}\{\int_0^u\frac{z(t)}{t}\,dt\}$, for all $u\geq 0$,
where $z\in C[0,\infty)$ is non-negative
such that $\lim_{t\to 0^+}z(t)/t\in [0,\infty)$ and
$\lim_{u\to \infty}z(u)=0$. Clearly,
$f(u)=u^p$, $f(u)=u^p\ln (u+1)$, and
$f(u)=u^p\arctan u$ ($p>1$) fall into this category.

Lemma \ref{lema2cras} provides a practical method to find functions $k$ which can be
considered in the statement of Theorem \ref{1karamata}. Here are some examples:
$k(t)=-1/\ln t$, $k(t)=t^{\alpha}$,
$k(t)={\rm exp}\left\{-1/t^{\alpha}\right\}$,
$k(t)={\rm exp}\left\{-\ln (1+\frac{1}{t})/t^{\alpha}\right\}$,
$k(t)={\rm exp}\left\{-\left[\arctan
\left(\frac{1}{t}\right)\right]/t^{\alpha}\right\}$,
 $ k(t)=t^{\alpha}/\ln (1+\frac{1}{t})$,
for some $\alpha>0$.

As we shall see, the uniqueness lies upon the crucial observation
\eq{po}, which shows that all  explosive solutions have the same
boundary behaviour. Note that the only case of Theorem
\ref{1karamata} studied so far is $f(u)=u^p$ ($p>1$) and
$k(t)=t^\alpha$ ($\alpha>0$) (see Garc\'{i}a-Meli\'{a}n,
Letelier-Albornoz, and Sabina de Lis \cite{gls}). For related
results on the uniqueness of explosive solutions (mainly in the
cases $b\equiv 1$ and $a=0$) we refer to Bandle and Marcus
\cite{bm}, Loewner and Nirenberg \cite{ln}, Marcus and V\'eron
\cite{mv1}.

\smallskip
{\it Proof of Lemma \ref{lema1cras}}.
From Property \ref{prop1cras} (iv) and Remark \ref{rem2cras}
$(i)$ we deduce
$a)\Longrightarrow b)$ and $\vartheta=\rho+1$.
Conversely,
$b)\Longrightarrow a)$ follows
by \ref{prop1cras} (iii) since $\vartheta\geq 1$ cf. $(A_1)$.

$b) \Longrightarrow c)$. Indeed, $\lim_{u\to \infty}
\frac{uf(u)}{F(u)}=1+\vartheta $, which yields
$\frac{\vartheta}{1+\vartheta}=
\lim_{u\to
\infty}\left[1-\left(\frac{F}{f}\right)'(u)\right]=1-\gamma$.

$c) \Longrightarrow b)$. Choose $s_1>0$ such that
$\left(\frac{F}{f}
\right)'(u)\geq \frac{\gamma}{2}$, $\forall u\geq
s_1$.
So, $\left(\frac{F}{f}\right)(u)\geq
\frac{(u-s_1)\gamma}{2}+
\left(\frac{F}{f}\right)(s_1)$, $\forall u\geq s_1$.
Passing to the
limit $u\to \infty$, we find
$\lim_{u\to \infty}\frac{F(u)}{f(u)}=\infty$. Thus,
$\lim_{u\to
\infty}\frac{uf(u)}{F(u)}=\frac{1}{\gamma}$. Since
$1-\gamma:=
\lim_{u\to \infty}\frac{F(u)f'(u)}{f^2(u)}$, we obtain
$\lim_{u\to
\infty}
\frac{uf'(u)}{f(u)}=\frac{1-\gamma}{\gamma}$. \qed
\smallskip

{\it Proof of Lemma \ref{lema2cras}}.   Since $\lim_{u\to
\infty}uS'(u)=\infty$
(cf. Property \ref{prop1cras} (i)), from Karamata Theorem we deduce
$\lim_{u\to \infty}
\frac{uS'(u)}{S(u)}=q+1>0$. Therefore, in any of the
cases $a)$, $b)$,
$c)$,
$\lim_{t\to 0^+}k(t)=0$ and $k$ is an increasing
$C^1$-function
on $(0,\nu)$, for $\nu>0$ sufficiently small.

$a)$ It is clear that $\lim_{t\to
0^+}\frac{tk'(t)}{k(t)\ln k(t)}=
\lim_{t\to 0^+}\frac{-S'(1/t)}{tS(1/t)}=-(q+1)$. By
l'Hospital's rule,
$\ell_0=\lim_{t\to 0^+}\frac{k(t)}{k'(t)}=0$ and
$\lim_{t\to 0^+}\frac{\left(\int_0^t
k(s)\,ds\right)\ln
k(t)}{tk(t)}=-\frac{1}{q+1}$. So,
$1-\ell_1:=\lim_{t\to 0^+}\frac{\left(\int_0^t
k(s)\,ds\right)
k'(t)}{k^2(t)}=1$.

$b)$ We see that $\lim_{t\to
0^+}\frac{tk'(t)}{k(t)}=\lim_{t\to 0^+}
\frac{S'(1/t)}{t S(1/t)}=q+1$. By l'Hospital's rule,
$\ell_0=0$ and
$\lim_{t\to 0^+}\frac{\int_0^t
k(s)\,ds}{tk(t)}=\frac{1}{q+2}$. So,
$ \ell_1=1-\lim_{t\to 0^+}\frac{\int_0^t
k(s)\,ds}{tk(t)}\,\frac{tk'(t)}{k(t)}=
\frac{1}{q+2}$.

$c)$ We have $\lim_{t\to 0^+}\frac{tk'(t)}{k^2(t)}=
\lim_{t\to 0^+}\frac{S'(1/t)}{tS(1/t)}=q+1$. By
l'Hospital's rule,
$\lim_{t\to 0^+}\frac{\int_0^t k(s)\,ds}{tk(t)}=1$.
Thus, $\ell_0=0$
and $\ell_1=1-\lim_{t\to 0^+}\frac{\int_0^t
k(s)\,ds}{t}\,\frac{tk'(t)}{k^2(t)}=1$.
\qed

\medskip
{\it Proof of Theorem \ref{1karamata}}.  Fix $a\in
(-\infty,\lambda_{\infty,1})$.
By Theorem \ref{teo0}, problem \eq{sep} has at least a large
solution.

If we prove that \eq{po} holds for an {\it arbitrary}
large solution
$u_a$ of \eq{sep}, then the uniqueness follows easily.
Indeed, if $u_1$ and $u_2$ are two
arbitrary large solutions of \eq{sep}, then \eq{po}
yields
$ \lim_{d(x)\to 0^+}\frac{u_1(x)}{u_2(x)}=1 $. Hence,
for
any $\ep\in (0,1)$, there exists
$\delta=\delta(\ep)>0$
such that
\neweq{uq} (1-\ep)u_2(x)\leq u_1(x)\leq (1+\ep)u_2(x),
\quad \forall x\in \Omega\ \mbox{with}\ 0<d(x)\leq
\delta.\endeq
Choosing eventually a smaller $\delta>0$, we can
assume that
$\ov{\Omega}_0\subset C_\delta$, where
$ C_\delta:=\{x\in \Omega:\ d(x)>\delta\}$.

It is clear that $u_1$ is a positive solution of
the boundary value problem
\neweq{mcm}
\Delta \phi+a\phi=b(x)f(\phi) \quad \mbox{in}\
C_\delta, \qquad
\phi=u_1 \quad \mbox{on}\ \partial C_\delta.
\endeq
By $(A_1)$ and \eq{uq}, we see that $\phi^-= (1-\ep)u_2$ (resp.,
$\phi^+=(1+\ep)u_2$) is a positive sub-solution (resp.,
super-solution) of \eq{mcm}. By the sub and super-solutions
method, \eq{mcm} has a positive solution $\phi_1$ satisfying
$\phi^-\leq \phi_1\leq \phi^+$ in $C_\delta$. Since $b>0$ on
$\ov{C}_\delta\setminus \ov{\Omega}_0$,
 we deduce that \eq{mcm} has a {\it
unique} positive
solution, that is, $u_1\equiv \phi_1$ in $C_\delta$. This
yields
$ (1-\ep)u_2(x)\leq u_1(x)\leq (1+\ep)u_2(x)$ in
$C_\delta$,
so that \eq{uq} holds in $\Omega$. Passing to the
limit
$\ep\to 0^+$, we conclude that $u_1\equiv u_2$.

In order to prove \eq{po} we state some useful
properties about $h$:

\noindent $(h_1)\ h\in C^2(0,\nu)$, $\lim_{t\to
0^+}h(t)=\infty$
(straightforward from \eq{mbcras}).

\noindent $(h_2)\ \lim_{t\to
0^+}\di\frac{h''(t)}{k^2(t)f(h(t)\xi)}=
\frac{1}{\xi^{\rho+1}}\,\frac{2+\rho \ell_1}{2+\rho}$,
$\forall \xi>0$
(so, $h''>0$ on $(0,2\delta)$, for $\delta>0$ small
enough).

\noindent $(h_3)\ \lim_{t\to 0^+}h(t)/h''(t)=
\lim_{t\to 0^+}h'(t)/h''(t)=0$.

\medskip
\noindent We check $(h_2)$ for $\xi=1$ only, since
$f\in RV_{\rho+1}$. Clearly,
$h'(t)=-k(t)\sqrt{2F(h(t))}$ and
\neweq{xm}
h''(t)=k^2(t)f(h(t))\left(1-2\frac{k'(t)\left(\int_0^t k(s)\,ds\right)}
{k^2(t)}\,\frac{ \sqrt{F(h(t))} } { f(h(t))\int_{h(t)}^\infty
[F(s)]^{-1/2}ds
}\right)\quad \forall t\in (0,\nu).\endeq
We see that $\lim_{u\to \infty}\sqrt{F(u)}/f(u)=0$.
Thus, from
l'Hospital's rule and Lemma \ref{lema1cras} we infer that
\neweq{xn} \lim_{u\to
\infty}\frac{\sqrt{F(u)}}{f(u)\int_u^\infty
[F(s)]^{-1/2}ds}
=\frac{1}{2}-\gamma=\frac{\rho}{2(\rho+2)}.\endeq
Using \eq{xm} and \eq{xn} we derive $(h_2)$ and also
\neweq{la} \lim_{t\to 0^+}
\frac{h'(t)}{h''(t)}=\frac{-2(2+\rho)}{2+\ell_1\rho}\,
\lim_{t\to 0^+}\frac{\int_0^t
k(s)\,ds}{k(t)}\,\lim_{u\to \infty}
\frac{\sqrt{F(u)}}{f(u)\int_u^\infty
[F(s)]^{-1/2}ds}=\frac{-\rho \ell_0}
{2+\ell_1\rho}=0.
\endeq
From $(h_1)$ and $(h_2)$, $\lim_{t\to 0^+}h'(t)=-\infty$. So,
l'Hospital's rule and \eq{la} yield $\lim_{t\to
0^+}\frac{h(t)}{h'(t)}=0$. This and \eq{la} lead to $\lim_{t\to
0^+}\frac{h(t)}{h''(t)}=0$ which proves $(h_3)$.

\medskip
{\it Proof of \eq{po}}. Fix $\ep\in (0,c/2)$. Since
$b\equiv 0$ on $\partial\Omega$ and $(B)$ holds, we
take $\delta>0$ so
that

\noindent $(i)\ \ \,\,d(x)$ is a $C^2$-function on the
set
$\{x\in \RR^N:\ d(x)<2\delta \}$;

\noindent $(ii)\ \,\,k^2$ is increasing on
$(0,2\delta)$;

\noindent $(iii)\
(c-\ep)k^2(d(x))<b(x)<(c+\ep)k^2(d(x))$,
$\forall x\in \Omega$ with $0<d(x)<2\delta$;

\noindent $(iv)\ \,\,h''(t)>0$ $\forall t\in
(0,2\delta)$ (from
$(h_2)$).

Let $\sigma\in (0,\delta)$ be arbitrary. We define
$\xi^\pm=
\left[\frac{2+\ell_1\rho}{(c\mp
2\ep)(2+\rho)}\right]^{1/\rho}$
and $v_\sigma^-(x)=h(d(x)+\sigma)\xi^-$, for all $x$
with
$d(x)+\sigma<2\delta$ resp.,
$v_\sigma^+(x)=h(d(x)-\sigma)\xi^+$, for all $x$ with
$\sigma<d(x)<2\delta$.

Using $(i)$-$(iv)$, when $\sigma<d(x)<2\delta$ we
obtain (since
$|\nabla d(x)|\equiv 1$)
$$ \begin{array}{lll}
&\di \Delta v_\sigma^++av_\sigma^+-b(x)f(v_{\sigma}^+)
\leq
\xi^+h''(d(x)-\sigma)\left(\frac{h'(d(x)-\sigma)}{h''(d(x)-\sigma)}
\Delta d(x)+a\frac{h(d(x)-\sigma)}{h''(d(x)-\sigma)}+1
\right.\\
&\di \qquad \qquad \quad
\hspace{5cm}\left.-(c-\ep)\frac{k^2(d(x)-\sigma)f(h(d(x)-\sigma)\xi^+)}
{h''(d(x)-\sigma)\xi^+}\right).
\end{array} $$

Similarly, when $d(x)+\sigma<2\delta $ we find
$$ \begin{array}{lll}
&\di \Delta v_\sigma^-+av_\sigma^--b(x)f(v_{\sigma}^-)
\geq
\xi^-h''(d(x)+\sigma)\left(\frac{h'(d(x)+\sigma)}{h''(d(x)+\sigma)}
\Delta d(x)+a\frac{h(d(x)+\sigma)}{h''(d(x)+\sigma)}+1
\right.\\
&\di \qquad \qquad \quad \hspace{5cm}\left.-
(c+\ep)\frac{k^2(d(x)+\sigma)f(h(d(x)+\sigma)\xi^-)}
{h''(d(x)+\sigma)\xi^-}\right).
\end{array} $$
Using $(h_2)$ and $(h_3)$ we see that, by diminishing
$\delta$, we can
assume
$$ \Delta
v_\sigma^+(x)+av_\sigma^+(x)-b(x)f(v_{\sigma}^+(x))\leq
0\quad
\forall x\ \mbox{with}\ \sigma<d(x)<2\delta; $$
$$ \Delta
v_\sigma^-(x)+av_\sigma^-(x)-b(x)f(v_{\sigma}^-(x))\geq
0\quad
\forall x\ \mbox{with}\ d(x)+\sigma<2\delta.
$$

Let $\Omega_1$ and $\Omega_2$ be smooth bounded
domains
such that $\Omega\subset\subset
\Omega_1\subset\subset\Omega_2$
and the first Dirichlet eigenvalue of $(-\Delta)$ in
the domain
$\Omega_1\setminus\ov{\Omega}$ is greater than $a$.
Let $p\in C^{0,\mu}(\ov{\Omega}_2)$ satisfy
$0<p(x)\leq b(x)$ for
$x\in \Omega\setminus C_{2\delta}$, $p=0$ on
$\ov{\Omega}_1\setminus
\Omega$ and $p>0$ on $\Omega_2\setminus
\ov{\Omega}_1$.
Denote by $w$ a positive large solution of
$$ \Delta w+aw=p(x)f(w)\quad \mbox{in}\ \Omega_2
\setminus
\ov{C}_{2\delta}.$$
The existence of $w$ is ensured by Theorem \ref{teo0}.

Suppose that $u_a$ is an arbitrary large solution of
\eq{sep} and let
$v:=u_a+w$. Then $v$ satisfies
$$ \Delta v+ a v-b(x)f(v)\leq 0\quad \mbox{in}\
\Omega\setminus
\ov{C}_{2\delta}.$$
Since
$v_{|\partial\Omega}=\infty>v^-_{\sigma|\partial\Omega}$
and $v_{|\partial
C_{2\delta}}=\infty>v^-_{\sigma|\partial
C_{2\delta}}$,
the maximum principle implies
\neweq{wo} u_a+w\geq v_\sigma^-\quad \mbox{on}\
\Omega\setminus
\ov{C}_{2\delta}.\endeq
Similarly,
\neweq{we} v_\sigma^++w\geq u_a\quad \mbox{on}\
C_\sigma\setminus
\ov{C}_{2\delta}. \endeq
Letting $\sigma\to 0$ in \eq{wo} and \eq{we}, we
deduce
$ h(d(x))\xi^++2w\geq u_a+w\geq h(d(x))\xi^-$, for all
$x\in \Omega
\setminus \ov{C}_{2\delta}$.
Since $w$ is uniformly bounded on $\partial\Omega$,
we have $ \di\xi^-\leq
\liminf_{d(x)\to 0}\frac{u_a(x)}{h(d(x))}\leq
\limsup_{d(x)\to 0}\frac{u_a(x)}{h(d(x))}\leq \xi^+$.
Letting $\ep\to 0^+$ we obtain \eq{po}. This concludes
the proof of
Theorem \ref{1karamata}. \qed

\medskip
 Bandle and  Marcus proved in \cite{bm1} that the blow-up rate of the
 unique large solution of \eq{sep}
depends on the curvature of the boundary of $\Omega$.
Our purpose in what follows is to refine the blow-up
rate of $u_a$ near $\partial\Omega$ by giving the second term in
its expansion near the boundary. This is a more subtle
question which represents the goal of more recent
literature (see Garc\'{i}a-Meli\'{a}n, Letelier-Albornoz, and Sabina de Lis \cite{gls}
and the references therein). The following
is very general and, as a novelty, it relies on the Karamata regular variation theory.

Recall that ${\cal K}$ denotes the set of all positive increasing
$C^1$-functions $k$ defined on $(0,\nu)$, for some $\nu>0$, which
satisfy $\lim_{t\searrow 0}(\int_0^t k(s)\,ds/k(t))^{(i)}
:=\ell_i$, $i\in \ov{0,1}$.
We also recall that $RV_q$ ($q\in \RR$) is the set of all positive measurable functions
$Z:[A,\infty)\ri\RR$ (for some $A>0$) satisfying
$ \lim_{u\to \infty}Z(\xi u)/Z(u)=\xi^q$,
$\forall\xi>0$.
Define by $NRV_q$
the class of functions $f$
in the form
$f(u)=Cu^{q}{\rm exp}\left\{\int_B^u \phi(t)/t\,dt\right\}$,
$\forall u\geq B>0$, where $C>0$ is a constant and $\phi\in
C[B,\infty)$ satisfies $\lim_{t\ri\infty}\phi (t)=0$.
The Karamata Representation Theorem shows that $NRV_q\subset RV_q$.

For any $\zeta>0$, set ${\cal K}_{0,\zeta}$ the subset of
${\cal K}$ with $\ell_1=0$ and
$\lim_{t\searrow 0}t^{-\zeta}
(\int_0^t k(s)\,ds/k(t))':=L_\star\in \RR$.
It can be proven that
${\cal K}_{0,\zeta}\equiv {\cal R}_{0,\zeta}$, where
$$ {\cal R}_{0,\zeta}=
\left\{\begin{array}{lll}
k: &
k({u}^{-1})=d_0 u\left[\Lambda(u)\right]^{-1}\,{\rm
exp}\left[-\int_{d_1}^u
\left(s\Lambda(s)\right)^{-1}ds \right]\ (u\geq d_1), \
0<\Lambda\in C^1[d_1,\infty),\\
& \lim_{u\to \infty}\Lambda(u)=\lim_{u\to \infty}
u\Lambda'(u)=0,\ \lim_{u\to \infty}u^{\zeta+1}\Lambda'(u)=\ell_\star\in
\RR,\
d_0,\, d_1>0
\end{array}\right\}.$$
Define
\begin{eqnarray*}
& & {\cal F}_{\rho\eta}=\left\{f\in NRV_{\rho+1}\
(\rho>0):\ \phi\in RV_\eta\ \mbox{or}\
-\phi\in RV_\eta\right\},\ \eta\in (-\rho-2,0];\\
& & {\cal F}_{\rho0,\tau}=
\{f\in {\cal F}_{\rho 0}:\ \lim_{u\to \infty}(\ln u)^\tau
\phi(u)=\ell^\star\in
\RR\},\ \tau\in (0,\infty).
\end{eqnarray*}

The following result establishes a precise asymptotic estimate in the neighbourhood
of the boundary.

\begin{thm}\label{t2asym} Assume that
\begin{equation}\label{b}
b(x)=k^2(d)(1+\tilde c d^\theta+o(d^\theta))\ \ \mbox{if}\ d(x)\to 0,
\ \ \mbox{where}\ k\in {\cal R}_{0,\zeta},\ \theta>0,\ \tilde c\in \RR.
\end{equation}

Suppose that $f$ fulfills
$(A_1)$ and one of the following growth conditions at infinity:

{\rm (i)} $f(u)=C u^{\rho+1}$ in a neighbourhood of infinity;

{\rm (ii)} $f\in {\cal F}_{\rho\eta}$ with $\eta\not=0$;

{\rm (iii)} $f\in {\cal F}_{\rho0,\tau_1}$ with
$\tau_1=\varpi/\zeta$, where $\varpi=\min\{\theta,\zeta\}$.

Then, for any $a\in
(-\infty,\lambda_{\infty,1})$, the unique positive solution $u_a$ of
\eq{sep}
satisfies
\begin{equation}\label{expan}
u_a(x) =\xi_0 h(d)(1+\chi d^\varpi+o(d^\varpi))\quad
\mbox{if }\ d(x)\ri 0,\quad \mbox{where}\
\xi_0=[2(2+\rho)^{-1}]^{1/\rho}\end{equation}
and $h$ is defined by
$ \int_{h(t)}^\infty [2F(s)]^{-1/2}ds=\int_0^t
k(s)\,ds$, for $t>0$ small enough.
The expression of $\chi$ is
$$
\chi=\left\{\begin{array}{lll}
& \di-(1+\zeta)\ell_\star(2\zeta)^{-1}\,
{\rm Heaviside}\,(\theta-\zeta)-
\tilde c\rho^{-1}\,
{\rm Heaviside}\,(\zeta-\theta):=\chi_1\ \ \mbox{if}\ {\rm (i)}\
\mbox{or}\
{\rm (ii)}\ \mbox{holds}\\
& \di\chi_1-
\ell^\star\rho^{-1} ( -\rho
\ell_\star/2)^{\tau_1} [1/(\rho+2)+\ln \xi_0]
\ \ \mbox{if}\ f\ \mbox{obeys } {\rm (iii)}.
\end{array} \right.$$
\end{thm}

Note that the only case related, in same way, to our Theorem \ref{t2asym}
corresponds
to $\Omega_0=\emptyset$, $f(u)=u^{\rho+1}$ on $[0,\infty)$,
$k(t)=ct^{\alpha}\in {\cal K}$ (where $c,\alpha>0$), $\theta=1$ in
\eq{b},
being studied in \cite{gls}. There, the two-term asymptotic expansion
of $u_a$ near $\partial\Omega$
($a\in \RR$ since $\lambda_{\infty,1}=\infty$)
involves both the distance function
$d(x)$ and the mean curvature $H$ of $\partial\Omega$.
However, the blow-up rate of $u_a$ we present in Theorem \ref{t2asym} is of a
different
nature since the class ${\cal R}_{0,\zeta}$ does not include
$k(t)=ct^\alpha$.

Our main result contributes to the knowledge in some new directions. More
precisely,
the blow-up rate of the unique
positive solution $u_a$ of \eq{sep} is refined as follows in the above result:

(a) on the maximal interval
$(-\infty,\lambda_{\infty,1})$
for the parameter $a$, which is in
connection with an appropriate semilinear
eigenvalue
problem; thus,
the condition $b>0$ in $\Omega$
is removed by defining the set
$\Omega_0$, but we maintain $b\equiv 0$ on $\partial\Omega$
since this is a {\em natural} restriction inherited from the logistic
problem.

(b) when $b$ satisfies \eq{b}, where $\theta$ is {\em any}
positive number and $k$ belongs to a very rich class of
functions, namely ${\cal R}_{0,\zeta}$ . The equivalence
${\cal R}_{0,\zeta}\equiv {\cal K}_{0,\zeta}$ shows the connection to
the larger class ${\cal K}$  for which the uniqueness of
$u_a$ holds.
In addition, the explicit form of $k\in {\cal R}_{0,\zeta}$
shows us how to built $k\in {\cal K}_{0,\zeta}$.

(c) for a wide class of
functions
$f\in NRV_{\rho+1}$ where either $\phi\equiv 0$ (case (i))
or $\phi$ (resp., $-\phi$) belongs to $RV_{\eta}$ with $\eta\in
(-\rho-2,0]$  (cases (ii) and (iii)).
Therefore,
the theory of regular variation
plays a key role in understanding the general framework and the approach
as well.

\medskip
{\it Proof of Theorem \ref{t2asym}}.
We first state two auxiliary results. Their proofs are straightforward and we shall omit them.

\begin{lem}\label{lema1asym}
Assume \eq{b} and $f\in NRV_{\rho+1}$ satisfies
$(A_1)$. Then $h$ has the
following properties:

{\rm (i)} $h\in
C^2(0,\nu)$, $\lim_{t\searrow 0}h(t)=\infty$ and $\lim_{t\searrow
0}h'(t)=-\infty$;

{\rm (ii)} $\lim_{t\searrow 0}h''(t)/[k^2(t)f(h(t)\xi)]=
(2+\rho\ell_1)/[\xi^{\rho+1}(2+\rho)]$, $\forall \xi>0$;

{\rm (iii)} $\lim_{t\searrow 0}h(t)/h''(t)= \lim_{t\searrow
0}h'(t)/h''(t)=\lim_{t\searrow 0} h(t)/h'(t)=0$;

{\rm (iv)} $\lim_{t\searrow 0}h'(t)/[th''(t)]=-\rho
\ell_1/(2+\rho\ell_1)$ and
$\lim_{t\searrow 0}h(t)/[t^2 h''(t)]=\rho^2
\ell_1^2/[2(2+\rho \ell_1)]$;

{\rm (v)} $\lim_{t\searrow 0}h(t)/[th'(t)]=\lim_{t\searrow 0}[\ln
t]/[\ln h(t)]
=-\rho\ell_1/2$;

{\rm (vi)} If $\ell_1=0$, then
$\lim_{t\searrow 0}t^j h(t)=\infty,$ for all $j>0$;

{\rm (vii)} $ \lim_{t\searrow 0} 1/[t^\zeta \ln h(t)]=
-\rho \ell_\star/2$ and
$\lim_{t\searrow 0}h'(t)/[t^{\zeta+1}h''(t)]=
\rho \ell_\star/(2\zeta)$, $ \forall k\in {\cal R}_{0,\zeta}$.
\end{lem}

Let $\tau>0$ be arbitrary. For any $u>0$, define
$ T_{1,\tau}(u)=\{
\rho/[2(\rho+2)]-\Xi(u)\}(\ln u)^\tau $
and $T_{2,\tau}(u)=
\{f(\xi_0u)/[\xi_0 f(u)]-\xi_0^\rho\}(\ln u)^\tau $.
Note that if $f(u)=C u^{\rho+1}$, for $u$ in a neighbourhood $V_\infty$
of infinity, then $T_{1,\tau}(u)=T_{2,\tau}(u)=0$ for each $u\in
V_\infty$.

\begin{lem}\label{lema2asym}
Assume $(A_1)$ and $f\in {\cal F}_{\rho\eta}$. The following hold:

{\rm (i)} If $f\in {\cal F}_{\rho 0,\tau}$, then
$\lim_{u\to \infty}T_{1,\tau}(u)=-\ell^\star/(\rho+2)^2 $ and
$\lim_{u\to \infty}T_{2,\tau}(u)=\xi_0^\rho\ell^\star \ln \xi_0$.

{\rm (ii)} If $f\in {\cal F}_{\rho\eta}$ with $\eta\not=0$, then
$\lim_{u\to \infty}T_{1,\tau}(u)=\lim_{u\to \infty}T_{2,\tau}(u)=0$.\end{lem}

Fix $\ep\in (0,1/2)$. We can find $\delta>0$ such that
$d(x)$ is of class $C^2$ on $\{x\in \RR^N:\
d(x)<\delta\}$,
$ k$ is nondecreasing on $(0,\delta)$, and
$h'(t)<0<h''(t)$
for all $t\in (0,\delta)$. A straightforward computation shows that
$\lim_{t\searrow 0}t^{1-\theta}k'(t)/k(t)=\infty$,
for every $\theta>0$.
Using now \eq{b}, it follows that
we can diminish $\delta>0$ such that
$k^2(t)\left[1+(\tilde
c-\ep)t^\theta\right]$ is increasing on $(0,\delta)$ and
\begin{equation}\label{3}
1+(\tilde c-\ep)d^\theta<b(x)/k^2(d)<
1+(\tilde c+ \ep)d^\theta, \quad\forall x\in
\Omega\ \mbox{with}\ d\in (0,\delta).
\end{equation}
We define
$u^\pm (x)=\xi_0h(d)(1+\chi_{\ep}^\pm d^\varpi)$,  with $d\in
(0,\delta)$, where
$ \chi_{\ep}^\pm=
\chi \pm \ep\,[1+{\rm
Heaviside}\,(\zeta-\theta)]/\rho.$
Take $\delta>0$ small enough such that $u^\pm(x)>0$, for each $x\in
\Omega$
with $d\in (0,\delta)$.
By the Lagrange mean value theorem, we obtain
$ f(u^\pm (x))=f(\xi_0 h(d))+\xi_0
\chi_{\ep}^\pm d^\varpi h(d)
f'(\Upsilon^ \pm(d)),
$
where $\Upsilon^\pm (d)=\xi_0
h(d)(1+\lambda^\pm(d) \chi_{\ep}^\pm d^\varpi)$, for some
$\lambda^\pm(d)\in [0,1]$. We claim that
\begin{equation}\label{5}
\lim_{d\searrow 0}f(\Upsilon^\pm(d))/f(\xi_0 h(d))=1.
\end{equation}
Fix $\sigma\in (0,1)$ and $M>0$ such that
$|\chi_{\ep}^\pm|<M$. Choose $\mu^\star>0$ so that
$ |(1\pm Mt)^{\rho+1}-1|<\sigma/2$, for all
$t\in (0,2\mu^\star).$ Let $\mu_\star\in
(0,(\mu^\star)^{1/\varpi})$ be such that, for every $x\in \Omega$
with $d\in (0,\mu_\star)$
$$
\left|f(\xi_0 h(d)(1\pm M\mu^\star))/f(\xi_0 h(d))
-(1\pm M\mu^\star)^{\rho+1}\right|<\sigma/2.
$$
Hence,
$ 1-\sigma<(1-M\mu^\star)^{\rho+1}-\sigma/2<
f(\Upsilon^\pm(d))/f(\xi_0
h(d))<(1+M\mu^\star)^{\rho+1}+\sigma/2<1+ \sigma $,
for every $x\in \Omega$ with $d\in
(0,\mu_\star)$. This proves \eq{5}.

{\em Step} 1. There exists $\delta_1\in (0,\delta)$ so that
$\Delta u^+ +au^{+}-k^2(d)[1+(\tilde c-\ep)d^\theta]f(u^+) \leq
0$, $\forall x\in \Omega$ with $d\in (0,\delta_1)$
and
$\Delta u^- +au^{-}-k^2(d)[1+(\tilde c+\ep)d^\theta]f(u^-) \geq
0$, $\forall x\in \Omega$ with $ d\in (0,\delta_1)$.

Indeed,  for every $x\in \Omega$ with
$d\in (0,\delta)$, we have
\begin{equation}\label{7}\begin{array}{lll}
&\Delta u^{\pm}+ au^{\pm}-k^2(d)\left[1+(\tilde c\mp\ep)d^\theta
\right]f(u^\pm)
\qquad \qquad \qquad & \\
   =\xi_0 d^\varpi h''(d)&\left[
a \chi_{\ep}^{\pm}\frac{h(d)}{h''(d)}+
\chi_{\ep}^{\pm}\Delta d\,\frac{h'(d)}{h''(d)}+ 2\varpi
\chi_{\ep}^{\pm}\,\frac{h'(d)}{dh''(d)}+
\varpi \chi_{\ep}^{\pm}\Delta d \,\frac{h(d)}{dh''(d)}
\right.&\nonumber \\
  &\left.
+\varpi(\varpi-1)\chi_{\ep}^{\pm}\,\frac{h(d)}{d^2 h''(d)}
+ \Delta d\,\frac{h'(d)}{d^\varpi h''(d)}
+\frac{a\,h(d)}{d^\varpi h''(d)}
+\sum_{j=1}^{4}{\cal S}_{j}^{\pm}(d)\right]& \nonumber
\end{array}
\end{equation}
where, for any $t\in (0,\delta)$, we denote
\begin{eqnarray*}
& & {\cal S}_{1}^{\pm}(t)=(-\tilde c\pm \ep)t^{\theta-\varpi}\,
k^2(t)f(\xi_0 h(t))/[\xi_0 h''(t)],\ \,
{\cal S}_{2}^{\pm}(t)=
\chi_{\ep}^{\pm}(1-k^2(t)h(t)f'(\Upsilon^{\pm}(t))
/h''(t)),\\
& & {\cal S}_{3}^{\pm}(t)=(-\tilde c\pm
\ep)\chi_{\ep}^{\pm}t^\theta\,
k^2(t)h(t)f'(\Upsilon^{\pm}(t))/h''(t),\ \,
{\cal
S}_{4}^{\pm}(t)=
t^{-\varpi}
(1-k^2(t)f(\xi_0 h(t))/[\xi_0 h''(t)]).
\end{eqnarray*}

By Lemma \ref{lema1asym} (ii), we find
$ \lim_{t\searrow 0}
k^2(t)f(\xi_0 h(t))[\xi_0 h''(t)]^{-1}=1,
$
which
yields $\lim_{t\searrow 0}{\cal S}_{1}^{\pm}(t)= (-\tilde c\pm
\ep){\rm Heaviside}\,(\zeta-\theta)$.
Using \eq{5}, we obtain
$ \lim_{t\searrow 0}k^2(t)h(t)f'(\Upsilon^{\pm}(t))/h''(t)
=\rho+1$. Hence, $\lim_{t\searrow 0}{\cal S}_{2}^\pm (t)=-\rho
\chi_{\ep}^\pm $ and $\lim_{t\searrow 0}{\cal S}_{3}^\pm (t)=0$.

Using the expression of $h''$, we derive $ {\cal S}_{4}^{\pm}(t)  =
\frac{k^2(t)f(h(t))}{h''(t)}\,\sum_{i=1}^{3}{\cal S}_{4,i}(t)$,
$\forall t\in (0,\delta)$, where we denote
$ {\cal S}_{4,1}(t)=
2\frac{\Xi(h(t))}{t^{\varpi}}
(\int_0^t k(s)\,ds/k(t))'$,  ${\cal
S}_{4,2}(t) = 2\frac{T_{1,\tau_1}(h(t))}
{[t^\zeta \ln h(t)]^{\tau_1}}$ and
${\cal S}_{4,3}(t)=
- \frac{T_{2,\tau_1}(h(t))}{[t^\zeta \ln h(t)]^{\tau_1}}$.

Since ${\cal R}_{0,\zeta}\equiv {\cal K}_{0,\zeta}$, we find
$ \lim_{t\searrow
0}{\cal S}_{4,1}(t)=-(1+\zeta)\rho \ell_\star
\zeta^{-1}(\rho+2)^{-1}\, {\rm
Heaviside}\,(\theta-\zeta)$.

{\em Cases} (i), (ii). By Lemma \ref{lema1asym} (vii) and Lemma \ref{lema2asym} (ii), we find
$ \lim_{t\searrow 0}S_{4,2}(t)=\lim_{t\searrow 0}S_{4,3}(t)=0.$
In view of Lemma \ref{lema1asym} (ii), we derive that
$ \lim_{t\searrow 0}S_4^\pm
(t)=-(1+\zeta)\rho\ell_\star(2\zeta)^{-1}\,
{\rm Heaviside}\,(\theta-\zeta)$.

{\em Case} (iii). By Lemma \ref{lema1asym} (vii) and Lemma \ref{lema2asym} (i),
$ \lim_{t\searrow 0}S_{4,2}(t)=-2\ell^\star(\rho+2)^{-2}
(-\rho\ell_\star/2)^{\tau_1}$ and
$\lim_{t\searrow 0}S_{4,3}(t)=-2\ell^\star(\rho+2)^{-1}
(-\rho\ell_\star/2)^{\tau_1}\ln \xi_0$.
Using Lemma \ref{lema1asym} (ii) once more, we arrive at
$ \lim_{t\searrow 0}S_{4}^\pm (t)=-(1+\zeta)\rho \ell_\star
(2\zeta)^{-1}\,{\rm
Heaviside}\,(\theta-\zeta)-\ell^\star
(-\rho\ell_\star/2)^{\tau_1}
[1/(\rho+2)+\ln \xi_0]$.

Note that in each of the cases (i)--(iii),
the definition of $\chi_{\ep}^\pm$
yields
$\lim_{t\searrow 0}\sum_{j=1}^4{\cal S}_j^+ (t)=-\ep<0$
and
$\lim_{t\searrow 0}\sum_{j=1}^4{\cal S}_j^- (t)=\ep>0$.
By Lemma \ref{lema1asym} (vii), $\lim_{t\searrow 0}\frac{h'(t)}{t^\varpi h''(t)}=0$.
But $\lim_{t\searrow 0}\frac{h(t)}{h'(t)}=0$, so
$\lim_{t\searrow 0}\frac{h(t)}{t^\varpi h''(t)}=0$.
Thus, using Lemma \ref{lema1asym} [(iii), (iv)], relation \eq{7} concludes
our Step 1.

{\em Step} 2. There exists $M^+$, $\delta^+>0$ such that
$u_a(x)\leq u^+(x)+M^+$,
for all $x\in \Omega$ with $0<d<\delta^+.$

Define $(0,\infty)\ni u\longmapsto \Psi_x(u)=au-b(x)f(u)$,
$\forall x$
with $d\in (0,\delta_1)$. Clearly, $\Psi_x(u)$ is decreasing when
$a\leq 0$.
Suppose $a\in (0,\lambda_{\infty,1})$.
Obviously, $f(t)/t:(0,\infty)\to (f'(0),\infty)$ is bijective. Let
$\delta_2\in (0,\delta_1)$ be such that $b(x)<1$, $\forall x$
with $d\in (0,\delta_2)$. Let $u_x$ define the unique positive solution
of $b(x)f(u)/u=a+f'(0)$, $\forall x$ with $d\in (0,\delta_2)$.
Hence, for any $x$ with $d\in (0,\delta_2)$, $u\to \Psi_x(u)$
is decreasing on $(u_x,\infty)$. But
$ \lim_{d(x)\searrow 0}\frac{b(x)f(u^+(x))}{u^+(x)}=
+\infty$ (use
$\lim_{d(x)\searrow 0}u^+(x)/h(d)=\xi_0$,
$(A_1)$ and Lemma \ref{lema1asym} [(ii) and (iii)]). So, for $\delta_2$ small enough,
$u^+(x)>u_x$, $\forall x$ with $d\in (0,\delta_2)$.

Fix $\sigma\in (0,\delta_2/4)$ and set ${\cal N}_\sigma:=
\{x\in \Omega:\ \sigma<d(x)<\delta_2/2\}$. We define
$u^*_\sigma(x)=u^+(d-\sigma,s)+M^+$, where $(d,s)$ are the local
coordinates
of $x\in {\cal N}_\sigma$. We choose $M^+>0$ large enough
to have
$u^*_\sigma(\delta_2/2,s)\geq u_a(\delta_2/2,s)$, $\forall \sigma
\in (0,\delta_2/4)$ and $\forall s\in \partial\Omega$.
Using \eq{3} and Step 1, we find
\begin{eqnarray*}
& & -\Delta u^*_\sigma (x)\geq
a u^+(d-\sigma,s)-[1+(\tilde c-\ep)(d-\sigma)^\theta]
k^2(d-\sigma)f(u^+(d-\sigma,s)) \\
& & \qquad \qquad \ \geq
a u^+(d-\sigma,s)-[1+(\tilde c-\ep)d^\theta]
k^2(d)f(u^+(d-\sigma,s))
\geq \Psi_x(u^+(d-\sigma,s))\\
& & \qquad \ \qquad \geq \Psi_x(u^*_\sigma)
= au^*_\sigma(x)-b(x)f(u^*_\sigma(x))\quad \mbox{in } {\cal
N}_{\sigma}.
\end{eqnarray*}
Thus, by the maximum principle, $u_a\leq u_\sigma^*$ in
${\cal N}_\sigma$, $\forall \sigma\in (0,\delta_2/4)$. Letting
$\sigma\to 0$, we have proved Step~2.

{\em Step} 3. There exists $M^-$, $\delta^->0$ such that
$ u_a(x)\geq u^-(x)-M^-$, for all $x\in \Omega$
with $0<d<\delta^-.$

For every $r\in (0,\delta)$, define $\Omega_r= \{x\in \Omega:\
0<d(x)<r\}$. We will prove that for $\lambda>0$ sufficiently small,
$\lambda u^-(x)\leq u_a(x)$, $\forall x\in \Omega_{\delta_2/4}$.
Indeed,
fix arbitrarily $\sigma\in (0,\delta_2/4)$.
Define $v^*_\sigma(x)=
\lambda u^-(d+\sigma,s)$, for $x=(d,s)\in \Omega_{\delta_2/2}$.
We choose $\lambda\in (0,1)$ small enough such that
$v_\sigma^*(\delta_2/4,s)\leq u_a(\delta_2/4,s)$, $\forall \sigma
\in (0,\delta_2/4)$, $\forall s\in \partial\Omega$.
Using \eq{3}, Step 1 and $(A_1)$, we find
\begin{eqnarray*}
& & \Delta v^*_\sigma (x)+av_\sigma^*(x)
\geq
\lambda k^2(d+\sigma)[1+(\tilde c+\ep)(d+\sigma)^\theta]
f(u^-(d+\sigma,s))\\
& & \qquad \qquad \quad \qquad \ \,
\geq k^2(d)[1+(\tilde c+\ep)d^\theta]f(\lambda u^-(d+\sigma,s))
\geq bf(v^*_\sigma),
\end{eqnarray*}
for all $x=(d,s)\in \Omega_{\delta_2/4}$, that is $v_\sigma^*$ is
a sub-solution of $\Delta u+au=b(x)f(u)$ in $\Omega_{\delta_2/4}$.
By the maximum principle, we conclude that $v_\sigma^*\leq u_a$ in
$\Omega_{\delta_2/4}$. Letting $\sigma\to 0$, we find $\lambda
u^-(x)\leq u_a(x)$, $\forall x\in \Omega_{\delta_2/4}$.

Since $\lim_{d\searrow 0}u^-(x)/h(d)=\xi_0$, by using $(A_1)$
and Lemma \ref{lema1asym} [(ii), (iii)], we can easily obtain
$\lim_{d\searrow 0}k^2(d)f(\lambda^2 u^-(x))/u^-(x)=\infty$.
So, there exists $\tilde \delta\in (0,\delta_2/4)$ such that
\begin{equation}\label{12}
  k^2(d)[1+(\tilde c+\ep)d^\theta]f(\lambda^2 u^-)/u^-
\geq \lambda^2|a|,\quad \forall x\in \Omega\ \mbox{with}\
0<d\leq \tilde\delta.
\end{equation}
By Lemma \ref{lema1asym} [(i) and (v)],
we deduce that $u^-(x)$ decreases with $d$ when $d\in (0,
\tilde \delta)$ (if necessary, $\tilde \delta>0$ is diminished).
Choose $\delta_*\in (0,\tilde
\delta)$, close enough to
$\tilde\delta$, such that
\begin{equation}\label{13}
h(\delta_*)(1+\chi_\ep^-\delta_*^\varpi)/
[h(\tilde\delta)(1+\chi_\ep^-\tilde
\delta^\varpi)]<1+\lambda.
\end{equation}
For each $\sigma\in (0,\tilde\delta-\delta_*)$, we define
$z_\sigma(x)=u^-(d+\sigma,s)-(1-\lambda)u^-(\delta_*,s)$. We prove
that $z_\sigma$ is a sub-solution of $\Delta u+au=b(x)f(u)$ in
$\Omega_{\delta_*}$. Using \eq{13}, $ z_\sigma(x) \geq
u^-(\tilde\delta,s)-(1-\lambda)u^-(\delta_*,s)>0$ $\forall
x=(d,s)\in \Omega_{\delta_*}$. By \eq{3} and Step 1, $z_\sigma$ is
a sub-solution of $\Delta u+au=b(x)f(u)$ in $\Omega_{\delta_*}$ if
\begin{equation}\label{14}
k^2(d +\sigma)[1+(\tilde c+\ep)(d+\sigma)^\theta]
\left[f(u^-(d+\sigma,s))-f(z_\sigma(d,s))\right]
\geq  a(1-\lambda)u^-(\delta_*,s),
\end{equation}
for all $(d,s)\in
\Omega_{\delta_*}$.
Applying the Lagrange mean
value theorem and $(A_1)$, we infer that \eq{14} is a consequence
of
$ k^2(d+\sigma)[1+(\tilde c+\ep)(d+\sigma)^\theta]
\,f(z_\sigma(d,s))/z_\sigma(d,s)
\geq |a|,\quad \forall (d,s)\in
\Omega_{\delta_*}$. This inequality holds
by virtue of
\eq{12}, \eq{13} and the decreasing character of $u^-$ with $d$.

On the other hand, $z_\sigma(\delta_*,s)\leq \lambda
u^-(\delta_*,s)\leq
u_a(x)$, $\forall x=(\delta_*,s)\in \Omega$.
Clearly, $\limsup_{d\to 0}(z_\sigma-u_a)(x)=-\infty$
and $b>0$ in $\Omega_{\delta_*}$. Thus, by the maximum principle,
$z_\sigma\leq u_a$ in $\Omega_{\delta_*}$,
$\forall\sigma\in (0,\tilde \delta-\delta_*)$. Letting $\sigma\to 0$,
we conclude the assertion of Step 3.

By Steps~2 and 3, $\chi_\ep^+\geq \{-1+u_a(x)/[\xi_0 h(d)]\}d^{-\varpi}
-M^+/[\xi_0 d^\varpi h(d)] $  $\forall x\in \Omega$
with $
d\in (0,\delta^+)$ and
$ \chi_\ep^-\leq \{-1+u_a(x)/[\xi_0 h(d)]\}d^{-\varpi}
+M^-/[\xi_0 d^\varpi h(d)]$ $ \forall x\in \Omega$
with $d\in (0,\delta^-)$.
Passing to the limit as $d\ri 0$ and using Lemma \ref{lema1asym} (vi), we
obtain
$  \chi_\ep^-\leq \liminf_{d\to 0}
\{-1+u_a(x)/[\xi_0 h(d)]\}d^{-\varpi}$ and
$\limsup_{d\to 0}
\{-1+u_a(x)/[\xi_0 h(d)]\}d^{-\varpi} \leq \chi_\ep^+ $.
Letting $\ep\to 0$, we conclude our proof.\qed

\section{Entire solutions blowing up at infinity of semilinear elliptic systems}
In this section we are concerned with the existence of solutions that blow up at infinity
for a class of semilinear elliptic systems defined on the whole space.

Consider the following semilinear elliptic system
    \neweq{s1}
    \left\{\begin{array}{lll}
    & \Delta u=p(x)g(v) &\quad \mbox{in}\ \RR^N,\\
    & \Delta v=q(x)f(u) &\quad \mbox{in}\ \RR^N,
    \end{array} \right.
    \endeq
    where $N\geq 3$ and $p,q\in C_{\rm loc}
    ^{0,\alpha}(\RR^N)$ ($0<\alpha<1$) are non-negative
    and radially
    symmetric functions.
    Throughout this paper we assume
    that $f,g\in C^{0,\beta}_{\rm loc}[0,\infty)$
    ($0<\beta<1$) are positive and non-decreasing on
    $(0,\infty)$.

    We are concerned here with the existence of
    positive {\it entire large solutions} of \eq{s1}, that
    is
    positive classical solutions which satisfy $u(x)\to
    \infty$ and $v(x)\to \infty$ as
    $|x|\to \infty$. Set $\RR^+=(0,\infty)$ and define

    \smallskip
    \noindent \quad $ {\cal G}=\{(a,b)\in \RR^+\times
    \RR^+;\
    (\exists)\mbox{ an entire radial solution of}\
    \eq{s1}\
    \mbox{so that}\ (u(0),v(0))=(a,b)\} $.

    \smallskip
    The case of pure powers in the non-linearities was
    treated
    by Lair and Shaker in \cite{aw2}. They proved that
    ${\cal G}=\RR^+\times\RR^+$ if $f(t)=t^{\gamma}$ and
    $g(t)=t^\theta$
    for $t\geq 0$ with $0<\gamma, \theta\leq 1$. Moreover,
    they established that
    all positive entire radial solutions of \eq{s1} are
    {\it large} provided that
    \neweq{s2} \int_{0}^{\infty}tp(t)\,dt=\infty,\qquad
    \int_{0}^{\infty}tq(t)\,dt=\infty. \endeq
    If, in turn
    \neweq{s3} \int_{0}^{\infty}tp(t)\,dt<\infty,\qquad
    \int_{0}^{\infty}tq(t)\,dt<\infty \endeq
    then all positive entire radial solutions of \eq{s1}
    are {\it bounded}.

    In what follows we generalize the above results to a
    larger class of systems. Theorems \ref{ts1} and \ref{teo2}
    are due to C\^{\i}rstea and R\u adulescu \cite{crjmpa}.

    \begin{thm}\label{ts1} Assume that
    \neweq{h3} \lim_{t\to \infty}
    \frac{g(cf(t))}{t}=0\quad \mbox{for all}\ c>0.\endeq
    Then ${\cal G}=\RR^+\times \RR^+$.
    Moreover, the following hold:

     (i) If $p$ and $q$ satisfy \eq{s2}, then
    all positive entire radial solutions of \eq{s1} are
    large.

    (ii) If $p$ and $q$ satisfy \eq{s3},
    then all positive entire radial solutions of \eq{s1}
    are bounded.
    Furthermore, if $f,g$ are locally Lipschitz continuous
    on $(0,\infty)$
    and $(u,v)$, $(\tilde u, \tilde v)$ denote two
    positive entire radial solutions
    of \eq{s1},
    then there exists a positive constant $C$ such that
    for all $r\in [0,\infty)$
    $$
     \max\,\{|u(r)-\tilde u(r)|,|v(r)-\tilde v(r)|\}
     \leq C\,\max\,\{|u(0)-\tilde u(0)|,
    |v(0)-\tilde v(0)|\}.$$
    \end{thm}

\proof
We start with the following auxiliary results.

\begin{lem}\label{eq} Condition \eq{s2} holds if and
    only if
    $\lim_{r\to \infty}A(r)=\lim_{r\to \infty}B(r)=\infty$
    where
    $$ A(r)\equiv
    \int_{0}^{r}t^{1-N}\int_{0}^{t}s^{N-1}p(s)\,ds\,dt,\qquad
    B(r)\equiv
    \int_{0}^{r}t^{1-N}\int_{0}^{t}s^{N-1}q(s)\,ds\,dt,
    \qquad \forall r>0.$$
    \end{lem}

  \proof Indeed, for any $r>0$
    \neweq{s4} A(r)=\frac{1}{N-2}\left[\int_0^r tp(t)\,dt-
    \frac{1}{r^{N-2}}\int_0^r t^{N-1}p(t)\,dt \right]
    \leq \frac{1}{N-2}\,\int_0^r tp(t)\,dt. \endeq
    On the other hand,
    $$ \begin{array}{lll}
    \di \int_0^r tp(t)\,dt-
    \frac{1}{r^{N-2}}\int_0^r t^{N-1}p(t)\,dt & \di=
    \frac{1}{r^{N-2}}\,\int_0^r\left(r^{N-2}-t^{N-2}\right)tp(t)\,dt\\
    & \di \geq

\frac{1}{r^{N-2}}\left[r^{N-2}-\left(\frac{r}{2}\right)^{N-2}\right]
    \int_{0}^{\frac{r}{2}}tp(t)\,dt.
    \end{array} $$
    This combined with \eq{s4} yields
    $$ \frac{1}{N-2}\,\int_0^r tp(t)\,dt\geq A(r)\geq
    \frac{1}{N-2}\left[1-\left(\frac{1}{2}\right)^{N-2}\right]
    \int_{0}^{\frac{r}{2}}tp(t)\,dt.$$
    Our conclusion follows now by letting $r\to \infty$.
    \qed

    \begin{lem}\label{dep} Assume that condition \eq{s3}
    holds. Let
    $f$ and $g$ be locally Lipschitz continuous functions
    on
    $(0,\infty)$. If $(u,v)$ and $(\tilde u, \tilde v)$
    denote two bounded positive entire radial solutions of
    \eq{s1},
    then there exists a positive constant $C$ such that
    for all $r\in [0,\infty)$
    $$
    \max\,\{ |u(r)-\tilde u(r)|,|v(r)-\tilde v(r)|\}
    \leq C\,\max\,\{|u(0)-\tilde u(0)|,
    |v(0)-\tilde v(0)|\}. $$
    \end{lem}

    \proof We first see that radial
    solutions of \eq{s1} are solutions of the ordinary
    differential equations system
    \neweq{s90} \left\{\begin{array}{lll}
    & \di u''(r)+\frac{N-1}{r}\,u'(r)=p(r)\,g(v(r)),\qquad
    r>0\\
    & \di v''(r)+\frac{N-1}{r}\,v'(r)=q(r)\,f(u(r)),\qquad
    r>0.
    \end{array} \right.\endeq
    Define $K=\max\,\{|u(0)-\tilde u(0)|,|v(0)-\tilde
    v(0)|\}$. Integrating
    the first equation of \eq{s90}, we get
    $$ u'(r)-\tilde
    u'(r)=r^{1-N}\int_{0}^{r}s^{N-1}p(s)(g(v(s))-
    g(\tilde v(s)))\,ds. $$
    Hence
    \neweq{s91} |u(r)-\tilde u(r)|\leq
    K+\int_{0}^{r}t^{1-N}
    \int_{0}^{t}s^{N-1}p(s)|g(v(s))-g(\tilde
    v(s))|\,ds\,dt.\endeq
    Since $(u,v)$ and $(\tilde u,\tilde v)$ are bounded
    entire radial solutions of \eq{s1}
    we have
    $$ \begin{array}{lll}
    & \di |g(v(r))-g(\tilde v(r))|\leq m|v(r)-\tilde
    v(r)|\quad
    &\mbox{for any}\ r\in
    [0,\infty) \\
    & \di |f(u(r))-f(\tilde u(r))|\leq m |u(r)-\tilde
    u(r)|\quad
    &\mbox{for any}\ r\in [0,\infty),
    \end{array} $$
    where $m$ denotes a Lipschitz constant for both
    functions $f$ and $g$.
    Therefore, using \eq{s91} we find
    \neweq{s92} |u(r)-\tilde u(r)|\leq K+
    m\int_{0}^{r}t^{1-N}\int_{0}^{t}s^{N-1}p(s)|v(s)-\tilde
    v(s)|\,ds\,dt.\endeq
    Arguing as above, but now with the second equation of
    \eq{s90}, we obtain
    \neweq{s93} |v(r)-\tilde v(r)|\leq K+
    m\int_{0}^{r}t^{1-N}\int_{0}^{t}s^{N-1}q(s)|u(s)-\tilde
    u(s)|\,ds\,dt.\endeq
    Define
    $$
    X(r)=K+m\int_{0}^{r}t^{1-N}\int_{0}^{t}s^{N-1}p(s)|v(s)-\tilde
    v(s)|
    \,ds\,dt.$$
    $$
    Y(r)=K+m\int_{0}^{r}t^{1-N}\int_{0}^{t}s^{N-1}q(s)|u(s)-\tilde
    u(s)|
    \,ds\,dt.$$
    It is clear that $X$ and $Y$ are non-decreasing
    functions with $X(0)=Y(0)=K$.
    By a simple calculation together with \eq{s92} and
    \eq{s93} we obtain
    \neweq{s94} \begin{array}{lll}
    & \di (r^{N-1}X')'(r)=mr^{N-1}p(r)|v(r)-\tilde
    v(r)|\leq mr^{N-1}p(r)Y(r)\\
    & \di (r^{N-1}Y')'(r)=mr^{N-1}q(r)|u(r)-\tilde
    u(r)|\leq mr^{N-1}q(r)X(r).
    \end{array}\endeq
    Since $Y$ is non-decreasing, we have
    \neweq{s95} X(r)\leq K+mY(r)A(r)\leq
    K+\frac{m}{N-2}\,Y(r)\int_{0}^{r}tp(t)\,dt
    \leq K+mC_{p}Y(r)\endeq
    where
    $C_p=\left(1/(N-2)\right)\int_{0}^{\infty}tp(t)\,dt$.
    Using \eq{s95}
    in the second inequality of \eq{s94} we find
    $$ (r^{N-1}Y')'(r)\leq mr^{N-1}q(r)(K+mC_pY(r)).$$
    Integrating twice this inequality from 0 to $r$, we
    obtain
    $$ Y(r)\leq
    K(1+mC_q)+\frac{m^2}{N-2}C_p\int_{0}^{r}tq(t)Y(t)\,dt,$$
    where
    $C_q=\left(1/(N-2)\right)\int_{0}^{\infty}tq(t)\,dt$.
    From
    Gronwall's inequality, we deduce
    $$ Y(r)\leq
    K(1+mC_q)e^{\frac{m^2}{N-2}C_p\int_{0}^{r}tq(t)\,dt}\leq
    K(1+mC_q)e^{m^2C_pC_q} $$
    and similarly for $X$. The conclusion follows now from
    the above
    inequality, \eq{s92} and \eq{s93}. \qed

    \medskip
    {\it Proof of Theorem {ts1} completed}.
Since the radial
    solutions of \eq{s1} are solutions of the ordinary
    differential equations system \eq{s90}
    it follows that the radial solutions of \eq{s1} with
    $u(0)=a>0$, $v(0)=b>0$
    satisfy
    \neweq{s5} u(r)=a+\int_{0}^{r}t^{1-N}\int_0^t
    s^{N-1}p(s)\,g(v(s))\,ds\,dt,
    \qquad r\geq 0. \endeq
    \neweq{s6} v(r)=b+\int_{0}^{r}t^{1-N}\int_0^t
    s^{N-1}q(s)\,f(u(s))\,ds\,dt,
    \qquad r\geq 0. \endeq
    Define $v_0(r)=b$ for all $r\geq 0$. Let $(u_k)_{k\geq
    1}$ and
    $(v_k)_{k\geq 1}$ be two sequences of functions given
    by
    $$ u_k(r)=a+\int_{0}^{r}t^{1-N}\int_0^t
    s^{N-1}p(s)\,g(v_{k-1}(s))\,ds\,dt,
    \qquad r\geq 0.$$
    $$ v_k(r)=b+\int_{0}^{r}t^{1-N}\int_0^t
    s^{N-1}q(s)\,f(u_k(s))\,ds\,dt,
    \qquad r\geq 0.$$
    Since $v_1(r)\geq b$, we find $u_2(r)\geq u_1(r)$
    for all $r\geq 0$. This implies $v_2(r)\geq v_1(r)$
    which further produces
    $u_3(r)\geq u_2(r)$ for all $r\geq 0$. Proceeding at
    the same manner
    we conclude that
    $$ u_k(r)\leq u_{k+1}(r)\quad \mbox{and}\quad
    v_k(r)\leq v_{k+1}(r),\qquad \forall r\geq 0\
    \mbox{and}\ k\geq 1.$$

    We now prove that the non-decreasing sequences
    $(u_k(r))
    _{k\geq 1}$ and $(v_k(r))_{k\geq 1}$ are bounded from
    above on bounded sets.
    Indeed, we have
    \neweq{s7} u_k(r)\leq u_{k+1}(r)\leq
    a+g(v_k(r))A(r),\qquad \forall r\geq 0\endeq
    and
    \neweq{s8} v_k(r)\leq
    b+f\left(u_k(r)\right)B(r),\qquad \forall r\geq
    0.\endeq
    Let $R>0$ be arbitrary. By \eq{s7} and \eq{s8} we find
    $$ u_k(R)\leq
    a+g\left(b+f\left(u_k(R)\right)B(R)\right)A(R), \qquad
    \forall k\geq 1 $$
    or, equivalently,
    \neweq{s9} 1\leq \frac{a}{u_k(R)}+
\frac{g\left(b+f\left(u_k(R)\right)B(R)\right)}{u_k(R)}\,A(R),\qquad
    \forall k\geq 1.\endeq
    By the monotonicity of $(u_k(R))_{k\geq 1}$,
    there exists $\lim_{k\to \infty}u_k(R):=L(R)$. We
    claim that $L(R)$ is finite.
    Assume the contrary. Then, by taking $k\to \infty$ in
    \eq{s9} and using
    \eq{h3} we obtain a contradiction.
    Since $u_{k}'(r)$, $v_k'(r)\geq 0$ we get that the map
    $(0,\infty)\ni R\to L(R)$ is non-decreasing on
    $(0,\infty)$ and
    \neweq{s10} u_k(r)\leq u_k(R)\leq L(R),\qquad \forall
    r\in [0,R],
    \ \forall k\geq 1.
    \endeq
    \neweq{s11} v_{k}(r)\leq
    b+f\left(L(R)\right)B(R),\qquad \forall
    r\in [0,R],\ \forall k\geq 1.\endeq
    It follows that there exists $\lim_{R\to
    \infty}L(R)=\ov{L}\in (0,\infty]$ and
    the sequences $(u_k(r))_{k\geq 1}$, $(v_k(r))_{k\geq
    1}$ are bounded
    above on bounded sets. Therefore, we can define
    $u(r):=\lim_{k\to \infty}
    u_k(r)$ and $v(r):=\lim_{k\to\infty}v_k(r)$ for all
    $r\geq 0$.
    By standard elliptic regularity theory we obtain that
    $(u,v)$ is a positive entire solution of \eq{s1} with
    $u(0)=a$ and $v(0)=b$.

    We now assume that, in addition, condition \eq{s3} is
    fulfilled. According
    to Lemma \ref{eq} we have that $\lim_{r\to
    \infty}A(r)=\ov{A}<\infty$ and
    $\lim_{r\to \infty}B(r)=\ov{B}<\infty$. Passing to the
    limit as $k\to \infty$
    in \eq{s9} we find
    $$ 1\leq
\frac{a}{L(R)}+\frac{g\left(b+f\left(L(R)\right)B(R)\right)}{L(R)}\,A(R)
    \leq
\frac{a}{L(R)}+\frac{g(b+f\left(L(R)\right)\ov{B})}{L(R)}\,\ov{A}.$$
    Letting $R\to \infty$ and using \eq{h3}
    we deduce $\ov{L}<\infty$. Thus, taking into account
    \eq{s10}
    and \eq{s11}, we obtain
    $$ u_k(r)\leq \ov{L}\quad \mbox{and}\quad
    v_k(r)\leq b+f(\ov{L})\ov{B},\qquad \forall r\geq 0,\
    \forall k\geq 1.$$
    So, we have found upper bounds for $(u_k(r))_{k\geq
    1}$
    and $(v_k(r))_{k\geq 1}$ which are independent of $r$.
    Thus, the solution
    $(u,v)$ is bounded from above. This shows that any
    solution of \eq{s5} and \eq{s6}
    will be bounded from above provided \eq{s3} holds.
    Thus,
    we can apply Lemma \ref{dep} to achieve the
    second assertion of $ii)$.

    Let us now drop the condition \eq{s3} and assume that
    \eq{s2} is fulfilled.
    In this case, Lemma \ref{eq} tells us that $\lim_{r\to
    \infty}A(r)=
    \lim_{r\to \infty}B(r)=\infty$. Let $(u,v)$ be an
    entire positive
    radial solution of \eq{s1}. Using \eq{s5} and \eq{s6}
    we obtain
    $$ u(r)\geq a+g(b)\,A(r),\qquad \forall r\geq 0.$$
    $$ v(r)\geq b+f(a)\,B(r),\qquad \forall r\geq 0.$$
    Taking $r\to\infty$ we get that $(u,v)$ is an entire
    large solution.
    This concludes the proof of Theorem~\ref{ts1}. \qed

\medskip
    If $f$ and $g$ satisfy the stronger regularity
    $f,g\in C^1[0,\infty)$, then we drop
    the assumption \eq{h3} and require, in turn,
    \smallskip

    \noindent $({\bf H_1})\quad f(0)=g(0)=0,\
    \liminf_{u\to \infty}
    \frac{f(u)}{g(u)}=:\sigma>0 $
    \smallskip

    \noindent and the Keller-Osserman condition
    \smallskip

    \noindent $({\bf H_2})\quad \di
    \int_{1}^{\infty}\frac{dt}{\sqrt{G(t)}}<\infty,\
    \mbox{where}\ G(t)=\int_{0}^{t}g(s)\,ds$.
    \smallskip

    Observe that assumptions $({\bf H_1})$ and $({\bf
    H_2})$
    imply that $f$ satisfies condition $({\bf H_2})$, too.

    Set $\eta=\min\,\{p,q\}$.

    Our main result in this case is

    \begin{thm}\label{teo2} Let $f,g\in C^1[0,\infty)$
    satisfy
    $({\bf H_1})$ and $({\bf H_2})$. Assume that \eq{s3} holds,
    $\eta$ is not identically zero at infinity and
    $\nu:=\max\,\{p(0),q(0)\}>0$.

    Then any entire radial solution $(u,v)$ of \eq{s1}
    with
    $(u(0),v(0))\in F({\cal G})$ is large.
    \end{thm}

\proof Under the assumptions of Theorem \ref{teo2} we prove
    the following auxiliary results.

    \begin{lem}\label{pr1} ${\cal G}\not=\emptyset$.
    \end{lem}

   \proof
    Cf. C\^{\i}rstea and R\u adulescu \cite{crjmpa},
    the problem
    $$\Delta \psi=(p+q)(x)(f+g)(\psi)\quad \mbox{in}\
    \RR^N,$$
    has a positive radial entire large solution. Since
    $\psi$ is radial,
    we have
    $$
    \psi(r)=\psi(0)+\int_{0}^{r}t^{1-N}\int_{0}^{t}s^{N-1}
    (p+q)(s)(f+g)(\psi(s))\,ds\,dt,\quad \forall r\geq
    0.$$
    We claim that $(0,\psi(0)]\times (0,\psi(0)]\subseteq
    {\cal G}$.
    To prove this, fix $0<a,b\leq\psi(0)$ and let
    $v_0(r)\equiv b$
    for all $r\geq 0$. Define the sequences $(u_k)_{k\geq
    1}$
    and $(v_k)_{k\geq 1}$ by
    \neweq{o} u_k(r)=
    a+\int_0^r t^{1-N}\int_0^t
    s^{N-1}p(s)g(v_{k-1}(s))\,ds\,dt,
    \quad \forall r\in[0,\infty),\quad \forall k\geq
    1,\endeq
    \neweq{p}
    v_{k}(r)=b+\int_0^r t^{1-N}\int_0^t
    s^{N-1}q(s)f(u_k(s))\,ds\,dt,
    \quad \forall r\in [0,\infty),\quad \forall k\geq
    1.\endeq
    We first see that $v_0\leq v_1$ which produces
    $u_1\leq u_2$.
    Consequently, $v_1\leq v_2$ which further yields
    $u_2\leq u_3$. With the
    same arguments, we obtain that $(u_k)$ and $(v_k)$
    are non-decreasing sequences. Since $\psi'(r)\geq 0$
    and
    $b=v_0\leq \psi(0)\leq \psi(r)$ for all $r\geq 0$ we
    find
    $$ \begin{array}{lll}
    u_1(r) &\di \leq a+
    \int_0^r t^{1-N}\int_0^t
    s^{N-1}p(s)g(\psi(s))\,ds\,dt\\
    & \di \leq \psi(0)+\int_0^r t^{1-N}
    \int_0^t
    s^{N-1}(p+q)(s)(f+g)(\psi(s))\,ds\,dt=\psi(r).
    \end{array} $$
    Thus $u_1\leq \psi$. It follows that
    $$ \begin{array}{lll}
    v_1(r) &\di \leq b+
    \int_0^r t^{1-N}\int_0^t
    s^{N-1}q(s)f(\psi(s))\,ds\,dt\\
    & \di \leq \psi(0)+\int_0^r t^{1-N}
    \int_0^t
    s^{N-1}(p+q)(s)(f+g)(\psi(s))\,ds\,dt=\psi(r).
    \end{array} $$
    Similar arguments show that
    $$ u_k(r)\leq \psi(r)\quad \mbox{and}\quad
    v_k(r)\leq \psi(r)\quad \forall r\in [0,\infty),\quad
    \forall k\geq 1.$$
    Thus, $(u_k)$ and $(v_k)$ converge and
    $(u,v)=\lim_{k\to\infty}(u_k,v_k)$
    is an entire radial solution of \eq{s1} such that
    $(u(0),v(0))=(a,b)$.
    This completes the proof. \qed

    An easy consequence of the above result is

    \begin{cor}\label{c3} If $(a,b)\in {\cal G}$, then
    $(0,a]\times
    (0,b]\subseteq {\cal G}$. \end{cor}
    \proof Indeed, the process used before can be repeated
    by taking
    $$ u_k(r)=
    a_0+\int_0^r t^{1-N}\int_0^t
    s^{N-1}p(s)g(v_{k-1}(s))\,ds\,dt,
    \quad \forall r\in[0,\infty),\quad \forall k\geq 1,$$
    $$ v_{k}(r)=b_0+\int_0^r t^{1-N}\int_0^t
    s^{N-1}q(s)f(u_k(s))\,ds\,dt,
    \quad \forall r\in [0,\infty),\quad \forall k\geq 1,$$
    where $0<a_0\leq a$, $0<b_0\leq b$ and $v_0(r)\equiv
    b_0$
    for all $r\geq 0$.

    Letting $(U,V)$ be the entire radial solution of
    \eq{s1} with central
    values $(a,b)$ we obtain as in Lemma~\ref{pr1},
    $$ u_k(r)\leq u_{k+1}(r)\leq U(r),\quad \forall r\in
    [0,\infty),\quad
    \forall k\geq 1,$$
    $$ v_k(r)\leq v_{k+1}(r)\leq V(r),\quad \forall r\in
    [0,\infty),
    \quad \forall k\geq 1.$$
    Set $(u,v)=\lim_{k\to \infty}(u_k,v_k)$. We see that
    $u\leq U$, $v\leq V$ on $[0,\infty)$ and $(u,v)$ is an
    entire
    radial solution of \eq{s1} with central values
    $(a_0,b_0)$.
    This shows that $(a_0,b_0)\in {\cal G}$, so that our
    assertion is proved.\qed

    \begin{lem}\label{pr2} ${\cal G}$ is
    bounded.\end{lem}

     \proof
    Set $0<\lambda <\min\,\{\sigma,1\}$ and let
    $\delta=\delta(\lambda)$
    be large enough so that
    \neweq{in} f(t)\geq \lambda g(t),\qquad \forall t\geq
    \delta.\endeq
    Since $\eta$ is radially symmetric and
    not identically zero at infinity, we can assume
    $\eta>0$ on $\partial B(0,R)$ for some $R>0$.
    Let $\zeta$ be
    a positive large solution $\zeta$
    of the problem
    $$ \Delta \zeta =\lambda
    \eta(x)g\left(\frac{\zeta}{2}\right)\qquad \mbox{in}\
    B(0,R).$$
    Arguing by contradiction, we assume that ${\cal
    G}$ is not bounded.
    Then, there exists $(a,b)\in {\cal G}$ such that
    $a+b>\max\,\{2\delta,\zeta(0)\}$. Let $(u,v)$ be the
    entire radial solution
    of \eq{s1} such that $(u(0),v(0))=(a,b)$. Since
    $u(x)+v(x)\geq
    a+b> 2\delta$ for all $x\in \RR^N$, by \eq{in}, we
    find
    $$ f(u(x))\geq f\left(\frac{u(x)+v(x)}{2}\right)\geq
    \lambda
    g\left(\frac{u(x)+v(x)}{2}\right) \qquad \mbox{if}\
    u(x)\geq v(x)$$
    and
    $$ g(v(x))\geq g\left(\frac{u(x)+v(x)}{2}\right)\geq
    \lambda
    g\left(\frac{u(x)+v(x)}{2}\right) \qquad \mbox{if}\
    v(x)\geq u(x).$$
    It follows that
    $$ \Delta (u+v)=p(x)g(v)+q(x)f(u)\geq
    \eta(x)(g(v)+f(u))\geq
    \lambda \eta(x)g\left(\frac{u+v}{2}\right)\qquad
    \mbox{in}\ \RR^N.$$
    On the other hand, $\zeta(x)\to \infty$ as $|x|\to R$
    and $u,v\in C^2(\overline{B(0,R)})$. Thus, by the
    maximum principle,
    we conclude that $u+v\leq \zeta$ in $B(0,R)$. But this
    is impossible
    since $u(0)+v(0)=a+b>\zeta(0)$.\qed

    \begin{lem}\label{pr4} $F({\cal G})\subset {\cal
    G}$. \end{lem}

    \proof
    Let $(a,b)\in F({\cal G})$. We claim that
    $(a-1/n_0,b-1/n_0)\in {\cal G}$
    provided $n_0\geq 1$ is large enough so that
    $\min\,\{a,b\}>1/n_0$.
    Indeed, if this is not true, by Corollary~\ref{c3}
    $$
    D:=\left[\left.a-\frac{1}{n_0},\infty\right.\right)\times
    \left[\left.b-\frac{1}{n_0},\infty\right.\right)\subseteq
    (\RR^+\times
    \RR^+)\setminus {\cal G}.$$
    So, we can find a small ball $B$ centered in $(a,b)$
    such that
    $B\subset\subset D$, i.e., $B\cap {\cal G}=\emptyset$.
    But this will
    contradict the choice of $(a,b)$. Consequently, there
    exists
    $(u_{n_0},v_{n_0})$ an entire radial solution of
    \eq{s1} such that
    $(u_{n_0}(0),
    v_{n_0}(0))=(a-1/n_0,b-1/n_0)$. Thus, for any $n\geq
    n_0$, we can define
    $$ u_n(r)=a-\frac{1}{n}+\int_0^r t^{1-N}\int_0^t
    s^{N-1}p(s)g(v_n(s))
    \,ds\,dt,\qquad
    r\geq 0,$$
    $$ v_n(r)=b-\frac{1}{n}+
    \int_0^r t^{1-N}\int_0^t
    s^{N-1}q(s)f(u_n(s))\,ds\,dt,\qquad r\geq 0.$$
    Using Corollary~\ref{c3} once more, we conclude that
    $(u_n)_{n\geq n_0}$
    and $(v_n)_{n\geq n_0}$ are non-decreasing sequences.
    We now prove that $(u_n)$ and $(v_n)$ converge on
    $\RR^N$. To this aim,
    let $x_0\in \RR^N$ be arbitrary. But $\eta$ is
    not identically zero at infinity so that, for some
    $R_0>0$, we have
    $\eta>0$ on
    $\partial B(0,R_0)$ and $x_0\in B(0,R_0)$.

    Since $\sigma=\liminf_{u\to
    \infty}\frac{f(u)}{g(u)}>0$, we find
    $\tau\in (0,1)$ such that
    $$ f(t)\geq \tau g(t),\qquad \forall t\geq
    \frac{a+b}{2}-\frac{1}{n_0}.$$
    Therefore, on the set where $u_n\geq v_n$, we have
    $$ f(u_n)\geq f\left(\frac{u_n+v_n}{2}\right)\geq \tau
    g\left(
    \frac{u_n+v_n}{2}\right).$$
    Similarly, on the set where $u_n\leq v_n$, we have
    $$ g(v_n)\geq g\left(\frac{u_n+v_n}{2}\right)\geq \tau
    g
    \left(\frac{u_n+v_n}{2}\right).$$
    It follows that, for any $x\in \RR^N$,
    $$ \Delta (u_n+v_n)=p(x)g(v_n)+q(x)f(u_n)\geq
    \eta(x)[g(v_n)+f(u_n)]\geq \tau
    \eta(x)g\left(\frac{u_n+v_n}{2}\right).$$
    On the other hand,  there exists
    a positive large
    solution of
    $$ \Delta \zeta=\tau
    \eta(x)g\left(\frac{\zeta}{2}\right)\qquad \mbox{in}\
    B(0,R_0).
    $$ The maximum principle yields $u_n+v_n\leq \zeta$ in
    $B(0,R_0)$.
    So, it makes sense to define $(u(x_0),v(x_0))=
    \lim_{n\to \infty}(u_n(x_0),v_n(x_0))$. Since $x_0$ is
    arbitrary,
    the functions $u$, $v$ exist on $\RR^N$. Hence $(u,v)$
    is an entire
    radial solution of \eq{s1} with central values
    $(a,b)$, i.e., $(a,b)\in
    {\cal G}$. \qed

    \medskip
    For $(c,d)\in (\RR^+\times\RR^+)\setminus
    {\cal G}$, define
    \neweq{rc} R_{c,d}=\sup\,\{r>0\,|\, \mbox{there exists
    a radial solution}\ \mbox{of}\ \eq{s1}\ \mbox{in}\
    B(0,r)\ \mbox{so that}\ (u(0),v(0))=(c,d)\}.
    \endeq

    \begin{lem}\label{pr5} If, in addition,
    $\nu=\max\,\{p(0),q(0)\}>0$,
    then $0<R_{c,d}<\infty$ where $R_{c,d}$ is defined by
    \eq{rc}.\end{lem}

     \proof Since $\nu>0$ and $p,q\in
    C[0,\infty)$, there exists
    $\epsilon>0$ such that
    $(p+q)(r)>0$ for all $0\leq r<\epsilon$. Let
    $0<R<\epsilon$
    be arbitrary. There exists a
    positive
    radial large solution of the problem
    $$ \Delta \psi_{R}=(p+q)(x)(f+g)(\psi_R)\qquad
    \mbox{in}\ B(0,R).$$
    Moreover, for any $0\leq r<R$,
    $$ \psi_R(r)=\psi_R(0)+
    \int_0^r t^{1-N}\int_0^t
    s^{N-1}(p+q)(s)(f+g)(\psi_R(s))\,ds\,dt.$$
    It is clear that $\psi_R'(r)\geq 0$. Thus, we find
    $$ \psi_R'(r)=r^{1-N}\int_0^r
    s^{N-1}(p+q)(s)(f+g)(\psi_R(s))\,ds\leq
    C(f+g)(\psi_R(r))$$
    where $C>0$ is a positive constant such that
    $\int_0^\epsilon
    (p+q)(s)\,ds\leq C$.

    Since $f+g$ satisfies $({\bf A_1})$ and
    $({\bf A_2})$, we may
     invoke Remark \ref{remarka} in Section~2 to conclude that
    $$ \int_1^\infty \frac{dt}{(f+g)(t)}<\infty.$$
    Therefore, we obtain
    $$
    -\frac{d}{dr}\int_{\psi_R(r)}^{\infty}\frac{ds}{(f+g)(s)}=
    \frac{\psi_R'(r)}{(f+g)(\psi_R(r))}\leq C \qquad
    \mbox{for any}\
    0<r<R.$$
    Integrating from $0$ to $R$ and recalling that
    $\psi_R(r)\to
    \infty$ as $r\nearrow R$, we obtain
    $$ \int_{\psi_R(0)}^{\infty}\frac{ds}{(f+g)(s)}\leq
    CR.$$
    Letting $R\searrow 0$ we conclude that
    $$ \lim_{R\searrow
    0}\int_{\psi_R(0)}^{\infty}\frac{ds}{(f+g)(s)}=0.$$
    This implies that $\psi_R(0)\to \infty$ as $R\searrow
    0$. So,
    there exists $0<\tilde R<\epsilon$
    such that $0<c,d\leq \psi_{\tilde R}(0)$. Set
    \neweq{d1} u_k(r)=c+\int_0^r t^{1-N}\int_0^t
    s^{N-1}p(s)g(v_{k-1}(s))
    \,ds\,dt,\qquad \forall r\in [0,\infty),\ \forall
    k\geq 1,\endeq
    \neweq{d2} v_k(r)=d+
    \int_0^r t^{1-N}\int_0^t
    s^{N-1}q(s)f(u_{k}(s))\,ds\,dt,
    \qquad \forall r\in [0,\infty),\ \forall k\geq
    1,\endeq
    where $v_0(r)=d$ for all $r\in [0,\infty)$. As in
    Lemma~\ref{pr1},
    we find that $(u_k)$ resp., $(v_k)$ are non-decreasing
    and
    $$ u_k(r)\leq \psi_{\tilde R}(r)\quad \mbox{and}\qquad
    v_k(r)\leq \psi_{\tilde R}(r),\qquad \forall r\in
    [0,\tilde R),
    \ \forall k\geq 1.$$
    Thus, for any $r\in [0,\tilde R)$, there exists
    $(u(r),v(r))=
    \lim_{k\to \infty}(u_k(r),v_k(r))$ which
    is, moreover, a radial solution of \eq{s1} in
    $B(0,\tilde R)$
    such that $(u(0),v(0))=(c,d)$. This shows that
    $R_{c,d}\geq \tilde R>0$.
    By the definition of $R_{c,d}$ we also derive
    \neweq{ls} \lim_{r\nearrow R_{c,d}}u(r)=\infty\quad
    \mbox{and}
    \quad  \lim_{r\nearrow R_{c,d}}v(r)=\infty.\endeq
    On the other hand, since $(c,d)\not\in {\cal G}$, we
    conclude that
    $R_{c,d}$ is finite. \qed
    \medskip

    {\it Proof of Theorem~\ref{teo2} completed.}
     Let $(a,b)\in F({\cal G})$ be arbitrary. By
    Lemma~\ref{pr4}, $(a,b)\in
    {\cal G}$ so that we can define $(U,V)$ an entire
    radial solution
    of \eq{s1} with $(U(0),V(0))=(a,b)$. Obviously, for
    any
    $n\geq 1$, $(a+1/n,b+1/n)\in (\RR^+\times
    \RR^+)\setminus
    {\cal G}$. By Lemma~\ref{pr5}, $R_{a+1/n,b+1/n}$ (in
    short, $R_n$)
    defined by \eq{rc} is a positive number. Let
    $(U_n,V_n)$ be the radial solution of \eq{s1} in
    $B(0,R_n)$
    with the central values $(a+1/n,b+1/n)$. Thus,
    \neweq{y} U_n(r)=a+\frac{1}{n}+
    \int_0^r t^{1-N}\int_0^t
    s^{N-1}p(s)g(V_n(s))\,ds\,dt,\qquad
    \forall r\in [0,R_n),\endeq
    \neweq{z} V_n(r)=b+\frac{1}{n}+
    \int_0^r t^{1-N}\int_0^t
    s^{N-1}q(s)f(U_n(s))\,ds\,dt,\qquad
    \forall r\in [0,R_n).\endeq
    In view of \eq{ls} we have
    $$\lim_{r\nearrow R_n}U_n(r)=\infty\quad
    \mbox{and}\quad
    \lim_{r\nearrow R_n}V_n(r)=\infty,\quad \forall n\geq
    1.$$
    We claim that $(R_n)_{n\geq 1}$ is a non-decreasing
    sequence. Indeed,
    if $(u_k)$, $(v_k)$ denote the sequences of functions
    defined by
    \eq{d1} and \eq{d2} with $c=a+1/(n+1)$ and
    $d=b+1/(n+1)$, then
    \neweq{di} u_k(r)\leq u_{k+1}(r)\leq U_n(r),\qquad
    v_k(r)\leq v_{k+1}(r)\leq V_n(r),\qquad \forall r\in
    [0,R_n),\
    \forall k\geq 1.\endeq
    This implies that $(u_k(r))_{k\geq 1}$ and
    $(v_k(r))_{k\geq 1}$
    converge for any $r\in [0,R_n)$. Moreover,
    $(U_{n+1},V_{n+1})=
    \lim_{k\to \infty}(u_k,v_k)$ is a radial solution
    of \eq{s1} in $B(0,R_n)$ with central values
    $(a+1/(n+1),b+1/(n+1))$.
    By the definition of $R_{n+1}$, it follows that
    $R_{n+1}\geq R_n$ for any $n\geq 1$.

    Set $R:=\lim_{n\to \infty}R_n$ and let $0\leq r<R$ be
    arbitrary. Then,
    there exists $n_1=n_1(r)$ such that $r<R_n$ for all
    $n\geq n_1$. From
    \eq{di} we see that $U_{n+1}\leq U_n$ (resp.,
    $V_{n+1}\leq V_n$)
    on $[0,R_n)$ for all $n\geq 1$. So, there exists
    $\lim_{n\to
    \infty}(U_n(r), V_{n}(r))$ which, by \eq{y} and
    \eq{z},
    is a radial solution of \eq{s1} in $B(0,R)$ with
    central values $(a,b)$.
    Consequently,
    \neweq{f} \lim_{n\to \infty} U_n(r)=U(r)\quad
    \mbox{and}\quad
    \lim_{n\to \infty}V_n(r)=V(r)\quad \mbox{for any}\
    r\in [0,R).\endeq
    Since $U_n'(r)\geq 0$, from \eq{z} we find
    $$ V_{n}(r)\leq b+\frac{1}{n}+
    f(U_n(r))\int_0^\infty t^{1-N}\int_0^t
    s^{N-1}q(s)\,ds\,dt.$$
    This yields
    \neweq{q1} V_n(r)\leq C_1 U_n(r)+C_2 f(U_n(r))\endeq
    where $C_1$ is an upper bound of
    $(V(0)+1/n)/(U(0)+1/n)$
    and $$C_2=\int_0^\infty t^{1-N}\int_0^t
    s^{N-1}q(s)\,ds\,dt
    \leq \frac{1}{N-2}\int_0^\infty sq(s)\,ds<\infty.$$
    Define $h(t)=g(C_1 t+C_2 f(t))$ for $t\geq 0$. It is
    easy to check that
    $h$ satisfies $({\bf A_1})$ and $({\bf A_2})$. Define
    $$ \Gamma(s)=\int_{s}^{\infty}\frac{dt}{h(t)},\qquad
    \mbox{for all}\ s>0.$$
    But $U_n$ verifies
    $$ \Delta U_n=p(x)g(V_n)$$ which combined with \eq{q1}
    implies
    $$ \Delta U_n\leq p(x)h(U_n).$$ A simple calculation
    shows that
    $$ \begin{array}{lll}
    \di \Delta \Gamma(U_n)&\di =\Gamma'(U_n)\Delta
    U_n+\Gamma''(U_n)|\nabla U_n|^2
    =\frac{-1}{h(U_n)}\Delta
    U_n+\frac{h'(U_n)}{[h(U_n)]^2}|\nabla U_n|^2\\
    & \di \geq
    \frac{-1}{h(U_n)}p(r)h(U_n)=-p(r)\end{array}$$
    which we rewrite as
    $$ \left(r^{N-1}\frac{d}{dr}\Gamma(U_n)\right)'\geq
    -r^{N-1}p(r)\qquad \mbox{for any}\ 0<r<R_n. $$
    Fix $0<r<R$. Then $r<R_n$ for all $n\geq n_1$ provided
    $n_1$ is large enough.
    Integrating the above inequality over $[0,r]$, we get
    $$ \frac{d}{dr}\Gamma(U_n)\geq -r^{1-N}\int_0^r
    s^{N-1}p(s)\,ds.$$
    Integrating this new inequality over $[r,R_n]$ we
    obtain
    $$ -\Gamma(U_n(r))\geq -\int_{r}^{R_n}t^{1-N}
    \int_0^t s^{N-1}p(s)\,ds\,dt,\qquad \forall n\geq n_1,
    $$
    since $U_n(r)\to \infty$ as $r\nearrow R_n$ implies
    $\Gamma(U_n(r))\to 0$
    as $r\nearrow R_n$. Therefore,
    $$ \Gamma(U_n(r))\leq \int_r^{R_n}t^{1-N}\int_0^t
    s^{N-1}p(s)\,ds\,dt,
    \qquad \forall n\geq n_1.$$
    Letting $n\to \infty$ and using \eq{f} we find
    $$ \Gamma(U(r))\leq \int_r^R t^{1-N}\int_0^t
    s^{N-1}p(s)\,ds\,dt,$$
    or, equivalently
    $$ U(r)\geq \Gamma^{-1}\left(\int_r^R t^{1-N}\int_0^t
    s^{N-1}p(s)
    \,ds\,dt\right).$$
    Passing to the limit as $r\nearrow R$ and using the
    fact that $\lim_{s\searrow
    0}\Gamma^{-1}(s)=\infty$ we deduce
    $$ \lim_{r\nearrow R}U(r)\geq
    \lim_{r\nearrow R}\Gamma^{-1}\left(\int_r^R
    t^{1-N}\int_0^t s^{N-1}p(s)
    \,ds\,dt\right)=\infty. $$
    But $(U,V)$ is an entire solution so that we conclude
    $R=\infty$
    and $\lim_{r\to \infty}U(r)=\infty$. Since \eq{s3}
    holds and $V'(r)\geq 0$ we find
    $$ \begin{array}{lll}
    U(r)& \di\leq a+g(V(r))\int_0^\infty t^{1-N}\int_0^t
    s^{N-1}p(s)\,ds\,dt\\
        & \di\leq a+g(V(r))\frac{1}{N-2}\int_0^\infty
    tp(t)\,dt, \qquad \forall r\geq 0.
    \end{array}$$
    We deduce $\lim_{r\to \infty}V(r)=\infty$, otherwise
    we obtain
    that $\lim_{r\to \infty}U(r)$ is finite, a
    contradiction. Consequently,
    $(U,V)$ is an entire large solution of \eq{s1}. This
    concludes our proof.\qed

    \section{Bifurcation problems for singular Lane-Emden-Fowler equations}
In this section we study the bifurcation  problem
\begin{equation} \tag{$P_{\,\la}$}
\left\{\begin{aligned} & -\Delta u=\la f(u)+a(x)g(u) &&
{\rm in}\ \Omega,\\
& u>0 && {\rm in}\ \Omega,\\
& u=0 && {\rm on}\ \partial\Omega,\\
\end{aligned} \right.
\end{equation}
where $\lambda\in \RR$ is a parameter and $\Omega\subset \RR^N$
($N\geq 2$) is a bounded domain with smooth boundary
$\partial\Omega$.
The main feature of this boundary value problem is the presence of the ``smooth" nonlinearity
$f$ combined with the ``singular" nonlinearity $g$. More exactly, we assume that
 $0<f\in C^{0,\beta}[0,\infty)$ and $0\leq
g\in C^{0,\beta}(0,\infty)$ $(0<\beta<1)$ fulfill the hypotheses
\begin{enumerate}
\item[($f1)$]\ \ $f$ is nondecreasing on $(0,\infty)$ while
$f(s)/s$ is nonincreasing for $s>0$; \item[($g1$)] \ \ $g$ is
nonincreasing on $(0,\infty)$ with $\lim_{s\searrow
0}g(s)=+\infty$; \item[($g2$)] \ \ there exists $C_0,\eta_0>0$ and
$\alpha\in(0,1)$ so that $g(s)\leq C_0s^{-\alpha}$, $\forall
s\in(0,\eta_0)$.
\end{enumerate}

The assumption ($g2$) implies
the following Keller-Osserman-type growth condition around the
origin
\neweq{KellerOser}\int_0^1\left(\int_0^tg(s)ds\right)^{-1/2}dt<+\infty.\endeq
As proved by B\'enilan, Brezis and Crandall in \cite{bbc},
condition \eq{KellerOser} is equivalent to the {\it property of
compact support}, that is, for any $h\in L^1(\RR^N)$ with compact
support, there exists a unique $u\in W^{1,1}(\RR^N)$ with compact
support such that $\Delta u\in L^1(\RR^N)$ and
$$-\Delta u+g(u)=h\qquad\mbox{a.e. in}\ \RR^N.$$

In many papers (see, e.g., Dalmasso \cite{dalma}, Kusano and
Swanson  \cite{kusano}) the potential $a(x)$ is assumed to depend
``almost'' radially on $x$, in the sense that $C_1\,p(|x|)\leq
a(x)\leq C_2\,p(|x|),$ where $C_1$, $C_2$ are positive constants
and $p(|x|)$ is a positive function satisfying some integrability
condition. We do not impose any growth assumption on $a$, but we
suppose throughout this paper that the variable potential $a(x)$
satisfies $a\in C^{0,\beta}(\overline{\Omega})$ and $a>0$ in
$\Omega$.

If $\lambda=0$ this equation is called the Lane-Emden-Fowler
equation and arises in the boundary-layer theory of viscous fluids
(see Wong \cite{wong}). Problems of this type, as well as the
associated evolution equations, describe naturally  certain
physical phenomena. For example, super-diffusivity equations of
this type have been proposed by de~Gennes \cite{gennes} as a model
for long range Van der Waals interactions in thin films spreading
on solid surfaces.

Our purpose is to study the effect of the asymptotically linear
perturbation $f(u)$ in $(P_\lambda)$, as well as to describe the
set of values of the positive parameter $\la$ such that problem
$(P_\la)$ admits a solution. In this case, we also prove a
uniqueness result. Due to the singular character of $(P_{\la})$,
we can not expect to find solutions in $C^2(\overline\Omega)$.
However, under the above assumptions we will show that $(P_{\la})$
has solutions in the class
$$\displaystyle{\mathcal E}:=\{\,u\in C^2(\Omega)\cap
C^{1,1-\alpha}({\overline{\Omega}});\,\,\Delta u\in
L^1(\Omega)\}.$$ We first observe that, in view of the assumption
$(f1)$, there exists
\[ m:=\lim_{s\ri\infty}\frac{f(s)}{s}\in[0,\infty).\]
This number plays a crucial role in our analysis. More precisely,
the existence of the solutions to $(P_{\la})$ will be separately
discussed for $m>0$ and $m=0.$  Let $\di
a_*=\min_{x\in\overline{\Omega}}a(x)$.

Theorems \ref{th1}--\ref{th2} have been established by
C\^{\i}rstea, Ghergu, and R\u adulescu \cite{cgr}.

\begin{thm} \label{th1} Assume $(f1),$ $(g1),$
$(g2)$ and $m=0.$ If $a_*>0$ (resp., $a_*=0$), then $(P_{\la})$
has a unique solution $u_{\la}\in{\mathcal E}$ for all $\la\in\RR$
(resp., $\la\geq 0$) with the properties:
\begin{enumerate}
\item[{\rm (i)}] $u_{\la}$ is strictly increasing with respect to
$\la$.
\item[{\rm (ii)}] there exist two positive constant $c_1,c_2>0$
depending on $\la$ such that $c_1d(x)\leq u_{\la}\leq c_2d(x)$ in
$\Omega.$
\end{enumerate}
\end{thm}

The bifurcation diagram in the ``sublinear" case $m=0$ is depicted
in Figure~\ref{eps:fig1}. \epspic{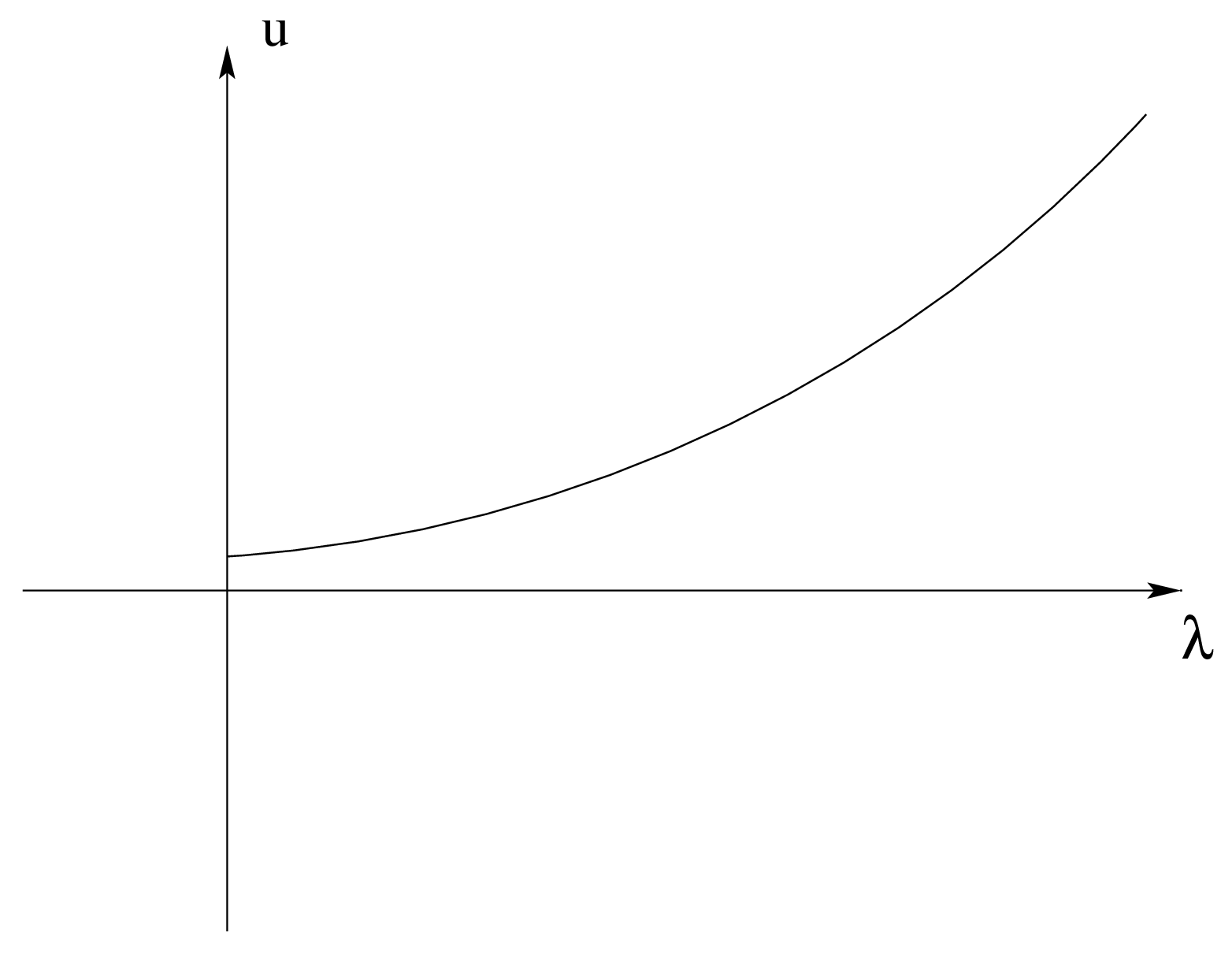}{6cm}{0}{The ``sublinear"
case $m=0$.}

\medskip
\proof
We first recall some auxiliary results that we need in the proof.

\begin{lem}\label{l1jmpa2} \emph{(Shi and Yao \cite{shi}).}
Let $\,F:\overline{\Omega}\times(0,\infty)\rightarrow\RR\,$ be a
H\"{o}lder continuous function with exponent $\beta\in(0,1)\,$ on
each compact subset of $\,\overline{\Omega}\times(0,\infty)\,$
which satisfies
\begin{enumerate}
\item[{\rm (F1)}] $\limsup _{s\rightarrow +\infty}
\left(s^{-1}\max_{x\in\overline{\Omega}}F(x,s)\right)<\la_1$;
\item[{\rm (F2)}] for each $\,t>0$, there exists a constant
$\,D(t)>0\,$ such that
$$F(x,r)-F(x,s)\geq -D(t)(r-s),\quad\mbox{for}\;\;x\in {\overline{\Omega}}\,
\;\;\mbox{and}\;\;\;r\geq s\geq t;$$ \item[{\rm (F3)}] there
exists  $\,\eta_0>0\,$ and an open subset
$\,\Omega_0\subset\Omega\,$ such that
$$\di\min_{x\in\overline{\Omega}}F(x,s)\geq
0\quad\mbox{for}\;s\in(0,\eta_0),$$ and
$$\di \lim_{s\searrow 0}\frac{F(x,s)}{s}=+\infty\quad\mbox{uniformly for }
x\in\Omega_0.$$  \end{enumerate} Then for any nonnegative function
$\,\phi_0\in C^{2,\beta}(\partial\Omega)$, the problem
$$\left\{\begin{aligned}
& -\Delta u=F(x,u) && {\rm in}\ \Omega,\,\\
& u> 0 && {\rm in}\ \Omega\,,\\
& u=\phi_0 && {\rm on}\ \partial\Omega,\,\\
\end{aligned} \right.$$
has at least one positive solution $\,u\in C^{2,\beta}(G)\cap
C(\overline{\Omega})$, for any compact set
$G\subset\Omega\cup\{x\in\partial\Omega;\,\phi_0(x)>0\}$.
\end{lem}

\begin{lem}\label{l2}\emph{(Shi and Yao \cite{shi}).}
Let $\,F:\overline{\Omega}\times(0,\infty)\rightarrow\RR\,$ be a
continuous function such that the mapping $\di(0,\infty)\ni
s\longmapsto\frac{F(x,s)}{s}\,$ is strictly decreasing at each
$\,x\in\Omega.$ Assume that there exists  $v$, $w\in
C^2(\Omega)\cap C({\overline{\Omega}})$ such that
\begin{enumerate}
\item[{\rm (a)}] $\Delta w+F(x,w)\leq 0\leq \Delta v+F(x,v)$ in $
\Omega$; \item[{\rm (b)}] $v,w>0$ in $\Omega$ and $v\leq w $ on
$\partial\Omega$; \item[{\rm (c)}] $\Delta v\in L^1(\Omega)$.
\end{enumerate}

Then $v\leq w$ in $\Omega$.\end{lem}

 Now, we are ready to give the proof of Theorem \ref{th1}.
This will be divided into four steps.

\smallskip {\sc Step 1.} {\bf Existence of  solutions to problem
$(P_{\la})$.}

For any $\la\in\RR,$ define the function
\neweq{bun}
\Phi_{\la}(x,s)=\la f(s)+a(x)g(s),
\quad(x,s)\in\overline\Omega\times (0,\infty).\endeq Taking into
account the assumptions of Theorem \ref{th1}, it follows that
$\Phi_\lambda$ verifies the hypotheses of Lemma~\ref{l1jmpa2} for
$\lambda\in \RR$ if $a_*>0$ and $\lambda\geq 0$ if $a_*=0$. Hence,
for $\lambda$ in the above range, $(P_{\la})$ has at least one
solution $u_{\la}\in C^{2,\beta}(\Omega)\cap
C(\overline{\Omega})$.

\smallskip {\sc Step 2.} {\bf Uniqueness of solution.}

Fix $\la\in\RR$ (resp., $\lambda\geq 0$) if $a_*>0$ (resp.,
$a_*=0$). Let $u_{\la}$ be a solution of $(P_{\la})$. Denote
$\la^{-}=\min\{0,\la\}$ and $\la^{+}=\max\{0,\la\}$. We claim that
$\Delta u_{\la}\in L^1(\Omega)$. Since $a\in
C^{0,\beta}(\overline{\Omega})$, by \cite[Theorem 6.14]{gt}, there
exists a unique nonnegative solution $\zeta\in
C^{2,\beta}(\overline{\Omega})$ of
$$\left\{\begin{aligned}
& -\Delta \zeta=a(x) && {\rm in}\ \Omega,\,\\
& \zeta=0 && {\rm on}\ \partial\Omega\,.\\
\end{aligned} \right.$$
By the weak maximum principle (see e.g., \cite[Theorem 2.2]{gt}),
$\zeta>0$ in $\Omega$. Moreover, we are going to prove that
\begin{enumerate}
\item[(a)] $z(x):=c\zeta(x)$ is a sub-solution of $(P_\la)$, for
$c>0$ small enough; \item[(b)] $z(x)\geq c_1 d(x)$ in
$\overline{\Omega}$, for some positive constant $c_1>0$;
\item[(c)] $u_\la\geq z$ in $\overline{\Omega}$.
\end{enumerate}

Therefore, by (b) and (c), $u_\la\geq c_1 d(x)$ in
$\overline{\Omega}$. Using $(g2)$, we obtain $g(u_\la)\leq C
d^{-\alpha}(x)$ in $\Omega$, where $C>0$ is a constant. So,
$g(u_\la)\in L^1(\Omega)$. This implies $$\Delta u_\la\in
L^1(\Omega). $$ {\em Proof of} (a). Using $(f1)$ and $(g1)$, we
have
$$\begin{aligned}
\Delta z(x)+\Phi_{\la}(x,z)& =-ca(x)+\la f(c\zeta)+
a(x)g(c\zeta)\\
& \geq -ca(x)+\la^-f(c\|\zeta\|_{\infty})+a(x)g(c\|\zeta\|_{\infty})\\
& \geq ca(x)\left[\frac{g(c\|\zeta\|_{\infty})}{2c}-1\right]
+f(c\|\zeta\|_{\infty})
\left[a_*\frac{g(c\|\zeta\|_{\infty})}{2f(c\|\zeta\|_{\infty})}+\la^-\right]\\
\end{aligned}$$
for each $x\in \Omega$. Since $\lambda<0$ corresponds to $a_*>0$,
using $\lim_{t\searrow 0}g(t)=+\infty$ and $\lim_{t\to 0}f(t)\in
(0,\infty)$, we can find $c>0$ small such that $$\Delta
z+\Phi_{\la}(x,z)\geq 0, \quad \forall x\in \Omega. $$ This
concludes (a).

{\em Proof of} (b). Since $\zeta\in
C^{2,\beta}(\overline{\Omega})$, $\zeta>0$ in $\Omega$ and
$\zeta=0$ on $\partial\Omega$, by Lemma 3.4 in Gilbarg and
Trudinger  \cite{gt}, we have
$$ \frac{\partial \zeta}{\partial \nu}(y)<0,\quad \forall y\in \partial\Omega.$$
Therefore, there exists a positive constant $c_0$ such that
$$ \frac{\partial \zeta}{\partial \nu}(y):=\lim_{x\in \Omega, x\to
y}\frac{\zeta(y)-\zeta(x)}{|x-y|}\leq -c_0, \quad \forall y\in
\partial\Omega. $$
So, for each $y\in \Omega$, there exists $r_y>0$  such that
\begin{equation} \label{hyt}
\frac{\zeta(x)}{|x-y|}\geq \frac{c_0}{2}, \quad \forall x\in
B_{r_y}(y)\cap \Omega. \end{equation} Using the compactness of
$\partial\Omega$, we can find a finite number $k$ of balls
$B_{r_{y_i}}(y_i)$ such that $\partial\Omega\subset \cup_{i=1}^k
B_{r_{y_i}}(y_i)$. Moreover, we can assume that for small $d_0>0$,
$$ \{x\in \Omega:\ \ d(x)<d_0\}\subset \cup_{i=1}^k B_{r_{y_i}}(y_i).
$$ Therefore, by (\ref{hyt}) we obtain
$$ \zeta(x)\geq \frac{c_0}{2}\,d(x),\quad \forall x\in \Omega
\ \mbox{with } d(x)<d_0.$$ This fact, combined with $\zeta>0$ in
$\Omega$, shows that for some constant $\tilde c>0$
$$ \zeta(x)\geq \tilde c d(x),\quad \forall x\in \Omega. $$
Thus, (b) follows by the definition of $z$.

{\em Proof of} (c). We distinguish two cases:

{\sc Case 1.} $\lambda\geq 0$.
We see that $\Phi_{\la}$ verifies
the hypotheses in Lemma~\ref{l2}. Since
$$ \begin{aligned}
& \Delta u_{\la}+\Phi_{\la}(x,u_{\la})\leq 0\leq \Delta
z+\Phi_{\la}(x,z) \quad\mbox{in}\;\;\Omega,\\
& \qquad \qquad u_{\la},\,z>0\quad\mbox{in}\;\;\Omega,\\
& \qquad \qquad u_{\la}=z\quad\mbox{on}\;\;\partial\Omega,\\
& \qquad \qquad \Delta z\in L^1(\Omega), \end{aligned} $$ by
Lemma~\ref{l2} it follows that $u_{\la}\geq z$ in
$\overline{\Omega}$.

Now, if $u_1$ and $u_2$ are two solutions of $(P_{\la})$, we can
use Lemma \ref{l2} in order to deduce that $u_1=u_2$.

{\sc Case 2.} $\lambda<0$ {\rm (corresponding to $a_*>0$).}
 Let $\ep>0$ be fixed. We prove that
\begin{equation} \label{supa}
z\leq u_\lambda+\ep (1+|x|^2)^{\tau} \quad \mbox{in } \Omega,
\end{equation}
where $\tau<0$ is chosen such that $\tau|x|^2+1>0$, $\forall x\in
\Omega$. This is always possible since $\Omega\subset \RR^N$
($N\geq 2$) is bounded.

We argue by contradiction. Suppose that there exists
$x_0\in\Omega$ such that $u_{\la}(x_0)+\ep(1+|x_0|)^\tau<z(x_0)$.
Then
$\min_{x\in\overline{\Omega}}\{u_{\la}(x)+\ep(1+|x|^2)^\tau-z(x)\}<0$
is achieved at some point $x_1\in\Omega$. Since $\Phi_{\la}(x,z)$
is nonincreasing in $z$, we have
$$\begin{aligned}
0&\geq -\Delta [u_{\la}(x)-z(x)+\ep(1+|x|^2)^\tau]|_{x=x_1}\\
&=\Phi_{\la}(x_1,u_{\la}(x_1))-\Phi_{\la}(x_1,z(x_1))-\ep
\Delta[(1+|x|^2)^\tau]|_{x=x_1}\\
& \geq -\ep \Delta[(1+|x|^2)^\tau]|_{x=x_1}
=-2\ep \tau (1+|x_1|^2)^{\tau-2}[(N+2\tau-2)|x_1|^2+N]\\
&\geq -4\ep\tau (1+|x_1|^2)^{\tau-2}(\tau|x_1|^2+1)>0.
\end{aligned} $$ This contradiction proves \eq{supa}. Passing to
the limit $\ep\to 0$, we obtain (c).

In a similar way we can prove that $(P_{\la})$ has a unique solution.

\smallskip {\sc Step 3.} {\bf Dependence on $\la$.}

We fix $\la_1<\la_2$, where $\la_1$, $\la_2\in\RR$ if $a_*>0$
resp., $\la_1$, $\la_2\in [0,\infty)$ if $a_*=0$. Let $u_{\la_1}$,
$u_{\la_2}$ be the corresponding solutions of $(P_{\la_1})$ and
$(P_{\la_2})$ respectively.

If $\la_1\geq 0$, then $\Phi_{\la_1}$ verifies the hypotheses in
Lemma~\ref{l2}. Furthermore, we have
$$ \begin{aligned}
& \Delta u_{\la_2}+\Phi_{\la_1}(x,u_{\la_2})\leq 0\leq \Delta
 u_{\la_1}+\Phi_{\la_1}(x,u_{\la_1}) \quad  \mbox{in } \Omega,\\
& \qquad \qquad u_{\la_1},u_{\la_2}>0 \quad \mbox{in } \Omega,\\
& \qquad \qquad u_{\la_1}=u_{\la_2} \quad \mbox{on } \partial\Omega,\\
& \qquad \qquad \Delta u_{\la_1}\in L^1(\Omega). \end{aligned}$$
Again by Lemma~\ref{l2}, we conclude that $\di u_{\la_1}\leq
u_{\la_2}$ in $\overline\Omega$. Moreover, by the maximum
principle,  $\di u_{\la_1}< u_{\la_2}$ in $\Omega$.

Let $\la_2\leq 0$; we show that $\di u_{\la_1}\leq u_{\la_2}$ in
$\overline\Omega$. Indeed, supposing the contrary, there exists
$x_0\in\Omega$ such that $\di u_{\la_1}(x_0)>u_{\la_2}(x_0)$. We
conclude now that $\di
\max_{x\in\overline{\Omega}}\{u_{\la_1}(x)-u_{\la_2}(x)\}>0$ is
achieved at some point in $\Omega$. At that point, say $\bar x$,
we have
$$\di 0\leq -\Delta (u_{\la_1}-u_{\la_2})(\bar x)=
\Phi_{\la_1}(\bar x,u_{\la_1}(\bar x))-\Phi_{\la_2}(\bar
x,u_{\la_2}(\bar x))<0,$$ which is a contradiction. It follows
that $u_{\la_1}\leq u_{\la_2}$ in $\overline{\Omega},$ and by
maximum principle we have $\di u_{\la_1}< u_{\la_2}$ in $\Omega$.

If $\la_1<0<\la_2$, then $u_{\la_1}<u_0<u_{\la_2}$ in $\Omega$.
This finishes the proof of Step~3.
\smallskip

 {\sc Step 4.} {\bf Regularity of the solution.}

Fix $\la\in\RR$ and let $\di u_{\la}\in C^2(\Omega)\cap
C(\overline\Omega)$ be the unique solution of  $(P_{\la}).$ An
important result in our approach is the following estimate
\begin{equation}\label{bunu}
\di c_1d(x)\leq u_{\la}(x)\leq c_2d(x),\quad\mbox{ for all
}\;x\in\Omega,\end{equation} where $c_1,c_2$ are positive
constants. The first inequality in \eq{bunu} was established in
Step 2. For the second one, we apply an idea found in Gui and Lin
\cite{gl}.

Using the smoothness of $\partial\Omega$, we can find
$\delta\in(0,1)$ such that for all
$x_0\in\Omega_{\delta}:=\{x\in\Omega\,;\,d(x)\leq \delta\}$, there
exists $y\in\RR^N\setminus\overline\Omega$ with
$d(y,\partial\Omega)=\delta$ and $d(x_0)=|x_0-y|-\delta$.

Let $K>1$ be such that diam$\,(\Omega)<(K-1)\delta$ and let $w$ be
the unique solution of the Dirichlet problem
\neweq{bdoi}
\left\{\begin{aligned} & -\Delta w=\la^+f(w)+g(w) &&
{\rm in}\ B_K(0)\setminus \overline{B_1}(0),\\
& w>0 && {\rm in}\ B_K(0)\setminus \overline{B_1}(0),\\
& w=0 && {\rm on}\ \partial(B_K(0)\setminus \overline{B_1}(0)),
\end{aligned} \right.
\end{equation}
where $B_r(0)$ is the open ball in $\RR^N$ of radius $r$ and
centered at the origin. By uniqueness, $w$ is radially symmetric.
Hence $w(x)=\tilde w(|x|)$ and
\neweq{btrei}
\left\{\begin{aligned}
 & \tilde w''+\frac{N-1}{r}\tilde w'+\la^+
f(\tilde w)+g(\tilde w)=0 && \mbox{for } r\in (1,K),\\
& \tilde w>0 && {\rm in}\ (1,K),\\
& \tilde w(1)=\tilde w(K)=0. &&
\end{aligned} \right.
\end{equation}
Integrating in \eq{btrei} we have
$$ \begin{aligned}
\tilde w'(t)& =\tilde w'(a)a^{N-1}t^{1-N}-t^{1-N}\int_a^t r^{N-1}
\left[\la^+f(\tilde w(r))+g(\tilde w(r))\right]dr,\\
&=\tilde w'(b)b^{N-1}t^{1-N}+t^{1-N}\int_t^b r^{N-1}
\left[\la^+f(\tilde w(r))+g(\tilde w(r))\right]dr,
\end{aligned} $$
where $1<a<t<b<K$. Since $g(\tilde w)\in L^1(1,K)$, we deduce that
both $\tilde w'(1)$ and $\tilde w'(K)$ are finite, so $\tilde w\in
C^2(1,K)\cap C^1[1,K]$. Furthermore,
\neweq{bpatru}
w(x)\leq C\min\{K-|x|,|x|-1\}, \quad\mbox{ for any } \;\;x\in
B_K(0)\setminus B_1(0). \endeq Let us fix $x_0\in\Omega_{\delta}.$
Then we can find $y_0\in\RR^N\setminus\overline\Omega$ with
$d(y_0,\partial\Omega)=\delta$ and $d(x_0)=|x_0-y|-\delta.$ Thus,
$\Omega\subset B_{K\delta}(y_0)\setminus B_{\delta}(y_0)$. Define
$\di v(x)=cw\left(\frac{x-y_0}{\delta}\right)$,
$x\in\overline\Omega$. We show that $v$ is a super-solution of
$(P_{\la})$, provided that $c$ is large enough. Indeed, if
$c>\max\{1,\delta^2\|a\|_{\infty}\}$, then for all $x\in\Omega$ we
have $$ \begin{aligned} \Delta v+\la f(v)+a(x)g(v) &
\leq\frac{c}{\delta^2}\left(\tilde w''(r)+ \frac{N-1}{r}\tilde
w'(r)\right)\\
& \quad +\la^+f(c\tilde w(r))+a(x)g(c\tilde w(r)), \end{aligned}
$$ where $ r=\di\frac{|x-y_0|}{\delta}\in(1,K)$. Using the
assumption $(f1)$ we get $f(c\tilde w)\leq cf(\tilde w)$ in
$(1,K)$. The above relations lead us to
$$\begin{aligned}
 \Delta v+\la f(v)+a(x)g(v)&
\leq \frac{c}{\delta^2}\left(\tilde w''+ \frac{N-1}{r}\tilde
w'\right)+\la^+ cf(\tilde w)+\|a\|_{\infty}g(\tilde w)\\
&\leq\frac{c}{\delta^2}\left(\tilde w''+
\frac{N-1}{r}\tilde w'+\la^+ f(\tilde w)+g(\tilde w)\right)\\
&=0.
\end{aligned}$$ Since $\Delta u_{\la}\in L^1(\Omega),$ with a
similar proof as in Step 2 we get $\di u_{\la}\leq v$ in $\Omega.$
This combined with \eq{bpatru} yields
$$ \di u_{\la}(x_0)\leq v(x_0)\leq
C\min\{K-\frac{|x_0-y_0|}{\delta},
\frac{|x_0-y_0|}{\delta}-1\}\leq\frac{C}{\delta}d(x_0).$$ Hence
$u_{\la}\leq \frac{C}{\delta}d(x)$ in $\Omega _{\delta}$ and the
last inequality in \eq{bunu} follows.

Let $G$ be the Green's function associated with the Laplace
operator in $\Omega$. Then, for all $x\in\Omega$ we have
$$\di u_{\la}(x)=-\int_{\Omega} G(x,y)\left[\la
f(u_{\la}(y))+a(y)g(u_{\la}(y))\right]dy,$$ and
$$\di \nabla u_{\la}(x)=-\int_{\Omega} G_x(x,y)\left[\la
f(u_{\la}(y))+a(y)g(u_{\la}(y))\right]dy.$$ If $x_1,x_2\in\Omega,$
using $(g2)$ we obtain
$$\begin{aligned}
|\nabla u_{\la}(x_1)-\nabla u_{\la}(x_2)|& \leq |\la|\int_{\Omega}
|G_x(x_1,y)-G_x(x_2,y)|\cdot
f(u_{\la}(y))dy\\
&\quad +\tilde c\int_{\Omega} |G_x(x_1,y)-G_x(x_2,y)|\cdot
u_{\la}^{-\alpha}(y)dy.
\end{aligned} $$
Now, taking into account that $u_{\la}\in C(\overline{\Omega})$,
by the standard regularity theory (see Gilbarg and Trudinger
\cite{gt}) we get
$$\di \int_{\Omega} |G_x(x_1,y)-G_x(x_2,y)|\cdot
f(u_{\la}(y))\leq \tilde c_1|x_1-x_2|.$$ On the other hand, with
the same proof as in \cite[Theorem 1]{gl}, we deduce
$$\di \int_{\Omega} |G_x(x_1,y)-G_x(x_2,y)|\cdot
u_{\la}^{-\alpha}(y)\leq \tilde c_2 |x_1-x_2|^{1-\alpha}.$$ The
above inequalities imply $u_{\la}\in C^2(\Omega)\cap
C^{1,1-\alpha}({\overline{\Omega}})$. The proof of
Theorem~\ref{th1} is now complete. \qed

\medskip
Next, consider the case $m>0.$ The results in this case are
different from those presented in Theorem \ref{th1}. A careful
examination of $(P_{\la})$ reveals the fact that the singular term
$g(u)$ is not significant. Actually, the conclusions are close to
those established in Mironescu and R\u adulescu \cite[Theorem A]{mr}, where an elliptic
problem associated to an asymptotically linear function is
studied.

Let $\la_1$ be the first Dirichlet eigenvalue of $(-\Delta)$ in
$\Omega$ and $\di \la^*=\frac{\la_1}{m}$. Our result in this case
is the following.

\begin{thm} \label{th2} Assume $(f1)$, $(g1)$, $(g2)$ and $m>0.$
Then the following hold.

\begin{enumerate}
\item[{\rm (i)}] If $\la\geq\la^*$, then $(P_{\la})$ has no
solutions in  ${\mathcal E}$. \item[{\rm (ii)}] If $a_*>0$ (resp.
$a_*= 0$) then $(P_{\la})$ has a unique solution
$u_{\la}\in{\mathcal E}$ for all $ -\infty<\la<\la^*$ (resp. $
0<\la< \la^*$) with the properties:

\noindent {\rm (ii1)} $u_{\la}$ is strictly increasing with
respect to $\la$;

\noindent {\rm (ii2)} there exists two positive constants $c_1,c_2>0$
depending on $\la$ such that $c_1d(x)\leq u_{\la}\leq c_2d(x)$ in
$\Omega$;

\noindent {\rm (ii3)} $ \lim\limits_{\la\nearrow
\la^*}u_{\la}=+\infty$, uniformly on compact subsets of $\Omega$.
\end{enumerate}
\end{thm}

The bifurcation diagram in the ``linear" case $m>0$ is depicted in
Figure~\ref{eps:fig2}.

\epspic{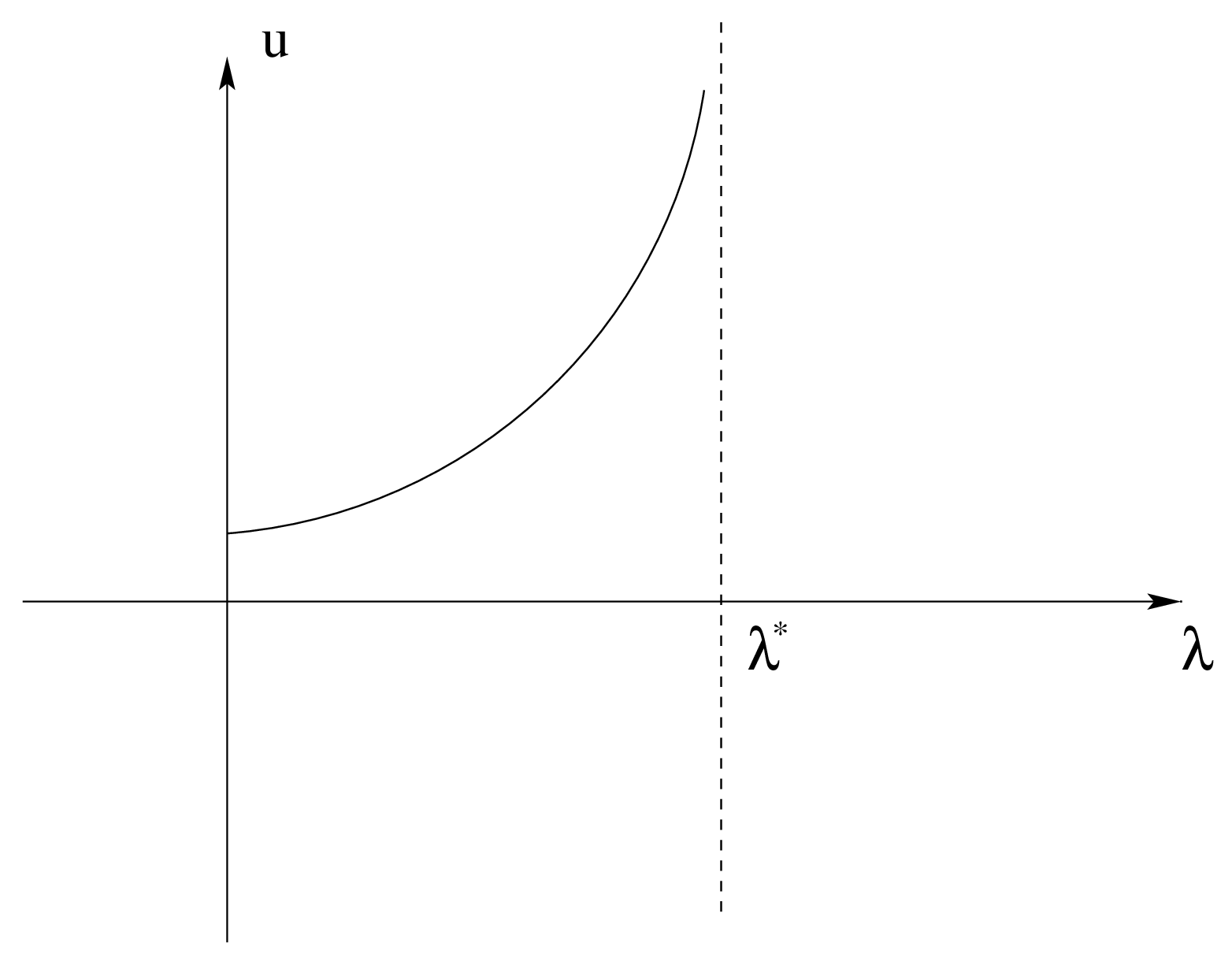}{6cm}{0}{The ``linear" case $m>0$.}

\medskip
\proof (i)  Let $\phi_1$ be the first eigenfunction of the Laplace
operator in $\Omega$ with  Dirichlet boundary condition. Arguing
by contradiction, let us suppose that there exists $\la\geq \la^*$
such that $(P_{\la})$ has a solution $u_{\la}\in {\mathcal E}$.

Multiplying by $\phi_1$ in $(P_{\la})$ and then integrating over
$\Omega$ we get \begin{equation} \label{cunu} \di
-\int_{\Omega}\phi_1\,\Delta
u_{\la}=\la\int_{\Omega}f(u_{\la})\phi_1+\int_{\Omega}
a(x)g(u_{\la})\phi_1 \end{equation} Since $\di \la\geq
\frac{\la_1}{m},$ in view of the assumption $(f1)$ we get $\la
f(u_{\la})\geq \la_1 u_{\la}$ in $\Omega$. Using this fact in
\eq{cunu} we obtain
$$-\int_{\Omega}\phi_1\,\Delta
u_{\la}>\la_1\int_{\Omega}u_{\la}\phi_1.$$ The regularity of
$u_{\la}$ yields $\di -\int_{\Omega}u_{\la}\Delta
\phi_1>\la_1\int_{\Omega}u_{\la}\phi_1$. This is clearly a
contradiction since $-\Delta \phi_1=\la_1\phi_1$ in $\Omega$.
Hence $(P_{\la})$ has no solutions in ${\mathcal E}$ for any
$\la\geq \la^*$.

(ii)  From now on, the proof of the existence, uniqueness and
regularity of solution is the same as in Theorem \ref{th1}.

(ii3) In what follows we shall apply some ideas developed in
Mironescu and R\u adulescu \cite{mr}. Due to the special character
of our problem, we will be able to prove that, in certain cases,
$L^2$--boundedness implies $H^1_0$--boundedness!

Let $u_{\la}\in{\mathcal E}$ be the unique solution of $(P_{\la})$
for $0<\la<\la^*$. We prove that $\di \lim_{\la\nearrow
\la^*}u_{\la}=+\infty$, uniformly on compact subsets of $\Omega$.
Suppose the contrary. Since $(u_{\la})_{0<\la<\la^*}$ is a
sequence of nonnegative super-harmonic functions in $\Omega,$ by
Theorem~4.1.9 in H\"ormander \cite{h}, there exists a subsequence
of $(u_{\la})_{\la<\la^*}$ (still denoted by
$(u_{\la})_{\la<\la^*}$ ) which is convergent in $L^1_{\rm
loc}(\Omega)$.

We first  prove that $(u_{\la})_{\la<\la^*}$ is bounded in
$L^2(\Omega)$. We argue by contradiction. Suppose that
$(u_{\la})_{\la<\la^*}$ is not bounded in $L^2(\Omega).$ Thus,
passing eventually at a subsequence we have
$u_{\la}=M(\la)w_{\la},$ where \begin{equation} \label{cdoi}
M(\la)=||u_{\la}||_{L^2(\Omega)}\ri\infty \mbox{ as }
\la\nearrow\la^* \mbox{  and  } w_{\la}\in L^2(\Omega),\;\;
\|w_{\la}\|_{L^2(\Omega)}=1.
\end{equation}
Using $(f1)$, $(g2)$ and the monotonicity assumption on $g$, we
deduce the existence of $A$, $B$, $C$, $D>0$ $(A>m)$ such that
\begin{equation} \label{ctrei}
f(t)\leq At+B,\quad g(t)\leq Ct^{-\alpha}+D,\quad \mbox{ for
all}\;t>0. \end{equation}
This implies
$$ \frac{1}{M(\la)}\left(\la
f(u_{\la})+a(x)g(u_{\la})\right)\ri 0 \quad \mbox{ in } L^1_{\rm
loc}(\Omega)\ \mbox{as }\la\nearrow \la^*$$ that is,
\begin{equation} \label{cpatru} -\Delta
w_{\la}\ri 0 \quad\mbox{in } L^1_{\rm loc}(\Omega)\ \mbox{as
}\la\nearrow \la^*.\end{equation} By Green's first identity, we
have
\begin{equation} \label{car1}
\int_\Omega \nabla w_\lambda\cdot \nabla \phi\,dx=-\int_\Omega
\phi\, \Delta w_\lambda\,dx=-\int_{ {\rm Supp}\, \phi} \phi\,
\Delta w_\lambda\,dx\quad \forall \phi\in C_0^\infty(\Omega).
\end{equation}
Using (\ref{cpatru}) we derive that
\begin{equation} \label{car2}
\begin{aligned}
\left|\int_{{\rm Supp}\,\phi} \phi\, \Delta w_\la\,dx\right| &
\leq
\int_{{\rm Supp}\,\phi}|\phi||\Delta w_\lambda|\,dx\\
& \leq \|\phi\|_{L^\infty}\int_{{\rm Supp}\,\phi}|\Delta
w_\lambda|\,dx\to 0\quad \mbox{as }\la\nearrow \la^*.
\end{aligned}
\end{equation}
Combining (\ref{car1}) and (\ref{car2}), we arrive at
\begin{equation} \label{car3}
\int_{\Omega} \nabla w_\lambda\cdot \nabla \phi\,dx \to 0 \ \
\mbox{as }\la\nearrow \la^*,\quad \forall \phi\in
C^\infty_0(\Omega). \end{equation} By definition, the sequence
$(w_{\la})_{0<\la<\la^*}$ is bounded in $L^2(\Omega)$.

We claim that $(w_{\la})_{\la<\la^*}$ is bounded in
$H^1_0(\Omega)$. Indeed, using \eq{ctrei} and H\"older's
inequality, we have
$$\begin{aligned}
\int_{\Omega}|\nabla w_{\la}|^2& =-\int_{\Omega}w_{\la}\Delta
w_{\la} =\frac{-1}{M(\la)}\int_{\Omega} w_{\la}
\Delta u_{\la}\\
& =\frac{1}{M(\la)}\int_{\Omega}\left[\la
w_{\la}f(u_{\la})+a(x)g(u_{\la})w_{\la}\right]\\
& \leq \frac{\la}{M(\la)}\int_{\Omega}
w_{\la}(Au_{\la}+B)+\frac{||a||_{\infty}}{M(\la)}
\int_{\Omega}w_{\la}(Cu^{-\alpha}_{\la}+D)\\
&=\la A\int_{\Omega}
w^2_{\la}+\frac{||a||_{\infty}C}{M(\la)^{1+\alpha}}
\int_{\Omega}w^{1-\alpha}_{\la}+\frac{\la
B+\|a\|_{\infty}D}{M(\la)}\int_{\Omega}w_{\la}\\
& \leq\la^*A+\frac{||a||_{\infty}C}{M(\la)^{1+\alpha}}
|\Omega|^{(1+\alpha)/2}+\frac{\la
B+\|a\|_{\infty}D}{M(\la)}|\Omega|^{1/2}.
\end{aligned} $$
From the above estimates, it is easy to see that
$(w_{\la})_{\la<\la^*}$ is bounded in $H^1_0(\Omega)$, so the
claim is proved. Then, there exists $w\in H^1_0(\Omega)$ such that
(up to a subsequence)
\neweq{ccinci}
w_{\la}\ \weak\ w \quad\mbox{ weakly in }\;\;H^1_0(\Omega) \
\mbox{ as } \la\nearrow \la^*
\end{equation} and, because $H_0^1(\Omega)$ is compactly embedded
in $L^2(\Omega)$,
\neweq{csase}
w_{\la}\ri w \quad\mbox{ strongly in }\;\;L^2(\Omega) \ \mbox{ as
} \la\nearrow \la^*.\end{equation} On the one hand, by
(\ref{cdoi}) and \eq{csase}, we derive that
$\|w\|_{L^2(\Omega)}=1$. Furthermore, using \eq{car3} and
\eq{ccinci}, we infer that
$$ \int_{\Omega} \nabla w\cdot \nabla \phi\,dx=0,\quad \forall
\phi\in C^\infty_0 (\Omega). $$ Since $w\in H_0^1(\Omega)$, using
the above relation and the definition of $H_0^1(\Omega)$, we get
$w=0$. This contradiction shows that $(u_{\la})_{\la<\la^*}$ is
bounded in $L^2(\Omega)$. As above for $w_\la$, we can derive that
$u_\la$ is bounded in $H_0^1(\Omega)$. So, there exists $u^*\in
H_0^1(\Omega)$ such that, up to a subsequence,
\neweq{car4}
\left\{\begin{aligned} & u_\la\ \weak\ u^*\ \ \mbox{weakly in }
H_0^1(\Omega)\
\ \mbox{as } \la\nearrow \la^*,\\
& u_\la \to u^*\ \ \mbox{strongly in } L^2(\Omega)\ \ \mbox{as }
\la\nearrow \la^*,\\
& u_\la\to u^* \ \ \mbox{a.e. in }\Omega \ \mbox{as }\la\nearrow
\la^* .
\end{aligned} \right.
\endeq

Now we can proceed to get a contradiction. Multiplying by $\phi_1$
in $(P_{\lambda})$ and integrating over $\Omega$ we have
\neweq{csapte}
\di -\int_{\Omega} \varphi_1\, \Delta
u_{\la}=\la\int_{\Omega}f(u_{\la}) \varphi_1+\int_{\Omega}
a(x)g(u_{\la})\varphi_1,\quad\mbox{ for all }\;0<\la<\la^*.
\end{equation}
On the other hand, by $(f1)$ it follows that $f(u_{\la})\geq
mu_{\la}$ in $\Omega,$ for all $0<\la<\la^*$. Combining this with
\eq{csapte} we obtain
\neweq{copt}
\di \la_1\int_{\Omega}u_{\la}\varphi_1\geq \la
m\int_{\Omega}u_{\la}\varphi_1+\int_{\Omega}
a(x)g(u_{\la})\varphi_1,\quad\mbox{ for all }\;0<\la<\la^*.
\endeq
Notice that by $(g1)$, \eq{car4} and the monotonicity of $u_\la$
with respect to $\la$ we can apply the Lebesgue convergence
theorem to find
$$ \int_\Omega a(x)g(u_\la)\varphi_1\,dx\to \int_\Omega
a(x)g(u^*)\varphi_1\,dx\ \ \mbox{as } \la\nearrow \lambda_1. $$
Passing to the limit in \eq{copt} as $\la\nearrow\la^*,$ and using
\eq{car4}, we get
\begin{equation}
\di \la_1\int_{\Omega}u^*\varphi_1\geq
\la_1\int_{\Omega}u^*\varphi_1+\int_{\Omega} a(x)g(u^*)\varphi_1.
\end{equation}
Hence $\di \int_{\Omega}a(x)g(u^*)\varphi_1=0,$ which is a
contradiction. This fact shows that $\di \lim_{\la\nearrow
\la^*}u_{\la}=+\infty$, uniformly on compact subsets of $\Omega$.
This ends the proof. \qed

\section{Sublinear singular elliptic problems with two bifurcation parameters}
Let $\Omega$ be a smooth bounded domain in
$\RR^N$ ($N\geq 2$).
 In this section we study the
existence or the nonexistence of solutions to the
following boundary value problem
$$\left\{\begin{tabular}{ll}
$-\Delta u+K(x)g(u)=\la f(x,u)+\mu h(x)$ \quad & ${\rm
in}\ \Omega,$\\
$u>0$ \quad & ${\rm in}\ \Omega,$\\
$u=0$ \quad & ${\rm on}\ \partial\Omega.$\\
\end{tabular} \right. \eqno(P_{\,\la,\,\mu})$$
Here $K,h\in C^{0,\gamma}(\overline\Omega ),$ with
$\,h>0\,$ on $\,\Omega\,$ and
$\,\la,\,\mu\,$ are positive real numbers.
We suppose that
$\,f:\overline{\Omega}\times[0,\infty)\rightarrow[0,\infty)\,$
is a H\"{o}lder continuous function
which is positive on
$\,\overline{\Omega}\times(0,\infty)$.
We also assume that $\,f\,$ is nondecreasing with
respect to the second
variable and is sublinear, that is,

\medskip


\noindent $\di(f1)\qquad$ the mapping $\di
(0,\infty)\ni
s\longmapsto\frac{f(x,s)}{s}\quad\mbox{is
nonincreasing for all}\;\, x\in\overline{\Omega};$
\medskip

\noindent $\di(f2)\qquad \lim_{s\downarrow
0}\frac{f(x,s)}{s}=+\infty\quad\mbox{and}\;\;
\lim_{s\rightarrow\infty}\frac{f(x,s)}{s}=0,\;\;\mbox{uniformly
for}\;\,x\in\overline{\Omega}.$
\medskip

\noindent We assume that $\,g\in
C^{0,\gamma}(0,\infty)\,$ is a nonnegative and
nonincreasing function satisfying
\medskip

\noindent $\di(g1)\qquad \lim_{s\downarrow
0}g(s)=+\infty;$
\medskip

\noindent $\di (g2)\qquad \mbox{there exists}\;\;C,\delta_0>0$ and
$\,\alpha\in(0,1)\,$ such that $\,g(s)\leq Cs^{-\alpha}\;\,$ for
all $\;s\in(0,\delta_0).$

\medskip
Our framework includes the Emden--Fowler equation that corresponds
to $g(s)=s^{-\gamma}$, $\gamma>0$ (see Wong \cite{wong}).

Denote
$\,\di{\cal E}=\{\,u\in C^2(\Omega)\cap
C({\overline{\Omega}});\,\,g(u)\in L^1(\Omega)\}.$

We  show in this section  that $(P_{\,\la,\,\mu})$ has at least
one solution in $\,{\cal E}\,$ for $\,\la,\mu\,$ belonging to a
certain range. We also prove that in some cases
$\,(P_{\,\la,\,\mu})\,$ has no solutions in $\,{\cal E},$ provided
that $\la$ and $\mu$ are sufficiently small.

\begin{rem}\label{r1}{\em $(i)\;$ If $\,u\in{\cal E},$
$\,v\in C^2(\Omega)\cap C({\overline{\Omega}})\,$
and $\,0<u<v\,$ in $\,\Omega,$   then $\,v\in{\cal
E}.$\\
$\,(ii)\;$ Let $\,u\in C^2(\Omega)\cap
C({\overline{\Omega}})\,$ be a
solution of  $\,(P_{\,\la,\,\mu}).$ Then $\,u\in{\cal
E}\,$ if and only if
$\;\Delta u\in L^1(\Omega).$}\end{rem}

A fundamental role will be played in our analysis by
the numbers
$$K^*=\max_{x\in{\overline{\Omega}}}K(x),\;K_*=\min_{x\in{\overline{\Omega}}}K(x).$$

Our main results (see Ghergu and R\u adulescu \cite{gr2}) are the
following.

\begin{thm}\label{th1jde}
Assume that $\,K_*>0\,$ and $\,f\,$ satisfies
$\,(f1)-(f2)$.\\
If $\di \,\int^1_0 g(s)ds=+\infty,\,$ then
$\,(P_{\,\la,\,\mu})\,$ has no solution in $\,{\cal
E}$
for any $\,\la,\mu>0.$
\end{thm}

\begin{thm}\label{th2jde}
Assume that $\,K_*>0,$ $\,f\,$ satisfies
$\,(f1)-(f2)\,$ and $\,g\,$ satisfies
$\,(g1)-(g2).$\\
Then there exists $\,\la_*,\mu_*>0\,$ such that

$(P_{\,\la,\,\mu})\,$ has at least one solution in
$\,{\cal E}\,$ if
$\,\la>\la_*\,$ or $\,\mu>\mu_*.$

$(P_{\,\la,\,\mu})\,$ has no solution in $\,{\cal
E}\,$ if
$\;\la<\la_*\,$ and $\,\mu<\mu_*.$

Moreover, if $\,\la>\la_*\,$ or $\,\mu>\mu_*,\,$ then
$\,(P_{\,\la,\,\mu})\,$ has a
maximal solution in $\,{\cal E}\,$ which is increasing
with respect to
$\,\la\,$ and $\,\mu.$
\end{thm}

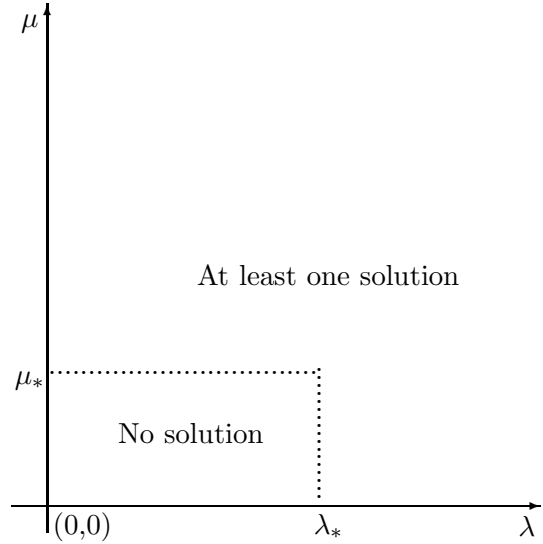
\begin{figure}[h]
\begin{center}
\begin{picture}(200,200)\setlength{\unitlength}{1pt}
\put(7,-15){\makebox(40,14){(0,0)}}
\put(175,-15){\makebox(40,15){$\lambda$}}
\put(100,-15){\makebox(40,15){$\lambda_*$}}
\vector(1,0){200}
\put(-190,-10){
\vector(0,20){200}}
\put(-213,175){\makebox(40,15){$\mu$}}
\put(-213,40){\makebox(40,15){$\mu_*$}}
\multiput(-186,50)(3,0){34}{.}
\multiput(-85,3)(0,3){17}{.}
\put(-152,20){\makebox(40,15){No solution}}
\put(-100,80){\makebox(40,15){At least one solution}}
\end{picture}
\end{center}
\caption{The dependence on $\,\lambda\,$ and $\,\mu$ in Theorem
\ref{th2jde}}
\end{figure}

\begin{thm}\label{th3jde}
Assume that $\,K^*\leq 0,$ $\,f\,$ satisfies
$\,(f1)-(f2)\,$ and $\,g\,$ satisfies
$\,(g1)-(g2).$\\
Then $\,(P_{\,\la,\,\mu})\,$ has a unique solution
$\,u_{\la,\mu}\in{\cal E}\,$
for any $\la,\, \mu>0.$
Moreover, $\,u_{\la,\mu}\,$ is increasing with respect
to
$\,\la\,$ and $\,\mu.$
\end{thm}

Theorems \ref{th2jde} and \ref{th3jde} also show the role played
by the sublinear term $f$ and the sign of $K(x)$. Indeed, if $f$
becomes linear then the situation changes radically. First, by the
results established by Crandall, Rabinowitz, and Tartar
\cite{crt}, the problem
$$
 \left\{\begin{tabular}{ll}
$-\Delta u-u^{-\alpha}= -u$ \quad & ${\rm in}\
\Omega,$\\
$u>0$ \quad & ${\rm in}\ \Omega,$\\
$u=0$ \quad & ${\rm on}\ \partial\Omega$\\
\end{tabular} \right.
$$
has a solution, for any $\alpha>0$. Next, as showed in Chen
\cite{chen}, the problem
$$
 \left\{\begin{tabular}{ll}
$-\Delta u+u^{-\alpha}= u$ \quad & ${\rm in}\
\Omega,$\\
$u>0$ \quad & ${\rm in}\ \Omega,$\\
$u=0$ \quad & ${\rm on}\ \partial\Omega$\\
\end{tabular} \right.
$$
has no solution, provided $0<\alpha <1$ and
$\lambda_1\geq 1$ (that is, if $\Omega$ is ``small''),
where
$\lambda_1$ denotes the first eigenvalue of
$(-\Delta)$ in $H^1_0(\Omega)$.

\begin{thm}\label{th4jde}
Assume that $K^*>0>K_*,$ $f$ satisfies
$(f1)-(f2)$ and $g$ verifies
$(g1)-(g2).$\\
Then there exists $\la_*,\, \mu_*>0$ such that
$(P_{\,\la,\,\mu})\,$ has at least one solution
$\,u_{\la,\mu}\in{\cal E}\,$ if
$\,\la>\la_*\,$ or $\,\mu>\mu_*.$
Moreover, for $\la>\la_*$ or $\mu>\mu_*,$
$u_{\la,\mu}$ is increasing with
respect to $\la$ and $\mu.$
\end{thm}

Before giving the proofs, we state some auxiliary results.

Let $\phi_1$ be the normalized positive eigenfunction
corresponding to the first eigenvalue $\la_1$ of the problem
\begin{equation}\label {doiunu}
 \left\{\begin{tabular}{ll}
$-\Delta u=\la u$ \quad & ${\rm in}\ \Omega,\,$\\
$u=0$ \quad & ${\rm on}\ \partial\Omega\,.$\\
\end{tabular} \right.
\end{equation}

\begin{lem}\label{l1jde}\emph{(Lazer and McKenna \cite{lm1}).}
$\di \int_{\Omega}\phi_1^{-s}dx<+\infty\,$ if and only
if $\,s<1.$
\end{lem}

Next, we observe that the hypotheses of Lemmas \ref{l1jmpa2} and
\ref{l2} are fulfilled for
\begin{equation}\label{doizero}
\Phi_{\la,\mu}(x,s)=\la f(x,s)+\mu h(x),\end{equation}
\begin{equation}\label{doizerozero}
\Psi_{\la,\mu}(x,s)=\la f(x,s)-K(x)g(s)+\mu
h(x),\quad\mbox{provided}\;\;K^*\leq 0.
\end{equation}

\begin{lem}\label{l4} Let $\,f\,$ satisfying
$\,(f1)-(f2)\,$ and $\,g\,$ satisfying
$\,(g1)-(g2).$ Then there exists
$\,\overline{\la}>0\,$ such that the problem
 \begin{equation}\label {doitrei}
 \left\{\begin{tabular}{ll}
$-\Delta v+g(v)=\la f(x,v)+\mu h(x)$ \quad & ${\rm
in}\ \Omega,\,$\\
$v>0$ \quad & ${\rm in}\ \Omega,\,$\\
$v=0$ \quad & ${\rm on}\ \partial\Omega.$\\
\end{tabular} \right.
\end{equation}
has at least one solution $\,v_{\la,\mu}\in{\cal E}\,$
for all $\,\la>\overline{\la}\,$
and for any  $\mu>0.$
\end{lem}

\proof Let $\,\la,\,\mu>0$. According to Lemmas \ref{l1jmpa2} and
\ref{l2}, the boundary value problem
\begin{equation}\label {doipatru}
 \left\{\begin{tabular}{ll}
$-\Delta U=\la f(x,U)+\mu h(x)$ \quad & ${\rm in}\
\Omega,\,$\\
$U>0$ \quad & ${\rm in}\ \Omega,\,$\\
$U=0$ \quad & ${\rm on}\ \partial\Omega$\\
\end{tabular} \right.
\end{equation}
has a unique solution $\,U_{\la,\mu}\in C^{2,\gamma}(\Omega)\cap
C(\overline{\Omega}).$ Then
$\,\overline{v}_{\la,\mu}=U_{\la,\mu}\,$ is a super-solution of
(\ref{doitrei}). The main point is to find a sub-solution of
(\ref{doitrei}). For this purpose, let
$\,H:[0,\infty)\rightarrow[0,\infty)\,$ be such that
\begin{equation}\label{doicinci}
\left\{\begin{tabular}{ll}
$H''(t)=g(H(t)),\quad \mbox{ for all }t>0,\,$\\
$H'(0)=H(0)=0.$\\
\end{tabular} \right.
\end{equation}
Obviously, $H\in C^2(0,\infty)\cap C^1[0,\infty)$
exists by our assumption
$(g2)$. From (\ref{doicinci}) it
follows that $\,H''\,$ is
nonincreasing, while $H$ and $H'$ are nondecreasing on
$\,(0,\infty).$ Using this fact and applying the mean
value
 theorem, we deduce that for all $\,t>0\,$ there exists
$\,\xi^1_t,\,\xi^2_t\in(0,t)\,$
such that
$$\di
\frac{H(t)}{t}=\frac{H(t)-H(0)}{t-0}=H'(\xi^1_t)\leq
H'(t);$$
$$\di
\frac{H'(t)}{t}=\frac{H'(t)-H'(0)}{t-0}=H''(\xi^2_t)\geq
H''(t).$$
The above inequalities imply
$$\di H(t)\leq tH'(t)\leq 2H(t),\qquad\mbox{for
all}\;\;t>0.$$
Hence
\begin{equation}\label{doisase}
\di 1\leq\frac{tH'(t)}{H(t)}\leq 2,\qquad\mbox{for
all}\;\;t>0.
\end{equation}

On the other hand, by $\,(g2)\,$ and
(\ref{doicinci}), there exists
$\,\eta>0\,$ such that
\begin{equation}\label{doisapte}
\left\{\begin{tabular}{ll}
$H(t)\leq\delta_0,\quad\mbox{for
all}\;\;t\in(0,\eta),$\\
$H''(t)\leq CH^{-\alpha}(t),\quad\mbox{for
all}\;\;t\in(0,\eta),$\\
\end{tabular} \right.
\end{equation}
which yields
\begin{equation}\label{doiopt}
\di H(t)\leq c\,t^{2/(\alpha+1)},
\quad\mbox{for all}\;\;t\in(0,\eta),
\end{equation}
where $\,c>0\,$ is a constant.

Now we look for a sub-solution of the form
$\,\underline{v}_{\la,\mu}=MH(\phi_1),$ for some constant $\,M>0.$
We have
\begin{equation}\label{doinoua}
\di -\Delta
\underline{v}_{\la,\mu}+g(\underline{v}_{\la,\mu})
=\la_1MH'(\phi_1)\phi_1+g(MH(\phi_1))-Mg(H(\phi_1))|\nabla\phi_1|^2\quad\mbox{in}\;\;
\Omega.
\end{equation}
Take $\,M\geq 1.$ The monotonicity  of
$\,g\,$ leads to
$$\di g(MH(\phi_1))\leq
g(H(\phi_1))\quad\mbox{in}\;\;\Omega,$$
and, by (\ref{doinoua}),
\begin{equation}\label{doizece}
\di -\Delta
\underline{v}_{\la,\mu}+g(\underline{v}_{\la,\mu})
\leq\la_1MH'(\phi_1)\phi_1+g(H(\phi_1))\left(1-M|\nabla\phi_1|^2\right)\quad\mbox{in}\;\;
\Omega.
\end{equation}
We claim that
\begin{equation}\label{doiunspe}
\di -\Delta
\underline{v}_{\la,\mu}+g(\underline{v}_{\la,\mu})
\leq 2\la_1MH'(\phi_1)\phi_1\quad\mbox{in}\;\;\Omega.
\end{equation}
Indeed, by Hopf's maximum principle, there exists
$\,\delta>0\,$ and
$\,\omega\subset\subset\Omega\,$ such that
$$\di |\nabla
\phi_1|\geq\delta\quad\mbox{in}\;\;\Omega\setminus\omega,$$
$$\di \phi_1\geq\delta\quad\mbox{in}\;\;\omega.$$
On $\,\Omega\setminus\omega\,$ we choose
$M\geq
M_1=\max\{1,\delta^{-2}\}.$ Then,
by (\ref{doizece})
we obtain
\begin{equation}\label{doidoispe}
\di -\Delta
\underline{v}_{\la,\mu}+g(\underline{v}_{\la,\mu})
\leq
\la_1MH'(\phi_1)\phi_1\quad\mbox{in}\;\;\Omega\setminus\omega.
\end{equation}
Fix
$M\geq
\max\left\{M_1,\frac{g(H(\delta))}{\la_1H'(\delta)\delta}\right\}.$
Then
$$\di g(H(\phi_1))\leq g(H(\delta))\leq
\la_1MH'(\delta)\delta
\leq\la_1MH'(\phi_1)\phi_1\quad\mbox{in}\;\;\omega.$$
From (\ref{doizece}) we deduce
\begin{equation}\label{doitreispe}
\di -\Delta
\underline{v}_{\la,\mu}+g(\underline{v}_{\la,\mu})
\leq 2\la_1MH'(\phi_1)\phi_1\quad\mbox{in}\;\;\omega.
\end{equation}
Hence our  claim (\ref{doiunspe})  follows from
(\ref{doidoispe})
and (\ref{doitreispe}).

Since $\,\phi_1>0\,$ in $\,\Omega,$ from
(\ref{doisase})
we have
\begin{equation}\label{doipaispe}
\di 1\leq\frac{H'(\phi_1)\phi_1}{H(\phi_1)}\leq
2\qquad\mbox{in}\;\;\Omega.
\end{equation}
Thus, (\ref{doiunspe}) and (\ref{doipaispe}) yield
\begin{equation}\label{doicincispe}
\di -\Delta
\underline{v}_{\la,\mu}+g(\underline{v}_{\la,\mu})
\leq
4\la_1MH(\phi_1)=4\la_1\underline{v}_{\la,\mu}\quad\mbox{in}\;\;\Omega.
\end{equation}
Take
$\di
\overline{\la}=4\la_1c^{-1}|\underline{v}_{\la,\mu}|_{\infty},$
where
$\,c=\inf\limits_{x\in\overline{\Omega}}f(x,|\underline{v}_{\la,\mu}|_{\infty})>0$.
If $\,\la>\overline{\la},$ the assumption $\,(f1)\,$
produces
$$\di\la\,\frac{f(x,\underline{v}_{\la,\mu})}{\underline{v}_{\la,\mu}}\geq
\overline{\la}\,\frac{f(x,|\underline{v}_{\la,\mu}|_{\infty})}
{|\underline{v}_{\la,\mu}|_{\infty}}\geq
4\la_1,\quad\mbox{for all}\;\;x\in\Omega.$$
This combined with (\ref{doicincispe}) gives
$$\di -\Delta
\underline{v}_{\la,\mu}+g(\underline{v}_{\la,\mu})\leq
\la\,f(x,\underline{v}_{\la,\mu})\quad\mbox{in}\;\;\Omega.$$ Hence
$\,\underline{v}_{\la,\mu}\,$ is a sub-solution of
(\ref{doitrei}), for all $\,\la>\overline{\la}\,$ and $\,\mu>0$.

We now prove  that
$\,\underline{v}_{\la,\mu}\in{\cal E}\!,$ that is
$\,g(\underline{v}_{\la,\mu})\in L^1(\Omega).$
Denote $\,\Omega_0=\{x\in\Omega;\;\phi_1(x)<\eta\}.$
By  (\ref{doisapte}) and
(\ref{doiopt}) it follows that
$$\di g(\underline{v}_{\la,\mu})=g(MH(\phi_1))\leq
g(H(\phi_1))\leq CH^{-\alpha}(\phi_1)
\leq
C_0\,\phi^{-2\alpha/(1+\alpha)}_1\qquad\mbox{in}\;\;\Omega_0,$$
$$\di g(\underline{v}_{\la,\mu})\leq
g(MH(\eta))\qquad\mbox{in}\;\;\Omega\setminus\Omega_0.$$ These
estimates combined with  Lemma \ref{l1jde} yield
$\,g(\underline{v}_{\la,\mu})\in L^1(\Omega)\,$ and so
$\,\Delta\underline{v}_{\la,\mu}\in L^1(\Omega).$ Hence
$$\di
\Delta\overline{v}_{\la,\mu}+\Phi_{\la,\mu}(x,\overline{v}_{\la,\mu})\leq
0\leq
\Delta\underline{v}_{\la,\mu}+\Phi_{\la,\mu}(x,\underline{v}_{\la,\mu})
\quad\mbox{in}\;\;\Omega,$$
$$\di
\underline{v}_{\la,\mu},\,\overline{v}_{\la,\mu}>0\quad\mbox{in}\;\;\Omega,$$
$$\di\underline{v}_{\la,\mu}=\overline{v}_{\la,\mu}\quad\mbox{on}\;\;\partial\Omega,$$
$$\di \Delta\underline{v}_{\la,\mu}\in L^1(\Omega).$$
By Lemma \ref{l2}, it follows that
$\,\underline{v}_{\la,\mu}\leq\overline{v}_{\la,\mu}\,$ on
$\,\overline{\Omega}.$  Now, standard elliptic arguments guarantee
the existence of a solution $\,v_{\la,\mu}\in C^2(\Omega)\cap
C(\overline{\Omega})\,$ for (\ref{doitrei}) such that
$\,\underline{v}_{\la,\mu}\leq
v_{\la,\mu}\leq\overline{v}_{\la,\mu}\,$ in $\,\overline{\Omega}.$
Since $\,\underline{v}_{\la,\mu}\in {\cal E}\!,$ by Remark
\ref{r1} we deduce that $\,v_{\la,\mu}\in {\cal E}.$ Hence, for
all $\,\la>\overline\la\,$ and $\,\mu>0,$
 problem (\ref{doitrei}) has at least a solution in
$\,{\cal E}.$ The proof of Lemma
\ref{l4} is now complete.
\qed

\medskip
We shall often refer in what follows to the following
approaching problem of
$(P_{\,\la,\,\mu}):$
$$\left\{\begin{tabular}{ll}
$-\Delta u+K(x)g(u)=\la f(x,u)+\mu h(x)$ \quad & ${\rm
in}\ \Omega,\,$\\
$u>0$ \quad & ${\rm in}\ \Omega,\,$\\
$\di u=\frac{1}{k}$ \quad & ${\rm on}\
\partial\Omega.\,$\\
\end{tabular} \right. \eqno(P^{\,k}_{\,\la,\,\mu\,})$$
where $\,k\,$ is a positive integer. We observe that any solution
of $\,(P_{\,\la,\,\mu})\,$ is a sub-solution of
$\,(P^{\,k}_{\,\la,\,\mu\,})$.

\medskip
{\it Proof of Theorem \ref{th1jde}}. Suppose to the contrary that
there exists $\la$ and $\mu$ such that $\,(P_{\,\la,\,\mu})\,$ has
a solution $\,u_{\la,\mu}\in{\cal E}$ and let $\,U_{\la,\mu}\,$ be
the solution of (\ref{doipatru}). Since
$$\di\Delta
U_{\la,\mu}+\Phi_{\la,\mu}(x,U_{\la,\mu})\leq 0\leq \Delta
u_{\la,\mu}+\Phi_{\la,\mu}(x,u_{\la,\mu})\quad\mbox{in}\;\;\Omega,$$
by Lemma \ref{l2} we get $\,u_{\la,\mu}\leq U_{\la,\mu}\,$ in
$\,\overline\Omega.$

Consider the perturbed problem
\begin{equation}\label {treiunu}
 \left\{\begin{tabular}{ll}
$-\Delta u+K_*g(u+\ep)=\la f(x,u)+\mu h(x)$ \quad &
${\rm in}\ \Omega,$\\
$u>0$ \quad & ${\rm in}\ \Omega,$\\
$u=0$ \quad & ${\rm on}\ \partial\Omega.$\\
\end{tabular} \right.
\end{equation}

Since $\,K_*>0,$ it follows that $\,u_{\la,\mu}\,$ and
$\,U_{\la,\mu}\,$ are sub and super-solution for (\ref{treiunu}),
respectively. So, by elliptic regularity,
 there exists
$\,u_{\ep}\in C^{2,\gamma}(\overline{\Omega})\,$ a
solution of
(\ref{treiunu}) such that
\begin{equation}\label{treidoi}
u_{\la,\mu}\leq u_{\ep}\leq U_{\la,\mu}
\quad\mbox{in}\;\;\Omega.
\end{equation}
Integrating in (\ref{treiunu}) we deduce
$$\di -\int_{\Omega}\Delta
u_{\ep}dx+K_*\int_{\Omega}g(u_{\ep}+\ep)dx=
\int_{\Omega}[\la f(x,u_{\ep})+\mu h(x)]dx.$$
Hence
\begin{equation}\label{treitrei}
\di-\int_{\partial\Omega}\frac{\partial
u_{\ep}}{\partial n}ds+
K_*\int_{\Omega}g(u_{\ep}+\ep)dx\leq M,
\end{equation}
where $\,M>0\,$ is a constant. Since
$\,\di\frac{\partial
u_{\ep}}{\partial n}\leq 0$ on $\partial\Omega$,
relation
(\ref{treitrei}) yields
$\,\di K_*\int_{\Omega}g(u_{\ep}+\ep)dx\leq M,$
and so
$\,\di K_*\int_{\Omega}g(U_{\la,\mu}+\ep)dx\leq M.$
Thus, for any compact subset
$\,\omega\subset\subset\Omega\,$ we have
$$\di K_*\int_{\omega}g(U_{\la,\mu}+\ep)dx\leq M.$$
Letting $\,\ep\rightarrow 0,$ the above relation leads
to
$\,\di K_*\int_{\omega}g(U_{\la,\mu})dx\leq M.$
Therefore
\begin{equation}\label{treitpatru}
\di K_*\int_{\Omega}g(U_{\la,\mu})dx\leq M.
\end{equation}

Choose $\,\delta>0\,$ sufficiently small and define
$\,\Omega_{\delta}:=\{x\in\Omega;\;\mbox{dist}(x,\partial\Omega)\leq\delta\}.$
Taking into account the regularity of domain, there
exists $\,k>0\,$ such that
$$U_{\la,\mu}\leq
k\,\mbox{dist}(x,\partial\Omega)\quad\mbox{for all}
\;x\in\Omega_{\delta}.$$
Then
$$\int_{\Omega}g(U_{\la,\mu})dx\geq\int_{\Omega_{\delta}}g(U_{\la,\mu})dx\geq
\int_{\Omega_{\delta}}g\left(k\,\mbox{dist}(x,\partial\Omega)\right)dx=+\infty,$$
which contradicts (\ref{treitpatru}). It follows that the problem
$\,(P_{\,\la,\,\mu})\,$ has no solutions in $\,{\cal E}\,$ and the
proof of Theorem \ref{th1jde} is now complete. \qed

\smallskip
 Using the same method as in Zhang
\cite[Theorem 2]{zhang}, we can prove that $\,(P_{\,\la,\,\mu})\,$
has no solution in $\,C^2(\Omega)\cap C^1({\overline{\Omega}}),$
as it was pointed out in Choi, Lazer, and McKenna \cite[Remark
2]{clm}.

\medskip
{\it Proof of Theorem \ref{th2jde}}. We split the proof into
several steps.

{\bf Step I.} {\sc Existence of the solutions of
$\,(P_{\,\la,\,\mu})\,$ for $\,\la\,$ large}.
By Lemma \ref{l4},
there exists $\,\overline{\la}\,$ such that for all
$\,\la>\overline{\la}\,$ and $\,\mu>0\,$ the problem
$$\left\{\begin{tabular}{ll}
$-\Delta v+K^*g(v)=\la f(x,v)+\mu h(x)$ \quad & ${\rm
in}\ \Omega,$\\
$v>0$ \quad & ${\rm in}\ \Omega,$\\
$v=0$ \quad & ${\rm on}\ \partial\Omega,$\\
\end{tabular} \right.$$
has at least one solution $\,v_{\la,\mu}\in{\cal E}.$ Then
$\di\,v_k=v_{\la,\mu}+\frac{1}{k}\,$ is  a sub-solution of
$\,(P^k_{\,\la,\,\mu})\,$ for all positive integers $\,k\geq 1.$

From Lemma \ref{l1jmpa2}, let $\,w\in
C^{2,\gamma}(\overline{\Omega})\,$ be the solution of
$$\left\{\begin{tabular}{ll}
$-\Delta w=\la f(x,w)+\mu h(x)$ \quad & ${\rm in}\
\Omega,$\\
$w>0$ \quad & ${\rm in}\ \Omega,$\\
$w=1$ \quad & ${\rm on}\ \partial\Omega.$\\
\end{tabular} \right.$$
It follows that $\,w\,$ is a super-solution of
$\,(P^k_{\,\la,\,\mu})\,$ for all $\,k\geq 1\,$ and
$$\di \Delta w+\Phi_{\la,\mu}(x,w)\leq 0\leq\Delta
v_1+\Phi_{\la,\mu}(x,v_1)
\quad\mbox{in}\;\;\Omega,$$
$$w,\;v_1>0\quad\mbox{in}\;\;\Omega,$$
$$w=v_1\quad\mbox{on}\;\;\partial\Omega,$$
$$\di \Delta v_1\in L^1(\Omega).$$
Therefore, by Lemma \ref{l2},  $1\leq v_1\leq w$ in
$\,\overline{\Omega}.$ Standard elliptic arguments  imply that
there exists a solution $\,u^1_{\la,\mu}\in
C^{2,\gamma}(\overline{\Omega})$ of $\,(P^1_{\,\la,\,\mu})\,$ such
that $\, v_1\leq u^1_{\la,\mu}\leq w\,$ in $\,\overline{\Omega}$.
Now, taking $\,u^1_{\la,\mu}\,$ and $\,v_2\,$ as a pair of super
and sub-solutions for $\,(P^2_{\,\la,\,\mu}),$ we obtain a
solution $\,u^2_{\la,\mu}\in C^{2,\gamma}(\overline{\Omega})$ of
$\,(P^2_{\,\la,\,\mu})\,$ such that $\, v_2\leq u^2_{\la,\mu}\leq
u^1_{\la,\mu}\,$ in $\,\overline{\Omega}.$ In this manner
 we find a sequence
$\di\,\{u^n_{\la,\mu}\}\,$ such that
\begin{equation}\label{treicinci}
\di v_n\leq u^n_{\la,\mu}\leq u^{n-1}_{\la,\mu}\leq
w\quad\mbox{in}\;\;
\overline{\Omega}. \end{equation}
Define
$\,u_{\la,\mu}(x)=\lim\limits_{n\rightarrow\infty}u^n_{\la,\mu}(x)\,$
for all
$\,x\in\overline{\Omega}.$ Standard bootstrap
arguments
 imply that
$\,u_{\la,\mu}\,$ is a solution of $\,(P_{\,\la,\,\mu}).$ From
(\ref{treicinci}) we have $\,v_{\la,\mu}\leq u_{\la,\mu}\leq w$ in
$\,\overline{\Omega}.$ Since $\,v_{\la,\mu}\in{\cal E},$ by Remark
\ref{r1} it follows that $\,u_{\la,\mu}\in{\cal E}.$ Consequently,
problem $\,(P_{\,\la,\,\mu})\,$ has at least a solution in
$\,{\cal E}\,$ for all $\,\la>\overline{\la}\,$ and $\,\mu>0.$

\smallskip
{\bf Step II.} {\sc Existence of the solutions of
$\,(P_{\,\la,\,\mu})\,$ for $\,\mu\,$ large.} Let us first notice
that $\,g\,$ verifies the hypotheses of Theorem 2 in D\'iaz,
Morel, and Oswald \cite{dmo}. We also remark that the assumption
$\,(g2)\,$ and Lemma \ref{l1jde} is essential to find a
sub-solution in the proof of Theorem~2 in D\'iaz, Morel, and
Oswald \cite{dmo}.

According to this result, there exists
$\,\overline\mu>0\,$ such that the problem
$$\left\{\begin{tabular}{ll}
$-\Delta v+K^*g(v)=\mu h(x)$ \quad & ${\rm in}\
\Omega,\,$\\
$v>0$ \quad & ${\rm in}\ \Omega,\,$\\
$v=0$ \quad & ${\rm on}\ \partial\Omega,\,$\\
\end{tabular} \right.$$
has at least a solution $\,v_{\mu}\in{\cal E}\,$ provided that
$\,\mu>\overline\mu.$ Fix $\,\la>0\,$ and denote
$\di\,v_k=v_{\mu}+\frac{1}{k},\;\,k\geq 1.$ Hence $v_k\,$ is a
sub-solution of $(P^k_{\,\la,\,\mu}),$ for all $\,k\geq 1.$
Similarly to the previous step we obtain a solution
$\,u_{\la,\mu}\in{\cal E}\,$ for all $\,\la>0\,$ and
$\,\mu>\overline\mu.$

\smallskip
{\bf Step III.} {\sc Nonexistence for $\,\la,\mu\,$ small}. Let
$\,\la,\,\mu>0.\;$ Since $\,K_*>0,$ the assumption $\,(g1)\,$
implies $\di\lim\limits_{s\downarrow
0}\Psi_{\la,\mu}(x,s)=-\infty,$ uniformly for
$\,x\!\in\overline{\Omega}.$ So,  there exists $\,c>0\,$ such that
\begin{equation}\label{treisase}
\di\Psi_{\la,\mu}(x,s)<0\quad\mbox{for
all}\;\;(x,s)\in\overline{\Omega}\times(0,c).
\end{equation}
Let $\,s\geq c.$ From $\,(f1)\,$ we deduce
$$\di\frac{\Psi_{\la,\mu}(x,s)}{s}\leq\la\,\frac{f(x,s)}{s}+\mu\,\frac{h(x)}{s}
\leq\la\,\frac{f(x,c)}{c}+\mu\,\frac{|h|_{\infty}}{s},$$
for all $\,x\in\overline\Omega.$
Fix $\di\mu<\frac{c\la_1}{2|h|_{\infty}}$
and let
$\,\di
M=\sup\limits_{x\in\overline\Omega}\frac{f(x,c)}{c}>0$.
From the above inequality
we have
\begin{equation}\label{treisapte}
\di\frac{\Psi_{\la,\mu}(x,s)}{s}\leq\,\la
M+\frac{\la_1}{2},\quad\mbox{for all}
\;\;(x,s)\in\overline{\Omega}\times[c,+\infty).
\end{equation}
Thus, (\ref{treisase}) and (\ref{treisapte}) yield
\begin{equation}\label{treiopt}
\di\Psi_{\la,\mu}(x,s)\leq
a(\la)s+\frac{\la_1}{2}s,\quad\mbox{for all}
\;\;(x,s)\in\overline{\Omega}\times(0,+\infty).
\end{equation}
Moreover, $\,a(\la)\ri 0\,$ as $\,\la\ri 0$.
If $\,(P_{\,\la,\,\mu})\,$ has a solution
$\,u_{\la,\mu},$ then
$$\begin{tabular}{ll}
$\di\la_1\int_{\Omega}u^2_{\la,\mu}(x)dx$&$\di\leq\int_{\Omega}
|\nabla
u_{\la,\mu}|^2dx=-\int_{\Omega}u_{\la,\mu}(x)\Delta
u_{\la,\mu}(x)dx$\\
&$\di\leq\int_{\Omega}u_{\la,\mu}(x)\Psi(x,u_{\la,\mu}(x))dx.$\\
\end{tabular} $$
Using (\ref{treiopt}), we get
$$\di\la_1\int_{\Omega}u^2_{\la,\mu}(x)dx\leq\left[a(\la)+
\frac{\la_1}{2}\right]\int_{\Omega}u^2_{\la,\mu}(x)dx.
$$
Since $\,a(\la)\ri 0\,$ as $\,\la\ri 0,$ the above
relation leads to a
contradiction for $\,\la,\mu>0\,$ sufficiently small.

\smallskip
{\bf Step IV.} {\sc Existence of a maximal solution of
$\,(P_{\,\la,\,\mu})\,$.} We  show that if $\,(P_{\,\la,\,\mu})\,$
has a solution $\,u_{\la,\mu}\in{\cal E},$ then  it has a maximal
solution. Let $\,\la,\mu>0\,$ be such that $\,(P_{\,\la,\,\mu})\,$
has a solution $\,u_{\la,\mu}\in{\cal E}$. If $\,U_{\la,\mu}\,$ is
the solution of (\ref{doipatru}), by Lemma \ref{l2} we have
$\,u_{\la,\mu}\leq U_{\la,\mu}\,$ in $\,\overline\Omega.$ For any
$j\geq 1$, denote
$$\di\Omega
_j=\left\{x\in\Omega;\;\mbox{dist}(x,\partial\Omega)>\frac
1j\right\}.$$
Let $\,U_0=U_{\la,\mu}\,$ and $\,U_j\,$ be the
solution of
$$\left\{\begin{tabular}{ll}
$-\Delta \zeta+K(x)g(U_{j-1})=\la f(x,U_{j-1})+\mu
h(x)$ \quad & ${\rm in}\ \Omega_j,$\\
$\zeta=U_{j-1}$ \quad & ${\rm in}\
\Omega\setminus\Omega_j.$\\
\end{tabular} \right.$$
Using the fact that $\,\Psi_{\la,\mu}\,$ is
nondecreasing with respect to the second
variable, we get
$$\di u_{\la,\mu}\leq U_j\leq U_{j-1}\leq
U_0\quad\mbox{in}\;\;\overline\Omega.$$ If $\,\overline
u_{\la,\mu}(x)=\lim\limits_{j\rightarrow\infty}U_j(x)\,$ for all
$\,x\in\overline{\Omega},$ by standard elliptic arguments (see
Gilbarg and Trudinger \cite{gt}) it follows that $\,\overline
u_{\la,\mu}\,$ is a solution of $\,(P_{\,\la,\,\mu})$. Since
$\,u_{\la,\mu}\leq \overline u_{\la,\mu}\,$ in
$\,\overline{\Omega},$ by Remark \ref{r1} we have $\,\overline
u_{\la,\mu}\in{\cal E}.$ Moreover, $\,\overline u_{\la,\mu}\,$
 is
a maximal solution of $\,(P_{\,\la,\,\mu}).$

\smallskip
{\bf Step V.} {\sc Dependence on $\,\la\,$ and $\,\mu\,$.} We
first  show  the dependence on $\,\la\,$ of the maximal solution
$\,\overline u_{\la,\mu}\in{\cal E}\,$ of $\,(P_{\,\la,\,\mu})$.
For this purpose, fix $\,\mu>0\,$ and define
$$A:=\{\la>0;\;(P_{\,\la,\,\mu})\;\,\mbox{has at least
a solution} \;\,u_{\la,\mu}\in{\cal E}\}.$$ Let $\,\la_*=\inf A$.
From the previous steps we have $\,A\not=\emptyset\,$ and
$\,\la_*>0$. Let $\,\la_1\in A\,$ and $\,\overline
u_{\la_1,\mu}\,$ be the maximal solution of
$\,(P_{\,\la_1,\,\mu})$. We prove that $\,(\la_1,+\infty)\subset
A$. If $\,\la_2>\la_1\,$ then  $\,\overline u_{\la_1,\mu}\,$ is a
sub-solution of $\,(P_{\,\la_2,\,\mu})\,$. On the other hand,
$$\di\Delta
U_{\la_2,\mu}+\Phi_{\la_2,\mu}(x,U_{\la_2,\mu})\leq
0\leq
\Delta \overline
u_{\la_1,\mu}+\Phi_{\la_2,\mu}(x,\overline
u_{\la_1,\mu})
\quad\mbox{in}\;\;\Omega,$$
$$\di U_{\la_2,\mu},\,\overline
u_{\la_1,\mu}>0\quad\mbox{in}\;\;\Omega,$$
$$\di U_{\la_2,\mu}\geq \overline
u_{\la_1,\mu}\quad\mbox{on}\;\;\partial\Omega,$$
$$\di \Delta\overline u_{\la_1,\mu}\in  L^1(\Omega).$$
By Lemma \ref{l2},  $\,\overline{u}_{\la_1,\mu}\leq
U_{\la_2,\mu}\,$ in $\,\overline\Omega$. In the same way as in
Step IV we find
 a solution
$\,u_{\la_2,\mu}\in{\cal E}\,$ of
$\,(P_{\,\la_2,\,\mu})$ such that
$$\di \overline u_{\la_1,\mu}\leq u_{\la_2,\mu}\leq
U_{\la_2,\mu}
\quad\mbox{in}\;\;\overline\Omega.$$
Hence $\,\la_2\in A\,$ and so
$\,(\la_*,+\infty)\subset
A.$
If $\,\overline u_{\la_2,\mu}\in{\cal E}\,$ is the
maximal solution of
$\,(P_{\,\la_2,\,\mu}),$ the above relation implies
$\,\di \overline u_{\la_1,\mu}\leq \overline
u_{\la_2,\mu}\,$
in $\,\overline\Omega$. By the maximum principle, it
follows that
$\,\di \overline u_{\la_1,\mu}< \overline
u_{\la_2,\mu}\,$
in $\,\Omega$. So, $\,\overline u_{\la,\mu}\,$ is
increasing with respect to $\,\la.$

To prove the dependence on $\,\mu,$ we fix $\,\la>0\,$
and define
$$B:=\{\mu>0;\;(P_{\,\la,\,\mu})\;\,\mbox{has at least
one solution} \;\,u_{\la,\mu}\in{\cal E}\}.$$ Let $\,\mu_*=\inf
B$. The conclusion follows in the same manner as above. The proof
of Theorem \ref{th2jde} is now complete. \qed

\medskip
{\it Proof of Theorem \ref{th3jde}}. Let $\,\la,\mu>0.$ We recall
that the function $\,\Psi_{\la,\mu}\,$ defined in
(\ref{doizerozero}) verifies the hypotheses of Lemma
\ref{l1jmpa2}, since $\,K^*\leq 0$. So, there exists
$\,u_{\la,\mu}\in C^{2,\gamma}(\Omega) \cap C(\overline{\Omega})$
a solution of $\,(P_{\,\la,\,\mu})$. If $\,U_{\la,\mu}\,$ is the
solution of (\ref{doipatru}), then
$$\di\Delta
u_{\la,\mu}+\Phi_{\la,\mu}(x,u_{\la,\mu})\leq 0\leq
\Delta U_{\la,\mu}+\Phi_{\la,\mu}(x,U_{\la,\mu})
\quad\mbox{in}\;\;\Omega,$$
$$\di u_{\la,\mu},\,
U_{\la,\mu}>0\quad\mbox{in}\;\;\Omega,$$
$$\di
u_{\la,\mu}=U_{\la,\mu}=0\quad\mbox{on}\;\;\partial\Omega.$$ Since
$\di \Delta U_{\la,\mu}\in L^1(\Omega),$  by Lemma \ref{l2} we get
$\,\di u_{\la,\mu}\geq U_{\la,\mu}\,$ in $\,\overline{\Omega}.$

We claim that there exists $\,c>0\,$ such that
\begin{equation}\label{treinoua}
\di U_{\la,\mu}\geq c\phi_1\quad\mbox{in}\;\;\Omega.
\end{equation}
Indeed, if not, there exists $\,\{x_n\}\subset\Omega\,$
and
$\,\ep_n\rightarrow 0\,$ such
that
\begin{equation}\label{treizece}
\di \left(U_{\la,\mu}-\ep_n\phi_1\right)(x_n)<0.
\end{equation}
Moreover, we can choose the sequence $\,\{x_n\}\,$
with the additional
property
\begin{equation}\label{treiunspe}
\di \nabla\left(U_{\la,\mu}-\ep_n\phi_1\right)(x_n)=0.
\end{equation}
Passing eventually at a subsequence, we can assume
that
$\,x_n\ri x_0\in\overline{\Omega}$. From
(\ref{treizece}) it follows that
$\,U_{\la,\mu}(x_0)\leq 0\,$ which implies
$\,U_{\la,\mu}(x_0)=0,$ that is
$\,x_0\in\partial\Omega$. Furthermore, from
(\ref{treiunspe}) we have
$\,\nabla U_{\la,\mu}(x_0)=0$. This is a contradiction
since
$\,\di\frac{\partial U_{\la,\mu}}{\partial n}(x_0)<0,$
by  Hopf's
strong maximum principle.
Our claim follows and so
\begin{equation}\label{treidoispe}
\di u_{\la,\mu}\geq U_{\la,\mu}\geq
c\phi_1\quad\mbox{in}\;\;\Omega.
\end{equation}
Then, $\,g(u_{\la,\mu})\leq g(U_{\la,\mu})\leq
g(c\phi_1)$ in $\,\Omega$. From the
assumption $\,(g2)\,$ and Lemma \ref{l1} (using the
same method as in the
proof of Lemma \ref{l4}) it follows that
$\,g(c\phi_1)\in L^1(\Omega)$. Hence
$\,u_{\la,\mu}\in{\cal E}$.

Let us now assume that $\,u^1_{\la,\mu},\,u^2_{\la,\mu}\in{\cal
E}\,$ are two solutions of $\,(P_{\,\la,\,\mu})$. In order to
prove the uniqueness, it is enough to show that
$\,u^1_{\la,\mu}\geq u^2_{\la,\mu}\,$ in $\,\overline{\Omega}$.
This follows by Lemma \ref{l2}.

Let us show now the dependence on $\,\la\,$ of the
solution of $\,(P_{\,\la,\,\mu})$.
For this purpose, let $\,0<\la_1<\la_2\,$ and
$\,u_{\la_1,\mu},\,u_{\la_2,\mu}\,$ be
the unique solutions of $\,(P_{\,\la_1,\,\mu})\,$ and
$\,(P_{\,\la_2,\,\mu})\,$
respectively, with $\,\mu>0\,$ fixed.
Since $\,u_{\la_1,\mu},\,u_{\la_2,\mu}\in{\cal E}\,$
and
$$\di \di\Delta
u_{\la_2,\mu}+\Phi_{\la_2,\mu}(x,u_{\la_2,\mu})\leq 0\leq \Delta
u_{\la_1,\mu}+\Phi_{\la_2,\mu}(x,u_{\la_1,\mu})
\quad\mbox{in}\;\;\Omega,$$ in virtue of Lemma \ref{l2} we find
$\di u_{\la_1,\mu}\leq u_{\la_2,\mu}\,$ in $\,\overline{\Omega}$.
So, by the maximum principle, $\di u_{\la_1,\mu}<u_{\la_2,\mu}\,$
in $\,\Omega$.

The dependence on $\,\mu\,$ follows similarly. The proof of
Theorem \ref{th3jde} is now complete. \qed

\medskip
{\it Proof of Theorem \ref{th4jde}}. {\bf Step I.} {\sc
Existence.} Using the fact that $\,K^*>0,$ from Theorem
\ref{th2jde} it follows that there exists $\,\la_*,\mu_*>0\,$ such
that the problem
$$\left\{\begin{tabular}{ll}
$-\Delta v+K^*g(v)=\la f(x,v)+\mu h(x)$ \quad & ${\rm
in}\ \Omega,\,$\\
$v>0$ \quad & ${\rm in}\ \Omega,\,$\\
$v=0$ \quad & ${\rm on}\ \partial\Omega.\,$\\
\end{tabular} \right.$$
has a maximal solution $\,v_{\la,\mu}\in{\cal E},$ provided
$\,\la>\la_*\,$ or $\,\mu>\mu_*$. Moreover, $\,v_{\la,\mu}\,$ is
increasing with respect to $\,\la\,$ and $\,\mu$. Then
$\di\,v_k=v_{\la,\mu}+\frac{1}{k}\,$ is a sub-solution of
$\,(P^k_{\,\la,\,\mu}),$ for all $\,k\geq 1$. On the other hand,
by Lemma \ref{l1jmpa2}, the boundary value problem
$$\left\{\begin{tabular}{ll}
$-\Delta w+K_*g(w)=\la f(x,w)+\mu h(x)$ \quad & ${\rm
in}\ \Omega,\,$\\
$w>0$ \quad & ${\rm in}\ \Omega,\,$\\
$\di w=\frac{1}{k}$ \quad & ${\rm on}\
\partial\Omega.\,$\\
\end{tabular} \right.$$
has a solution $\,w_k\in C^{2,\gamma}(\overline{\Omega})$.
Obviously, $\,w_k\,$ is a super-solution of
$\,(P^k_{\,\la,\,\mu})$.

Since $\,K^*>0>K_*,$ we have
$$\di\Delta w_k+\Phi_{\la,\mu}(x,w_k)\leq 0\leq
\Delta v_k+\Phi_{\la,\mu}(x,v_k)
\quad\mbox{in}\;\;\Omega,$$
and
$$\di w_k,\, v_k>0\quad\mbox{in}\;\;\Omega,$$
$$\di w_k=v_k\quad\mbox{on}\;\;\partial\Omega,$$
$$\,\di \Delta v_k\in  L^1(\Omega).$$
From Lemma \ref{l2} it follows that $\,v_k\leq w_k\,$ in
$\,\overline{\Omega}$. By standard super and sub-solution
argument, there exists a minimal solution $\,u^1_{\la,\mu}\in
C^{2,\gamma}(\overline{\Omega})$ of $\,(P^1_{\,\la,\,\mu})\,$ such
that $\, v_1\leq u^1_{\la,\mu}\leq w_1\,$ in
$\,\overline{\Omega}$. Now, taking $\,u^1_{\la,\mu}\,$ and
$\,v_2\,$ as a pair of super and sub-solutions for
$\,(P^2_{\,\la,\,\mu}),$ we deduce that there exists a minimal
solution $\,u^2_{\la,\mu}\in C^{2,\gamma}(\overline{\Omega})$ of
$\,(P^2_{\,\la,\,\mu})\,$ such that $\, v_2\leq u^2_{\la,\mu}\leq
u^1_{\la,\mu}\,$ in $\,\overline{\Omega}.$ Arguing in the same
manner, we obtain a sequence $\di\,\{u^k_{\la,\mu}\}\,$ such that
\begin{equation}\label{treitreispe}
\di v_k\leq u^k_{\la,\mu}\leq u^{k-1}_{\la,\mu}\leq
w_1\quad\mbox{in}\;\; \overline{\Omega}. \end{equation} Define
$\,u_{\la,\mu}(x)=\lim\limits_{k\rightarrow\infty}u^k_{\la,\mu}(x)\,$
for all $\,x\in\overline{\Omega}$. With a similar argument to that
used in  the proof of Theorem \ref{th2jde}, we find that
$\,u_{\la,\mu}\in{\cal E}\,$ is a solution of
$\,(P_{\,\la,\,\mu})$. Hence, problem  $\,(P_{\,\la,\,\mu})\,$ has
at least a solution in $\,{\cal E},$ provided that $\,\la>\la_*\,$
or $\,\mu>\mu_*$.

{\bf Step II.} {\sc Dependence on $\,\la\,$ and $\,\mu\,$}. As
above, it is enough to justify only the dependence on $\la$. Fix
$\la_*<\la_1<\la_2,$ $\mu>0$ and let $u_{\la_1,\mu}$,
$u_{\la_2,\mu}\in{\cal E}$ be the solutions of
$\,(P_{\,\la_1,\,\mu})\,$ and $\,(P_{\,\la_2,\,\mu})\,$
respectively that we have obtained in Step I. It follows that
$\,u^k_{\la_2,\mu}\,$ is a super-solution of
$\,(P^k_{\,\la_1,\,\mu})$. So, Lemma \ref{l2} combined with the
fact that $\,v_{\la,\mu}\,$ is increasing with respect to
$\,\la>\la_*\,$ yield
$$\di u^k_{\la_2,\mu}\geq
v_{\la_2,\mu}+\frac{1}{k}\geq v_{\la_1,\mu}+\frac{1}{k}
\quad\mbox{in}\;\;\overline{\Omega}.$$ Thus,
$\,u^k_{\la_2,\mu}\geq u^k_{\la_1,\mu}\,$ in
$\,\overline{\Omega}\,$ since $\,u^k_{\la_1,\mu}\,$ is the minimal
solution of $\,(P^k_{\,\la_1,\,\mu})\,$ which satisfies
$\,u^k_{\la_1,\mu}\geq v_{\la_1,\mu}+1/k$ in
$\,\overline{\Omega}$. It follows that $\,u_{\la_2,\mu}\geq
u_{\la_1,\mu}\,$ in $\,\overline{\Omega}$. By the maximum
principle we deduce that $\,u_{\la_2,\mu}>u_{\la_1,\mu}\,$ in
$\,\Omega$. This concludes the proof. \qed

\section{Bifurcation and asymptotics for the
singular Lane-Emden-Fowler equation with a convection term} Let
$\Omega\subset \RR^N$ $(N\geq 2)$ be a bounded domain with a
smooth boundary. In this section we are concerned with singular
elliptic problems of the following type
\neweq{P}
\left\{\begin{tabular}{ll}
$-\Delta u=g(u)+\la|\nabla u|^p+\mu f(x,u)$ \quad & $\mbox{\rm in}\ \Omega,$\\
$u>0$ \quad & $\mbox{\rm in}\ \Omega,$\\
$u=0$ \quad & $\mbox{\rm on}\ \partial\Omega,$\\
\end{tabular} \right.
\endeq
where $0<p\leq 2$ and $\la,\mu\geq 0.$ As remarked by
Choquet-Bruhat and Leray \cite{leray} and by Kazdan and Warner
\cite{kazdan}, the requirement that the nonlinearity grows at most
quadratically in $|\nabla u|$ is natural in order to apply the
maximum principle.

Throughout this section we suppose that
$f:\overline{\Omega}\times[0,\infty)\rightarrow[0,\infty)$ is a
H\"{o}lder continuous function which is nondecreasing with respect
to the second variable and is positive on $\overline\Omega\times
(0,\infty).$ We assume that $g:(0,\infty)\ri(0,\infty)$ is a
H\"{o}lder continuous function which is nonincreasing and
$\lim_{s\searrow 0}g(s)=+\infty.$

 Many papers have been devoted to the case $\la=0,$ where the
problem \eq{P} becomes
\neweq{laedin}
\left\{\begin{tabular}{ll}
$-\Delta u=g(u)+\mu f(x,u)$ \quad & $\mbox{\rm in}\ \Omega,$\\
$u>0$ \quad & $\mbox{\rm in}\  \Omega,$\\
$u=0$ \quad & $\mbox{\rm on}\ \partial\Omega,$\\
\end{tabular} \right.
\endeq
If $\mu=0,$ then \eq{laedin} has a unique solution (see Crandall,
Rabinowitz, and Tartar \cite{crt}, Lazer and McKenna \cite{lm1}).
When $\mu>0,$ the study of \eq{laedin} emphasizes the role played
by the nonlinear term $f(x,u).$  For instance, if one of the
following assumptions are fulfilled

$(f1)\;$ there exists $c>0$ such that $f(x,s)\geq cs$ for all
$(x,s)\in \overline\Omega\times [0,\infty);$

$(f2)\;$ the mapping $(0,\infty)\ni s\longmapsto \frac{f(x,s)}{s}$
is nondecreasing for all $x\in\overline\Omega,$\\
then problem \eq{laedin} has solutions only if $\mu>0$ is small
enough (see Coclite and Palmieri \cite{cp}). In turn, when $f$
satisfies the following assumptions

$\di(f3)\;$ the mapping $ (0,\infty)\ni
s\longmapsto\frac{f(x,s)}{s}\quad\mbox{is nonincreasing for
all}\;\, x\in\overline{\Omega};$

$(f4)\; \lim_{s\rightarrow\infty}\frac{f(x,s)}{s}=0,\;\;
\mbox{uniformly for}\;\,x\in\overline{\Omega},$\\
then problem \eq{laedin} has at least one solutions for all
$\mu>0$ (see Coclite and Palmieri \cite{cp}, Shi and Yao
\cite{shi} and the references therein). The same assumptions will
be used in the study of \eq{P}.

By the monotonicity of $g,$ there exists
$$\di a=\lim_{s\ri\infty}g(s)\in[0,\infty).$$

The main results in this section have been obtained by Ghergu and
R\u adulescu \cite{gr4,grjmaa}.

We are first concerned with the case $\la=1$ and $1<p\leq 2.$ In
the statement of the following result we do not need assumptions
$(f1)-(f4);$ we just require that $f$ is a H\"{o}lder continuous
function which is nondecreasing with respect to the second
variable and is positive on $\overline \Omega\times(0,\infty).$

\begin{thm}\label{th1edin} Assume $\la=1$ and $1<p\leq 2.$\\
{\rm (i) } If $p=2$ and $a\geq \la_1,$ then \eq{P} has no solutions;\\
{\rm (ii) } If $p=2$ and $a<\la_1$ or $1<p<2,$ then there exists
$\mu^*>0$ such that \eq{P} has at least one classical solution for
$\mu<\mu^*$ and no solutions exist if $\;\mu>\mu^*.$
\end{thm}

If $\la=1$ and $0<p\leq 1$ the study of existence is close related
to the asymptotic behaviour of the nonlinear term $f(x,u).$ In
this case we prove

\begin{thm}\label{th2edin} Assume $\la=1$ and $0<p\leq 1.$\\
{\rm (i) } If $f$ satisfies $(f1)$ or $(f2),$ then there exists
$\mu^*>0$ such that \eq{P} has at least one classical solution for
$\mu<\mu^*$ and no solutions exist if $\;\mu>\mu^*;$\\
{\rm (ii) } If $\;0<p<1$ and $f$ satisfies $(f3)-(f4),$ then
\eq{P} has at least one solution for all $\mu\geq 0.$
\end{thm}

 Next we are concerned with the case $\mu=1.$ Our result is the
following

\begin{thm}\label{th3edin} Assume $\mu=1$ and $f$ satisfies
assumptions $(f3)$ and $(f4).$ Then the following properties hold true.\\
{\rm (i) } If $\;0<p<1,$ then \eq{P} has at least one classical
solution for all $\la\geq 0$;\\
{\rm (ii) } If $\;1\leq p\leq 2,$ then there exists
$\la^*\in(0,\infty]$ such that \eq{P} has at least one classical
solution for $\la<\la^*$ and no solution exists if $\;\la>\la^*.$
Moreover, if $\,1<p\leq 2,$ then $\la^*$ is finite.
\end{thm}
Related to the above result we raise the following {\bf open
problem:} if $p=1$ and $\mu=1,$ is $\la^*$ a finite number?

Theorem \ref{th3edin} shows the importance of the convection term
$\lambda |\nabla u|^p$ in \eq{P}. Indeed, according to Theorem
\ref{th3jde} and for any $\mu>0$, the boundary value problem
\neweq{PQ}
\left\{\begin{tabular}{ll}
$-\Delta u=u^{-\alpha}+\la|\nabla u|^p+\mu u^\beta$ \quad & $\mbox{\rm in}\ \Omega,$\\
$u>0$ \quad & $\mbox{\rm in}\ \Omega,$\\
$u=0$ \quad & $\mbox{\rm on}\ \partial\Omega$\\
\end{tabular} \right.
\endeq
has a unique solution, provided $\lambda=0$, $\alpha$, $\beta\in
(0,1)$. The above theorem shows that if $\la$ is not necessarily
0, then the following situations may occur : (i) problem \eq{PQ}
has solutions if $p\in (0,1)$ and for all $\lambda\geq 0$; (ii) if
$p\in (1,2)$ then there exists $\la^*>0$ such that problem \eq{PQ}
has a solution for any $\la<\la^*$ and no solution exists if
$\la>\la^*.$

To see the dependence between $\la$ and $\mu$ in \eq{P}, we
consider the special case $f\equiv 1$ and $p=2.$ In this case we
can say more about the problem \eq{P}. More precisely we have

\begin{thm}\label{th4edin}
Assume that $p=2$ and $f\equiv 1.$\\
{\rm (i) } The problem \eq{P} has solution if and only if $\la(a+\mu)<\la_1;$\\
{\rm (ii) } Assume $\mu>0$ is fixed, $g$ is decreasing and let
$\di \la^*=\frac{\la_1}{a+\mu}.$ Then \eq{P} has a unique solution
$u_\la$ for all $\la<\la^*$ and the sequence $(u_\la)_{\la<\la^*}$
is increasing with respect to $\la.$\\
Moreover, if $\,\di \limsup_{s\searrow 0}s^\alpha g(s)<+\infty,$
for some $\alpha\in(0,1),$ then the sequence of solutions
$(u_\la)_{0<\la<\la^*}$ has the following properties

\qquad{\rm (ii1)} For all $0<\la<\la^*$ there exist two positive
constants $c_1,c_2$ depending on $\la$ such that $c_1\,{\rm
dist}(x,\partial \Omega)\leq u_\la\leq c_2\,{\rm dist}(x,\partial
\Omega)$ in $\Omega;$

$\qquad${\rm (ii2)} $u_\la\in C^{1,1-\alpha}(\overline\Omega)\cap
C^2(\Omega);$

$\qquad${\rm (ii3)} $u_{\la}\longrightarrow +\infty$ as
$\la\nearrow \la^*$, uniformly on compact subsets of $\Omega.$
\end{thm}

As regards the uniqueness of the solutions to problem \eq{P}, we
may say that this does not seem to be a feature easy to achieve.
Only when $f(x,u)$ is constant in $u$ we can use classical methods
in order to prove the uniqueness. It is worth pointing out here
that the uniqueness of the solution is a delicate issue even for
the simpler problem \eq{laedin}. We have already observed  that if
$f$ fulfills $(f3)-(f4)$ and $g$ satisfies the same growth
condition as in Theorem \ref{th4edin}, then  this solution is
unique, provided that problem \eq{laedin} has a solution. On the
other hand, if $f$ satisfies $(f2),$ the uniqueness generally does
not occur. In that sense we refer the interested reader to Haitao
\cite{hai}. In the case $f(x,u)=u^q,$ $g(u)=u^{-\gamma},$
$0<\gamma<\frac{1}{N}$ and $1<q<\frac{N+2}{N-2},$ we learn from
\cite{hai} that problem \eq{laedin} has at least two classical
solutions provided $\mu$ belongs to a certain range.

\smallskip
Our approach relies on finding of appropriate sub- and
super-solutions of \eq{P}. This will allows us to enlarge the
study of bifurcation to a class of problems more generally to that
studied in Zhang and Yu \cite{zy}. However, neither the method
used in \cite{zy}, nor our method gives a precise answer if
$\la^*$ is finite or not in the case $p=1$ and $\mu=1.$

We start with some auxiliary results.

 Let $\varphi_1$ be the
normalized positive eigenfunction corresponding to the first
eigenvalue $\la_1$ of $(-\Delta )$ in $H^1_0(\Omega)$. As it is
well known $\la_1>0,$ $\varphi_1\in C^2({\overline{\Omega}})$ and
\neweq{fi}
C_1\,{\rm dist}(x,\partial \Omega)\leq\varphi_1\leq C_2\,{\rm
dist}(x,\partial \Omega)\quad\mbox{ in }\;\Omega,
\endeq
for some positive constants $C_1,C_2>0.$
 From the characterization
of $\la_1$ and $\varphi_1$ we state the following elementary
result. For the convenience of the reader we shall give a complete
proof.

\begin{lem}\label{l01} Let $F:\overline{\Omega}\times(0,\infty)\rightarrow\RR$ be a
continuous function such that $F(x,s)\geq \la_1s+b$ for some $b>0$
and for all $(x,s)\in \overline\Omega\times (0,\infty).$ Then the
problem
\neweq{la1}
\left\{\begin{tabular}{ll}
$-\Delta u=F(x,u)$ \quad & $\mbox{\rm in}\ \Omega,$\\
$u> 0$ \quad & $\mbox{\rm in}\ \Omega,$\\
$u=0$\quad & $\mbox{\rm on}\ \partial\Omega,$\\
\end{tabular} \right.\endeq
has no solutions.
\end{lem}

 \proof By contradiction, suppose that \eq{la1} admits a
solution. This will provide a super-solution of the problem
\neweq{la2}
\left\{\begin{tabular}{ll}
$-\Delta u=\la_1 u+b$ \quad & $\mbox{\rm in}\ \Omega,$\\
$u> 0$ \quad & $\mbox{\rm in}\ \Omega,$\\
$u=0$\quad & $\mbox{\rm on}\ \partial\Omega,$\\
\end{tabular} \right.\endeq
Since 0 is a sub-solution, by the sub and super-solution method
and classical regularity theory it follows that \eq{la1} has a
solution $u\in C^2(\overline\Omega).$ Multiplying by $\varphi_1$
in \eq{la2} and then integrating over $\Omega,$ we get
$$\di -\int_{\Omega}\varphi_1\Delta u =\la_1\int_\Omega \varphi_1u+b\int_{\Omega}\varphi_1,$$
that is $\di \la_1\int_{\Omega}\varphi_1 u =\la_1\int_\Omega
\varphi_1u+b\int_{\Omega}\varphi_1,$ which implies $\di
\int_\Omega\varphi_1=0.$ This is clearly a contradiction since
$\varphi_1>0$ in $\Omega.$ Hence \eq{la1} has no solutions. \qed

According to Lemma \ref{l1jmpa2}, there exists $\zeta\in
C^2({\overline \Omega})$ a solution of the problem
\neweq{treidoiedin}
\left\{\begin{tabular}{ll}
$\di -\Delta \zeta=g(\zeta)$&$\mbox {\rm in}\ \Omega,$\\
$\di \zeta>0$&$\mbox {\rm in}\ \Omega,$\\
$\di \zeta=0$&$\mbox{\rm on}\ \partial\Omega.$\\
\end{tabular}\right.\endeq
Clearly $\zeta$ is a sub-solution of \eq{P} for all $\la\geq 0.$
It is worth pointing out here that the sub-super solution method
still works for the problem \eq{P}. With the same proof as in
Zhang and Yu \cite[Lemmma 2.8]{zy} that goes back to the
pioneering work of Amann \cite{amann} we state the following
result.

\begin{lem}\label{sub-sup}
Let $\la,\mu\geq 0.$ If \eq{P} has a super-solution $\overline
u\in C^2(\Omega) \cap C(\overline \Omega)$ such that
$\zeta\leq\overline u$ in $\Omega,$ then \eq{P} has at least a
solution.
\end{lem}

\begin{lem}\label{l4edin}\emph{(Alaa and Pierre \cite{ap}).} If $\;p>1,$ then there exists
a real number $\bar \sigma>0$ such that the problem
\begin{equation}\label {doitreiedin}
 \left\{\begin{tabular}{ll}
$-\Delta u=|\nabla u|^p+\sigma$ \quad & $\mbox{\rm in}\ \Omega,$\\
$u=0$ \quad & $\mbox{\rm on}\ \partial\Omega,$\\
\end{tabular} \right.
\end{equation}
has no solutions for $\sigma>\bar\sigma.$
\end{lem}

\begin{lem}\label{lp2}
Let $F:\overline\Omega\times(0,\infty)\ri[0,\infty)$ and
$G:(0,\infty)\ri(0,\infty)$ be two H\"{o}lder continuous functions
that verify

$(A1)\;$ $F(x,s)>0,$ for all
$\;(x,s)\in\overline\Omega\times(0,\infty);$

$(A2)\;$ The mapping $[0,\infty)\ni s\longmapsto F(x,s)$ is
nondecreasing for all $x\in\overline\Omega;$

$(A3)\;$ $G$ is nonincreasing and $\;\lim_{s\searrow
0}G(s)=+\infty.$

Assume that $\tau>0$ is a positive real number. Then the following holds.\\
{\rm (i) } If $\;\tau \lim_{s\ri\infty}G(s)\geq \la_1,$ then the
problem
\neweq{trunu}
\left\{\begin{tabular}{ll}
$\di -\Delta u=G(u)+\tau|\nabla u|^2+\mu F(x,u)$&$\mbox{\rm in}\ \Omega,$\\
$\di u>0$&$\mbox{\rm in}\ \Omega,$\\
$\di u=0$&$\mbox{\rm on}\ \partial\Omega,$\\
\end{tabular}\right.
\endeq
has no solutions.\\
{\rm (ii) } If $\;\tau\lim_{s\ri\infty}G(s)<\la_1,$ then there
exists $\bar \mu>0$ such that the problem \eq{trunu} has at least
one solution for all $0\leq \mu<\bar\mu.$
\end{lem}
\proof (i) With the change of variable $\di v= e^{\tau u}-1,$ the
problem \eq{trunu} takes the form
\neweq{trdoi}
\left\{\begin{tabular}{ll}
$-\Delta v=\Psi_\mu(x,u)$\quad & $\mbox{\rm in}\ \Omega,$\\
$v>0$ \quad & $\mbox{\rm in}\ \Omega,$\\
$v=0$ \quad & $\mbox{\rm on}\ \partial\Omega,$\\
\end{tabular} \right. \endeq
where
$$\di \Psi_\mu(x,s)=\tau(s+1)G\left(\frac{1}{\tau}
\ln(s+1)\right)+\mu \tau(s+1)
F\left(x,\frac{1}{\tau}\ln(s+1)\right),$$
for all $(x,s)\in\overline\Omega\times (0,\infty).$\\
Taking into account the fact that $G$ is nonincreasing and
$\tau\lim_{s\ri\infty}G(s)\geq \la_1,$ we get
$$\Psi_\mu(x,s)\geq \la_1(s+1)\quad\mbox{ in }\;
\overline\Omega\times(0,\infty), \mbox{ for all }\; \mu\geq 0.$$
By Lemma \ref{l01} we conclude that \eq{trdoi} has no solutions.
Hence \eq{trunu} has no solutions.\\
(ii) Since
$$\di \lim _{s\rightarrow+\infty}
\frac{\tau(s+1)G\left(\frac{1}{\tau}\ln(s+1)\right)+1}{s}<\la_1$$
and
$$\di \lim _{s\searrow 0}
\frac{\tau(s+1)G\left(\frac{1}{\tau}\ln(s+1)\right)+1}{s}=+\infty,$$
we deduce that the mapping $(0,\infty)\ni s\longmapsto
\tau(s+1)G\left(\frac{1}{\tau}\ln(s+1)\right)+1$ fulfills the
hypotheses in Lemma \ref{l1jmpa2}. According to this one, there
exists $\overline v\in C^2(\Omega)\cap C(\overline\Omega)$ a
solution of the problem
$$\left\{\begin{tabular}{ll}
$-\Delta v=\di
\tau(v+1)G\left(\frac{1}{\tau}\ln(v+1)\right)+1$\quad & $
\mbox{\rm in}\ \Omega,$\\
$v>0$ \quad & $\mbox{\rm in}\ \Omega,$\\
$v=0$ \quad & $\mbox{\rm in}\ \partial\Omega.$\\
\end{tabular} \right.$$
Define
$$\di \bar \mu:=\frac{1}{\tau (\|\overline v\|_{\infty}+1)}\cdot
\frac{1}{\di \max\limits_{x\in\overline\Omega}
F\left(x,\frac{1}{\tau}\ln(\|\overline v\|_{\infty}+1)\right)}.$$
It follows that $\overline v$ is a super-solution of \eq{trdoi}
for all
$0\leq \mu<\bar \mu.$\\
Next we provide a sub-solution $\underline v$ of \eq{trdoi} such
that $\underline v\leq \overline v$ in $\Omega.$ To this aim, we
apply Lemma \ref{l1jmpa2} to get that there exists $\underline
v\in C^2(\Omega)\cap C(\overline\Omega)$ a solution of the problem
$$\left\{\begin{tabular}{ll}
$-\Delta v=\di \tau G\left(\frac{1}{\tau}\ln(v+1)\right)$\quad &
$\mbox{\rm in}\ \Omega,$\\
$v>0$ \quad & $\mbox{\rm in}\ \Omega,$\\
$v=0$ \quad & $\mbox{\rm on}\ \partial\Omega.$\\
\end{tabular} \right.$$
Clearly, $\underline v$ is a sub-solution of \eq{trdoi} for all
$0\leq \mu<\bar\mu.$ Let us prove now that $\underline v\leq
\overline v$ in $\Omega.$ Assuming the contrary, it follows that
$\max_{x\in\overline\Omega} \{\underline v-\overline v\}>0$ is
achieved in $\Omega.$ At that point, say $x_0,$ we have
$$\begin{aligned}
\di 0&\di \leq -\Delta (\underline v-\overline v)(x_0)\\
&\di \leq \tau \left[G\left(\frac{1}{\tau}\ln(\underline
v(x_0)+1)\right)-
G\left(\frac{1}{\tau}\ln(\overline v(x_0)+1)\right)\right]-1<0,\\
\end{aligned}$$
which is a contradiction. Thus, $\underline v\leq \overline v$ in
$\Omega.$ We have proved that $(\underline v,\overline v)$ is an
ordered pair of sub-super solutions of \eq{trdoi} provided $0\leq
\mu<\bar \mu.$ It follows that \eq{trunu} has at least one
classical solution for all $0\leq \mu<\bar \mu$ and the proof of
Lemma \ref{lp2} is now complete. \qed

\medskip
{\it Proof of Theorem \ref{th1edin}}. According to Lemma
\ref{lp2}(i) we deduce that \eq{P} has no solutions if $p=2$ and
$a\geq \la_1.$ Furthermore, if $p=2$ and $a<\la_1,$ in view of
Lemma \ref{lp2}(ii), we deduce that \eq{P} has at least one
classical solution if $\mu$ is small enough. Assume now $1<p<2$
and let us fix $C>0$ such that
\neweq{CC}
aC^{p/2}+C^{p-1}<\la_1.
\endeq
Define
$$\psi:[0,\infty)\ri[0,\infty),\quad \psi(s)=\frac{s^p}{s^2+C}.$$
A careful examination reveals the fact that $\psi$ attains its
maximum at $\bar s=\left(\frac{Cp}{2-p}\right)^{2-p}.$ Hence
$$\di \psi(s)\leq \psi(\bar s)=\frac{p^{p/2}(2-p)^{(2-p)/2}}{2C^{1-p/2}},
\quad\mbox{ for all }\;s\geq 0.$$ By the classical Young's
inequality we deduce
$$\di p^{p/2}(2-p)^{(2-p)/2}\leq 2,$$
which yields $\psi(s)\leq C^{p/2-1},$ for all $s\geq 0.$ Thus, we
have proved
\neweq{inq}
\di s^p\leq C^{p/2}s^2+C^{p/2-1},\quad \mbox{ for all }\;s\geq 0.
\endeq
Consider the problem
\neweq{trtrei}
\left\{\begin{tabular}{ll}
$-\Delta u=g(u)+C^{p/2-1}+C^{p/2}|\nabla u|^2+\mu f(x,u)$&$\mbox {in }\Omega$\\
$u>0$ \quad & $\mbox{\rm in}\ \Omega,$\\
$u=0$ \quad & $\mbox{\rm on}\ \partial\Omega,$\\
\end{tabular} \right. \endeq
By virtue of \eq{inq}, any solution of \eq{trtrei} is a
super-solution of \eq{P}.

 Using now \eq{CC} we get
$$\lim_{s\ri\infty}C^{p/2}(g(u)+C^{p/2-1})<\la_1.$$
The above relation enables us to apply Lemma \ref{lp2}(ii) with
$G(s)=g(s)+C^{p/2-1}$ and $\tau=C^{p/2}.$ It follows that there
exists $\bar\mu>0$ such that \eq{trtrei} has at least a solution
$u.$ With a similar argument to that used in the proof of Lemma
\ref{lp2}, we obtain $\zeta \leq u$ in $\Omega,$ where $\zeta$ is
defined in \eq{treidoiedin}. By Lemma \ref{sub-sup} we get that
\eq{P} has at least one solution if $0\leq \mu<\bar\mu.$

We have proved that \eq{P} has at least one classical solution for
both cases $p=2$ and $a<\la_1$ or $1<p<2,$ provided $\mu$ is
nonnegative small enough. Define next
$$A=\{\mu\geq 0;\;\mbox{ problem }\eq{P}\mbox{ has at least one solution}\}.$$
The above arguments implies that $A$ is nonempty. Let $\mu^*=\sup
A.$ We first show that $[0,\mu^*)\subseteq A.$ For this purpose,
let $\mu_1\in A$ and $0\leq \mu_2<\mu_1.$ If $u_{\mu_1}$ is a
solution of \eq{P} with $\mu=\mu_1,$ then $u_{\mu_1}$ is a
super-solution of \eq{P} with $\mu=\mu_2.$ It is easy to prove
that $\zeta\leq u_{\mu_1}$ in $\Omega$ and by virtue of Lemma
\ref{sub-sup} we conclude that the problem \eq{P} with $\mu=\mu_2$
has at least one solution.

Thus we have proved $[0,\mu^*)\subseteq A.$ Next we show $\mu^*<+\infty.$\\
Since $\lim_{s\searrow 0}g(s)=+\infty,$ we can choose $s_0>0$ such
that $g(s)>\bar \sigma$ for all $s\leq s_0.$ Let
$$\di \mu_0=\frac{\bar \sigma}
{\min_{x\in\overline\Omega}f(x,s_0)}.$$ Using the monotonicity of
$f$ with respect to the second argument, the above relations
yields
$$\di g(s)+\mu f(x,s)\geq \bar\sigma,\quad \mbox{ for all }\;
(x,s)\in\overline\Omega\times (0,\infty)\;\,\mbox{ and }\;
\mu>\mu_0.$$ If \eq{P} has a solution for $\mu>\mu_0,$ this would
be a super-solution of the problem
\begin{equation}\label {treipatru}
\left\{\begin{tabular}{ll}
$-\Delta u=|\nabla u|^p+\bar\sigma$ \quad & $\mbox{\rm in}\ \Omega,$\\
$u=0$ \quad & $\mbox{\rm on}\ \partial\Omega.$\\
\end{tabular} \right.
\end{equation}
Since 0 is a sub-solution, we deduce that \eq{treipatru} has at
least one solution. According to Lemma \ref{l4edin}, this is a
contradiction. Hence $\mu^*\leq \mu_0<+\infty.$ This concludes the
proof of Theorem \ref{th1edin}. \qed

\medskip
{\it Proof of Theorem \ref{th2edin}} (i) We fix $p\in(0,1]$ and
define
$$q=q(p)=\left\{\begin{tabular}{ll}
$p+1$\quad&${\rm if\; } 0<p<1,$\\
$3/2$\quad&${\rm if\; } p=1.$\\
\end{tabular}\right.$$
Consider the problem
\neweq{treisup}
\left\{\begin{tabular}{ll}
$\di -\Delta u=g(u)+1+|\nabla u|^{q}+\mu f(x,u)$&$\mbox {\rm in}\ \Omega,$\\
$\di u>0$&$\mbox{\rm in}\ \Omega,$\\
$\di u=0$&$\mbox{\rm on}\ \partial\Omega.$\\
\end{tabular}\right.
\endeq
Since $s^p\leq s^q+1,$ for all $s\geq 0,$ we deduce that any
solution of \eq{treisup} is a super-solution of \eq{P}.
Furthermore, taking into account the fact that $1<q<2,$ we can
apply Theorem \ref{th1edin}(ii) in order to get that \eq{treisup}
has at least one solution if $\mu$ is small enough. Thus, by Lemma
\ref{sub-sup} we deduce that \eq{P} has at least one classical
solution. Following the method used in the proof of Theorem
\ref{th1edin}, we set
$$\di A=\{\mu\geq 0;\;\mbox{ problem }\eq{P}\mbox{ has at least one solution}\}$$
and let $\mu^*=\sup A.$ With the same arguments we prove that
$[0,\mu^*)\subseteq A.$ It remains only to show that
$\mu^*<+\infty.$ \smallskip

Let us assume first that $f$ satisfies $(f1).$ Since
$\lim_{s\searrow 0}g(s)=+\infty,$ we can choose
$\mu_0>\frac{2\la_1}{c}\,$ such that $\frac{1}{2}\mu_0cs+g(s)\geq
1$ for all $s>0.$ Then
$$\di g(s)+\mu f(x,s)\geq \la_1 s+1,\quad\mbox{ for all }\;
(x,s)\in\overline\Omega\times(0,\infty)\;\mbox{ and }\; \mu\geq
\mu_0.$$ By virtue of Lemma \ref{l01} we obtain that \eq{P} has no
classical solutions if $\mu\geq\mu_0,$ so $\mu^*$ is finite.

Assume now that $f$ satisfies $(f2).$ Since $\lim_{s\searrow
0}g(s)=+\infty,$ there exists $s_0>0$ such that
\neweq{f21}
\di g(s)\geq \la_1(s+1)\quad\mbox{ for all }\; 0<s<s_0.
\endeq
On the other hand, the assumption $(f2)$ and the fact that
$\Omega$ is bounded implies that the mapping
$$\di(0,\infty)\ni s\longmapsto
\frac{\min_{x\in\overline\Omega}f(x,s)}{s+1}$$ is nondecreasing,
so we can choose $\tilde\mu>0$ with the property
\neweq{f22}
\tilde\mu\cdot \frac{\min_{x\in\overline\Omega}f(x,s)}{s+1}\geq
\la_1\quad \mbox{ for all }\;s\geq s_0.
\endeq
Now \eq{f21} combined with \eq{f22} yields
$$\di g(s)+\mu f(x,s)\geq \la_1 (s+1),\quad\mbox{ for all }\;(x,s)\in
\overline\Omega\times(0,\infty)\;\mbox { and
}\;\mu\geq\tilde\mu.$$ Using Lemma \ref{l01}, we deduce that
\eq{P} has no solutions if
$\mu>\tilde\mu,$ that is, $\mu^*$ is finite.\\
The first part in Theorem \ref{th2edin} is therefore established.
\smallskip

 (ii) The strategy is to find a super-solution $\overline
u_\mu\in C^2(\Omega)\cap C(\overline\Omega)$ of \eq{P} such that
$\zeta\leq \overline u_\mu$ in $\Omega.$ To this aim, let $h\in
C^2(0,\eta]\cap C[0,\eta]$ be such that
\begin{equation}\label{dunu}
\left\{\begin{tabular}{ll}
$h''(t)=-g(h(t)),\quad \mbox{ for all } 0<t<\eta,$\\
$h(0)=0,$\\
$h>0\quad\mbox{\rm in}\  (0,\eta].$
\end{tabular} \right.
\end{equation}
The existence of $h$ follows by classical arguments of ODE. Since
$h$ is concave, there exists $h'(0+)\in$$(0,+\infty].$ By taking
$\eta>0$ small enough, we can assume that $h'>0$ in $(0,\eta],$ so
$h$ is increasing on $[0,\eta].$
\begin{lem}\label{lh}
{\rm (i)} $h\in C^1[0,\eta]$ if and only if $\di\int_0^1g(s)ds<+\infty;$\\
{\rm (ii)} If $0<p\leq 2,$ then there exist $c_1,c_2>0$ such that
$$\di (h')^p(t)\leq c_1g(h(t))+c_2,\quad\mbox{ for all }\;0<t<\eta.$$
\end{lem}
\noindent\proof (i) Multiplying by $h'$ in \eq{dunu} and then
integrating on $[t,\eta],$ $0<t<\eta,$ we get
\neweq{hprim}
\di (h')^2(t)-(h')^2(\eta)=2\int_t^\eta
g(h(s))h'(s)ds=2\int_{h(t)}^{h(\eta)} g(\tau)d\tau.
\endeq
This gives
\neweq{hpr}
\di (h')^2(t)= 2G(h(t))+(h')^2(\eta)\quad\mbox{ for all }\;
0<t<\eta,
\endeq
where $G(t)=\di\int_t^{h(\eta)}g(s)ds.$ From \eq{hpr} we deduce
that $h'(0+)$ is finite
if and only if $G(0+)$ is finite, so (i) follows.\\
(ii) Let $p\in(0,2].$ Taking into account the fact that $g$ is
nonincreasing, the inequality \eq{hpr} leads to
\neweq{ddoibis}
\di (h')^2(t)\leq 2h(\eta)g(h(t))+(h')^2(\eta),\quad\mbox{ for all
}0<t<\eta.
\endeq
Since $s^p\leq s^2+1,$ for all $s\geq 0,$ from \eq{ddoibis} we
have
\neweq{ddoi}
\di (h')^p(t)\leq c_1g(h(t))+c_2,\quad\mbox{ for all }0<t<\eta
\endeq
where $c_1=2h(\eta)$ and $c_2=(h')^2(\eta)+1.$ This completes the
proof of our Lemma. \qed \smallskip

{\it Proof of Theorem \ref{th2edin} completed.} Let $p\in(0,1)$
and $\mu\geq 0$ be fixed. We also fix $c>0$ such that
$c\|\varphi_1\|_\infty<\eta.$ By Hopf's maximum principle, there
exist $\delta>0$ small enough and $\theta_1>0$ such that
\neweq{dtrei}
\di |\nabla \varphi_1|>\theta_1 \quad\mbox{ in }\;\Omega_\delta,
\endeq
where $\Omega_\delta:=\{x\in\Omega;\,\mbox{dist}(x,\partial\Omega)
\leq \delta\}.$\\
Moreover, since $\lim_{s\searrow 0}g(h(s))=+\infty,$ we can pick
$\delta$ with the property
\neweq{dpatru}
\di (c\theta_1)^2g(h(c\varphi_1))-3\mu
f(x,h(c\varphi_1))>0\quad\mbox{ in }\;\Omega_\delta.
\endeq
Let $\di
\theta_2:=\inf\limits_{\Omega\setminus\Omega_\delta}\varphi_1>0.$
We choose $M>1$ with
\neweq{dcinci}
M(c\theta_1)^2>3,
\endeq
\neweq{dsase}
Mc\la_1\theta_2h'(c\|\varphi_1\|_{\infty})>3g(h(c\theta_2)).
\endeq
Since $p<1,$ we also may assume
\neweq{dnouap}
(Mc)^{1-p}\la_1(h')^{1-p}(c\|\varphi_1\|_{\infty})\geq
3\|\nabla\varphi_1\|_\infty^p.
\endeq
On the other hand, by Lemma \ref{lh}(ii) we can choose $M>1$ such
that
\neweq{dsaptebis}
3(h'(c\varphi_1))^p\leq M^{1-p}(c\theta_1)^{2-p}g(h(c\varphi_1))
\quad\mbox{ in }\;\Omega_\delta.
\endeq
The assumption $(f4)$ yields
$$\di
\lim_{s\ri\infty}\frac{3\mu
f(x,sh(c\|\varphi_1\|_\infty))}{sh(c\|\varphi_1\|_\infty)}=0.$$ So
we can choose $M>1$ large enough such that
$$\di \frac{3\mu f(x,Mh(c\|\varphi_1\|_\infty))}{Mh(c\|\varphi_1\|_\infty)}<
\frac{c\la_1\theta_2h'(c\|\varphi_1\|_\infty)}{h(c\|\varphi_1\|_\infty)},$$
uniformly in $\Omega.$ This leads us to
\neweq{dsaptep}
3\mu
f(x,Mh(c\|\varphi_1\|_\infty))<Mc\la_1\theta_2h'(c\|\varphi_1\|_{\infty}),
\quad \mbox{ for all }\;x\in\Omega.
\endeq

For $M$ satisfying \eq{dcinci}-\eq{dsaptep}, we prove that
$\overline u_\mu=Mh(c\varphi_1)$ is a super-solution of \eq{P}. We
have
\neweq{calculp}
\di -\Delta \overline
u_\la=Mc^2g(h(c\varphi_1))|\nabla\varphi_1|^2+
Mc\la_1\varphi_1h'(c\varphi_1) \quad\mbox{ in }\;\Omega.
\endeq
First we prove that
\neweq{dzecep}
Mc^2g(h(c\varphi_1))|\nabla\varphi_1|^2\geq g(\overline
u_\mu)+|\nabla \overline u_\mu|^p+\mu f(x,\overline
u_\mu)\quad\mbox{ in }\;\Omega_\delta.
\endeq
From \eq{dtrei} and \eq{dcinci} we get
\neweq{dunspep}
\di \frac{1}{3}Mc^2g(h(c\varphi_1))|\nabla\varphi_1|^2\geq
g(h(c\varphi_1))\geq g(Mh(c\varphi_1))=g(\overline u_\mu)
\quad\mbox{ in }\;\Omega_\delta.
\endeq
By \eq{dtrei} and \eq{dsaptebis} we also have
\neweq{ddoispep}
\di \frac{1}{3}Mc^2g(h(c\varphi_1))|\nabla\varphi_1|^2\geq
(Mc)^p(h')^p(c\varphi_1))|\nabla\varphi_1|^p=|\nabla \overline
u_\mu|^p \quad\mbox{ in }\;\Omega_\delta.
\endeq
The assumption $(f3)$ and \eq{dpatru} produce
\neweq{dtreispep}
 \frac{1}{3}Mc^2g(h(c\varphi_1))|\nabla\varphi_1|^2\geq
\mu Mf(x,h(c\varphi_1))\geq \mu f(x,Mh(c\varphi_1)) \quad\mbox{ in
}\;\Omega_\delta.
\endeq
Now, by \eq{dunspep}, \eq{ddoispep} and \eq{dtreispep} we conclude
that \eq{dzecep} is fulfilled.
\medskip

Next we prove
\neweq{dpaispep}
Mc\la_1\varphi_1h'(c\varphi_1)\geq g(\overline u_\mu)+|\nabla
\overline u_\mu|^p+\mu f(x,\overline u_\mu)\quad\mbox{ in }\;
\Omega\setminus\Omega_\delta.
\endeq
From \eq{dsase} we obtain
\neweq{dcincispep}
\di \frac{1}{3}Mc\la_1\varphi_1h'(c\varphi_1)\geq
g(h(c\varphi_1))\geq g(Mh(c\varphi_1))=g(\overline
u_\mu)\quad\mbox{ in }\; \Omega\setminus\Omega_\delta.
\endeq
From \eq{dnouap} we get
\neweq{dsaispep}
\di \frac{1}{3}Mc\la_1\varphi_1h'(c\varphi_1)\geq
(Mc)^p(h')^p(c\varphi_1)|\nabla\varphi_1|^p=|\nabla \overline
u_\mu|^p \quad\mbox{ in }\; \Omega\setminus\Omega_\delta.
\endeq
By \eq{dsaptep} we deduce
\neweq{dsaptispep}
\di \frac{1}{3}Mc\la_1\varphi_1h'(c\varphi_1)\geq \mu
f(x,Mh(c\varphi_1))= \mu f(x,\overline u_\mu) \quad\mbox{ in }\;
\Omega\setminus\Omega_\delta.
\endeq
Obviously, \eq{dpaispep} follows now by
\eq{dcincispep}, \eq{dsaispep} and \eq{dsaptispep}.\\
Combining \eq{calculp} with \eq{dzecep} and \eq{dpaispep} we find
that $\overline u_\mu$ is a super-solution of \eq{P}. Moreover,
$\zeta\leq \overline u_\mu$ in $\Omega.$ Applying Lemma
\ref{sub-sup}, we deduce that \eq{P} has at least one solution for
all $\mu\geq 0.$ This finishes the proof of Theorem \ref{th2edin}.
\qed

\medskip
{\it Proof of Theorem \ref{th3edin}} The proof case relies on the
same arguments used in the proof of Theorem \ref{th2edin}. In
fact, the main point is to find a super-solution $\overline
u_\la\in C^2(\Omega)\cap(\overline\Omega)$ of \eq{P}, while
$\zeta$ defined in \eq{treidoiedin} is a sub-solution. Since $g$
is nonincreasing, the inequality $\zeta\leq \overline u_\la$ in
$\Omega$ can be proved easily and the existence of solutions to
\eq{P} follows by Lemma \ref{sub-sup}.
\medskip

Define $c,\delta$ and $\theta_1,\theta_2$ as in the proof of
Theorem \ref{th2edin}. Let M satisfying \eq{dcinci} and
\eq{dsase}. Since $g(h(s))\ri+\infty$ as $s\searrow 0,$ we can
choose $\delta>0$ such that
\neweq{dpa}
\di (c\theta_1)^2g(h(c\varphi_1))-3f(x,h(c\varphi_1))>0\quad\mbox{
in }\;\Omega_\delta.
\endeq

The assumption $(f4)$ produces
$$\di
\lim_{s\ri\infty}\frac{f(x,sh(c\|\varphi_1\|_\infty))}{sh(c\|\varphi_1\|_\infty)}=0,
\quad\mbox{ uniformly for }\;\,x\in\overline\Omega.$$ Thus, we can
take $M>3$ large enough, such that
$$\di \frac{f(x,Mh(c\|\varphi_1\|_\infty))}{Mh(c\|\varphi_1\|_\infty)}<
\frac{c\la_1\theta_2h'(c\|\varphi_1\|_\infty)}{3h(c\|\varphi_1\|_\infty)}.$$
The above relation yields
\neweq{dsapte}
3f(x,Mh(c\|\varphi_1\|_\infty))<
Mc\la_1\theta_2h'(c\|\varphi_1\|_{\infty}),\quad\mbox{for all }
\;\;x\in\overline\Omega.
\endeq
Using Lemma \ref{lh}(ii) we can take $\la>0$ small enough such
that the following inequalities hold
\neweq{dopt}
3\la M^{p-1}(h')^p(c\varphi_1)\leq
g(h(c\varphi_1))(c\theta_1)^{2-p}\quad\mbox{\ in } \;\Omega_\delta
\endeq
\neweq{dnoua}
\la_1\theta_2 h'(c\|\varphi_1\|_\infty)>
 3\la (Mc)^{p-1}
(h')^p(c\theta_2)\|\nabla\varphi_1\|_\infty^p.
\endeq
For $M$ and $\la$ satisfying \eq{dcinci}-\eq{dsase} and
\eq{dpa}-\eq{dnoua}, we claim that $\overline
u_\la=Mh(c\varphi_1)$ is a super-solution of \eq{P}. First we have
\neweq{calcul}
\di -\Delta \overline
u_\la=Mc^2g(h(c\varphi_1))|\nabla\varphi_1|^2+
Mc\la_1\varphi_1h'(c\varphi_1) \quad\mbox{ in }\;\Omega.
\endeq
Arguing as in the proof of Theorem \ref{th2edin}, from \eq{dtrei},
\eq{dcinci}, \eq{dpa}, \eq{dopt} and the assumption $(f3)$ we
obtain
\neweq{dzece}
Mc^2g(h(c\varphi_1))|\nabla\varphi_1|^2\geq g(\overline
u_\la)+\la|\nabla \overline u_\la|^p+f(x,\overline
u_\la)\quad\mbox{ in }\;\Omega_\delta.
\endeq
On the other hand, \eq{dsase}, \eq{dsapte} and \eq{dnoua} gives
\neweq{dpaispe}
Mc\la_1\varphi_1h'(c\varphi_1)\geq g(\overline u_\la)+\la|\nabla
\overline u_\la|^p+f(x,\overline u_\la)\quad\mbox{ in }\;
\Omega\setminus\Omega_\delta.
\endeq
Using now \eq{calcul} and \eq{dzece}-\eq{dpaispe} we find that
$\overline u_\la$
is a super-solution of \eq{P} so our claim follows.\\
As we have already argued at the beginning of this case, we easily
get that $\zeta\leq\overline u_\la$ in $\Omega$ and by Lemma
\ref{sub-sup} we deduce that problem \eq{P} has at least one
solution if $\la>0$ is sufficiently small.

 Set
$$ \di A=\{\;\la\geq 0; \mbox{ problem }\eq{P}\mbox{ has at least one classical solution}\}.$$
From the above arguments, $A$ is nonempty. Let $\la^*=\sup A.$
First we claim that if $\la\in A,$ then $[0,\la)\subseteq A.$ For
this purpose, let $\la_1\in A$ and $0\leq \la_2<\la_1.$ If
$u_{\la_1}$ is a solution of \eq{P} with $\la=\la_1,$ then
$u_{\la_1}$ is a super-solution for \eq{P} with $\la=\la_2$ while
$\zeta$ defined in \eq{treidoiedin} is a sub-solution. Using Lemma
\ref{sub-sup} once more, we have that \eq{P} with $\la=\la_2$ has
at least one classical solution. This proves the claim. Since
$\la\in A$ was arbitrary chosen, we conclude that
$[0,\la^*)\subset A.$
\medskip

Let us assume now $p\in(1,2].$ We prove that $\la^*<+\infty.$ Set
$$\di m:=\inf_{(x,s)\in\overline\Omega\times(0,\infty)}\Big(g(s)+f(x,s)\Big).$$
Since $\lim_{s\searrow 0}g(s)=+\infty$ and the mapping $\di
(0,\infty)\ni s\longmapsto \min_{x\in\overline\Omega}f(x,s)$ is
positive and nondecreasing, we deduce that $m$ is a positive real
number. Let $\la>0$ be such that \eq{P} has a solution $u_\la.$ If
$v=\la^{1/(p-1)}u_\la,$ then $v$ verifies
\neweq{doptispe}
\left\{\begin{tabular}{ll} $\di -\Delta v\geq |\nabla
v|^p+\la^{1/(p-1)}m$ \quad & $\mbox{\rm in }\Omega,$\\
$v>0$ \quad & $\mbox{\rm in }\Omega,$\\
$v=0$ \quad & $\mbox{\rm on }\, \partial\Omega.$\\
\end{tabular} \right.
\endeq
It follows that $v$ is a super-solution of \eq{doitreiedin} for
$\di \sigma=\la^{1/(p-1)}m.$ Since 0 is a sub-solution, we obtain
that \eq{doitreiedin} has at least one classical solution for
$\sigma$ defined above. According to Lemma \ref{l4edin}, we have
$\sigma\leq \bar\sigma,$ and so $\di
\la\leq\left(\frac{\bar\sigma}{m}\right)^{p-1}.$ This means that
$\la^*$ is finite.
\medskip

Assume now $p\in(0,1)$ and let us prove that $\la^*=+\infty.$
Recall that $\zeta$ defined in \eq{treidoiedin} is a sub-solution.
To get a super-solution, we proceed in the same manner. Fix
$\la>0.$ Since $p<1$ we can find $M>1$ large enough such that
\eq{dcinci}-\eq{dsase} and \eq{dsapte}-\eq{dnoua} hold. From now
on, we follow the same steps as above. The proof of Theorem
\ref{th3edin} is now complete. \qed

\smallskip
 We remark that if $\int _0^1g(s)ds<\infty,$ then the above method can
be applied in order to extend the study of \eq{P} to the case
$\mu=1$ and $p>2.$ Indeed, by Lemma \ref{lh}(i) it follows $h\in
C^1[0,\eta].$ Using this fact, we can choose $c_1,c_2>0$ large
enough such that the conclusion of Lemma \ref{lh}(ii) holds.
Repeating the above arguments we prove that if $p>2$ then there
exists a real number $\la^*>0$ such that \eq{P} has at least one
solution if $\la<\la^*$ and no solutions exist if $\la>\la^*.$

\medskip
{\it Proof of Theorem \ref{th4edin}}. (i) If $\la=0,$ the
existence of the solution follows by using Lemma \ref{l1jmpa2}.
 Next we assume that $\la>0$
and let us fix $\mu\geq 0.$ With the change of variable $v=e^{\la
u}-1,$ the problem \eq{P} becomes
\neweq{Q}
\left\{\begin{tabular}{ll} $-\Delta v=\Phi_\la(v)$ \quad &
$\mbox{\rm
in } \Omega,$\\
$v>0$ \quad & $\mbox{\rm in } \Omega,$\\
$v=0$ \quad & $\mbox{\rm on } \partial\Omega,$\\
\end{tabular} \right.
\endeq
where
$$\di
\Phi_\la(s)=\la(s+1)g\left(\frac{1}{\la}\ln(s+1)\right)+\la\mu
(s+1),$$ for all $s\in (0,\infty).$ Obviously $\Phi_\la$ is not
monotone but we still have that the mapping $\di (0,\infty)\ni
s\mapsto \frac{\Phi_\la(s)}{s}\,$ is decreasing for all $\la>0$
and
$$\di \lim _{s\rightarrow+\infty}
\frac{\Phi_\la(s)}{s}=\la(a+\mu) \quad\mbox{ and }\quad
\lim_{s\searrow
0}\frac{\Phi_\la(s)}{s}=+\infty,$$ uniformly for $\la>0.$\\
We first remark that $\Phi_{\la}$ satisfies the hypotheses in
Lemma \ref{l1jmpa2} provided $\la(a+\mu)<\la_1.$ Hence \eq{Q} has
at least one solution.

On the other hand, since $g\geq a$ on $(0,\infty),$ we get
\neweq{Phi}
\di \Phi_\la(s)\geq\la(a+\mu)(s+1),\quad\mbox{ for all }\;
\la,s\in(0,\infty).
\endeq
Using now Lemma \ref{l01} we deduce that \eq{Q} has no solutions
if $\la(a+\mu)\geq \la_1.$ The proof of the first part in Theorem
\ref{th4edin} is therefore complete. \smallskip

(ii) We split the proof into several steps.

\noindent {\sc Step 1.} {\bf Existence of solutions.} This follows
directly from (i).

\smallskip
 {\sc Step 2.} {\bf Uniqueness of the solution.}\\
Fix $\la\geq 0.$ Let $u_1$ and $u_2$ be two classical solutions of
\eq{P} with $\la<\la^*.$ We show that $u_1\leq u_2$ in $\Omega.$
Supposing the contrary, we deduce that
$\max\limits_{\overline\Omega}\{u_1-u_2\}>0$ is achieved in a
point $x_0\in \Omega.$ This yields $\nabla (u_1-u_2)(x_0)=0$ and
$$\di 0\leq -\Delta
(u_1-u_2)(x_0)=g(u_1(x_0))-g(u_2(x_0))<0,$$ a contradiction. We
conclude that $u_1\leq u_2$ in $\Omega;$ similarly $u_2\leq u_1.$
Therefore $u_1= u_2$ in $\Omega$ and the uniqueness is proved.
\smallskip

 {\sc Step 3.} {\bf Dependence on $\la$.}
Fix $0\leq \la_1<\la_2<\la^*$ and let $u_{\la_1},$ $u_{\la_2}$ be
the unique solutions of \eq{P} with $\la=\la_1$ and $\la=\la_2$
respectively. If $\{x\in\Omega;u_{\la_1}>u_{\la_2}\}$ is nonempty,
then $\max\limits_{\overline\Omega}\{u_{\la_1}-u_{\la_2}\}>0$ is
achieved in $\Omega.$ At that point, say $\bar x,$ we have $\nabla
(u_{\la_1}-u_{\la_2})({\bar x})=0$ and
$$0\leq -\Delta(u_{\la_1}-u_{\la_2})(\bar x)=g(u_{\la_1}(\bar x))-g(u_{\la_2}(\bar x))
+(\la_1-\la_2)|\nabla u_{\la_1}|^p(\bar x)<0,$$
which is a contradiction.\\
Hence $u_{\la_1}\leq u_{\la_2}$ in $\overline\Omega.$ The maximum
principle also gives $u_{\la_1}< u_{\la_2}$ in $\Omega.$
\smallskip

 {\sc Step 4.} {\bf Regularity.} We fix $0<\la<\la^*,$ $\mu>0$
and assume that $\limsup_{s\searrow 0}s^{\alpha}g(s)<+\infty.$
This means that $g(s)\leq cs^{-\alpha}$ in a small positive
neighborhood of the origin. To prove the regularity, we will use
again the change of variable $v=e^{\la u}-1.$ Thus, if $u_\la$ is
the unique solution of \eq{P}, then $v_\la=e^{\la u_\la}-1$ is the
unique solution of \eq{Q}. Since $\di \lim_{s\searrow
0}\frac{e^{\la s}-1}{s}=\la,$ we conclude that (ii1) and (ii2) in
Theorem \ref{th4edin} are established if we prove

${\rm (a)}\quad \tilde c_1\,{\rm dist}(x,\partial \Omega) \leq
v_\la(x)\leq \tilde c_2\,{\rm dist} (x,\partial \Omega)$ in
$\Omega,$ for some positive constants $\tilde c_1,\tilde c_2>0.$

${\rm (b)}\quad v_\la\in C^{1,1-\alpha}(\overline \Omega).$\\
{\it Proof of} (a). By the monotonicity of $g$ and the fact that
$g(s)\leq cs^{-\alpha}$ near the origin, we deduce the existence
of $A,B,C>0$ such that
\neweq{ctreiedin}
\di \Phi_\la(s)\leq As+Bs^{-\alpha}+C,\quad \mbox{ for
all}\;0<\la<\la^*\mbox { and }\,s>0.
\endeq
Let us fix $m>0$ such that $m\la_1\|\varphi_1\|_{\infty}<\la\mu.$
Combining this with \eq{Phi} we deduce
\neweq{sssu}
\di -\Delta(v_\la-m\varphi_1)=\Phi_\la(v_\la)-m\la_1\varphi_1\geq
\la\mu-m\la_1\varphi_1\geq 0
\endeq
in $\Omega.$ Since $v_\la-m\varphi_1=0$ on $\partial\Omega,$ we
conclude
\neweq{fi2}
v_\la\geq m\varphi_1\quad\mbox{ in }\;\Omega.
\endeq
Now, \eq{fi2} and \eq{fi} imply $v_\la\geq \tilde c_1\,{\rm dist}
(x,\partial \Omega)$ in $\Omega,$ for some positive constant
$\tilde c_1>0.$ The first inequality in the statement of (a) is
therefore established. For the second one, we apply an idea found
in Gui and Lin \cite{gl}. Using \eq{fi2} and the estimate
\eq{ctreiedin}, by virtue of Lemma \ref{l1jde} we deduce
$\Phi_\la(v_\la)\in L^1(\Omega),$ that is, $\Delta v_\la\in
L^1(\Omega).$

Using the smoothness of $\partial\Omega,$ we can find
$\delta\in(0,1)$ such that for all
$x_0\in\Omega_{\delta}:=\{x\in\Omega\,;\,{\rm
dist}(x,\partial\Omega)\leq \delta\},$ there exists
$y\in\RR^N\setminus\overline\Omega$ with ${\rm
dist}(y,\partial\Omega)=\delta$ and ${\rm
dist}(x_0,\partial\Omega)=|x_0-y|-\delta.$

Let $K>1$ be such that diam$\,(\Omega)<(K-1)\delta$ and let $\xi$
be the unique solution of the Dirichlet problem
$$
\left\{\begin{tabular}{ll} $-\Delta \xi=\Phi_\la(\xi)$ \quad &
${\rm
in}\ B_K(0)\setminus B_1(0),$\\
$\xi>0$ \quad & ${\rm in}\ B_K(0)\setminus B_1(0),$\\
$\xi=0$ \quad & ${\rm on}\ \partial(B_K(0)\setminus B_1(0)).$\\
\end{tabular} \right.
$$
where $B_r(0)$ denotes the open ball in $\RR^N$ of radius $r$ and
centered at the origin. By uniqueness, $\xi$ is radially
symmetric. Hence $\xi(x)=\tilde \xi(|x|)$ and
\neweq{btreiedin}
\left\{\begin{tabular}{ll} $\di \tilde \xi''+\frac{N-1}{r}\tilde
\xi'+\Phi_\la(\tilde\xi)=0$
\quad & ${\rm in}\ (1,K),$\\
$\tilde \xi>0$ \quad & ${\rm in}\ (1,K),$\\
$\tilde \xi(1)=\tilde \xi(K)=0.$ \quad & \\
\end{tabular} \right.
\end{equation}
Integrating in \eq{btreiedin} we have
$$\begin{aligned}
\di \tilde \xi'(t)& =\tilde \xi'(a)a^{N-1}t^{1-N}-t^{1-N}\int_a^t
r^{N-1} \Phi_\la(\tilde \xi(r))dr\\
& =\tilde \xi'(b)b^{N-1}t^{1-N}+t^{1-N}\int_t^b
r^{N-1} \Phi_\la(\tilde \xi(r))dr,\\
\end{aligned}$$
where $1<a<t<b<K.$ With the same arguments as above we have
$\Phi_\la(\tilde \xi)\in L^1(1,K)$ which implies that both $\tilde
\xi(1)$ and $\tilde \xi(K)$ are finite. Hence $\tilde \xi\in
C^2(1,K)\cap C^1[1,K].$ Furthermore,
\neweq{bpatruedin}
\xi(x)\leq \tilde C\min\{K-|x|,|x|-1\}, \quad\mbox{ for any }
\;\;x\in B_K(0)\setminus B_1(0).
\endeq
Let us fix $x_0\in\Omega_{\delta}.$ Then we can find
$y_0\in\RR^N\setminus\overline\Omega$ with ${\rm
dist}(y_0,\partial\Omega)=\delta$ and ${\rm
dist}(x_0,\partial\Omega)=|x_0-y|-\delta.$ Thus, $\Omega\subset
B_{K\delta}(y_0)\setminus B_{\delta}(y_0).$ Define $\di \overline
v(x)=\xi\left(\frac{x-y_0}{\delta}\right),$ for all
$x\in\overline\Omega.$ We show that $\overline v$ is a
super-solution of \eq{Q}. Indeed, for all $x\in\Omega$ we have
$$\begin{tabular}{rl}
$\di\Delta \overline v+\Phi_\la(\overline v)$ &
$\di=\frac{1}{\delta^2}\left(\tilde \xi''+
\frac{N-1}{r}\tilde \xi'\right)+\Phi_\la(\tilde \xi)$\\
&$\di \leq\frac{1}{\delta^2}\left(\tilde \xi''+
\frac{N-1}{r}\tilde \xi'+\Phi_\la(\tilde \xi)\right)$\\
&$=0,$\\
\end{tabular}$$
where $\di r=\frac{|x-y_0|}{\delta}.$ We have obtained that
$$\di \Delta \overline v+\Phi_\la(\overline v)\leq 0\leq \Delta v_\la+\Phi_\la(v_\la)
\quad\mbox{ in }\;\Omega,$$
$$\di \overline v,v_\la>0\;\;\mbox{ in }\;\Omega,\;\overline v=v_\la
\;\;\mbox{ on }\;\partial\Omega $$
$$\Delta v_\la\in L^1(\Omega).$$
By Lemma \ref{l2} we get $\di v_{\la}\leq \overline v$ in
$\Omega.$ Combining this with \eq{bpatruedin} we obtain
$$\di v_{\la}(x_0)\leq \overline v(x_0)\leq
\tilde C\min\left\{K-\frac{|x_0-y_0|}{\delta},
\frac{|x_0-y_0|}{\delta}-1\right\}\leq\frac{\tilde C}{\delta}{\rm
dist}(x_0,\partial\Omega).$$
Hence $v_{\la}\leq \frac{\tilde
C}{\delta}{\rm dist}(x,\partial\Omega)$ in $\Omega_{\delta}$ and
the second inequality in the statement of (a) follows. \smallskip

 {\it Proof of} (b). Let $G$ be the Green's function associated
with the Laplace operator in $\Omega.$ Then, for all $x\in\Omega$
we have
$$\di v_{\la}(x)=-\int_{\Omega} G(x,y)\Phi_\la(v_{\la}(y))dy$$
and
$$\di \nabla v_{\la}(x)=-\int_{\Omega} G_x(x,y)\Phi_\la(v_{\la}(y))dy.$$
If $x_1,x_2\in\Omega,$ using \eq{ctreiedin} we obtain
$$\begin{tabular}{ll}
$\di |\nabla v_{\la}(x_1)-\nabla v_{\la}(x_2)|$&$\di
\leq \int_{\Omega} |G_x(x_1,y)-G_x(x_2,y)|\cdot(Av_\la+C)dy$\\
&$\di \;\;\,+B\int_{\Omega}
|G_x(x_1,y)-G_x(x_2,y)|\cdot v_{\la}^{-\alpha}(y)dy.$\\
\end{tabular}$$
Now, taking into account that $v_{\la}\in C(\overline\Omega),$ by
the standard regularity theory (see Gilbarg and Trudinger
\cite{gt}) we get
$$\di \int_{\Omega} |G_x(x_1,y)-G_x(x_2,y)|\cdot
(Av_\la+C)dy\leq \tilde c_1|x_1-x_2|.$$ On the other hand, with
the same proof as in \cite[Theorem 1]{gl}, we deduce
$$\di \int_{\Omega} |G_x(x_1,y)-G_x(x_2,y)|\cdot
v_{\la}^{-\alpha}(y)\leq \tilde c_2 |x_1-x_2|^{1-\alpha}.$$ The
above inequalities imply $u_{\la}\in C^2(\Omega)\cap
C^{1,1-\alpha}({\overline{\Omega}}).$

{\sc Step 5.} {\bf Asymptotic behaviour of the solution.} This
follows with the same lines as in the proof of Theorem \ref{th2}.
 \qed

\medskip
We are concerned in what follows with the closely related
Dirichlet problem
$$\left\{\begin{tabular}{ll}
$-\Delta u+K(x)g(u)+|\nabla u|^a=\lambda f(x,u)$ \quad & ${\rm
in}\ \Omega,$\\
$u>0$ \quad & ${\rm in}\ \Omega,$\\
$u=0$ \quad & ${\rm on}\ \partial\Omega,$\\
\end{tabular} \right. \eqno(1)_\lambda$$
where $\Omega$ is a smooth bounded domain in $\RR^N$ ($N\geq 2$),
$\lambda>0,$ $0<a\leq 2$ and $K\in C^{0,\gamma}(\overline\Omega)$,
$0<\gamma<1.$ We assume from now on that
$f:\overline{\Omega}\times[0,\infty)\rightarrow[0,\infty)$ is a
H\"{o}lder continuous function which is positive on
$\overline{\Omega}\times(0,\infty)$ such that $f$ is nondecreasing
with respect to the second variable and is sublinear, in the sense
that the mapping $$ (0,\infty)\ni
s\longmapsto\frac{f(x,s)}{s}\quad\mbox{is nonincreasing for
all}\;\, x\in\overline{\Omega}$$ and $$\lim_{s\ri
0^+}\frac{f(x,s)}{s}=+\infty\quad\mbox{and}\;\;
\lim_{s\rightarrow\infty}\frac{f(x,s)}{s}=0,\;\;\mbox{uniformly
for}\;\,x\in\overline{\Omega}.$$

 We also assume that $g\in C^{0,\gamma}(0,\infty)$ is a nonnegative
and nonincreasing function satisfying $$ \lim_{s\ri
0^+}g(s)=+\infty.$$

Problem $(1)_\lambda$ has been considered in Section~7 in the
absence of the gradient term $|\nabla u|^a$ and assuming that the
singular term $g(t)$ behaves like $t^{-\alpha}$ around the origin,
with $t\in (0,1)$. In this case it has been shown that the sign of
the extremal values of $K$ plays a crucial role. In this sense, we
have proved in Section~7 that if $K<0$ in $\overline\Omega,$ then
problem $(1)_\lambda$ (with $a=0$) has a unique solution in the
class ${\cal E}=\{u\in C^2(\Omega)\cap
C({\overline{\Omega}});\,\,g(u)\in L^1(\Omega)\},$ for all
$\lambda>0.$ On the other hand, if $K>0$ in $\overline\Omega$,
then there exists $\lambda^*$ such that problem $(1)_\lambda$ has
solutions in $\cal E$ if $\lambda>\lambda^*$ and no solution
exists if $\lambda<\lambda^*$. The case where $f$ is
asymptotically linear, $K\leq 0$, and $a=0$ has been discussed in
Section~6. In this framework, a major role is played by
$\lim_{s\ri\infty}f(s)/s=m>0.$ More precisely, there exists a
solution (which is unique) $u_\lambda\in C^2(\Omega)\cap
C^{1}(\overline\Omega)$
 if and only if $\lambda<\lambda^*:=\lambda_1/m.$ An additional result asserts that the mapping
$(0,\lambda^*)\longmapsto u_\lambda$  is increasing and
$\lim_{\lambda\nearrow \lambda^*} u_\lambda=+\infty$ uniformly on
compact subsets of $\Omega.$

Due to the singular character of our problem $(1)_\lambda,$ we
cannot expect to have solutions in $C^2(\overline\Omega).$ We are
seeking in this paper classical solutions of $(1)_\lambda,$ that
is, solutions $u\in C^2(\Omega)\cap C(\overline\Omega)$ that
verify $(1)_\lambda.$ Closely related to our problem is the
following one, which has been considered in the first part of this
Section:

\neweq{edinburgh}
\left\{\begin{tabular}{ll}
$-\Delta u=g(u)+|\nabla u|^a+\lambda f(x,u)$ \quad & $\mbox{\rm in}\ \Omega,$\\
$u>0$ \quad & $\mbox{\rm in}\ \Omega,$\\
$u=0$ \quad & $\mbox{\rm on}\ \partial\Omega,$\\
\end{tabular} \right.
\endeq
where $f$ and $g$ verifies the above assumptions. We recall that
we have proved that if $0<a<1$ then problem \eq{edinburgh} has at
least one classical solution for all $\lambda\geq 0.$ In turn, if
$1<a\leq 2,$ then problem \eq{edinburgh} has no solutions for
large values of $\lambda>0.$

The existence results for our problem $(1)_\lambda$ are quite
different to those of \eq{edinburgh} presented in the first part
of this Section. More exactly, we prove in what follows that
problem $(1)_\lambda$ has at least one solution only when
$\lambda>0$ is large enough and $g$ satisfies a naturally growth
condition around the origin. Thus, we extend the results in
Barles, G. D\'iaz, and J. I. D\'iaz \cite[Theorem 1]{bdd},
corresponding to $K\equiv 0,$ $f\equiv f(x)$ and $a\in[0,1).$

The main difficulty in the treatment of $(1)_\lambda$ is the lack
of the usual maximal principle between super and sub-solutions,
due to the singular character of the equation. To overcome it, we
state an improved comparison principle that fit to our problem
$(1)_\lambda$ (see Lemma \ref{l} below).

In our first result we assume that $K<0$ in $\Omega.$ Note that
$K$ may vanish on $\partial\Omega$ which leads us to a competition
on the boundary between the potential $K(x)$ and the singular term
$g(u).$ We prove the following result.

\begin{thm}\label{th0jmaa}
Assume that $K<0$ in $\Omega.$ Then, for all $\lambda>0$, problem
$(1)_\lambda$ has at least one classical solution.
\end{thm}

Next, we assume that $K>0$ in $\overline\Omega.$ In this case, the
existence of a solution to $(1)_\lambda$ is closely related to the
decay rate around its singularity. In this sense, we prove that
problem $(1)_\lambda$ has no solution, provided that $g$ has a
``strong" singularity at the origin. More precisely, we have

\begin{thm}\label{th1jmaa}
Assume that   $K>0$ in $\overline\Omega$ and $\int^1_0
g(s)ds=+\infty$. Then problem $(1)_\lambda$ has no classical
solutions.
\end{thm}

In the following result, assuming that $\int^1_0 g(s)ds<+\infty$,
we show that problem $(1)_\lambda$ has at least one solution,
provided that $\lambda>0$ is large enough. More precisely, we
prove

\begin{thm}\label{th2jmaa}
Assume that $K>0$ in $\overline\Omega$ and $\int^1_0
g(s)ds<+\infty.$ Then there exists $\lambda^*>0$ such that problem
$(1)_\lambda$ has at least one classical solution if
$\lambda>\lambda^*$ and no solution exists if $\lambda<\lambda^*.$
\end{thm}

A very useful auxiliary result in the proofs of the above theorems
is the following comparison principle that improves Lemma
\ref{l2}. Our proof uses some ideas from Shi and Yao \cite{shi},
that go back to the pioneering work by Brezis and Kamin \cite{bk}.

\begin{lem}\label{l}
Let $\Psi:\overline{\Omega}\times(0,\infty)\rightarrow\RR$ be a
continuous function such that the mapping  $\di(0,\infty)\ni
s\longmapsto\frac{\Psi(x,s)}{s}$ is strictly decreasing at each
$x\in\Omega.$ Assume that there exists  $v$, $w\in C^2(\Omega)\cap
C({\overline{\Omega}})$ such that

$(a)\qquad \Delta w+\Psi(x,w)\leq 0\leq \Delta v+\Psi(x,v)$ in
$\Omega;$

$(b)\qquad v,w>0$ in $\Omega$ and $v\leq w$ on $\partial\Omega;$

$(c)\qquad\Delta v\in L^1(\Omega)\;\mbox{ or }\;
\Delta w\in L^1(\Omega).$\\
Then $v\leq w$ in $\Omega.$
\end{lem}

 {\it Proof.} We argue by contradiction and assume that
$v\geq w$ is not true in $\Omega.$ Then, we can find
$\ep_0,\delta_0>0$ and a ball $B\subset\subset\Omega$ such that
$v-w\geq \ep_0$  in $B$ and
\neweq{shi2}
\di \int_Bvw\left(\frac{\Psi(x,w)}{w}-\frac{\Psi(x,v)}{v}\right)
dx\geq \delta_0.
\endeq
The case $\Delta v\in L^1(\Omega)$ was stated in Lemma \ref{l2}.
Let us assume now that $\Delta w\in L^1(\Omega)$ and set
$M=\max\{1,\|\Delta w\|_{L^1(\Omega)}\}$,
$\ep=\min\left\{1,\ep_0,2^{-2}\delta_0/M\right\}.$ Consider  a
nondecreasing function $\theta\in C^1(\RR)$ such that
$\theta(t)=0,$ if $t\leq 1/2,$ $\theta(t)=1,$ if $t\geq 1,$ and
$\theta(t)\in(0,1)$ if $t\in(1/2,1).$ Define
$$\di \theta_\ep(t)=\theta\left(\frac{t}{\ep}\right),\quad t\in\RR.$$
Since $w\geq v$ on $\partial\Omega,$ we can find a smooth
subdomain $\Omega^*\subset\subset\Omega$ such that
$$B\subset\Omega^*\quad\mbox{ and }\;
v-w<\frac{\ep}{2}\;\mbox{ in }\,\Omega\setminus\Omega^*.$$ Using
the hypotheses (a) and (b) we deduce
\neweq{shi3}
\di \int_{\Omega^*}(w\Delta v-v\Delta w)\theta_\ep(v-w)dx
\geq\int_{\Omega^*}
vw\left(\frac{\Psi(x,w)}{w}-\frac{\Psi(x,v)}{v}\right)
\theta_\ep(v-w)dx.
\endeq
By \eq{shi2} we have
$$\begin{tabular}{ll}
$\di\int_{\Omega^*}
vw\left(\frac{\Psi(x,w)}{w}-\frac{\Psi(x,v)}{v}\right)
\theta_\ep(v-w)dx$&$\di\geq \int_{B}
vw\left(\frac{\Psi(x,w)}{w}-\frac{\Psi(x,v)}{v}\right)
\theta_\ep(v-w)dx$\\
&$\di=\int_{B}
vw\left(\frac{\Psi(x,w)}{w}-\frac{\Psi(x,v)}{v}\right)dx\geq
\delta_0.$
\end{tabular}$$
To raise a contradiction, we need only to prove that the left-hand
side in \eq{shi3} is smaller than $\delta_0.$ For this purpose, we
define
$$\di \Theta_\ep(t)=\int_0^ts\theta'_\ep(s)ds,\quad t\in\RR.$$
It is easy to see that
\neweq{shi4}
\Theta_\ep(t)=0, \;\mbox{ if }\,t<\frac{\ep}{2}\quad\mbox{ and
}\,\; 0\leq \Theta_\ep(t)\leq 2\ep,\;\mbox{ for all }t\in\RR.
\endeq
Now, using the Green theorem, we evaluate the left-hand side of
\eq{shi3}:
$$\begin{tabular}{ll}
&$\di \int_{\Omega^*}(w\Delta v-v\Delta w)\theta_\ep(v-w)dx$\\
$=$&$\di \int_{\partial\Omega^*}w\theta_\ep(v-w)\frac{\partial v}
{\partial n}d\sigma-
\int_{\Omega^*}(\nabla w\cdot\nabla v)\theta_\ep(v-w)dx$\\
&$\di-\int_{\Omega^*}w\theta'_\ep(v-w)\nabla v \cdot\nabla(
v-w)dx- \int_{\partial\Omega^*}v\theta_\ep(v-w)\frac{\partial w}
{\partial n}d\sigma$\\
&$\di+\int_{\Omega^*}(\nabla w\cdot \nabla v)\theta_\ep(v-w)dx+
\int_{\Omega^*}v\theta'_\ep(v-w)\nabla w\cdot\nabla(v-w)dx$\\
$=$&$\di\int_{\Omega^*}\theta'_\ep(v-w)(v\nabla w-w\nabla v)
\cdot\nabla( v-w)dx.$\\
\end{tabular}$$
The above relation can also be rewritten as
$$\begin{tabular}{ll}
$\di \int_{\Omega^*}(w\Delta v-v\Delta w)\theta_\ep(v-w)dx=$&
$\di\int_{\Omega^*}w\theta'_\ep(v-w)\nabla(w-v)
\cdot\nabla(v-w)dx$\\
&$\di+\int_{\Omega^*}(v-w)\theta'_\ep(v-w)\nabla w
\cdot\nabla(v-w)dx.$\\
\end{tabular}$$
Since $\di \int_{\Omega^*}w\theta'_\ep(v-w)\nabla(w-v)
\cdot\nabla(v-w)dx\leq 0,$ the last equality yields
$$\di \int_{\Omega^*}(w\Delta v-v\Delta w)\theta_\ep(v-w)dx\leq
\di\int_{\Omega^*}(v-w)\theta'_\ep(v-w)\nabla w
\cdot\nabla(v-w)dx,$$ that is,
$$\di \int_{\Omega^*}(w\Delta v-v\Delta w)\theta_\ep(v-w)dx\leq
\di\int_{\Omega^*}\nabla w\cdot \nabla(\Theta_\ep(v-w))dx.$$ Again
by Green's first formula and by \eq{shi4} we have
$$\begin{tabular}{ll}
$\di \int_{\Omega^*}(w\Delta v-v\Delta
w)\theta_\ep(v-w)dx$&$\di\leq
\di\int_{\partial\Omega^*}\Theta_\ep(v-w) \frac{\partial
v}{\partial n}d\sigma-
\di\int_{\Omega^*}\Theta_\ep(v-w)\Delta wdx$\\
&$\di\leq  - \int_{\Omega^*}\Theta_\ep(v-w)\Delta wdx\leq
2\ep\int_{\Omega^*}|\Delta w|dx$\\
&$\di\leq 2\ep M<\frac{\delta_0}{2}.$\\
\end{tabular}$$
Thus, we have obtained a contradiction. Hence $v\leq w$ in
$\Omega$ and the proof of Lemma \ref{l} is now complete. \qed

\medskip We are now ready to prove our main results.

\medskip
{\it Proof of Theorem \ref{th0jmaa}}. Fix $\lambda>0.$ Obviously,
$\Psi(x,s)=\lambda f(x,s)-K(x)g(s)$ satisfies the hypotheses in
Lemma \ref{l1jmpa2} since $K<0$ in $\Omega.$ Hence, there exists a
solution $\overline u_\lambda$ of the problem
$$
\left\{\begin{tabular}{ll}
$-\Delta u=\lambda f(x,u)-K(x)g(u)$ \quad & ${\rm in}\ \Omega,$\\
$u>0$ \quad & ${\rm in}\ \Omega,$\\
$u=0$ \quad & ${\rm on}\ \partial\Omega.$\\
\end{tabular}\right.
$$
We observe that $\overline u_\lambda$ is a super-solution of
problem $(1)_\lambda.$ To find a sub-solution, let us denote
$$\di p(x)=\min\{\lambda f(x,1);-K(x)g(1)\}, \quad x\in\overline\Omega.$$
Using the monotonicity of $f$ and $g,$ we observe that $p(x)\leq
\lambda f(x,s)-K(x)g(s)$ for all $(x,s)\in\Omega\times
(0,\infty).$ We now consider the problem
\neweq{subsol}
\left\{\begin{tabular}{ll}
$-\Delta v+|\nabla v|^a=p(x)$ \quad & ${\rm in}\ \Omega,$\\
$v=0$ \quad & ${\rm on}\ \partial\Omega.$\\
\end{tabular}\right. \endeq
First, we observe that $v=0$ is a sub-solution of \eq{subsol}
while $w$ defined by
$$\left\{\begin{tabular}{ll}
$-\Delta w=p(x)$ \quad & ${\rm in}\ \Omega,$\\
$w=0$ \quad & ${\rm on}\ \partial\Omega,$\\
\end{tabular}\right. $$
is a super-solution. Since $p>0$ in $\Omega$ we deduce that $w\geq
0$ in $\Omega.$ Thus, the problem \eq{subsol} has at least one
classical solution $v.$ We claim that $v$ is positive in $\Omega.$
Indeed, if $v$ has a minimum in $\Omega,$ say at $x_0,$ then
$\nabla v(x_0)=0$ and $\Delta v(x_0)\geq 0.$ Therefore
$$\di 0\geq -\Delta v(x_0)+|\nabla v|^a(x_0)=p(x_0)>0,$$
which is a contradiction. Hence $\min_{x\in\overline\Omega} v=
\min_{x\in\partial\Omega} v=0,$ that is, $v>0$ in $\Omega.$ Now
$\underline u_\lambda=v$ is a sub-solution of $(1)_\lambda$ and we
have
$$\di -\Delta \underline u_\lambda=p(x)\leq \lambda f(x,\overline u_\lambda)-K(x)g(\overline u_\lambda)=
-\Delta \overline u_\lambda\quad\mbox{ in }\,\Omega.$$ Since
$\underline u_\lambda=\overline u_\lambda=0$ on $\partial\Omega,$
from the above relation we may conclude that $\underline
u_\lambda\leq \overline u_\lambda$ in $\Omega$ and so, there
exists at least one classical solution for $(1)_\lambda.$ The
proof of Theorem \ref{th0jmaa} is now complete. \qed

\medskip
{\it Proof of Theorem \ref{th1jmaa}}. We  give a direct proof,
without using any change of variable, as in Zhang \cite{z1}. Let
us assume that there exists $\lambda>0$ such that the problem
$(1)_\lambda$ has a classical solution $u_\la.$ By our hypotheses
on $f$, we deduce by Lemma \ref{l1jmpa2} that for all $\la>0$
there exists $U_\la\in C^2(\overline\Omega)$ such that
\begin{equation}\label {UU}
 \left\{\begin{tabular}{ll}
$-\Delta U_\la=\la f(x,U_\la)$ \quad & ${\rm in}\
\Omega,$\\
$U_\la>0$ \quad & ${\rm in}\ \Omega,$\\
$U_\la=0$ \quad & ${\rm on}\ \partial\Omega.$\\
\end{tabular} \right.
\end{equation}
Moreover, there exist $c_1,c_2>0$ such that
\neweq{udst}
c_1\,\mbox{dist}\,(x,\partial\Omega)\leq U_{\la}(x)\leq c_2\,
\mbox{dist}\,(x,\partial\Omega)\quad\mbox{for all} \;x\in\Omega.
\endeq
Consider the perturbed problem
\begin{equation}\label {epsilon}
 \left\{\begin{tabular}{ll}
$-\Delta u+K_*g(u+\ep)=\la f(x,u)$ \quad &
${\rm in}\ \Omega,$\\
$u>0$ \quad & ${\rm in}\ \Omega,$\\
$u=0$ \quad & ${\rm on}\ \partial\Omega,$\\
\end{tabular} \right.
\end{equation}
where $K_*=\min_{x\in\overline\Omega}K(x)>0.$ It is clear that
$u_{\la}$ and $U_{\la}$ are respectively sub and super-solution of
(\ref{epsilon}). Furthermore, we have
$$\di \Delta U_\la+f(x,U_\la)\leq 0\leq \Delta u_\la+f(x,u_\la)\quad\mbox{ in }\Omega,$$
$$\di U_\la,u_\la>0\quad\mbox{ in }\Omega,$$
$$U_\la=u_\la=0\quad\mbox{ on }\partial\Omega,$$
$$\Delta U_\la\in L^1(\Omega)\;\,(\mbox{ since }\,
U_\la\in C^2(\overline\Omega)).$$ In view of Lemma \ref{l} we get
$u_\la\leq U_\la$ in $\Omega.$ Thus, a standard bootstrap argument
(see Gilbarg and Trudinger \cite{gt}) implies that there exists a
solution $u_{\ep}\in C^{2}(\overline{\Omega})$ of (\ref{epsilon})
such that
$$
u_{\la}\leq u_{\ep}\leq U_{\la} \quad\mbox{in}\;\;\Omega.
$$
Integrating in (\ref{epsilon}) we obtain
$$\di -\int_{\Omega}\Delta
u_{\ep}dx+K_*\int_{\Omega}g(u_{\ep}+\ep)dx=\la
\int_{\Omega}f(x,u_{\ep})dx.$$ Hence
\begin{equation}\label{33jmaa}
\di-\int_{\partial\Omega}\frac{\partial u_{\ep}}{\partial n}ds+K_*
\int_{\Omega}g(u_{\ep}+\ep)dx\leq M,
\end{equation}
where $M>0$ is a positive constant. Taking into account the fact
that $\di\frac{\partial u_{\ep}}{\partial n}\leq 0$ on
$\partial\Omega,$ relation (\ref{33jmaa}) yields $\di
K_*\int_{\Omega}g(u_{\ep}+\ep)dx\leq M.$ Since $u_\ep\leq U_{\la}$
in $\overline\Omega,$ from the last inequality we can conclude
that $\di \int_{\Omega}g(U_{\la}+\ep)dx\leq C,$ for some $C>0.$
Thus, for any compact subset $\omega\subset\subset\Omega$ we have
$$\di \int_{\omega}g(U_{\la}+\ep)dx\leq C.$$
Letting $\ep\ri 0^+,$ the above relation produces $\di
\int_{\omega}g(U_{\la})dx\leq C.$ Therefore
\begin{equation}\label{34jmaa}
\di \int_{\Omega}g(U_{\la})dx\leq C.
\end{equation}
On the other hand, using \eq{udst} and the hypothesis
$\int_0^1g(s)ds=+\infty,$ it follows
$$\di \int_{\Omega}g(U_{\la})dx\geq
\int_{\Omega}g(c_2\mbox{dist}\,(x,\partial\Omega))dx=+\infty,$$
which contradicts (\ref{34jmaa}). Hence, $(1)_{\la}$ has no
classical solutions and the proof of Theorem \ref{th1jmaa} is now
complete. \qed

\medskip
{\it Proof of Theorem \ref{th2jmaa}}. Fix $\la>0.$ We first note
that $U_\la$ defined in \eq{UU} is a super-solution of $(1)_\la.$
We now focuss on finding a sub-solution $\underline u_\la$ such
that $\underline u_\la\leq U_\la$ in $\Omega.$

Let $h:[0,\infty)\rightarrow[0,\infty)$ be such that
\begin{equation}\label{25jmaa}
\left\{\begin{tabular}{ll}
$h''(t)=g(h(t)),\quad \mbox{ for all }t>0,$\\
$h>0,\quad \mbox{ in }(0,\infty),$\\
$h(0)=0.$\\
\end{tabular} \right.
\end{equation}
Multiplying by $h'$ in \eq{25jmaa} and then integrating over
$[s,t]$  we have
$$(h')^2(t)-(h')^2(s)=2\int^{h(t)}_{h(s)}g(\tau)d\tau,\quad
\mbox{ for all }\,t>s>0.$$ Since $\int^1_0g(\tau)d\tau<\infty,$
from the above equality we deduce that we can extend $h'$ in
origin by taking $h'(0)=0$ and so $h\in C^2(0,\infty)\cap
C^1[0,\infty).$ Taking into account the fact that $h'$ is
increasing and $h''$ is decreasing on $(0,\infty),$  the mean
value theorem implies that
$$\di\frac{h'(t)}{t}=\frac{h'(t)-h'(0)}{t-0}\geq h''(t),\quad\mbox{ for all }\,t>0.$$
Hence $h'(t)\geq th''(t),$ for all $t>0.$ Integrating in the last
inequality we get
\neweq{has}
th'(t)\leq 2h(t),\quad\mbox{ for all }\,t>0.
\endeq

Let $\phi_1$ be the normalized positive eigenfunction
corresponding to the first eigenvalue $\la_1$ of the problem
$$
 \left\{\begin{tabular}{ll}
$-\Delta u=\la u$ \quad & ${\rm in}\ \Omega,\,$\\
$u=0$ \quad & ${\rm on}\ \partial\Omega\,.$\\
\end{tabular} \right.
$$
It is well known that $\phi_1\in C^2(\overline\Omega).$
Furthermore, by Hopf's maximum principle there exist $\delta>0$
and $\Omega_0\subset\subset\Omega$ such that $\di |\nabla
\phi_1|\geq\delta$ in $\Omega\setminus\Omega_0.$ Let
$M=\max\{1,2K^*\delta^{-2}\},$ where
$K^*=\max_{x\in\overline\Omega}K(x).$ Since
$$\lim_{{\rm dist}\,(x,\partial\Omega)\ri 0^+}
\Big\{-K^*g(h(\phi_1))+M^a(h')^a(\phi_1)|\nabla\phi_1|^a\Big\}=-\infty,$$
by letting $\Omega_0$ close enough to the boundary of $\Omega$ we
can assume that
\neweq{omegaz}
-K^*g(h(\phi_1))+M^a(h')^a(\phi_1)|\nabla\phi_1|^a<0\quad\mbox{ in
}\;\Omega\setminus\Omega_0.
\endeq
We now are able to show that $\underline{u}_{\la}=Mh(\phi_1)$ is a
sub-solution of $(1)_\la$ provided $\la>0$ is sufficiently large.
Using the monotonicity of $g$ and \eq{has} we have
\begin{equation}\label{29jmaa}
\begin{tabular}{ll}
$\di -\Delta
\underline{u}_{\la}+K(x)g(\underline{u}_{\la})+|\nabla\underline u_\la|^a=$\\
$\quad\leq -Mg(h(\phi_1))|\nabla\phi_1|^2+\la_1Mh'(\phi_1)\phi_1+
K^*g(Mh(\phi_1))+M^a(h')^a(\phi_1)|\nabla\phi_1|^a$\\
$\quad\leq
g(h(\phi_1))(K^*-M|\nabla\phi_1|^2)+\la_1Mh'(\phi_1)\phi_1+
M^a(h')^a(\phi_1)|\nabla\phi_1|^a$\\
$\quad\leq g(h(\phi_1))(K^*-M|\nabla\phi_1|^2)+2\la_1Mh(\phi_1)+
M^a(h')^a(\phi_1)|\nabla\phi_1|^a.$\\
\end{tabular}
\end{equation}
The definition of $M$ and \eq{omegaz} yield
\neweq{bord}
\di -\Delta
\underline{u}_{\la}+K(x)g(\underline{u}_{\la})+|\nabla\underline
u_\la|^a \leq 2\la_1Mh(\phi_1)=2\la_1\underline u_\la\quad\mbox{
in }\;\Omega\setminus\Omega_0.\endeq Let us choose $\la>0$ such
that
\neweq{lambda1}
\di
\la\frac{\min_{x\in\overline\Omega_0}f(x,Mh(\|\phi_1\|_\infty))}
{M\|\phi_1\|_\infty}\geq 2\la_1.
\endeq
Then, by virtue of the assumptions on $f$ and using \eq{lambda1},
we have
$$\di \la\frac{f(x,\underline u_\la)}{\underline u_\la}\geq
\la\frac{f(x,Mh(\|\phi_1\|_\infty))}{M\|\phi_1\|_\infty}\geq
2\la_1 \quad\mbox{ in }\Omega\setminus\Omega_0.$$ The last
inequality combined with \eq{bord} yield
\neweq{bordfinal}
\di -\Delta
\underline{u}_{\la}+K(x)g(\underline{u}_{\la})+|\nabla\underline
u_\la|^a \leq 2\la_1 \underline u_\la\leq \la f(x,\underline
u_\la)\quad\mbox{ in }\Omega\setminus\Omega_0.\endeq On the other
hand, from \eq{29jmaa} we obtain
\neweq{interion}
\di -\Delta
\underline{u}_{\la}+K(x)g(\underline{u}_{\la})+|\nabla\underline
u_\la|^a \leq K^*g(h(\phi_1))+2\la_1Mh(\phi_1)+
M^a(h')^a(\phi_1)|\nabla\phi_1|^a\quad\mbox{ in }\Omega_0.
\end{equation}
Since $\phi_1>0$ in $\overline\Omega_0$ and $f$ is positive on
$\overline \Omega_0\times(0,\infty),$ we may choose $\la>0$ such
that
\neweq{lambda2}
\di \la\min_{x\in\overline\Omega_0} f(x,Mh(\phi_1))\geq
\max_{x\in\overline\Omega_0}\Big\{K^*g(h(\phi_1))+2\la_1Mh(\phi_1)+
M^a(h')^a(\phi_1)|\nabla\phi_1|^a\Big\}.
\endeq
From \eq{interion} and \eq{lambda2} we deduce
\neweq{intfinal}
\di -\Delta
\underline{u}_{\la}+K(x)g(\underline{u}_{\la})+|\nabla\underline
u_\la|^a \leq \la f(x,\underline u_\la)\quad\mbox{ in }\Omega_0.
\endeq
Now, \eq{bordfinal} together with \eq{intfinal} shows that
$\underline u_\la=Mh(\phi_1)$ is a sub-solution of $(1)_\la$
provided $\la>0$ satisfy \eq{lambda1} and \eq{lambda2}. With the
same arguments as in the proof of Theorem \ref{th1jmaa} and using
Lemma \ref{l}, one can prove that $\underline u_\la\leq U_\la$ in
$\Omega.$ By a standard bootstrap argument (see Gilbarg and
Trudinger \cite{gt}) we obtain a classical solution $u_\la$ such
that $\underline u_\la\leq u_\la\leq U_\la$ in $\Omega.$

We have proved that $(1)_\la$ has at least one classical solution
when $\la>0$ is large. Set
$$ \di A=\{\la>0; \mbox{ problem }(1)_\la
\mbox{ has at least one classical solution}\}.$$ From the above
arguments we deduce that $A$ is nonempty. Let $\la^*=\inf A.$ We
claim that if $\la\in A,$ then $(\la,+\infty)\subseteq A.$ To this
aim, let $\la_1\in A$ and $\la_2>\la_1.$ If $u_{\la_1}$ is a
solution of $(1)_{\la_1},$ then $u_{\la_1}$ is a sub-solution for
$(1)_{\la_2}$ while $U_{\la_2}$  defined in \eq{UU} for
$\la=\la_2$  is a super-solution. Moreover, we have
$$\di \Delta U_{\la_2}+\la_2f(x,U_{\la_2})\leq 0\leq \Delta u_{\la_1}+
\la_2f(x,u_{\la_1})\quad\mbox{ in }\Omega,$$
$$\di U_{\la_2},u_{\la_1}>0\quad\mbox{ in }\Omega,$$
$$U_{\la_2}=u_{\la_1}=0\quad\mbox{ on }\partial\Omega$$
$$\Delta U_{\la_2}\in L^1(\Omega).$$
Again by Lemma \ref{l} we get $u_{\la_1}\leq U_{\la_2}$ in
$\Omega.$ Therefore, the problem $(1)_{\la_2}$ has at least one
classical solution. This proves the claim. Since $\la\in A$ was
arbitrary chosen, we conclude that $(\la^*,+\infty)\subset A.$

To end the proof, it suffices to show that $\la^*>0.$ In that
sense, we will prove that there exists $\la>0$ small enough such
that $(1)_\la$ has no classical solutions. We first remark that
$$\di \lim_{s\ri 0^+}(f(x,s)-K(x)g(s))=-\infty \quad
\mbox{ uniformly for }\, x\in\Omega.$$ Hence, there exists $c>0$
such that
\neweq{n1}
\di f(x,s)-K(x)g(s)<0, \quad \mbox{ for all }\,
(x,s)\in\Omega\times(0,c).
\endeq
On the other hand, the assumptions on $f$ yield
\neweq{n2}
\di \frac{f(x,s)-K(x)g(s)}{s}\leq \frac{f(x,s)}{s}\leq
\frac{f(x,c)}{c} \quad \mbox{ for all }\,
(x,s)\in\Omega\times[c,+\infty).
\endeq
Let $m=\max_{x\in\overline\Omega}\frac{f(x,c)}{c}.$ Combining
\eq{n1} with \eq{n2} we find
\neweq{n3}
\di f(x,s)-K(x)g(s)<ms, \quad \mbox{ for all }\,
(x,s)\in\Omega\times(0,+\infty).
\endeq
Set $\la_0=\min\left\{1,\la_1/2m\right\}.$ We show that problem
$(1)_{\la_0}$ has no classical solution. Indeed, if $u_0$ would be
a classical solution of $(1)_{\la_0},$ then, according to \eq{n3},
$u_0$ is a sub-solution of
\neweq{n4}
\left\{\begin{tabular}{ll} $\di-\Delta u=\frac{\la_1}{2}u$ \quad &
${\rm in}\
\Omega,$\\
$u>0$ \quad & ${\rm in}\ \Omega,$\\
$u=0$ \quad & ${\rm on}\ \partial\Omega.$\\
\end{tabular} \right.
\endeq
Obviously, $\phi_1$ is a super-solution of \eq{n4} and by Lemma
\ref{l} we get $u_0\leq \phi_1$ in $\Omega.$ Thus, by standard
elliptic arguments, problem \eq{n4} has a solution $u\in
C^2(\overline\Omega).$ Multiplying by $\phi_1$ in \eq{n4} and then
integrating over $\Omega$ we have
$$\di -\int_\Omega\phi_1\Delta udx=\frac{\la_1}{2}\int_\Omega
u\phi_1dx,$$ that is,
$$\di -\int_\Omega u\Delta\phi_1dx=\frac{\la_1}{2}\int_\Omega
u\phi_1dx.$$ The above equality yields $\int_\Omega u\phi_1dx=0,$
which is clearly a contradiction, since $u$ and $\phi_1$ are
positive in $\Omega.$ If follows that problem $(1)_{\la_0}$ has no
classical solutions which means that $\la^*>0.$ This completes the
proof of Theorem \ref{th2jmaa}. \qed


\begin{thebibliography}{999}  {\footnotesize

\bibitem{ap} N. E. Alaa and M. Pierre, Weak solutions of
some quasilinear elliptic equations with data measures, {\it SIAM
J. Math. Anal.} {\bf 24} (1993), 23-35.

\bibitem{al_t1} S. Alama and G. Tarantello, On the solvability of a
semilinear elliptic equation via an associated
eigenvalue problem, {\it Math. Z.}, {\bf 221} (1996), 467-493.

\bibitem{amann} H. Amann, Existence and multiplicity theorems for
semilinear elliptic boundary value problems, {\it Math. Z.} {\bf
150} (1976), 567-597.

\bibitem {a} R. Aris, {\it The Mathematical Theory of
Diffusion and Reaction in Permeable Catalysts}, Clarendon Press,
Oxford, 1975.

\bibitem{bandle} C. Bandle, Asymptotic behaviour of large solutions of quasilinear
elliptic problems, {\it Z. Angew. Math. Phys.} {\bf 54} (2003), 731-738.

\bibitem{be} C. Bandle and M. Ess\`{e}n, On the solutions of quasilinear elliptic problems
with boundary blow-up, in {\it Partial differential equations of elliptic type}
(Cortona, 1992), Sympos. Math. {\bf 35}, Cambridge Univ. Press,
Cambridge, 1994, p. 93-111.

\bibitem{bangia} C. Bandle and E. Giarrusso, Boundary blow-up for semilinear
elliptic equations with nonlinear gradient terms,
{\it Advances in Differential Equations} {\bf 1} (1996), 133-150.

\bibitem{bm} C. Bandle and M. Marcus,  'Large' solutions of semilinear
elliptic equations: Existence, uniqueness, and asymptotic behaviour,
{\it J. Anal. Math.} {\bf 58} (1992), 9-24.

\bibitem{bm1} C. Bandle and M. Marcus,  Dependence of blowup rate of large solutions
of semilinear elliptic equations on the curvature of the boundary,
{\it Complex Variables, Theory Appl.} {\bf 49} (2004), 555-570.

\bibitem {bdd} G. Barles, G. D\'iaz, and J. I. D\'iaz, Uniqueness and continuum of
foliated solutions for a quasilinear elliptic equation with a non
lipschitz nonlinearity, {\it Comm. Partial Differential Equations}
{\bf 17} (1992), 1037-1050.

\bibitem{bbc} P. B\'enilan, H. Brezis, and M. Crandall,
A semilinear equation in $L^1(\RR^N),$ {\it Ann. Scuola Norm. Sup.
Pisa} \textbf{4} (1975), 523-555.

\bibitem{bi} L. Bieberbach,  $\Delta u=e^u$ und die automorphen
Funktionen, {\it Math. Ann.} {\bf 77} (1916), 173-212.

\bibitem{bgt} N. H. Bingham, C. M. Goldie, and J. L. Teugels,
{\it Regular Variation}, Cambridge University Press, Cambridge, 1987.

\bibitem{bk} H. Brezis and S. Kamin, Sublinear elliptic equations in
$\RR^N,$ {\it Manuscripta Math.} {\bf 74} (1992), 87-106.

\bibitem{bros} H. Brezis and L. Oswald, Remarks on sublinear
elliptic equations, {\it Nonlinear Anal., T.M.A.} {\bf 10}
(1986), 55-64.

\bibitem{caf} L. Caffarelli, R. Hardt, and L. Simon, Minimal surfaces with isolated singularities,
{\it Manuscripta Math.} {\bf 48} (1984), 1-18.

\bibitem{cn1} A. Callegari and A. Nachman, Some
singular nonlinear equations arising in boundary layer theory,
{\it J.~Math. Anal. Appl.} {\bf 64} (1978), 96-105.

\bibitem {cn} A. Callegari and A. Nachman, A nonlinear
singular boundary value problem in the theory of pseudoplastic
fluids, {\it SIAM J. Appl. Math.} {\bf 38} (1980), 275-281.

\bibitem{chen} H. Chen, On a singular nonlinear elliptic equation,
{\it Nonlinear Anal., T.M.A.} {\bf 29} (1997), 337-345.

\bibitem{chipot} M. Chipot, {\it Elements of Nonlinear Analysis},
Birkh\"auser Advanced Texts, Birkh\"auser Verlag, 2000.

\bibitem{clm} Y. S. Choi, A. C.  Lazer, and P. J.
McKenna, Some remarks on a singular elliptic boundary value
problem, {\it Nonlinear Anal., T.M.A.} {\bf 3} (1998), 305-314.

\bibitem{leray} Y. Choquet-Bruhat and J. Leray, Sur le probl\`eme
de Dirichlet quasilin\'eaire d'ordre 2, {\it C.~R. Acad. Sci.
Paris, Ser.~A} {\bf 274} (1972), 81-85.

\bibitem{cgr} F.-C. C\^{\i}rstea, M. Ghergu, and V. R\u adulescu, Combined effects
of asymptotically linear and singular nonlinearities in bifurcation problems of Lane-Emden-Fowler type,
{\it J.~Math. Pures Appl.} {\bf 84} (2005), 493-508.

\bibitem{crna} F.-C. C\^{\i}rstea and V. R\u adulescu, Blow-up solutions for semilinear
elliptic problems, {\it Nonlinear Analysis, T.M.A.} {\bf 48}
(2002), 541-554.

\bibitem{crasun} F.-C. C\^{\i}rstea and V. R\u adulescu,
Uniqueness of the blow-up boundary solution of logistic
equations with absorption, {\it C.~R. Acad. Sci. Paris, Ser.~I} {\bf 335} (2002), 447-452.

\bibitem{crjmpa} F.-C. C\^{\i}rstea and V. R\u adulescu,
Entire solutions blowing-up at infinity for semilinear
elliptic systems, {\it J.~Math. Pures Appliqu\'ees} {\bf 81} (2002), 827-846.

\bibitem{crccm} F.-C. C\^{\i}rstea and V. R\u adulescu,
Existence and uniqueness of blow-up solutions for a class of
logistic equations, {\it Commun. Contemp. Math.} {\bf 4} (2002), 559-586.

\bibitem{crasas} F.-C. C\^{\i}rstea and V. R\u adulescu,
Asymptotics for the blow-up boundary solution of the
logistic equation with absorption, {\it C. R. Acad. Sci. Paris,
Ser.~I} {\bf 336} (2003), 231-236.

\bibitem{crhouston} F.-C. C\^{\i}rstea and V. R\u adulescu,
Solutions with boundary blow-up for a class of nonlinear
elliptic problems, {\it Houston J.~Math.} {\bf 29} (2003), 821-829.

\bibitem{crcras04} F.-C. C\^{\i}rstea and V. R\u adulescu,
Extremal singular solutions for degenerate logistic-type
equations in anisotropic media, {\it C. R. Acad. Sci. Paris,
Ser.~I} {\bf 339} (2004), 119-124.

\bibitem{crasan} F.-C. C\^{\i}rstea and V. R\u adulescu,
Nonlinear problems with boundary blow-up: a Karamata regular
variation theory approach", {\it Asymptotic Analysis}, in press.

\bibitem{crtams} F.-C. C\^{\i}rstea and V. R\u adulescu, Boundary blow-up in nonlinear elliptic equations of
Bieberbach--Rademacher type", {\it Transactions Amer. Math. Soc.}, in press.

\bibitem {ck} D. S. Cohen and H. B. Keller, Some
positive problems suggested by nonlinear heat generators, {\it J.
Math. Mech.} {\bf 16} (1967), 1361-1376.

\bibitem{cp} M. Coclite and G. Palmieri, On a singular nonlinear Dirichlet
problem, {\it Commun. Partial Diff. Equations} \textbf{14} (1989),
1315-1327.

\bibitem{crt} M. G. Crandall, P. H. Rabinowitz,
and L. Tartar, On a Dirichlet problem with a singular
nonlinearity, {\it Commun. Partial Diff. Equations} \textbf{2}
(1977), 193-222.

\bibitem{dalma} R. Dalmasso, Solutions d'\'equations elliptiques semi-lin\'eaires
singuli\`eres, {\it Ann. Mat. Pura Appl.} {\bf 153} (1989), 191-201.

\bibitem{gennes} P. G. de Gennes, Wetting: statics and dynamics, {\it Review of
Modern Physics} \textbf{57} (1985), 827x-863.

\bibitem{diaz} J. I. D\'iaz,  {\it Nonlinear Partial Differential Equations and Free
Boundaries. Vol. I. Elliptic Equations}, Research Notes in
Mathematics, vol.~106, Pitman (Advanced Publishing Program),
Boston, MA, 1985.

\bibitem{dmo} J. I. D\'iaz, J. M. Morel, and L. Oswald,
An elliptic equation with singular nonlinearity, {\it Comm. Partial Differential
Equations} {\bf 12} (1987), 1333-1344.

\bibitem{dh} Y. Du and Q. Huang, Blow-up solutions for a class of semilinear elliptic
and parabolic equations, {\it SIAM J. Math. Anal.} {\bf 31} (1999), 1-18.

\bibitem{dgrad} L. Dupaigne, M. Ghergu, and V. R\u adulescu,
Singular elliptic problems with convection term in anisotropic
media, in preparation.

\bibitem{ful} W. Fulks and J. S. Maybee, A singular nonlinear equation, {\it Osaka J.~Math.}
{\bf 12} (1960), 1-19.

\bibitem{gal} V. Galaktionov and J.-L. V\'azquez, The problem of blow-up in nonlinear
parabolic equations, {\it Discrete Contin. Dynam. Systems, Ser. A}
 {\bf 8} (2002), 399-433.

\bibitem{gls} J. Garc\'{i}a-Meli\'{a}n, R. Letelier-Albornoz, and J. Sabina de Lis,
 Uniqueness and asymptotic behaviour for solutions
of semilinear problems with boundary blow-up, {\it Proc. Amer. Math. Soc.}
{\bf 129} (2001), 3593-3602.

\bibitem{gnr} M. Ghergu, C. Niculescu, and V. R\u adulescu, Explosive solutions of
elliptic equations with absorption and nonlinear gradient term, {\it Proc. Indian Acad. Sci. (Math. Sci.)}
{\bf 112} (2002), 441-451.

\bibitem{gr1} M. Ghergu and V. R\u adulescu, Bifurcation and asymptotics for the
Lane-Emden-Fowler equation, {\it C. R. Acad. Sci. Paris, Ser.~I} \textbf{337} (2003), 259-264.

\bibitem{gr2} M. Ghergu and V. R\u adulescu, Sublinear singular
elliptic problems with two parameters, {\it J. Differential Equations} \textbf{195} (2003), 520-536.

\bibitem{grracsam} M. Ghergu and V. R\u adulescu, Explosive solutions of
semilinear elliptic systems with gradient term, {\it RACSAM Rev. Real  Acad. Cienc. Exactas
F\'{\i}s. Nat. Ser.~A  Mat.} {\bf 97} (2003), 437-445.

\bibitem{graaa} M. Ghergu and V. R\u adulescu, Existence and non-existence of entire solutions to the
logistic differential equation, {\it Abstract and Applied Analysis} {\bf 17} (2003), 995-1003.

\bibitem{grcras03} M. Ghergu and V. R\u adulescu, Bifurcation and asymptotics for the Lane-Emden-Fowler
equation, {\it C.~R. Acad. Sci. Paris, Ser.~I} {\bf 337} (2003), 259-264.

\bibitem{grcpaa} M. Ghergu and V. R\u adulescu, Nonradial blow-up solutions of
sublinear elliptic equations with gradient term, {\it Commun. Pure Appl. Anal.} {\bf 3} (2004), 465-474.

\bibitem{gr3} M. Ghergu and V. R\u adulescu, Bifurcation for a class of singular elliptic
problems with quadratic convection term, {\it C. R. Acad. Sci. Paris, Ser.~I} \textbf{338} (2004), 831-836.

\bibitem{gr4} M. Ghergu and V. R\u adulescu, Multiparameter bifurcation and asymptotics for the
singular Lane-Emden-Fowler equation with a convection term, {\it Proc. Royal Soc. Edinburgh
Sect. A} {\bf 135} (2005), 61-84.

\bibitem {grjmaa} M. Ghergu and V. R\u adulescu, On a class of sublinear
singular elliptic problems with convection term, {\it J.~Math. Anal. Appl.} {\bf 311} (2005), 635-646.

\bibitem{gt} D. Gilbarg and N. S. Trudinger, {\it Elliptic Partial
Differential Equations of Second Order}, 2nd ed., Springer Verlag,
Berlin, 1983.

\bibitem{gomes} S. M. Gomes, On a singular nonlinear
elliptic problem, {\it SIAM J. Math. Anal.} {\bf 17} (1986)
1359-1369.

\bibitem{gl} C. Gui and F. H. Lin, Regularity of an elliptic problem with a singular nonlinearity, {\it Proc. Royal
Soc. Edinburgh Sect. A} {\bf 123} (1993), 1021-1029.

\bibitem{hai} Y. Haitao, Multiplicity and asymptotic behavior of positive solutions for a singular semilinear elliptic
problem, {\it J.~Differential Equations} {\bf 189} (2003), 487-512.

\bibitem{her} J. Hern\'andez, F. J. Mancebo, and J. M. Vega, On the linearization
of some singular nonlinear elliptic problems and applications,
{\it Ann. Inst. H.~Poincar\'e, Anal. Non Lin\'eaire} {\bf 19}
(2002), 777-813.

\bibitem{her1} J. Hern\'andez, F. J. Mancebo, and J. M. Vega, Nonlinear singular elliptic problems:
recent results and open problems, {\it Preprint}, 2005.

\bibitem{h} L. H\"ormander, {\it The Analysis of Linear
Partial Differential Operators I}. Springer Verlag, Berlin, 1983.

\bibitem{kazdan} J. Kazdan and F. W. Warner, Remarks on some
quasilinear elliptic equations, {\it Comm. Pure Appl. Math.} {\bf
28} (1975), 567-597.

\bibitem{ke} J. B. Keller,  On solutions of $\Delta u=f(u)$, {\it Comm. Pure Appl.
Math.} {\bf 10} (1957), 503-510.

\bibitem{kusano} T. Kusano and C. A. Swanson, Entire positive solutions of singular
 elliptic equations, {\it Japan J.~Math.} {\bf 11} (1985), 145-155.

\bibitem{aw2} A. V. Lair and A. W. Shaker, Existence of entire large positive solutions of
semilinear elliptic systems, {\it J. Differential Equations} {\bf 164} (2000), 380-394.

\bibitem{lw} A. V. Lair and A. W. Wood, Large solutions of semilinear elliptic equations
with nonlinear gradient terms, {\it Internat.~J. Math. Math. Sci.} {\bf 22} (1999), 869-883.

\bibitem{ll} J. M. Lasry and P.-L. Lions, Nonlinear elliptic equations with singular boundary conditions and
stochastic control with state constraints; the model problem, {\it
Math. Ann.} {\bf 283} (1989), 583-630.

\bibitem{lm1} A. C. Lazer and P. J. McKenna, On a
singular nonlinear elliptic boundary value problem, {\it Proc.
Amer. Math. Soc.} {\bf 3} (1991), 720-730.

\bibitem{lm} A. C. Lazer and P. J. McKenna, On a problem of Bieberbach and
Rademacher, {\it Nonlinear Anal., T.M.A.} {\bf 21} (1993), 327-335.

\bibitem{lm2} A. C. Lazer and P. J. McKenna, Asymptotic behaviour of solutions of boundary
blowup problems, {\it Differential Integral Equations} {\bf 7} (1994), 1001-1019.

\bibitem{gall} J. F. Le Gall, A path-valued Markov process
and its connections with partial differential equations, in {\it
First European Congress of Mathematics}, Vol. II (Paris, 1992),
185-212, Progr. Math., 120, Birkh\"auser Verlag, Basel, 1994.

\bibitem{karamata} J. Karamata, Sur un mode de
croissance r\'{e}guli\`ere de fonctions. Th\'eor\`emes fondamentaux,
{\it Bull. Soc. Math. France} {\bf 61} (1933), 55-62.

\bibitem{ln} C. Loewner and L. Nirenberg, Partial differential equations
invariant under conformal or projective transformations, in {\it Contribution to
Analysis}, Academic Press, New York, 1974, p. 245-272.

\bibitem{marcus7} M. Marcus, On solutions with blow-up at the boundary for a
class of semilinear elliptic equations, in {\it Developments in Partial
Differential Equations and Applications to Mathematical Physics} (G.~Buttazzo {\it et
al.}, Eds.),  Plenum Press, New York (1992), 65-77.

\bibitem{mv1} M. Marcus and L. V\'eron, Uniqueness and asymptotic behavior of solutions with boundary blow-up for a
class of nonlinear elliptic equations, {\it Ann. Inst. H. Poincar\'e, Anal. Non Lin\'eaire} {\bf 14} (1997), 237-274.

\bibitem{mv2} M. Marcus and L. V\'eron, Existence and uniqueness results for large solutions of general nonlinear
elliptic equations, {\it J.~Evol. Equations} {\bf 3} (2003),
637-652.

\bibitem{mea} A. Meadows, Stable and singular solutions of the equation $\Delta u=1/u$, {\it Indiana Univ. Math.~J.}
{\bf 53} (2004),  1681-1703.

\bibitem{mr} P. Mironescu and V. R\u adulescu, The study
of a bifurcation problem associated to an asymptotically linear
function, {\it Nonlinear Anal., T.M.A.} \textbf{26} (1996), 857-875.

\bibitem{os} R. Osserman, On the inequality $\Delta u\geq f(u)$,
{\it Pacific J. Math.} {\bf 7} (1957), 1641-1647.

\bibitem{pino} M. del Pino, A global estimate for the
gradient in a singular elliptic boundary value problem, {\it Proc.
Roy. Soc. Edinburgh Sect. A} {\bf 122} (1992), 341-352.

\bibitem{qu} P. Quittner, Blow-up for semilinear parabolic equations with a gradient term,
{\it Math. Meth. Appl. Sci.} {\bf 14} (1991), 413-417.

\bibitem{rad} H. Rademacher, Einige besondere Probleme der partiellen
Differentialgleichungen, in {\it Die Differential und Integralgleichungen der
Mechanik und Physik I}, 2nd. edition, (P.~Frank und R. von Mises, eds.), Rosenberg,
New York, 1943, p. 838-845.

\bibitem{rrv} A. Ratto, M. Rigoli, and L. V\'eron, Scalar curvature and conformal deformation
of hyperbolic space, {\it J. Funct. Anal.} {\bf 121} (1994), 15-77.

\bibitem{radbirkh} V. R\u adulescu, Bifurcation and asymptotics for elliptic problems
with singular nonlinearity, in {\it Studies in Nonlinear Partial
Differential Equations: In Honor of Haim Brezis, Fifth European
Conference on Elliptic and Parabolic Problems: A special tribute
to the work of Haim Brezis}, Gaeta, Italy, May 30--June 3, 2004
(C.~Bandle, H.~Berestycki, B.~Brighi, A.~Brillard, M.~Chipot,
J.-M.~Coron, C.~Sbordone, I.~Shafrir, V.~Valente, G.~Vergara
Caffarelli, Eds.), Birkh\"auser Verlag, 2005, pp. 349-362.

\bibitem{seneta}  E. Seneta, {\it Regularly Varying
Functions}, Lecture Notes in Mathematics 508, Springer Verlag,
Berlin Heidelberg, 1976.

\bibitem {shi} J. Shi and M. Yao, On a singular nonlinear semilinear elliptic problem,
{\it Proc. Royal Soc. Edinburgh, Sect. A} {\bf 128} (1998), 1389-1401.

\bibitem{sy2} J. Shi and M. Yao, Positive solutions for elliptic equations with singular nonlinearity, {\it Electronic
Journal of Differential Equations} {\bf 4} (2005), 1-11.

\bibitem{stu} C. A. Stuart, Existence and approximation of solutions of nonlinear elliptic
equations, {\it Math.~Z.} {\bf 147} (1976), 53-63.

\bibitem{stuart1} C. A. Stuart, Self-trapping of an
electromagnetic field and bifurcation from the essential spectrum,
{\it Arch. Rational Mech. Anal.} {\bf 113} (1991), 65-96.

\bibitem{stuart2} C. A. Stuart and H.--S. Zhou, A variational problem related to self-trapping
of an electromagnetic field, {\it Math. Methods Appl. Sci.} {\bf 19} (1996), 1397-1407.

\bibitem{ver} L. V\'eron, {\it Singularities of Solutions of Second Order Quasilinear
Equations}, Pitman Res. Notes Math. Ser., Vol. 353, Longman, Harlow, 1996.

\bibitem{wong} J. S. W. Wong, On the generalized Emden-Fowler
equation, {\it SIAM Rev.} {\bf 17} (1975), 339-360.

\bibitem {zhang} Z. Zhang, On a Dirichlet problem with a
singular nonlinearity, {\it J.~Math. Anal. Appl.} {\bf 194}
(1995), 103-113.

\bibitem {z1} Z. Zhang, Nonexistence of positive classical solutions
of a singular nonlinear Dirichlet problem with a convection term,
{\it Nonlinear Anal., T.M.A.} {\bf 8} (1996), 957-961.

\bibitem {zy} Z. Zhang and J. Yu, On a singular nonlinear Dirichlet
problem with a convection term, {\it SIAM J. Math. Anal.} {\bf 4}
(2000), 916-927.

}

\end{thebibliography}
\end{document}